\documentclass[11pt]{article}

\usepackage[top=2.5cm, bottom=2.5cm, left=2.5cm, right=2.5cm]{geometry}   %

\usepackage{preambule}

\usepackage{notations}
\graphicspath{{figures/}}
\pgfplotsset{table/search path={figures},}

\makeatletter
\def\pgfplots@install@path@replacements{%
    \ifpgfplots@path@replace@ellipse
        \let\tikz@do@circle=\pgfplots@path@@tikz@do@circle
        \let\tikz@do@ellipse=\pgfplots@path@@tikz@do@circle
        \expandafter\def\expandafter\pgfinterruptpicture\expandafter{\pgfinterruptpicture
            \let\tikz@do@circle=\pgfplots@path@@tikz@do@circle@orig
            \let\tikz@do@ellipse=\pgfplots@path@@tikz@do@ellipse@orig
        }%
    \fi
}%
\let\pgfplots@path@@tikz@do@circle@orig=\tikz@do@circle
\let\pgfplots@path@@tikz@do@ellipse@orig=\tikz@do@ellipse

\let\pgfplots@path@@tikz@do@circle@oldandbroken=\pgfplots@path@@tikz@do@circle
\def\pgfplots@path@@tikz@do@circle#1{\pgfplots@path@@tikz@do@circle@oldandbroken{#1}{#1}}
\def\pgfplots@path@@tikz@do@ellipse#1#2{\pgfplots@path@@tikz@do@circle@oldandbroken{#1}{#2}}
\makeatother

\newboolean{afficherGraphes}
\setboolean{afficherGraphes}{true}

\newboolean{afficherCommentaires}
\setboolean{afficherCommentaires}{false}

\newboolean{brouillon}
\setboolean{brouillon}{false}

\ifthenelse{\boolean{brouillon}}{%
  \includeonly{sections/C_functional_framework.tex}%
}{}

\title{Time-harmonic wave propagation in junctions of two periodic half-spaces}
\author{Pierre Amenoagbadji, Sonia Fliss, Patrick Joly$^*$}
\date{POEMS, CNRS, Inria, ENSTA Paris, Institut Polytechnique de Paris, 91120 Palaiseau, France}

\begin{document}
\allowdisplaybreaks

\newcounter{taggedeq}
\setcounter{taggedeq}{0}
\pretocmd{\equation}{\stepcounter{taggedeq}}{}{}
\renewcommand*{\theHequation}{\thetaggedeq.\theequation}

\maketitle

\begin{abstract}
  We are interested in the Helmholtz equation in a junction of two periodic half-spaces. When the overall medium is periodic in the direction of the interface, Fliss and Joly (2019) proposed a method which consists in applying a partial Floquet-Bloch transform along the interface, to obtain a family of waveguide problems parameterized by the Floquet variable. In this paper, we consider two model configurations where the medium is no longer periodic in the direction of the interface. Inspired by the works of Gérard-Varet and Masmoudi (2011, 2012), and Blanc, Le Bris, and Lions (2015), we use the fact that the overall medium has a so-called quasiperiodic structure, in the sense that it is the restriction of a higher dimensional periodic medium. Accordingly, the Helmholtz equation is lifted onto a higher dimensional problem with coefficients that are periodic along the interface. This periodicity property allows us to adapt the tools previously developed for periodic media. However, the augmented PDE is elliptically degenerate (in the sense of the principal part of its differential operator) and thus more delicate to analyse.  
\end{abstract}

\ifthenelse{\boolean{brouillon}}{%
  \section{Introduction and description of the method}
\subsection{Presentation of the model problem}
Initially motivated by solid state theory in the 1940s \cite{kittel1996introduction}, the study of periodic media has since then sparked some significant interest in the mathematical physics \surligner{literature}, especially with the recent advent of photonic crystals \cite{joannopoulos2008molding, johnson2001photonic, kuchment2001mathematics, sakoda2004optical}. When it comes to wave propagation, one of the most remarkable properties of such structures from a practical point of view is the possibility of having \textit{band gaps}, that is, frequency ranges for which waves cannot propagate (we refer to \cite{kuchment2004some, kuchment2001mathematics} for a mathematical explanation). This feature is of particular interest in micro or nano-technology, for building optical filters or wave localization devices. In order to understand the underlying phenomena and to design these applications, it is important to have efficient numerical tools to solve the time-harmonic wave equation in presence of periodic media.

\begin{figure}[ht!]
  \centering
  \makebox[\textwidth][c]{
    \includegraphics[page=1]{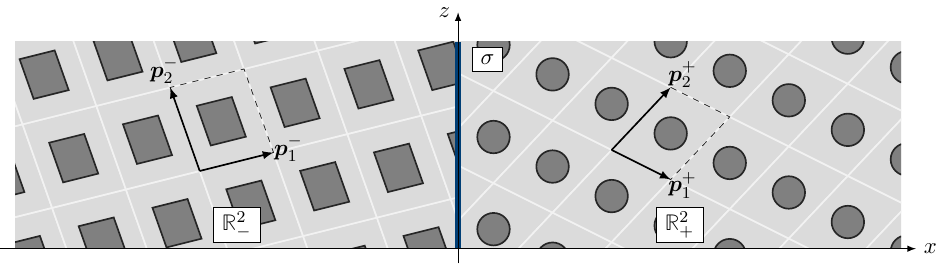}
  }
  \caption{Juxtaposition of arbitrary periodic half-spaces\label{H2Dmod:fig:general_configuration}}
\end{figure}

\vspace{1\baselineskip} \noindent
This paper is devoted to time-harmonic wave propagation in presence of two periodic half-spaces $\R^2_\pm := \{\xv = (\xvi, \zvi) \in \R^2\ /\ \pm \xvi > 0\}$ (see Figure \ref{H2Dmod:fig:general_configuration}). More precisely, let $\sigma := \{\xv = (\xvi, \zvi) \in \R^2\ /\ x = 0\}$ denote the interface between $\R^2_+$ and $\R^2_-$. The canonical basis of $\R^2$ is $\{\ev_\xvi, \ev_\zvi\}$. We are interested in \textit{finding $u \in H^1(\R^2)$ such that}
\begin{equation}
  \label{H2Dmod:eq:transmission_problem}
  \left\{
    \begin{array}{r@{\ =\ }l@{\quad}l}
      \displaystyle- \transp{\nabla}\, \aten\, \nabla u - \rho\, \omega^2\, u & 0 & \textnormal{in}\ \ \R^2_+ \cup \R^2_-,
      \\[8pt]
      \displaystyle \llbracket \aten\, \nabla u \cdot \ev_\xvi \rrbracket_{\sigma} & g & \textnormal{on}\ \ \sigma,
    \end{array}
  \right. \tag{$\mathscr{P}$} 
\end{equation}
where the presence of some (arbitrarily small) absorption is assumed (see Remark \ref{H2Dmod:rmk:justification_absorption}), namely
\begin{equation}\label{H2Dmod:eq:absorption}
  \text{the frequency }\omega \text{ satisfies }\Imag \omega > 0.
\end{equation}
The overall medium is represented by the tensor $\aten \in L^\infty(\R^2; \R^{2 \times 2})$ and the coefficient $\rho \in L^\infty(\R^2)$. We assume that $\aten$ is symmetric, and that
\begin{equation}
  \displaystyle
  \spexists a, r > 0, \quad \spforall \xv, \boldsymbol{\xi} \in \R^2, \qquad 0 < a \; |\boldsymbol{\xi}|^2 \leq \transp{\boldsymbol{\xi}}\; \aten(\xv)\; \boldsymbol{\xi} \quad \textnormal{and} \quad 0 < r \leq \rho(\xv).
  \label{H2Dmod:eq:ellipticity_assumption}
\end{equation}
Let $\pv^\pm_1, \pv^\pm_2 \in \R^2$ be linearly independent vectors. The tensor $\aten$ (\emph{resp.} the coefficient $\rho$) is assumed to coincide in $\R^2_\pm$ with a continuous function $\aten^\pm$ (\emph{resp.} $\rho^\pm$) which is $\Z \pv^\pm_1 + \Z \pv^\pm_2$--periodic, as illustrated in Figure \ref{H2Dmod:fig:general_configuration}. The exact structural properties of $(\aten, \rho)$ considered in this paper, namely the assumptions on $(\aten^\pm, \rho^\pm)$ and $(\pv^\pm_1, \pv^\pm_2)$, are detailed in Section \ref{H2Dmod:sec:specific_configs} and shown in Figures \ref{H2Dmod:fig:configuration_a} and \ref{H2Dmod:fig:configuration_b}.

\vspace{1\baselineskip} \noindent
In \eqref{H2Dmod:eq:transmission_problem}, we use the notation $\llbracket \boldsymbol{w} \cdot \ev_\xvi \rrbracket_{\sigma} \in H^{-1/2}(\sigma)$ for the jump accross $\sigma$ of the normal component of a function $\boldsymbol{w} \in [L^2(\R^2)]^2$ such that $\boldsymbol{w}|_{\R^2_\pm} \in H(\dive; \R^2_\pm)$:
\begin{equation}
  \llbracket \boldsymbol{w} \cdot \ev_\xvi \rrbracket_{\sigma} := (\boldsymbol{w}^- \cdot \ev_\xvi)|_{\sigma} - (\boldsymbol{w}^+ \cdot \ev_\xvi)|_{\sigma}, \quad \textnormal{with} \quad \boldsymbol{w}^\pm := \boldsymbol{w}|_{\R^2_\pm}.
\end{equation}
Finally $g \in H^{-1/2}(\sigma)$ is a given jump data \surligner{(see Remark \ref{H2Dmod:rmk:jump_data})}.

\begin{rmk}
  \begin{enumerate}[label=$(\textit{\alph*})$., ref=\thermk.\textit{\alph*}, wide = 0pt]
    \item\label{H2Dmod:rmk:jump_data} \surligner{We have considered a jump condition in \eqref{H2Dmod:eq:transmission_problem} for simplicity. This condition allows to reduce \eqref{H2Dmod:eq:transmission_problem} to an interface equation, as shown in Section \ref{H2Dmod:sec:DtN_waveguide} (see also Figure \ref{H2Dmod:fig:Uhat_from_Uplusminushat}). In Section \ref{H2Dmod:sec:extension_source_term}, we explain how the method in this paper extends to the case of a volume source term (that is, a right-hand side in the volume equation of \eqref{H2Dmod:eq:transmission_problem}).}
    \item\label{H2Dmod:rmk:justification_absorption}%
    For real-valued frequencies $\omega$, \surligner{one expects that \eqref{H2Dmod:eq:transmission_problem} does not admit any solution in $H^1(\R^2)$, due to a lack of sufficient decay in $|\xv|$ at infinity, caused by wave propagation up to infinity. Moreover, uniqueness of a solution in $H^1_{\textit{loc}}(\R^2)$ does not hold in general. One classical way to define the physical solution for wave propagation problems in unbounded media is to add a so-called \emph{radiation condition}, which describes the behavior of the solution at infinity. In the homogeneous case, that is, when $\aten$ and $\rho$ are constant, this is the so-called Sommerfeld radiation condition \cite{sommerfeld1912}.} However, to our knowledge, such a radiation condition is not known in presence of unbounded periodic media. 

    \hspace*{6mm} \surligner{The physical solution can also be defined using the \emph{limiting absorption principle}}, which consists in (\emph{1}) adding an imaginary part to $\omega$ (called absorption), and in (\emph{2}) studying the limit process as $\Imag \omega \to 0$. \surligner{If the limit in a certain sense exists, then it is the physical solution.} The limiting absorption principle is well-established for time-harmonic wave propagation problems in unbounded media that are homogeneous or stratified outside a bounded domain; see for instance \cite{agmon1975spectral, eidus1986limiting, wilcox1966wave,weder1990spectral}. More recently, it has been successfully applied to periodic closed waveguides \cite{fliss2009analyse, hoang2011limiting, fliss2016solutions, fliss2021dirichlet}, periodic layers \cite{kirsch2022scattering, kirsch2018radiation, kirsch2018limiting} (see also \cite{bonnet2002diffraction}), and to $n$--dimensional fully periodic media with $n \geq 2$ \cite{murata2006asymptotics,radosz2015new,mandel2019limiting}. However, as far as we know, there is no complete answer for transmission problems between periodic half-spaces. This is the reason why we assume the presence of some arbitrarily small but non-vanishing absorption (namely \eqref{H2Dmod:eq:absorption}), which is an essential step in understanding the non-absorbing case. \surligner{We also refer to Section \ref{H2Dmod:sec:extension_non_absorbing} for further discussion on the limiting absorption principle.}
  \end{enumerate}
\end{rmk}

\vspace{1\baselineskip} \noindent
Under Assumptions \eqref{H2Dmod:eq:absorption} and \eqref{H2Dmod:eq:ellipticity_assumption}, it follows from Lax-Milgram's theorem that Problem \eqref{H2Dmod:eq:transmission_problem} admits a unique solution $u \in H^1(\R^2)$. Moreover, if $g \in L^2(\sigma)$ is compactly supported, then a Combes-Thomas estimate \cite{combes1973asymptotic} \surligner{(see also Remark \ref{H2D:rmk:proof_exponential_decay})} allows to prove that this solution decays exponentially at infinity in all directions, namely
\begin{equation}\label{H2Dmod:eq:exp_decay_all_directions}
  \displaystyle
  \spexists c,\, \alpha > 0, \quad \big\|u\, \exp(\alpha \Imag \omega\, |\xv|) \big\|_{H^1(\R^2)} \leq c\; \|g\|_{L^2(\sigma)}.
\end{equation}
To solve \eqref{H2Dmod:eq:transmission_problem}, a naive method relying on this decay estimate would consist in truncating the computational domain at a certain distance related to $\Imag \omega$, with homogeneous Dirichlet boundary conditions for instance. However the accuracy of such a method is prone to deterioration as $\Imag \omega \to 0$. Worse, if $g$ is not compactly supported, then only the exponential decay in the direction normal to the interface is guaranteed:
\begin{equation}\label{H2Dmod:eq:exp_decay}
 \spexists c,\, \alpha > 0, \quad \big\|u\, \exp(\alpha \Imag \omega\, |\xvi|)\big\|_{H^1(\R^2)} \leq c\; \|g\|_{L^2(\sigma)}.
\end{equation}
Using the properties of $\aten$ and $\rho$, our goal is to develop a rigorous numerical method which allows to deal with the unboundedness of the domain, and which we hope will remain robust as $\Imag \omega$ tends to $0$. In the case where $\aten$ and $\rho$ are periodic in the direction of the interface $\sigma$, such a method has been proposed by \cite{flisscassanbernier2010}, in the spirit of \cite{besse:hal-00698916, fliss2009exact,fliss2012wave}. The underlying idea is to apply a Floquet-Bloch transform with respect to the variable $\zvi$ along the interface, to obtain a family of closed waveguide problems parameterized by the Floquet variable, and which can be solved using the Dirichlet-to-Neumann (DtN) approach proposed in \cite{fliss2009analyse,joly2006exact}.

\vspace{1\baselineskip} \noindent
Our goal is to extend the method in \cite{flisscassanbernier2010} to the case where the global medium represented by $\aten$ and $\rho$ is no longer periodic in the direction of the interface. We shall use the crucial (but non-obvious) observation that the medium has a \textit{quasiperiodic} structure along the interface, namely, that it is the restriction of a higher dimensional periodic structure. Accordingly, the idea is to interpret \eqref{H2Dmod:eq:transmission_problem} as the \textquote{restriction} of an \emph{augmented} boundary value problem in higher dimensions, where periodicity along the interface is recovered. This so-called \textit{lifting approach} allows one to adapt the ideas in \cite{flisscassanbernier2010}, but comes with the price that the augmented partial differential equation (PDE) is elliptically degenerate (in the sense of the principal part of its differential operator), and hence more complicated to analyse and approximate.

\vspace{1\baselineskip} \noindent
Although the lifting approach already appears in the works by Kozlov \cite{kozlov1979averaging}, the idea of using this approach for periodic half-spaces came to us mainly from the more recent papers \cite{gerard2011homogenization,gerard2012homogenization,blanc2015local}. In these papers, the lifting approach has been used in the homogenization setting (that is for instance, when the period of the structure is negligible with respect to the wavelength) in presence of a boundary or an interface. The lifting approach has also been used in \cite{bouchitte2010homogenization, wellander2019homogenization} for the homogenization of quasicrystals. However, to our knowledge, it seems that this idea has never been used for the analysis or simulation of wave propagation phenomena. The first approach in this direction has been developed in \cite{amenoagbadji2023wave} for wave propagation in $1$--dimensional quasiperiodic structures. 

\subsection{Formal description of the method}\label{H2Dmod:sec:formal_description_method}
We provide a brief and formal description of the method developed in this paper. In what follows, the generic $3$--dimensional space variable is denoted by $\xva = (\xvi, \zvi_1, \zvi_2)$. The canonical basis of $\R^3$ is denoted by $\{\eva_\xvi, \eva_1, \eva_2\}$. 

\vspace{1\baselineskip} \noindent
Depending on the restrictions $(\aten^\pm, \rho^\pm)$ and the periodicity vectors $(\pv^\pm_1, \pv^\pm_2)$, $\aten$ and $\rho$ may not be periodic with respect to $\zvi$. Nevertheless, under some suitable assumptions on $(\aten^\pm, \rho^\pm)$ and $(\pv^\pm_1, \pv^\pm_2)$ that we shall introduce in Section \ref{H2Dmod:sec:specific_configs} (see also Figures \ref{H2Dmod:fig:configuration_a} and \ref{H2Dmod:fig:configuration_b}), we have from Section \ref{H2Dmod:sec:hidden_quasiperiodicity} that
\begin{equation}
  \label{H2Dmod:eq:intro_general_quasiperiodic_expression}
  \aeforall \xv \in \R^2, \quad \aten(\xv) = \aten_p(\cutmat\, \xv) \quad \textnormal{and} \quad \rho(\xv) = \rho_p(\cutmat\, \xv),
\end{equation}
for some $\cutmat \in \R^{3 \times 2}$. The tensor $\aten_p \in L^\infty(\R^3; \R^{2 \times 2})$ (\emph{resp.} the coefficient $\rho_p \in L^\infty(\R^3)$) coincides in the half-space $\R^3_\pm := \{\xva = (\xvi, \zvi_1, \zvi_2) \in \R^3\ /\ \pm \xvi > 0\}$ with a $\Z^3$--periodic function $\aten^\pm_p$ (\emph{resp.} $\rho^\pm_p$), as Figures \ref{H2Dmod:fig:augmented_structure_config_a} and \ref{H2Dmod:fig:augmented_structure_config_b} show. As a consequence, contrary to $\aten$ and $\rho$, $\aten_p$ and $\rho_p$ are $1$--periodic with respect to the variables $(\zvi_1, \zvi_2)$ of the interface $\Sigma$ between $\R^3_+$ and $\R^3_-$:
\begin{equation}
  \label{H2Dmod:eq:def_Sigma}
  \Sigma := \{\xva = (\xvi, \zvi_1, \zvi_2) \in \R^3\ /\ \xvi = 0\}.
\end{equation}
Equation \eqref{H2Dmod:eq:intro_general_quasiperiodic_expression} shows the \emph{quasiperiodic} nature of $\aten$ and $\rho$, and suggests seeking the solution $u$ of \eqref{H2Dmod:eq:transmission_problem} as the restriction of a $3$D function $U$ along the hyperplane $\cutmat\, \R^2$, that is:
\begin{equation}
  \displaystyle
  \spforall \xv \in \R^2, \quad u(\xv) = U(\cutmat\, \xv).
  \label{H2Dmod:eq:lifting_ansatz}
\end{equation}
The extension $U$ shall be characterized as the solution of a 3D \textquote{augmented} problem with periodic coefficients $\aten_p$ and $\rho_p$. In order to construct such a problem for $U$, we formally use a chain rule, which links the partial derivatives of $u$ with those of $U$: if $\bsnabla := \transp{(\partial_\xvi\ \ \partial_{\zvi_1}\ \ \partial_{\zvi_2})}$ denotes the $3$D gradient operator, then given $F \in \mathscr{C}^1(\R^3)$, one has
\begin{equation}
  \displaystyle
  \label{H2Dmod:eq:chain_rule}
  \spforall \xv \in \R^2, \quad \big[\nabla F(\cutmat \cdot)\big](\xv) = \big[\transp{\cutmat} \bsnabla F\big](\cutmat\, \xv).
\end{equation}
Using the ansatz \eqref{H2Dmod:eq:lifting_ansatz} and the chain rule \eqref{H2Dmod:eq:chain_rule} in the volume equation satisfied by $u$ suggests to introduce
\begin{subequations}\label{H2Dmod:eq:formal_augmented_BVP}
\begin{equation}
  \displaystyle- \transp{\bsnabla} \cutmat\vts \aten_p\!  \transp{\cutmat} \bsnabla\, U - \rho_p\, \omega^2\, U = 0 \quad \textnormal{in}\ \  \R^3_+ \cup \R^3_-.
  \label{H2Dmod:eq:formal_augmented_PDE}
\end{equation}
In addition, the jump condition in \eqref{H2Dmod:eq:transmission_problem} on the line $\sigma$ may be formally lifted into a jump condition on the plane $\Sigma$, which is the interface defined by \eqref{H2Dmod:eq:def_Sigma}. This can be written as
\begin{equation}
  \displaystyle
  \llbracket (\cutmat\vts \aten_p \! \transp{\cutmat} \bsnabla\, U) \cdot \eva_\xvi \rrbracket_{\Sigma} = G.
  \label{H2Dmod:eq:formal_jump_condition}
\end{equation}
\end{subequations}
Here, $\llbracket \cutmat\, \boldsymbol{W} \cdot \eva_\xvi \rrbracket_{\Sigma} := (\cutmat\, \boldsymbol{W}^- \cdot \eva_\xvi)|_{\Sigma} - (\cutmat\, \boldsymbol{W}^+ \cdot \eva_\xvi)|_{\Sigma}$, with $\boldsymbol{W}^\pm := \boldsymbol{W}|_{\R^3_\pm}$ for any $\boldsymbol{W} : \R^3 \to \C^2$, and the data $G: \Sigma \to \C$ must formally satisfy the following condition
\begin{equation}
 \aeforall \xv = (0, \zvi) \in \sigma, \quad G(\cutmat\, \xv) = g(\xv),
  \label{H2Dmod:eq:compatibility}
\end{equation}
by consistency with the jump condition in \eqref{H2Dmod:eq:transmission_problem}. The appropriate functional framework to study the augmented problem \eqref{H2Dmod:eq:formal_augmented_BVP} is developed in Section \ref{H2Dmod:sec:functional_framework}.

\begin{figure}[ht!]
  \noindent\makebox[\textwidth]{%
    \begin{tikzpicture}[scale=0.5]
      \def\angleinterface{55}%
      \definecolor{coplancoupe}{RGB}{211, 227, 65}
      \tdplotsetmaincoords{77}{40}
      \def\coteCyl{1.5}
      \def\coteInf{6}
      \def\BBxmin{-4.5}
      \def\BBxmax{+4.5}
      \def\BBymin{-6}
      \def\BBymax{+6}
      \def\BBzmin{0}
      \def\BBzmax{+2.75}
      \def\coupevecY{1.414213562}%
      \def\coupevecZ{0.75}
      \def\ecart{2.5}
      \pgfmathsetmacro\origineStrip{\BBxmax+\ecart-\BBxmin}
      \pgfmathsetmacro\origineCyl{\origineStrip+\BBxmax+\ecart-\BBxmin}
      \pgfmathsetmacro\origineCell{\origineCyl+\BBxmax+\ecart+\coteCyl}
      \pgfmathsetmacro\xcoupemax{min(\BBymax/\coupevecY, \BBzmax/\coupevecZ)}
      \pgfmathsetmacro\xcoupemin{max(\BBymin/\coupevecY, \BBzmin/\coupevecZ)}
      \begin{scope}[tdplot_main_coords, scale=0.9]%
        %
        %
        \begin{scope}
          \fill[canvas is zx plane at y=\BBymax, fill=secondary!70, opacity=0.5] (0, \BBxmin) rectangle (\coteInf, 0);
          \fill[canvas is zx plane at y=\BBymax, fill=tertiary!70, opacity=0.25] (0, \BBxmax) rectangle (\coteInf, 0);
          \fill[canvas is yz plane at x=\BBxmin, fill=secondary!60, opacity=0.675] (\BBymin, 0) rectangle (\BBymax, \coteInf);
          \fill[secondary!80, opacity=0.675] (\BBxmin, \BBymin, 0) -- +(-\BBxmin, 0, 0) -- +(-\BBxmin, \BBymax-\BBymin, 0) -- +(0, \BBymax-\BBymin, 0) -- cycle;
          \fill[canvas is xy plane at z=\coteInf, secondary!85, opacity=0.675] (\BBxmin, \BBymin) -- +(-\BBxmin, 0) -- +(-\BBxmin, \BBymax-\BBymin) -- +(0, \BBymax-\BBymin) -- cycle;
          \filldraw[thick, draw=white, fill=secondary!95, opacity=0.875] (0, \BBymin, 0) -- +(0, \BBymax-\BBymin, 0) -- +(0, \BBymax-\BBymin, \coteInf) -- +(0, 0, \coteInf) -- cycle;
          \fill[canvas is zx plane at y=\BBymin, fill=secondary!80, opacity=0.675] (0, \BBxmin) rectangle (\coteInf, 0);
          \fill[canvas is xy plane at z=0, tertiary!60, opacity=0.675] (\BBxmax, \BBymin) -- +(-\BBxmax, 0) -- +(-\BBxmax, \BBymax-\BBymin) -- +(0, \BBymax-\BBymin) -- cycle;
          \fill[canvas is xy plane at z=\coteInf, tertiary!85, opacity=0.675] (\BBxmax, \BBymin) -- +(-\BBxmax, 0) -- +(-\BBxmax, \BBymax-\BBymin) -- +(0, \BBymax-\BBymin) -- cycle;
          \fill[canvas is zx plane at y=\BBymin, fill=tertiary!80, opacity=0.675] (0, \BBxmax) rectangle (\coteInf, 0);
          \fill[canvas is yz plane at x=\BBxmax, fill=tertiary!87, opacity=0.675] (\BBymin, 0) rectangle (\BBymax, \coteInf);
          \foreach \idI in {-3, ..., 3} {
            \draw[white, opacity=0.4, thin] (\BBxmax, \idI*\coteCyl, 0) -- (\BBxmax, \idI*\coteCyl, \coteInf);
            \draw[white, opacity=0.4, thin] (\BBxmax, \idI*\coteCyl, \coteInf) -- (\BBxmin, \idI*\coteCyl, \coteInf);
          }
          \foreach \idI in {1, ..., 4} {
            \draw[white, opacity=0.4, thin] (\BBxmax, \BBymin, \idI*\coteCyl) -- (\BBxmax, \BBymax, \idI*\coteCyl);
            \draw[white, opacity=0.4, thin] (\BBxmax, \BBymin, \idI*\coteCyl) -- (\BBxmin, \BBymin, \idI*\coteCyl);
          }
          \foreach \idI in {-3, ..., 3} {
            \draw[white, opacity=0.4, thin] (\idI*\coteCyl, \BBymin, \coteInf) -- (\idI*\coteCyl, \BBymax, \coteInf);
            \draw[white, opacity=0.4, thin] (\idI*\coteCyl, \BBymin, 0) -- (\idI*\coteCyl, \BBymin, \coteInf);
          }
          \draw (0, 0, -2.5) node {$\R^3$};
        \end{scope}

        \begin{scope}
          \fill[canvas is zx plane at y=\BBymax, fill=secondary!70, opacity=0.5] (0, \BBxmin+\origineStrip) rectangle (\coteCyl, \origineStrip);
          \fill[canvas is yz plane at x=\BBxmin+\origineStrip, fill=secondary!60, opacity=0.675] (\BBymin, 0) rectangle (\BBymax, \coteCyl);
          \fill[secondary!80, opacity=0.675] (\BBxmin+\origineStrip, \BBymin, 0) -- +(-\BBxmin, 0, 0) -- +(-\BBxmin, \BBymax-\BBymin, 0) -- +(0, \BBymax-\BBymin, 0) -- cycle;
          \fill[secondary!90!black, opacity=0.75] (\BBxmin+\origineStrip, \BBymin, \coteCyl) -- +(-\BBxmin, 0, 0) -- +(-\BBxmin, \BBymax-\BBymin, 0) -- +(0, \BBymax-\BBymin, 0) -- cycle;
          \filldraw[thick, draw=white, fill=secondary!95, opacity=0.875] (\origineStrip, \BBymin, 0) -- +(0, \BBymax-\BBymin, 0) -- +(0, \BBymax-\BBymin, \coteCyl) -- +(0, 0, \coteCyl) -- cycle;
          \fill[canvas is zx plane at y=\BBymin, fill=secondary!80, opacity=0.675] (0, \BBxmin+\origineStrip) rectangle (\coteCyl, \origineStrip);
          \fill[canvas is zx plane at y=\BBymax, fill=tertiary!70, opacity=0.25] (0, \BBxmax+\origineStrip) rectangle (\coteCyl, \origineStrip);
          \fill[canvas is xy plane at z=0, tertiary!60, opacity=0.675] (\BBxmax+\origineStrip, \BBymin) -- +(-\BBxmax, 0) -- +(-\BBxmax, \BBymax-\BBymin) -- +(0, \BBymax-\BBymin) -- cycle;
          \fill[canvas is xy plane at z=\coteCyl, tertiary!95!black, opacity=0.75] (\BBxmax+\origineStrip, \BBymin) -- +(-\BBxmax, 0) -- +(-\BBxmax, \BBymax-\BBymin) -- +(0, \BBymax-\BBymin) -- cycle;
          \fill[canvas is zx plane at y=\BBymin, fill=tertiary!80, opacity=0.675] (0, \BBxmax+\origineStrip) rectangle (\coteCyl, \origineStrip);
          \fill[canvas is yz plane at x=\BBxmax+\origineStrip, fill=tertiary!87, opacity=0.675] (\BBymin, 0) rectangle (\BBymax, \coteCyl);
          \foreach \idI in {-4, ..., 3} {
            \draw[white, opacity=0.4, thin] (\BBxmax+\origineStrip, \idI*\coteCyl, 0) -- (\BBxmax+\origineStrip, \idI*\coteCyl, \coteCyl);
            \draw[white, opacity=0.4, thin] (\BBxmax+\origineStrip, \idI*\coteCyl, \coteCyl) -- (\BBxmin+\origineStrip, \idI*\coteCyl, \coteCyl);
          }
          \foreach \idI in {-3, ..., 3} {
            \draw[white, opacity=0.4, thin] ({\idI*\coteCyl+\origineStrip}, \BBymin, \coteCyl) -- ({\idI*\coteCyl+\origineStrip}, \BBymax, \coteCyl);
            \draw[white, opacity=0.4, thin] ({\idI*\coteCyl+\origineStrip}, \BBymin, 0) -- ({\idI*\coteCyl+\origineStrip}, \BBymin, \coteCyl);
          }
          \draw (\origineStrip, 0, -2.5) node {$\R^2 \times (0, 1)$};
          \node (fl1) at (\BBxmax, 0, 0.75*\coteInf) {};
          \node (fl2) at (\origineStrip, 0, \coteCyl) {};
          \draw[-latex] (fl1) to[bend left=40] (fl2);
          \draw (\origineStrip+0.25, 0, \coteCyl+3.5) node {Step \eqref{H2Dmod:item:lifting_step_1}};
        \end{scope}%
        \begin{scope}
          \fill[secondary!80, opacity=0.5] (\BBxmin+\origineCyl, 0, 0) -- +(-\BBxmin, 0, 0) -- +(-\BBxmin, \coteCyl, 0) -- +(0, \coteCyl, 0) -- cycle;
          \fill[secondary!90, opacity=0.625] (\BBxmin+\origineCyl, \coteCyl, 0) -- +(-\BBxmin, 0, 0) -- +(-\BBxmin, 0, \coteCyl) -- +(0, 0, \coteCyl) -- cycle;
          \draw[dashed, white] (\BBxmin+\origineCyl, \coteCyl, 0) -- (\BBxmax+\origineCyl, \coteCyl, 0);
          \filldraw[draw=white, fill=secondary!90, opacity=0.625] (\BBxmin+\origineCyl, 0, 0) -- +(-\BBxmin, 0, 0) -- +(-\BBxmin, 0, \coteCyl) -- +(0, 0, \coteCyl) -- cycle;
          \filldraw[draw=white, fill=secondary!80, opacity=0.5] (\BBxmin+\origineCyl, 0, \coteCyl) -- +(-\BBxmin, 0, 0) -- +(-\BBxmin, \coteCyl, 0) -- +(0, \coteCyl, 0) -- cycle;
          \filldraw[draw=white, fill=secondary!95!black, opacity=0.875] (\origineCyl, 0, 0) -- (\origineCyl, \coteCyl, 0) -- (\origineCyl, \coteCyl, \coteCyl) -- (\origineCyl, 0, \coteCyl) -- cycle;
          \fill[tertiary!80, opacity=0.5] (\BBxmax+\origineCyl, 0, 0) -- +(-\BBxmax, 0, 0) -- +(-\BBxmax, \coteCyl, 0) -- +(0, \coteCyl, 0) -- cycle;
          \fill[tertiary!90, opacity=0.625] (\BBxmax+\origineCyl, \coteCyl, 0) -- +(-\BBxmax, 0, 0) -- +(-\BBxmax, 0, \coteCyl) -- +(0, 0, \coteCyl) -- cycle;
          \draw[dashed, white] (\BBxmax+\origineCyl, \coteCyl, 0) -- (\BBxmax+\origineCyl, \coteCyl, 0);
          \filldraw[draw=white, fill=tertiary!90, opacity=0.625] (\BBxmax+\origineCyl, 0, 0) -- +(-\BBxmax, 0, 0) -- +(-\BBxmax, 0, \coteCyl) -- +(0, 0, \coteCyl) -- cycle;
          \filldraw[draw=white, fill=tertiary!80, opacity=0.5] (\BBxmax+\origineCyl, 0, \coteCyl) -- +(-\BBxmax, 0, 0) -- +(-\BBxmax, \coteCyl, 0) -- +(0, \coteCyl, 0) -- cycle;
          \foreach \idI in {-3, ..., 3} {
            \draw[white, opacity=0.4, thin] ({\idI*\coteCyl+\origineCyl}, 0, \coteCyl) -- ({\idI*\coteCyl+\origineCyl}, \coteCyl, \coteCyl);
            \draw[white, opacity=0.4, thin] ({\idI*\coteCyl+\origineCyl}, 0, 0) -- ({\idI*\coteCyl+\origineCyl}, 0, \coteCyl);
          }
          \draw (\origineCyl, 0, -1.5) node {$\R \times (0, 1)^2$};
          \node (fl1) at (0.5*\BBxmax+\origineStrip, 0, \coteCyl) {};
          \node (fl2) at (\origineCyl, 0, \coteCyl) {};
          \draw[-latex] (fl1) to[bend left=60] (fl2);
          \draw (0.25*\BBxmax+0.5*\origineStrip+0.5*\origineCyl, 0, \coteCyl+3) node {Step \eqref{H2Dmod:item:lifting_step_2}}; 
        \end{scope}%
        \begin{scope}
          \fill[secondary!80, opacity=0.5] (\origineCell-0.5, 0, 0) -- +(-\coteCyl, 0, 0) -- +(-\coteCyl, \coteCyl, 0) -- +(0, \coteCyl, 0) -- cycle;
          \fill[secondary!90, opacity=0.625] (\origineCell-0.5, \coteCyl, 0) -- +(-\coteCyl, 0, 0) -- +(-\coteCyl, 0, \coteCyl) -- +(0, 0, \coteCyl) -- cycle;
          \draw[dashed, white] (\origineCell-0.5, \coteCyl, 0) -- (\origineCell-0.5, \coteCyl, 0);
          \filldraw[canvas is yz plane at x={\origineCell-0.5-\coteCyl}, secondary!90, draw=white, dashed, opacity=0.625] (0, 0) rectangle (\coteCyl, \coteCyl);
          \filldraw[canvas is yz plane at x=\origineCell-0.5, fill=secondary!90, draw=white, opacity=0.75] (0, 0) rectangle (\coteCyl, \coteCyl);
          \filldraw[draw=white, fill=secondary!90, opacity=0.625] (\origineCell-0.5, 0, 0) -- +(-\coteCyl, 0, 0) -- +(-\coteCyl, 0, \coteCyl) -- +(0, 0, \coteCyl) -- cycle;
          \filldraw[draw=white, fill=secondary!80, opacity=0.5] (\origineCell-0.5, 0, \coteCyl) -- +(-\coteCyl, 0, 0) -- +(-\coteCyl, \coteCyl, 0) -- +(0, \coteCyl, 0) -- cycle;
          \filldraw[canvas is yz plane at x=\origineCell, fill=teal!90, draw=white, opacity=1] (0, 0) rectangle (\coteCyl, \coteCyl);
          \fill[tertiary!80, opacity=0.5] (\coteCyl+\origineCell+0.5, 0, 0) -- +(-\coteCyl, 0, 0) -- +(-\coteCyl, \coteCyl, 0) -- +(0, \coteCyl, 0) -- cycle;
          \fill[tertiary!90, opacity=0.625] (\coteCyl+\origineCell+0.5, \coteCyl, 0) -- +(-\coteCyl, 0, 0) -- +(-\coteCyl, 0, \coteCyl) -- +(0, 0, \coteCyl) -- cycle;
          \draw[dashed, white] (\coteCyl+\origineCell+0.5, \coteCyl, 0) -- (\coteCyl+\origineCell+0.5, \coteCyl, 0);
          \filldraw[canvas is yz plane at x=\origineCell+0.5, tertiary!90, draw=white, dashed, opacity=0.625] (0, 0) rectangle (\coteCyl, \coteCyl);
          \filldraw[canvas is yz plane at x=\origineCell+0.5+\coteCyl, fill=tertiary!90, draw=white, opacity=0.75] (0, 0) rectangle (\coteCyl, \coteCyl);
          \filldraw[draw=white, fill=tertiary!90, opacity=0.625] (\coteCyl+\origineCell+0.5, 0, 0) -- +(-\coteCyl, 0, 0) -- +(-\coteCyl, 0, \coteCyl) -- +(0, 0, \coteCyl) -- cycle;
          \filldraw[draw=white, fill=tertiary!80, opacity=0.5] (\coteCyl+\origineCell+0.5, 0, \coteCyl) -- +(-\coteCyl, 0, 0) -- +(-\coteCyl, \coteCyl, 0) -- +(0, \coteCyl, 0) -- cycle;
          \draw (\origineCell+0.5*\coteCyl, 0, -1) node {$(-1, 1) \times (0, 1)^2$};
          \node (fl1) at (0.5*\coteCyl+\origineCyl, 0, \coteCyl) {};
          \node (fl2) at (\origineCell+0.5*\coteCyl, 0, \coteCyl) {};
          \draw[-latex] (fl1) to[bend left=40] (fl2);
          \draw (0.5*\origineCyl+0.5*\origineCell+0.5*\coteCyl, 0, \coteCyl+2.25) node {Step \eqref{H2Dmod:item:lifting_step_3}}; 
        \end{scope}%

        \pgfmathsetmacro\xbasisvec{\BBxmax}
        \pgfmathsetmacro\ybasisvec{\BBymin-8}
        \pgfmathsetmacro\zbasisvec{0.5}
        \node (evax) at (\xbasisvec+1.5, \ybasisvec, \zbasisvec) {$\eva_\xvi$}; 
        \draw[-latex] (\xbasisvec-0.5, \ybasisvec, \zbasisvec) -- (evax);
        \node (evay) at (\xbasisvec, \ybasisvec+3, \zbasisvec) {$\eva_1$}; 
        \draw[-latex] (\xbasisvec, \ybasisvec-0.5, \zbasisvec) -- (evay);
        \node (evaz) at (\xbasisvec, \ybasisvec, \zbasisvec+1.5) {$\eva_2$}; 
        \draw[-latex] (\xbasisvec, \ybasisvec, \zbasisvec-0.5) -- (evaz);
      \end{scope}%
    \end{tikzpicture}
  }
  \caption{Reduction of the domains through the steps \eqref{H2Dmod:item:lifting_step_1}, \eqref{H2Dmod:item:lifting_step_2}, and \eqref{H2Dmod:item:lifting_step_3}.\label{H2Dmod:fig:steps_resolution}}
\end{figure}

\vspace{1\baselineskip} \noindent
The advantage of Problem \eqref{H2Dmod:eq:formal_augmented_BVP} lies in the periodicity properties of $(\aten_p, \rho_p)$ along $\Sigma$, which are exploited in Section \ref{H2Dmod:sec:DtN_approach_augmented_problem} to reduce computations to a bounded domain. More precisely, the procedure is divided into three steps, each of which is devoted to bounding the domain in one direction, as illustrated in Figure \ref{H2Dmod:fig:steps_resolution}:
\begin{enumerate}[label=$(\textit{\alph*})$, ref=$\textit{\alph*}$]
  \item\label{H2Dmod:item:lifting_step_1} \textit{Periodic augmented jump data}. The condition \eqref{H2Dmod:eq:compatibility} offers some great latitude in choosing $G$. If $g$ is smooth enough, then one obvious pick would be a $G$ that is constant in the $\eva_2$--direction for instance, namely $G(0, \zvi_1, \zvi_2) := g(0, \zvi_1/\cuti_1)$ for $(\zvi_1, \zvi_2) \in \R^2$. Then, since $\aten_p$ and $\rho_p$ are $1$--periodic in the $\eva_2$--direction, one can expect $U$ to be also $1$--periodic in the $\eva_2$--direction, that is $U(\cdot + \vts \eva_2) = U$, so that \eqref{H2Dmod:eq:formal_augmented_BVP} reduces to a problem defined in the strip $\R^2 \times (0, 1)$. This property holds if one considers, for generality, extensions $G$ that are formally speaking $1$--periodic in the $\eva_2$--direction, that is $G(\cdot + \vts \eva_2) = G$. The corresponding strip problem is studied in Section \ref{H2Dmod:sec:analysis_augmented_problem}. 
  \item\label{H2Dmod:item:lifting_step_2} \textit{Partial Floquet-Bloch transform}. In Section \ref{H2Dmod:sec:FB_transform}, we introduce the partial Floquet-Bloch transform with respect to $\zvi_1$. Thanks to the periodicity of $(\aten_p, \rho_p)$ in the $\eva_1$--direction, this transform leads to a family of waveguide problems defined in the cylinder $\R \times (0, 1)^2$.
  \item\label{H2Dmod:item:lifting_step_3} \textit{DtN approach}. At last, using the periodicity of $(\aten^\pm_p, \rho^\pm_p)$ in the $\eva_\xvi$--direction, we resort in Section \ref{H2Dmod:sec:DtN_waveguide} to the DtN method developed in \cite{fliss2009exact, joly2006exact}. This method allows to reduce computations to the interface $\{0\} \times (0, 1)^2$ thanks to local Dirichlet problems defined in the cells $(0, \pm 1) \times (0, 1)^2$, and a propagation operator which satisfies a Riccati equation.
\end{enumerate}

\subsection{Outline of the paper}
The paper is structured as follows. In Section \ref{H2Dmod:sec:model_problem_and_augmented_structure}, two specific configurations of $2$D media represented by $\aten$ and $\rho$ are introduced. We prove \eqref{H2Dmod:eq:intro_general_quasiperiodic_expression}, namely that these $2$D media can be viewed as the restriction of $3$D augmented structures (represented by $\aten_p$ and $\rho_p$) along the hyperplane oriented along a particular matrix $\cutmat \in \R^{3 \times 2}$. Accordingly, as explained in Section \ref{H2Dmod:sec:formal_description_method}, we extend \eqref{H2Dmod:eq:transmission_problem} onto a $3$D augmented problem which is studied and solved in Section \ref{H2Dmod:sec:DtN_approach_augmented_problem} using the Floquet-Bloch transform and the DtN approach developed in \cite{joly2006exact}. But beforehand, we have to dedicate Section \ref{H2Dmod:sec:functional_framework} to setting up an appropriate functional framework for the analysis of the augmented 3D problem. In particular, we introduce anisotropic Sobolev spaces for which we establish trace theorems and Green's formulas. In Section \ref{H2Dmod:sec:algo_results}, we summarize the algorithm and present its discretization. A particular attention is given to the discretization of local Dirichlet cell problems involved in the procedure. Even though these problems can be solved directly using $3$D finite elements, we exploit their anisotropic structure to instead solve a family of $2$D cell problems. This so-called quasi $2$--dimensional (or quasi-$2$D) method is inspired from the quasi-$1$D developed in \cite{amenoagbadji2023wave}, although less obvious, as we demonstrate. Section \ref{H2Dmod:sec:numerical_results} provides numerical results to illustrate the efficiency of the method in different situations. Finally, the present work suggests some extensions and perspectives, which we highlight in Section \ref{H2Dmod:sec:extensions_perspectives}.

\subsection*{Notation}
\begin{itemize}
  \item We use the notation $\mathcal{O}_\idom \subset \R^n$, $\idom \in \{\varnothing, +, -\}$ to refer to any triple of the form $\{\mathcal{O}, \mathcal{O}_+, \mathcal{O}_-\}$, with the convention that $\mathcal{O}_\idom = \mathcal{O}$ for $\idom = \varnothing$.
  \item The indicator function of an open set $\mathcal{O} \subset \R^n$ is denoted by $\mathds{1}_{\mathcal{O}}$.
  \item Given $\pv_1, \pv_2 \in \R^2$, we denote by $\Z \pv_1 + \Z \pv_2$ the lattice $\Z \pv_1 + \Z \pv_2 := \{n_1\, \pv_1 + n_2\, \pv_2\ /\  n_1, n_2 \in \Z\}$. Similarly, for $\pva_1, \pva_2, \pva_3 \in \R^3$, let $\Z \pva_1 + \Z \pva_2 + \Z \pva_3 := \{n_1\, \pva_1 + n_2\, \pva_2 + n_3\, \pva_3 \ /\ n_1, n_2, n_3 \in \Z\}$.
  \item $\transp{\mathbb{A}} \in \C^{m \times n}$ denotes the transpose of $\mathbb{A} \in \C^{n \times m}$. The $n \times n$ identity matrix is denoted by $\mathbb{I}_n$.
  \item For $\mathcal{O} \subset \R^n$, the scalar product on $L^2(\mathcal{O})$ is denoted by $(\cdot, \cdot)_{L^2(\mathcal{O})}$. We denote by $\langle \cdot, \cdot \rangle_{\partial \mathcal{O}}$ the duality product between $H^{-1/2}(\partial \mathcal{O})$ and $H^{1/2}(\partial \mathcal{O})$.
  \item Given two Banach spaces $\mathscr{X}_1$ and $\mathscr{X}_2$, the space of bounded linear operators from $\mathscr{X}_1$ to $\mathscr{X}_2$ is denoted by $\mathscr{L}(\mathscr{X}_1, \mathscr{X}_2)$. We set $\mathscr{L}(\mathscr{X}) := \mathscr{L}(\mathscr{X}, \mathscr{X})$.
\end{itemize}

\subsection*{Acknowledgments}
\surligner{%
This research is supported in part by Simons Foundation Math + X Investigator Award \#376319 (Michael I. Weinstein). The authors would like to thank the Isaac Newton Institute for Mathematical Sciences, Cambridge, for support and hospitality during the programme \emph{Mathematical theory and applications of multiple wave scattering} where work on this paper was undertaken. This work was supported by EPSRC grant no EP/R014604/1.
}%
}{%
}

\ifthenelse{\boolean{brouillon}}{%
  \section{The model configurations and their hidden quasiperiodic nature}\label{H2Dmod:sec:model_problem_and_augmented_structure}
In this section, we present the properties of the functions $\aten \in L^\infty(\R^2; \R^{2 \times 2}), \rho \in L^\infty(\R^2)$ involved in \eqref{H2Dmod:eq:transmission_problem}, and we highlight their quasiperiodic nature as restrictions of $3$D functions along a hyperplane.

\subsection{Two specific configurations}\label{H2Dmod:sec:specific_configs}
We assume that $\aten$ and $\rho$ can be written as:
\begin{equation}
  \label{H2Dmod:eq:restriction_A_rho}
  \displaystyle
  \aeforall \xv = (\xvi, \zvi) \in \R^2, \quad %
  \aten(\xv) := \left\{\begin{array}{l@{\quad \textnormal{if}\ \ }l}
    \aten^-(\xv) & \xvi < 0,
    \\[6pt]
    \aten^+(\xv) & \xvi > 0,
  \end{array}
  \right. \quad \textnormal{and} \quad
  \rho(\xv) := \left\{\begin{array}{l@{\quad \textnormal{if}\ \ }l}
    \rho^-(\xv) & \xvi < 0,
    \\[6pt]
    \rho^+(\xv) & \xvi > 0,
  \end{array}
  \right.
\end{equation}
where $\aten^\pm \in \mathscr{C}^0(\R^2; \R^{2 \times 2})$, $\rho^\pm \in \mathscr{C}^0(\R^2)$. We consider two specific classes of functions $(\aten^\pm, \rho^\pm)$ (\surligner{see Section \ref{H2Dmod:sec:extension_general_config} for the extension to more general cases}).
\begin{enumerate}[ref = $\mathscr{\Alph*}$, wide = 0pt]
  \item\label{H2Dmod:item:config_a} \underline{Configuration \eqref{H2Dmod:item:config_a} --- \textit{The media are periodic along the interface $\sigma = \{\xvi = 0\}$}}
  \\[\topsep]
  This setting consists in $\aten^\pm$ and $\rho^\pm$ being $\Z \ev_\xvi + \Z (p^\pm_\zvi \ev_\zvi)$--periodic for some numbers $p^\pm_\zvi > 0$, that is,
  \begin{equation}
    \spforall \xv \in \R^2, \quad
    \left\{%
    \begin{array}{r@{\ =\ }l@{\quad \textnormal{and} \quad}r@{\ =\ }l}
      \aten^\pm(\xv + \ev_\xvi) & \aten^\pm(\xv), & \aten^\pm(\xv + p^\pm_\zvi \ev_\zvi) & \aten^\pm(\xv),
      \\[8pt]
      \rho^\pm(\xv + \ev_\xvi) & \rho^\pm(\xv), & \rho^\pm(\xv + p^\pm_\zvi \ev_\zvi) & \rho^\pm(\xv).
    \end{array}
    \right.
    \label{H2Dmod:eq:config_a}
  \end{equation}
  This is illustrated in Figure \ref{H2Dmod:fig:configuration_a}.

  \begin{rmk}\label{H2Dmod:rmk:periodicity_parameters_ratio}~
    \begin{enumerate}[label=$(\textit{\alph*})$., ref = \theprop.\textit{\alph*}]
      \setlength\itemsep{1em}
      \item If the ratio $p^+_\zvi/p^-_\zvi$ is a rational number that can be written as $a / b$ for some coprime $(a, b) \in \Z \times \N^*$, then $(\aten^+, \rho^+)$ and $(\aten^-, \rho^-)$ share a common period $\tau := b\, p^+_\zvi = a\, p^-_\zvi$ in the $\ev_\zvi$--direction. It follows that the overall medium represented by $(\aten, \rho)$ is periodic in the $\ev_\zvi$--direction. Hence, as done in \cite{flisscassanbernier2010}, a Floquet-Bloch transform can be applied with respect to the variable $\zvi$ along the interface, reducing \eqref{H2Dmod:eq:transmission_problem} to a family of waveguide problems set in $\R \times (0, \tau)$, and parameterized by the Floquet variable. However, this method becomes more costly as the denominator of $p^+_\zvi/p^-_\zvi$, and thus the period $\tau$ increases. It could then be relevant to use, in this case, the method described in this paper.   
      \item When $p^+_\zvi/p^-_\zvi$ is irrational, the method in \cite{flisscassanbernier2010} cannot be applied directly because $\aten$ and $\rho$ are no longer periodic along the interface. In this case, one might be tempted to construct a rational approximation $(a_n/b_n)_n$ of $p^+_\zvi/p^-_\zvi$, and to compute, using \cite{flisscassanbernier2010}, the solution $u_n$ of \eqref{H2Dmod:eq:transmission_problem}, obtained by replacing $p^+_\zvi/p^-_\zvi$ by $a_n/b_n$ for $n$ large enough. However, in addition to the theoretical questions that such a strategy raises (regarding for instance the convergence of $(u_n)_n$ to $u$), there are some numerical drawbacks. In fact, for a sequence of rationals $(a_n/b_n)_n$ to converge to an irrational $p^+_\zvi/p^-_\zvi$, the sequence of denominators $(b_n)_n$ must tend to infinity. Therefore, as explained in the first point, computational costs would inevitably increase with the periods $(\tau_n)_n$ as $(a_n/b_n)_n$ tends to $p^+_\zvi/p^-_\zvi$. Furthermore, the approximation quality would be strongly related to the irrationality measure of $p^+_\zvi/p^-_\zvi$, which indicates how efficiently it can be approximated by rational numbers. More detail about this aspect on rational approximation can be found in \cite[Chapter XI]{hardy1979introduction}.
    \end{enumerate}  
  \end{rmk}

  \begin{figure}[ht!]
    \centering
    \makebox[\textwidth][c]{
      \includegraphics[page=2]{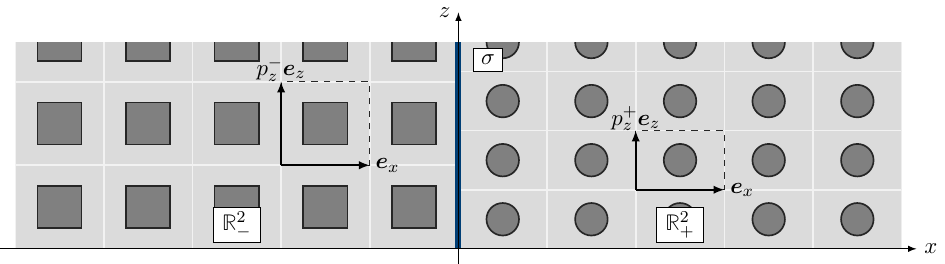}
    }
    \caption{Configuration \eqref{H2Dmod:item:config_a}: Juxtaposition of two media that are periodic along the interface\label{H2Dmod:fig:configuration_a}}
  \end{figure}

  \item\label{H2Dmod:item:config_b} \underline{Configuration \eqref{H2Dmod:item:config_b} --- Junction of a homogeneous medium and a periodic one}
  \\[\topsep]
  This corresponds to the case where $\aten^-$ and $\rho^-$ are constant while $\aten^+$ and $\rho^+$ are $\Z \ev_\xvi + \Z \pv^+$--periodic for some vector $\pv^+ = (p^+_\xvi, p^+_\zvi) \in \R^2$ such that $p^+_\zvi \neq 0$: for $\xv \in \R^2$,
  \begin{equation}
    \spforall \xv \in \R^2, \quad \left\{%
    \begin{array}{r@{\ \equiv\ }l@{\quad \textnormal{and} \quad}r@{\ =\ }l@{, \quad }r@{\ =\ }l}
      \aten^-(\xv) & \aten^- & \aten^+(\xv + \ev_\xvi) & \aten^+(\xv) & \aten^+(\xv + \pv^+) & \aten^+(\xv),
      \\[8pt]
      \rho^-(\xv) & \rho^- & \rho^+(\xv + \ev_\xvi) & \rho^+(\xv) & \rho^+(\xv + \pv^+) & \rho^+(\xv),
    \end{array}
    \right.
    \label{H2Dmod:eq:config_b}
  \end{equation}
  as illustrated in Figure \ref{H2Dmod:fig:configuration_b}.
  \begin{rmk}\label{H2Dmod:rmk:x_component_periodicity_vector}
    Note that $p^+_\zvi \ev_\zvi = \pv^+ - p^+_\xvi \ev_\xvi$. Therefore if $p^+_\xvi$ is rational with the irreducible form $a / b$, $(a, b) \in \Z \times \N^*$, then $b\, p^+_\zvi \ev_\zvi = b\vts \pv^+ - a\vts \ev_\xvi$ and thus $(\aten^+, \rho^+)$ are $\tau$--periodic in the $\ev_\zvi$--direction with $\tau := b\, p^+_\zvi$. Conversely, if $p^+_\xvi$ is irrational, then the medium is no longer periodic in the direction of the interface. In that regard, Remark \ref{H2Dmod:rmk:periodicity_parameters_ratio} translates to Configuration \eqref{H2Dmod:item:config_b} when $p^+_\zvi/p^-_\zvi$ is replaced by $p^+_\xvi$.
  \end{rmk}

  \begin{figure}[H]
    \makebox[\textwidth][c]{
      \includegraphics[page=3]{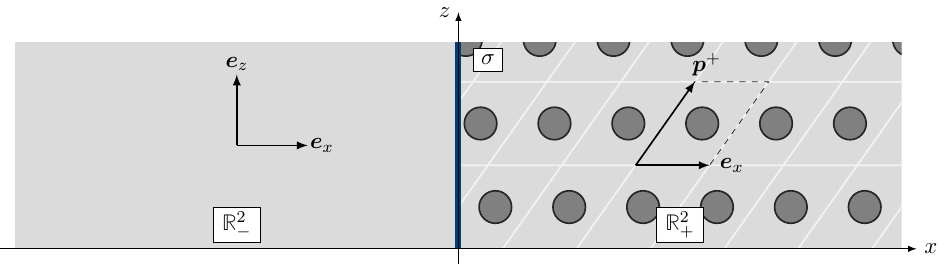}
    }
    \caption{Configuration \eqref{H2Dmod:item:config_b}: Juxtaposition of a periodic medium and a homogeneous one\label{H2Dmod:fig:configuration_b}}
  \end{figure}
\end{enumerate}

\subsection{A hidden quasiperiodicity along the interface}
\label{H2Dmod:sec:hidden_quasiperiodicity}
The goal of this section is to show that the functions $\aten$ and $\rho$ given by Configurations \eqref{H2Dmod:item:config_a} and \eqref{H2Dmod:item:config_b} can be viewed as quasiperiodic along the interface $\sigma$, in the sense that they are restrictions to a hyperplane of $3$D functions that are periodic along a $2$D interface containing $\sigma$. More precisely, we shall prove \eqref{H2Dmod:eq:intro_general_quasiperiodic_expression}, which is recalled below:
\begin{subequations}
  \label{H2Dmod:eq:general_quasiperiodic_form}
  \begin{equation}
    \label{H2Dmod:eq:general_quasiperiodic_expression}
    \aeforall \xv \in \R^2, \quad \aten(\xv) = \aten_p(\cutmat\, \xv) \quad \textnormal{and} \quad \rho(\xv) = \rho_p(\cutmat\, \xv),
  \end{equation}
  where $\aten_p$ and $\rho_p$ are defined by
  \begin{equation}
    \displaystyle
    \spforall \xva = (\xvi, \zvi_1, \zvi_2) \in \R^3, \quad  \aten_p(\xva) := %
    \left\{\begin{array}{l@{,\quad}l}
      \aten^+_p(\xva) & x > 0
      \\[6pt]
      \aten^-_p(\xva) & x < 0
    \end{array}
    \right. \quad \textnormal{and} \quad %
    \rho_p(\xva) := %
    \left\{\begin{array}{l@{,\quad}l}
      \rho^+_p(\xva) & x > 0
      \\[6pt]
      \rho^-_p(\xva) & x < 0,
    \end{array}
    \right.
    \label{H2Dmod:eq:def_Ap_rho_p}
  \end{equation}
  with $\aten^\pm_p, \rho^\pm_p \in \mathscr{C}^0(\R^3)$ which are $\Z^3$--periodic. Furthermore, the matrix $\cutmat \in  \R^{3 \times 2}$, which will be referred to as the \emph{cut matrix}, has the following form: 
  \begin{equation}
    \label{H2Dmod:eq:def_cut_matrix}
    \displaystyle
    \cutmat = %
    \begin{pmatrix}
     1 & 0
     \\
     0 & \cuti_1
     \\
     0 & \cuti_2
   \end{pmatrix} \quad \textnormal{with} \quad \cuti_1 \neq 0,
  \end{equation}
  the exact expression of $\cuti_1$ and $\cuti_2$ depending on the configuration.
\end{subequations}

\subsubsection{Quasiperiodicity of Configuration \texorpdfstring{\eqref{H2Dmod:item:config_a}}{(A)}} \label{H2Dmod:sec:augmented_structure_config_a}%
To find $(\aten_p, \rho_p)$, the formal idea is to split the tangential variable $\zvi$ into two variables $\zvi_1$ and $\zvi_2$, where $\zvi_1$ is associated to the periodicity of the medium in $\R^2_+$, and where $\zvi_2$ is associated to the periodicity of the medium in $\R^2_-$. In addition, to ensure that the medium is $1$--periodic with respect to $\zvi_1$ and $\zvi_2$, this splitting comes with a rescaling: $\zvi_1$ will correspond to $\zvi / p^+_\zvi$, while $\zvi_2$ will correspond to $\zvi / p^-_\zvi$. This suggests to introduce 
\begin{equation}
  \displaystyle
  \label{H2Dmod:eq:Apm_rhopm_config_a}
  \spforall \xva = (\xvi, \zvi_1, \zvi_2) \in \R^2, \quad \left\{%
  \begin{array}{r@{\ :=\ }l@{\qquad \textnormal{and} \qquad}r@{\ :=\ }l}
    \aten^+_p(\xva) & \aten^+(\xvi, \vts p^+_\zvi\, \zvi_1) & \rho^+_p(\xva) & \rho^+(\xvi, \vts p^+_\zvi\, \zvi_1)
    \\[10pt]
    \aten^-_p(\xva) & \aten^-(\xvi, \vts p^-_\zvi\, \zvi_2) & \rho^-_p(\xva) & \rho^-(\xvi, \vts p^-_\zvi\, \zvi_2),
  \end{array}
  \right.
\end{equation}
where $(\aten^+_p, \rho^+_p)$ are $1$--periodic in $\xvi$, $1$--periodic in $\zvi_1$, and independent of $\zvi_2$, while $(\aten^-_p, \rho^-_p)$ are $1$--periodic in $\xvi$, independent of $\zvi_1$, and $1$--periodic in $\zvi_2$, as illustrated in Figure \ref{H2Dmod:fig:augmented_structure_config_a}.

\vspace{1\baselineskip} \noindent
Now let $(\aten_p, \rho_p)$ be given by \eqref{H2Dmod:eq:def_Ap_rho_p} using $(\aten^\pm_p, \rho^\pm_p)$. Then, since both $(\aten^+_p, \rho^+_p)$ and $(\aten^-_p, \rho^-_p)$ are $1$--periodic with respect to $\zvi_1$ and $\zvi_2$, it follows that $(\aten_p, \rho_p)$ are periodic along the interface $\Sigma$.

\begin{figure}[ht!]
  \noindent\makebox[\textwidth]{%
    \centering
    \includegraphics[page=4]{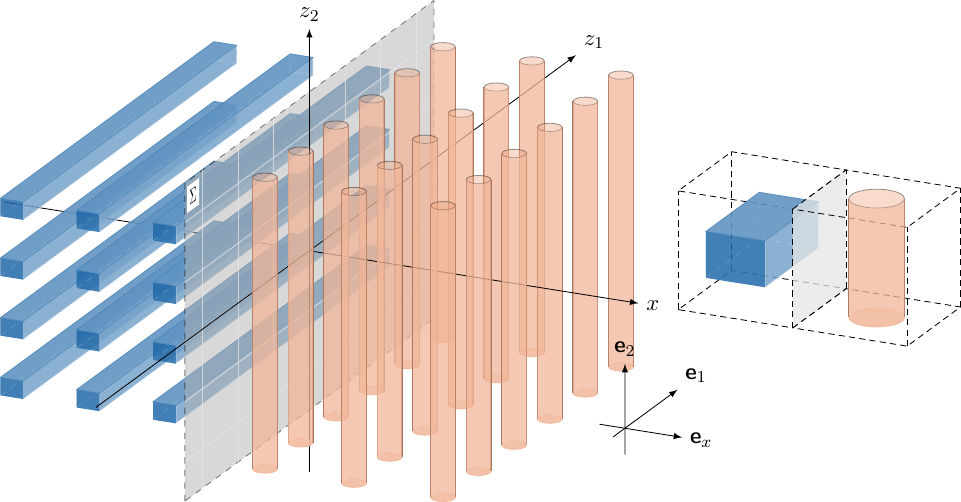}
  }
  \caption{\surligner{Left: Schematic representation of the augmented structure for the Configuration \eqref{H2Dmod:item:config_a}--type medium represented in Figure \ref{H2Dmod:fig:configuration_a}. Right: Close-up view of a periodicity cell.}\label{H2Dmod:fig:augmented_structure_config_a}}
\end{figure}

\vspace{1\baselineskip} \noindent
Moreover, by \textquote{inverting} \eqref{H2Dmod:eq:Apm_rhopm_config_a}, we deduce that
\begin{equation}
  \displaystyle
  \spforall \xv = (\xvi, \zvi) \in \R^2, \quad \aten^\pm(\xv) = \aten^\pm_p (\xvi,\vts \zvi / p^+_\zvi,\vts \zvi / p^-_\zvi) \quad \textnormal{and} \quad \rho^\pm(\xv) = \rho^\pm_p (\xvi,\vts \zvi / p^+_\zvi,\vts \zvi / p^-_\zvi).
\end{equation}
This corresponds to \eqref{H2Dmod:eq:general_quasiperiodic_expression} and \eqref{H2Dmod:eq:def_cut_matrix} with
\begin{equation}
  \displaystyle
  \cuti_1 := 1/p^+_\zvi \quad \textnormal{and} \quad \cuti_2 := 1/p^-_\zvi.
\end{equation}%
Finally, note that the half-plane $\cutmat\, \R^2_\pm$ (\emph{resp.} the line $\cutmat\, \sigma$) is included in the half-space $\R^3_\pm$ (\emph{resp.} the interface $\Sigma$ defined by \eqref{H2Dmod:eq:def_Sigma}).

\subsubsection{Quasiperiodicity of Configuration \texorpdfstring{\eqref{H2Dmod:item:config_b}}{(B)}} \label{H2Dmod:sec:augmented_structure_config_b}%
From its properties \eqref{H2Dmod:eq:config_b}, we shall see that $\aten^+$ admits the expression \eqref{H2Dmod:eq:general_quasiperiodic_expression}, although it is less obvious than in Configuration \eqref{H2Dmod:item:config_a}. It is useful to \textquote{look at $\aten^+$ in the basis $\{\ev_\xvi, \pv^+\}$} by defining
\begin{equation}
  \spforall (\accentset{\circ}{\xvi}, \accentset{\circ}{\zvi}) \in \R^2, \quad \accentset{\circ}{\aten}^+(\accentset{\circ}{\xvi}, \accentset{\circ}{\zvi}) := \aten^+ (\accentset{\circ}{\xvi}\, \ev_\xvi + \accentset{\circ}{\zvi}\, \pv^+) := \aten^+(\accentset{\circ}{\xvi} + \accentset{\circ}{\zvi}\, p^+_\xvi,\, \accentset{\circ}{\zvi}\, p^+_\zvi).
  \label{H2Dmod:eq:def_A0}
\end{equation}
From the periodicity properties \eqref{H2Dmod:eq:config_b} of $\aten^+$ it follows that $\accentset{\circ}{\aten}^+$ is 1--periodic with respect to its variables. Moreover, considering the change of variables
\begin{equation*}
  \displaystyle
  \left.
  \begin{array}{r@{\ =\ }l}
    \accentset{\circ}{\xvi} + p^+_\xvi \accentset{\circ}{\zvi} & \xvi 
    \\[5pt]
    p^+_\zvi \accentset{\circ}{\zvi} & \zvi 
  \end{array}
  \right\}
  \quad \Longleftrightarrow \quad
  \left\{
    \begin{array}{r@{\ =\ }l}
      \accentset{\circ}{\xvi} & \xvi - (p^+_\xvi / p^+_\zvi)\, \zvi
      \\[5pt]
      \accentset{\circ}{\zvi} & \zvi / p^+_\zvi,
    \end{array}
  \right.
\end{equation*}
which is well-defined because $p^+_\zvi \neq 0$, \eqref{H2Dmod:eq:def_A0} can be \textquote{inverted} as
\begin{equation}
  \spforall \xv = (\xvi, \zvi) \in \R^2, \quad \aten^+(\xv) = \accentset{\circ}{\aten}^+(\xvi - (p^+_\xvi / p^+_\zvi)\, \zvi, \zvi / p^+_\zvi).
  \label{H2Dmod:eq:A_from_A0}
\end{equation}
Thus, similarly to Configuration \eqref{H2Dmod:item:config_a}, we split the tangential variable $\zvi$ into a variable $\zvi_1$ corresponding to $\zvi / p^+_\zvi$, and a variable $\zvi_2$ corresponding to $- (p^+_\xvi / p^+_\zvi)\, \zvi$. Accordingly, it is natural to define the function
\begin{align}\label{H2Dmod:eq:Ap_from_A0}
  \spforall \xva = (\xvi, \zvi_1, \zvi_2) \in \R^3, \quad \aten^+_p(\xva) := \accentset{\circ}{\aten}^+ (x + \zvi_2, \zvi_1), 
\end{align}
where $\aten^+_p$ is $\Z^3$--periodic, and in particular periodic with respect to the variables $(\zvi_1, \zvi_2)$ of the interface $\Sigma$ defined by \eqref{H2Dmod:eq:def_Sigma}, as illustrated in Figure \ref{H2Dmod:fig:augmented_structure_config_b}. Then, \eqref{H2Dmod:eq:A_from_A0} becomes
\begin{equation}
  \label{H2Dmod:eq:A_from_Ap_config_b}
  \displaystyle
  \spforall \xv = (\xvi, \zvi) \in \R^2, \quad \aten^+(\xv) = \aten^+_p(\xvi, \zvi / p^+_\zvi, - (p^+_\xvi / p^+_\zvi)\, \zvi).
\end{equation}
The same arguments can be applied to $\rho$, to define $\rho^+_p$ similarly to \eqref{H2Dmod:eq:Ap_from_A0}. In addition, $\aten^-$ and $\rho^-$ are extended as constant functions over $\R^3$. More precisely, we set
\begin{equation}
  \label{H2Dmod:eq:def_Am_rhom_config_b}
  \spforall \xva = (\xvi, \zvi_1, \zvi_2) \in \R^3, \quad \displaystyle
  \aten^-_p(\xva) := \aten^- \quad \textnormal{and} \quad \rho^-_p(\xva) := \rho^-,
\end{equation}
and we define $\aten_p$ and $\rho_p$ as in \eqref{H2Dmod:eq:def_Ap_rho_p}, so that \eqref{H2Dmod:eq:general_quasiperiodic_expression} and \eqref{H2Dmod:eq:def_cut_matrix} hold with 
\begin{equation}
  \displaystyle
  \cuti_1 := 1/p^+_\zvi \neq 0 \quad \textnormal{and} \quad \cuti_2 := - p^+_\xvi / p^+_\zvi.
\end{equation}
Finally, note that the half-plane $\cutmat\, \R^2_\pm$ (\emph{resp.} the line $\cutmat\, \sigma$) is included in the half-space $\R^3_\pm$ (\emph{resp.} the interface $\Sigma$).

\begin{figure}[ht!]
  \noindent\makebox[\textwidth]{%
    \centering
    \includegraphics[page=5]{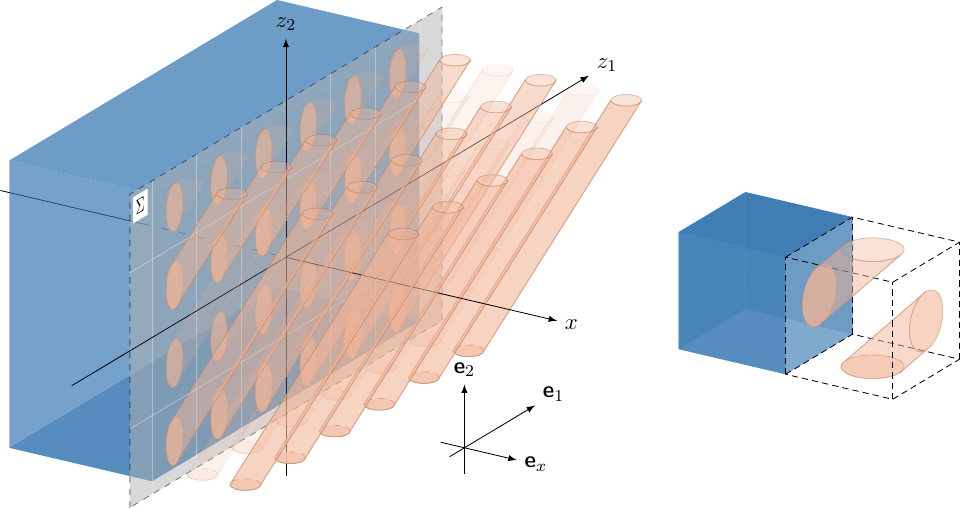}
  }
  \caption{\surligner{Left: Schematic representation of the augmented structure for the Configuration \eqref{H2Dmod:item:config_b}--type medium represented in Figure \ref{H2Dmod:fig:configuration_b} (some cylinders have been made more transparent for readability). Right: Close-up view of a periodicity cell.}\label{H2Dmod:fig:augmented_structure_config_b}}
\end{figure}

}{%
}

\ifthenelse{\boolean{brouillon}}{%
  \section{Functional framework}\label{H2Dmod:sec:functional_framework}
The formal procedure described in Section \ref{H2Dmod:sec:formal_description_method} consists in seeking the solution $u$ of \eqref{H2Dmod:eq:transmission_problem} under the form $u(\xv) = U(\cutmat \xv)$, where $U$ satisfies
\begin{equation}\label{H2Dmod:eq:formal_augmented_PDE_recalled}
  \left\{
    \begin{array}{r@{\ =\ }l}
      \displaystyle- \transp{\bsnabla} \cutmat\vts \aten_p\!  \transp{\cutmat} \bsnabla\, U - \rho_p\, \omega^2\, U & 0 \quad \textnormal{in}\ \  \R^3_+ \cup \R^3_-,
      \ret
    \llbracket (\cutmat\vts \aten_p \! \transp{\cutmat} \bsnabla\, U) \cdot \eva_\xvi \rrbracket_{\Sigma} & G.
    \end{array}
  \right.
\end{equation}
The advantage of the augmented equation \eqref{H2Dmod:eq:formal_augmented_PDE_recalled} lies in its periodic nature, which allows to use adapted tools for periodic PDEs, such as the Floquet-Bloch transform. Nevertheless, a difficulty is that the differential operator $-\transp{\bsnabla} \cutmat\vts \aten_p\!\! \transp{\cutmat} \bsnabla$ is \emph{elliptically degenerate} because the matrix $\cutmat$ given by \eqref{H2Dmod:eq:def_cut_matrix} is of rank $2$. More precisely,
\[
  \displaystyle
  \aeforall \xv \in \R^2, \quad \Ker \big[\cutmat\vts \aten_p (\xv) \transp{\cutmat}\big] = \vect \begin{pmatrix} 0 \\ -\cuti_2 \\ \phantom{-}\cuti_1 \end{pmatrix}.
\]
As a consequence, the properties of \eqref{H2Dmod:eq:formal_augmented_PDE_recalled} differ substantially from those of the classical Helmholtz equation given by $-\transp{\bsnabla}\vts \aten_p\vts \bsnabla U - \rho\, \omega^2\, U = 0$. In particular, one needs an appropriate functional framework to take the anisotropic nature of \eqref{H2Dmod:eq:formal_augmented_PDE_recalled} into account. This is the object of this section, which sets up the framework that will be used afterwards. 

\subsection{Anisotropic Sobolev spaces: motivation and outline}
To begin, let us define for any open set $\mathcal{O} \subset \R^3$ the anisotropic Sobolev spaces
\begin{equation}
  \displaystyle
  \label{H2Dmod:eq:def_H1C_HdivC}
  \begin{array}{r@{\ :=\ }l}
    H^1_\cutmat (\mathcal{O}) & \{U \in L^2(\mathcal{O})\ /\ \transp{\cutmat} \bsnabla U \in [L^2(\mathcal{O})]^2 \},
    \\[15pt]
    H_\cutmat(\mathbf{div}; \mathcal{O}) & \{\boldsymbol{W} \in [L^2(\mathcal{O})]^2\ /\ \transp{\bsnabla} \cutmat \boldsymbol{W} \in L^2(\mathcal{O})\}.
  \end{array}
\end{equation}
These are Hilbert spaces when equipped with the respective scalar products
\[
  \begin{array}{c@{\quad\ }r@{\ :=\ }l}
    \spforall U, V \in H^1_\cutmat (\mathcal{O}), & \displaystyle (U, V)_{H^1_\cutmat (\mathcal{O})} & \displaystyle \int_{\mathcal{O}} \big[\vts U\, \overline{V} + \big(\transp{\cutmat} \bsnabla U\big) \cdot \big(\transp{\cutmat} \bsnabla \overline{V}\big)\vts\big],
    \\[15pt]
    \displaystyle \spforall \boldsymbol{W}, \widetilde{\boldsymbol{W}} \in H_\cutmat(\mathbf{div}; \mathcal{O}), & \displaystyle (\boldsymbol{W}, \widetilde{\boldsymbol{W}})_{H_\cutmat(\mathbf{div}; \mathcal{O})} & \displaystyle \int_{\mathcal{O}} \big[\vts \boldsymbol{W} \cdot \overline{\widetilde{\boldsymbol{W}}} + (\transp{\bsnabla} \cutmat \boldsymbol{W})\; (\transp{\bsnabla} \cutmat \overline{\widetilde{\boldsymbol{W}}}) \vts\big].
  \end{array}
\]
We shall denote by $\|\cdot\|_{H^1_\cutmat (\mathcal{O})}$ and $\|\cdot\|_{H_\cutmat(\mathbf{div}; \mathcal{O})}$ the respective induced norms. 

\vspace{1\baselineskip} \noindent
Due to the domains introduced throughout the paper, a specific attention will be given to "rectangle" based cylindrical domains in what follows. Let $I_\xvi, I_1$ be intervals which do not need to be bounded, and such that $0 \in \overline{I_\xvi}$. We consider the $3$--dimensional domain $\genO$ and the $2$--dimensional transverse set $\genS^\tau$ given by
\begin{equation}
  \label{H2Dmod:eq:cylindrical_domain_infinite_along_e2}
  \displaystyle
  \begin{array}{r@{\ :=\ \big\{(\xvi, \zvi_1, \zvi_2) \in \R^3\ /\ }c@{\ \ \zvi_1 \in I_1,\ \ \zvi_2 \in \R}l@{\quad \equiv\ }l}
    \genO & \xvi \in I_\xvi, & \big\} & I_\xvi \times I_1 \times \R,
    \\[10pt]
    \spforall \tau \in \overline{I_\xvi}, \qquad \genS^\tau & \xvi = \tau, & \big\} & \{\tau\} \times I_1 \times \R,
  \end{array}
\end{equation}
with $\genS^\tau =: \genS$ for $\tau = 0$. The above domains are unbounded in the $\eva_2$--direction. We will be interested in the "$\zvi_2$--cell domains", bounded in $\zvi_2$:
\begin{equation}
  \label{H2Dmod:eq:cylindrical_domain}
  \displaystyle
  \begin{array}{r@{\ :=\ \big\{(\xvi, \zvi_1, \zvi_2)\ }c@{\ /\  \zvi_2 \in (0, 1)}l@{\quad \equiv\ }l}
    \genO_\periodic &  \in \genO & \big\} & I_\xvi \times I_1 \times (0, 1),
    \\[10pt]
    \spforall \tau \in \overline{I_\xvi}, \qquad \genS^\tau_\periodic &  \in \genS^\tau & \big\} & \{\tau\} \times I_1 \times (0, 1),
  \end{array}
\end{equation}
with $\genS^\tau_\periodic =: \genS_\periodic$ for $\tau = 0$. \surligner{The subscript "$\periodic$" refers to the fact that $\genO_\periodic$ and $\genS^\tau_\periodic$ are bounded in the $\eva_2$--direction. However, note that $\genO_\periodic$ and $\genS^\tau_\periodic$ may be bounded or unbounded in the $\eva_1$--direction, depending on whether $I_1$ is bounded or not. The same applies to $\genO_\periodic$ in the $\eva_\xvi$--direction.} Figure \ref{H2Dmod:fig:cylindrical_domain} (left) represents these domains for $I_\xvi = \R_-$ and $I_1 \subset \R_+$.

\begin{figure}[H]
  \def\BBxmin{-4.5}
  \def\BBxmax{+1.5}
  \def\BBymin{-0.25}
  \def\BBymax{+5.25}
  \def\BBzmin{-0.25}
  \def\BBzmax{+2.75}
  \def\coupevecY{1.414213562}%
  \def\coupevecZ{0.75}
  \def\planymin{0}
  \def\planymax{4}
  \begin{tikzpicture}
    \tdplotsetmaincoords{75}{40}
    \pgfmathsetmacro\xcoupemax{min(\BBymax/\coupevecY, \BBzmax/\coupevecZ)}
    \pgfmathsetmacro\xcoupemin{max(\BBymin/\coupevecY, \BBzmin/\coupevecZ)}
    \begin{scope}[tdplot_main_coords, scale=0.9, shift={(-4.15cm, 1.15cm)}]%
      \node (Omega) at (\BBxmin-1.25, {\planymin+0.125*(\planymax+\planymin)}, 1.5) {$\genO_\periodic$};
      \draw[-latex] (Omega) to[bend right=40] (\BBxmin+1, {\planymin+0.125*(\planymax+\planymin)}, 0.5);
      \node[tertiary!60!black] (SigmaTau) at (0.6*\BBxmin, {\planymax+0.5*(\planymax+\planymin)}, 2) {$\genS^\tau_\periodic$};

      \draw[-latex, tertiary!60!black] (SigmaTau) to[bend left=40] (0.6*\BBxmin, \planymax, 0.75);
      %
      \fill[secondary!80, opacity=0.5] (\BBxmin, \planymax, 0) -- +(-\BBxmin, 0, 0) -- +(-\BBxmin, 0, 1) -- +(0, 0, 1) -- cycle;
      \draw[white, dashed, thick] (\BBxmin, \planymax, 0) -- +(-\BBxmin, 0, 0);
      \fill[secondary!80, opacity=0.5] (\BBxmin, \planymin, 0) -- +(-\BBxmin, 0, 0) -- +(-\BBxmin, \planymax-\planymin, 0) -- +(0, \planymax-\planymin, 0) -- cycle;
      \filldraw[thick, draw=white, fill=tertiary!60, opacity=1] (0.6*\BBxmin, \planymin, 0) -- +(0, \planymax-\planymin, 0) -- +(0, \planymax-\planymin, 1) -- +(0, 0, 1) -- cycle;
      \fill[secondary!80, opacity=0.7] (\BBxmin, \planymin, 1) -- +(-\BBxmin, 0, 0) -- +(-\BBxmin, \planymax-\planymin, 0) -- +(0, \planymax-\planymin, 0) -- cycle;
      \filldraw[thick, draw=white, fill=tertiary!60, opacity=0.6] (0, \planymin, 0) -- +(0, \planymax-\planymin, 0) -- +(0, \planymax-\planymin, 1) -- +(0, 0, 1) -- cycle;
      \filldraw[thick, draw=white, fill=secondary!90, opacity=0.5] (\BBxmin, \planymin, 0) -- +(-\BBxmin, 0, 0) -- +(-\BBxmin, 0, 1) -- +(0, 0, 1) -- cycle;
      %
      \node (xaxis) at ({\BBxmax+0.25*(\BBxmax-\BBxmin)}, 0, 0) {$\xvi$};%
      \draw[-latex] ({\BBxmin-0.1*(\BBxmax-\BBxmin)}, 0, 0) -- (xaxis);%
      \node (yaxis) at (0, {\BBymax+0.25*(\BBymax-\BBymin)}, 0) {$\zvi_1$};%
      \draw[-latex] (0, {\BBymin-0.1*(\BBymax-\BBymin)}, 0) -- (yaxis);%
      \draw[-latex] (0, 0, {\BBzmin-0.25*(\BBzmax-\BBzmin)}) -- (0, 0, {\BBzmax+0.25*(\BBzmax-\BBzmin)}) node[above] {$\smash{\zvi_2}$};
      \draw (0.6*\BBxmin, \planymin, 0) node[below] {$\tau$};
      %
      \node[tertiary!60!black] (Sigma) at (\BBxmax+1.25, {\planymin+0.5*(\planymax+\planymin)}, 2.5) {$\genS_\periodic$};
      \draw[-latex, tertiary!60!black] (Sigma) to[bend left=40] (0, {\planymin+0.5*(\planymax+\planymin)}, 0.5);
      \begin{scope}[canvas is yz plane at x=0.1]
        \draw[decorate,decoration={brace, mirror}] (\planymin, 0) -- (\planymax, 0);
      \end{scope}
      \draw (1, 0.5*\planymin+0.5*\planymax, 0) node {$I_1$};
    \end{scope}%
  \end{tikzpicture}
  \hfill
  \begin{tikzpicture}
    \def\lezero{3.25}
    \fill[thick, fill=secondary!70, opacity=0.75] (\BBxmin, {\planymin/\coupevecY}) -- (\BBxmin, {\planymax/\coupevecY}) -- (0, {\planymax/\coupevecY}) -- (0, {\planymin/\coupevecY});
    \draw (0.5*\BBxmin, {0.5*(\planymin + \planymax)/\coupevecY}) node {$\Qcuti$};
    %
    %
    \draw[-latex] (\BBxmin-0.5, 0) -- (1, 0) node[right] {$\xvi$};
    \draw[-latex] (0, {\BBzmin-0.125*(\BBzmax-\BBzmin)}) -- (0, {0.115+\BBzmax+0.25*(\BBzmax-\BBzmin)}) node[above] {$\zvi$};
    \draw [tertiary, line width = 0.75mm] (0, {\planymin/\coupevecY}) -- (0, {\planymax/\coupevecY});
    \node[tertiary] (Stau) at (1.5, {0.5*(\planymin + \planymax)/\coupevecY}) {$\Scuti$};
    \draw[-latex, tertiary] (Stau) -- (0, {0.5*(\planymin + \planymax)/\coupevecY});
  \end{tikzpicture} 
  \caption{Left: the rectangle-based cylindrical domain $\genO_\periodic$ and the lateral set $\genS_\periodic$ given by \eqref{H2Dmod:eq:cylindrical_domain}. Right: the domains $\Qcuti$ and $\Scuti$ given by \eqref{H2Dmod:eq:2D_1D_cylindrical_domains}. $I_\xvi = \R_-$ and $I_1 \subset \R_+$. \label{H2Dmod:fig:cylindrical_domain}}
\end{figure}

\noindent
For smooth functions $V \in \mathscr{C}^1(\genO_\periodic)$, $\boldsymbol{W} \in [\mathscr{C}^1(\genO_\periodic)]^2$, one has the Green's formula
\begin{equation*}
  \displaystyle
  \int_{\genO_\periodic} \big[ (\divecut \boldsymbol{W})\, \overline{V} + \boldsymbol{W} \cdot \overline{\gradcut V} \big] = \int_{\partial\vts \genO_\periodic} (\cutmat\, \boldsymbol{W} \cdot \nva)\, \overline{V},
\end{equation*}
where $\nva$ denotes the unit normal vector on $\partial\vts \genO_\periodic$ which is outward with respect to $\genO_\periodic$. \surligner{In what follows, we shall be more} interested in two cases:
\begin{enumerate}[label=$(\textnormal{C}.\arabic*)$, ref=$\textnormal{C}.\arabic*$]
  \item\label{H2Dmod:item:case_strip} \emph{The strip}: $I_1 = \R$ and $V, \boldsymbol{W}$ are $1$--periodic with respect to $\zvi_2$;
  \item\label{H2Dmod:item:case_cylinder} \emph{The cylinder}: $I_1 = (0, 1)$ and $V, \boldsymbol{W}$ are $1$--periodic with respect to $\zvi_1$ and $\zvi_2$.
\end{enumerate}
In both cases, the above Green's formula becomes
\begin{equation}\label{H2Dmod:eq:Green_formula_smooth}
  \displaystyle
  \int_{\genO_\periodic} \big[ (\divecut \boldsymbol{W})\, \overline{V} + \boldsymbol{W} \cdot \overline{\gradcut V} \big] = \sum_{\tau \in \partial I_\xvi} \int_{\genS^\tau_\periodic} (\cutmat\, \boldsymbol{W} \cdot \nva)\, \overline{V},
\end{equation}
where the right-hand side is reduced to the surface integral on the transverse faces $\genS^\tau_\periodic$, $\tau \in \partial I_\xvi$. Our goal is to extend the above formula to $(V, \boldsymbol{W}) \in H^1_\cutmat (\genO_\periodic) \times H_\cutmat(\mathbf{div}; \genO_\periodic)$ (that is, to prove \eqref{H2Dmod:eq:Green_formula_strips} and \eqref{H2Dmod:eq:Green_cylinder} for the cases \eqref{H2Dmod:item:case_strip} and \eqref{H2Dmod:item:case_cylinder} defined above respectively). This requires to define properly the traces $V|_{\genS^\tau_\periodic}$ (Section \ref{H2Dmod:sec:gamma0}) and $(\cutmat\, \boldsymbol{W} \cdot \nva)|_{\genS^\tau_\periodic}$ (Section \ref{H2Dmod:sec:gamma1_Greens_formula_strip}), and to replace the integral in the right-hand side of \eqref{H2Dmod:eq:Green_formula_smooth} by a duality product between the appropriate trace spaces (Section \ref{H2Dmod:sec:gamma1_Greens_formula_strip} for Case \eqref{H2Dmod:item:case_strip} and Section \ref{H2Dmod:sec:e1_e2_periodic_functions} for Case \eqref{H2Dmod:item:case_cylinder}).

\vspace{1\baselineskip} \noindent
As a first step, it is useful to introduce a so-called shear map (see \eqref{H2Dmod:eq:def_shear}), to charaterize any function in $H^1_\cutmat (\genO_\periodic)$ (\emph{resp.} in $H_\cutmat (\mathbf{div}; \genO_\periodic)$) which is periodic with respect to $\zvi_2$ in terms of its $L^2$--regularity in the $\eva_2$--direction and its $H^1$--regularity in the $\eva_\xvi$ and $\cuti_1\vts \eva_1 + \cuti_2\vts \eva_2$--directions. This is the object of the next section.

\subsection{Anisotropic spaces of \texorpdfstring{$\Z\eva_2$}{Ze2}--periodic functions} \label{H2Dmod:sec:e2_periodic_functions}
Given numbers $-\infty \leq a < b \leq +\infty$ and a Banach space $(\mathscr{X}, \|\cdot\|_{\mathscr{X}})$, we recall that
\begin{equation}\label{H2Dmod:eq:def_L2per_Xvalued}
  L^2(a, b; \mathscr{X}) := \displaystyle \Big\{\surligner{\boldsymbol{V}}\ /\ \  \| \surligner{\boldsymbol{V}} \|^2_{L^2(a, b; \mathscr{X})} := \int_a^b \|\surligner{\boldsymbol{V}}(\cdot, s)\|^2_{\mathscr{X}} \; ds < +\infty \Big\}.
\end{equation}
Our starting point is the oblique transformation
\begin{equation}
  \label{H2Dmod:eq:oblique_mapping}
  \displaystyle
  \operatorname{T}_\cuti: (\xv, s) \in \R^3\ \longmapsto\ \cutmat\, \xv + s\, \eva_2 = (\xvi, \cuti_1\vts \zvi, \cuti_2\vts \zvi + s) \in \R^3,
\end{equation}
which relates the $3$D operator $\gradcut$ to the $2$D gradient $\nabla$ by means of the chain rule \eqref{H2Dmod:eq:chain_rule} recalled as $\nabla [F \circ \operatorname{T}_\cuti (\cdot, s)] (\xv) = (\gradcut F) \circ \operatorname{T}_\cuti (\xv, s)$, for any $\mathscr{C}^1$ function $F$. Consider the $2$D domain $\Qcuti$ and the edge $\Scuti^\tau$ given by
\begin{equation}
  \label{H2Dmod:eq:2D_1D_cylindrical_domains}
  \displaystyle
  \begin{array}{r@{\ :=\ \big\{(\xvi, \zvi) \in \R^2\ /\ }c@{\ \ \cuti_1\vts \zvi \in I_1}l}
    \Qcuti & \xvi \in I_\xvi, & \big\},
    \\[10pt]
    \Scuti^\tau & \xvi = \tau, & \big\},
  \end{array}
\end{equation}
with $\Scuti^\tau =: \Scuti$ for $\tau = 0$ (see Figure \ref{H2Dmod:fig:cylindrical_domain} (right) for $I_\xvi = \R_-$ and $I_1 \subset \R_+$). We note that $\Qcuti$ and $\Scuti^\tau$ are well-defined, since for Configurations \eqref{H2Dmod:item:config_a} and \eqref{H2Dmod:item:config_b}, we have $\cuti_1 \neq 0$. Moreover, $\operatorname{T}_\cuti$ is bijective from $\Qcuti \times \R$ (\emph{resp.} $\Scuti^\tau \times \R$) to $\genO$ (\emph{resp.} $\genS^\tau$). The image of $\Qcuti \times (0, 1)$ (\emph{resp.} $\Scuti^\tau \times (0, 1)$) by $\operatorname{T}_\cuti$ is the oblique domain $\genO_{\periodic, \cuti}$ (\emph{resp.} the oblique face $\genS^\tau_{\periodic, \cuti}$) represented in Figure \ref{H2Dmod:fig:cylindrical_domain_b}:
\begin{equation}
  \label{H2Dmod:eq:oblique_cylindrical_domain}
  \displaystyle
  \genO_{\periodic, \cuti} := \operatorname{T}_\cuti \big(\Qcuti \times (0, 1) \big) %
  \qquad \textnormal{and} \qquad
  \genS^\tau_{\periodic, \cuti} := \operatorname{T}_\cuti \big(\Scuti^\tau \times (0, 1) \big), %
\end{equation}
with $\genS^\tau_{\periodic, \cuti} =: \genS_{\periodic, \cuti}$ for $\tau = 0$. 

\vspace{1\baselineskip} \noindent
As a preliminary step, one can use the change of variables $(\xv, s) \mapsto \operatorname{T}_\cuti(\xv, s)$ and Fubini's theorem, to prove the charactization
\begin{equation}\label{H2Dmod:eq:fibered_characterization_H1cut_Hdivecut}
  \begin{array}{r@{\ =\ }l}
      H^1_\cutmat (\genO_{\periodic, \cuti}) & \Big\{V \in L^2(\genO_{\periodic, \cuti})\ /\ V \circ \operatorname{T}_\theta \in L^2\big(\R; H^1(\Qcuti)\big) \Big\},
      \\[15pt]
      H_\cutmat(\mathbf{div}; \genO_{\periodic, \cuti}) & \Big\{\boldsymbol{W} \in [L^2(\genO_{\periodic, \cuti})]^2 \ /\ \boldsymbol{W} \circ \operatorname{T}_\theta \in L^2\big(\R; H(\dive; \Qcuti)\big) \Big\},
    \end{array}
\end{equation}
which is useful to derive the properties of $H^1_\cutmat (\genO_{\periodic, \cuti})$ and $H_\cutmat(\mathbf{div}; \genO_{\periodic, \cuti})$ from the well-studied spaces $L^2(\R; H^1(\Qcuti))$, $L^2(\R; H(\dive; \Qcuti))$. We would like to provide similar expressions for the spaces of $\Z \eva_2$--periodic $H^1_\cutmat$ or $H_\cutmat (\mathbf{div})$ functions defined on the straight domain $\genO_\periodic$ (instead of $\genO_{\periodic, \cuti}$). To this end, we shall, roughly speaking, identify functions on $\genO_\periodic$ and $\genO_{\periodic, \cuti}$, using the notion of \emph{periodic extension in the $\eva_2$--direction}. In the sequel, we fix $d \in \{1, 2\}$, and we use the notation $\boldsymbol{V}$ (in bold) to refer to a $\C^d$--valued function.

\begin{figure}[ht!]
  \def\BBxmin{-4.5}
  \def\BBxmax{+1.5}
  \def\BBymin{-0.25}
  \def\BBymax{+4.75}
  \def\BBzmin{-0.25}
  \def\BBzmax{+2.75}
  \def\coupevecY{1.414213562}%
  \def\coupevecZ{0.75}
  \def\planymin{0}
  \def\planymax{4}
  \pgfmathsetmacro\planymaxYcut{4*sqrt(\coupevecY^2 + \coupevecZ^2)/\coupevecY}
  \begin{subfigure}{0.4\textwidth}
    \caption{\label{H2Dmod:fig:cylindrical_domain_a}}
    \begin{tikzpicture}
      \tdplotsetmaincoords{65}{40}
      \pgfmathsetmacro\xcoupemax{min(\BBymax/\coupevecY, \BBzmax/\coupevecZ)}
      \pgfmathsetmacro\xcoupemin{max(\BBymin/\coupevecY, \BBzmin/\coupevecZ)}
      \begin{scope}[tdplot_main_coords, scale=0.8]%
        %
        %
        \fill[secondary!80, opacity=0.5, pattern=north west lines, pattern color=tertiary!60!black] (\BBxmin, \planymin, 0) -- +(-\BBxmin, 0, 0) -- +(-\BBxmin, \planymaxYcut-\planymin, 0) -- +(0, \planymaxYcut-\planymin, 0) -- cycle;
        \draw[dashed, gray!80] (\BBxmin, \planymaxYcut, 0) -- +(-\BBxmin, 0, 0) -- +(-\BBxmin, 0, 1) -- +(0, 0, 1) -- cycle;
        \fill[dashed, gray!20, opacity=0.25] (\BBxmin, \planymaxYcut, 0) -- +(-\BBxmin, 0, 0) -- +(-\BBxmin, 0, 1) -- +(0, 0, 1) -- cycle;
        \draw[gray!80, dashed] (\BBxmin, \planymin, 1) -- +(-\BBxmin, 0, 0) -- +(-\BBxmin, \planymaxYcut-\planymin, 0) -- +(0, \planymaxYcut-\planymin, 0) -- cycle;
        \fill[tertiary!20, dashed, opacity=0.25] (\BBxmin, \planymin, 1) -- +(-\BBxmin, 0, 0) -- +(-\BBxmin, \planymaxYcut-\planymin, 0) -- +(0, \planymaxYcut-\planymin, 0) -- cycle;
        \draw[dashed, gray!80] (\BBxmin, \planymin, 0) -- +(-\BBxmin, 0, 0) -- +(-\BBxmin, 0, 1) -- +(0, 0, 1) -- cycle;
        \fill[dashed, tertiary!20, opacity=0.25] (\BBxmin, \planymin, 0) -- +(-\BBxmin, 0, 0) -- +(-\BBxmin, 0, 1) -- +(0, 0, 1) -- cycle;
        \fill[tertiary!20, opacity=0.25, canvas is yz plane at x=0] (0, \planymin) rectangle (\planymaxYcut, 1);
        %
        %
        %
        \node (xaxis) at ({\BBxmax+0.25*(\BBxmax-\BBxmin)}, 0, 0) {$\xvi$};%
        \draw[-latex] ({\BBxmin-0.1*(\BBxmax-\BBxmin)}, 0, 0) -- (xaxis);%
        \node (yaxis) at (0, {\BBymax+0.25*(\BBymax-\BBymin)}, 0) {$\zvi$};%
        \draw[-latex] (0, {\BBymin-0.1*(\BBymax-\BBymin)}, 0) -- (yaxis);%
        \draw[-latex] (0, 0, {\BBzmin-0.25*(\BBzmax-\BBzmin)}) -- (0, 0, {\BBzmax+0.5*(\BBzmax-\BBzmin)}) node[above] {$s$};
        \draw[thick, tertiary, line width = 0.75mm] (0, 0, 0) -- (0, \planymaxYcut, 0);
        \begin{scope}[canvas is yz plane at x=0.1]
          \draw[decorate,decoration={brace, mirror}] (\planymin, 0) -- (\planymaxYcut, 0);
        \end{scope}
        \draw (1, 0.5*\planymin+0.5*\planymaxYcut, 0) node {$\Scuti$};
        \node [align=center] at ({\BBxmin+0.65}, {\planymaxYcut-0.75}, 0) {$\Qcuti$};
      \end{scope}%
    \end{tikzpicture}
  \end{subfigure}
  \hfill
  \begin{tikzpicture}
    \node at (0, 0.6) {};
    \draw[-latex] (-1, 3) -- (1, 3);
    \draw (0, 3) node[above] {$\operatorname{T}_\cuti$};
  \end{tikzpicture}
  \hfill
  \begin{subfigure}{0.4\textwidth}
    \caption{\label{H2Dmod:fig:cylindrical_domain_b}}
    \begin{tikzpicture}
      \tdplotsetmaincoords{75}{40}
      \pgfmathsetmacro\xcoupemax{min(\BBymax/\coupevecY, \BBzmax/\coupevecZ)}
      \pgfmathsetmacro\xcoupemin{max(\BBymin/\coupevecY, \BBzmin/\coupevecZ)}
      \begin{scope}[tdplot_main_coords, scale=0.8, shift={(-4.15cm, 1.15cm)}]%
        \node[tertiary!60!black] (Omega) at (\BBxmin-2, {0.75*\planymax}, {1+0.75*\planymax*\coupevecZ/\coupevecY}) {$\genO_{\periodic, \cuti}$};
        \draw[tertiary!60!black, -latex] (Omega) to[bend right=40] (\BBxmin+0.25, {0.75*\planymax}, {0.75*\planymax*\coupevecZ/\coupevecY});
        %
        %
        %
        \fill[secondary!80, opacity=0.5] (\BBxmin, \planymax, 0) -- +(-\BBxmin, 0, 0) -- +(-\BBxmin, 0, 1) -- +(0, 0, 1) -- cycle;
        \draw[white, dashed, thick] (\BBxmin, \planymax, 0) -- +(-\BBxmin, 0, 0);
        \fill[secondary!80, opacity=0.5] (\BBxmin, \planymin, 0) -- +(-\BBxmin, 0, 0) -- +(-\BBxmin, \planymax-\planymin, 0) -- +(0, \planymax-\planymin, 0) -- cycle;
        \fill[tertiary!60, opacity=0.8, pattern=north west lines, pattern color=tertiary!50!black] (\BBxmin, \planymin, 0) -- +(0, {\coupevecY/\coupevecZ}, 1) -- +(-\BBxmin, {\coupevecY/\coupevecZ}, 1) -- +(-\BBxmin, 0, 0) -- cycle;
        \fill[secondary!80, opacity=0.5] (\BBxmin, \planymin, 1) -- +(-\BBxmin, 0, 0) -- +(-\BBxmin, \planymax-\planymin, 0) -- +(0, \planymax-\planymin, 0) -- cycle;
        \fill[tertiary!60, opacity=0.8, pattern=north west lines, pattern color=tertiary!50!black] (\BBxmin, {\planymin+(\coupevecY/\coupevecZ)}, 1) -- (\BBxmin, \planymax, {\planymax*\coupevecZ/\coupevecY}) -- (0, \planymax, {\planymax*\coupevecZ/\coupevecY}) -- (0, {\planymin+(\coupevecY/\coupevecZ)}, 1) -- cycle;
        \draw[white, dashed, thick] (\BBxmin, \planymax, {\planymax*\coupevecZ/\coupevecY}) -- +(-\BBxmin, 0, 0);
        \draw[thin, white, dashed, opacity=0.75] (\BBxmin, {\planymin+(\coupevecY/\coupevecZ)}, 1) -- +(-\BBxmin, 0, 0);
        \fill[tertiary!50, opacity=0.25] (\BBxmin, \planymax, {\planymax*\coupevecZ/\coupevecY}) -- +(-\BBxmin, 0, 0) -- +(-\BBxmin, 0, 1) -- +(0, 0, 1) -- cycle;
        \fill[tertiary!80, opacity=0.6] (\BBxmin, \planymin, 1) -- (\BBxmin, \planymax, {1 + \planymax*\coupevecZ/\coupevecY}) -- (0, \planymax, {1 + \planymax*\coupevecZ/\coupevecY}) -- (0, \planymin, 1) -- cycle;
        \filldraw[thick, draw=white, fill=secondary!80!black, opacity=0.5] (0, \planymin, 0) -- +(0, \planymax-\planymin, 0) -- +(0, \planymax-\planymin, 1) -- +(0, 0, 1) -- cycle;
        \filldraw[thick, draw=white, fill=secondary!90, opacity=0.5] (\BBxmin, \planymin, 0) -- +(-\BBxmin, 0, 0) -- +(-\BBxmin, 0, 1) -- +(0, 0, 1) -- cycle;
        \filldraw[thick, draw=white, fill=tertiary!90, opacity=0.6] (\BBxmin, \planymin, 0) -- +(-\BBxmin, 0, 0) -- +(-\BBxmin, 0, 1) -- +(0, 0, 1) -- cycle;
        \filldraw[thick, draw=white, fill=tertiary!90, opacity=0.6] (0, \planymin, 0) -- (0, \planymax, {\planymax*\coupevecZ/\coupevecY}) -- (0, \planymax, {1 + \planymax*\coupevecZ/\coupevecY}) -- (0, \planymin, 1) -- cycle;
        \draw[very thick, dashed, -latex, white] (0, \planymin+0.5, {0.5+(\planymin+0.5)*\coupevecZ/\coupevecY}) -- (0, \planymin+2, {0.5+(\planymin+2)*\coupevecZ/\coupevecY});
        \draw[very thick, dashed, -latex, white] (-0.5, \planymin, 0.5) -- +(-2, 0, 0);
        \draw[dashed] (0, \planymax-0.05, {\planymax*\coupevecZ/\coupevecY}) -- (0, \planymax-0.05, 1);
        %
        %
        \node (xaxis) at ({\BBxmax+0.25*(\BBxmax-\BBxmin)}, 0, 0) {$\xvi$};%
        \draw[-latex] ({\BBxmin-0.1*(\BBxmax-\BBxmin)}, 0, 0) -- (xaxis);%
        \node (yaxis) at (0, {\BBymax+0.25*(\BBymax-\BBymin)}, 0) {$\zvi_1$};%
        \draw[-latex] (0, {\BBymin-0.1*(\BBymax-\BBymin)}, 0) -- (yaxis);%
        \draw[-latex] (0, 0, {\BBzmin-0.25*(\BBzmax-\BBzmin)}) -- (0, 0, {\BBzmax+0.75*(\BBzmax-\BBzmin)}) node[above] {$\smash{\zvi_2}$};
        %
        \node[tertiary!60!black] (Sigma) at (\BBxmax, {0.9*(\planymin+\planymax)}, {2.5+0.9*\planymax*\coupevecZ/\coupevecY}) {$\genS_{\periodic, \cuti}$};
        \draw[-latex, tertiary!60!black] (Sigma) to[bend left=40] (0, {\planymin+0.9*(\planymax-\planymin)}, {0.5+0.9*\planymax*\coupevecZ/\coupevecY});
        \begin{scope}[canvas is yz plane at x=0.1]
          \draw[decorate,decoration={brace, mirror}] (\planymin, 0) -- (\planymax, 0);
        \end{scope}
        \draw (1, 0.5*\planymin+0.5*\planymax, 0) node {$I_1$};
      \end{scope}%
    \end{tikzpicture} 
  \end{subfigure}
  \begin{minipage}[b]{0.5\textwidth}
    \raggedright
    \caption*{%
      (a) The strip $\Qcuti \times (0, 1)$, where $\Qcuti$ is the hatched domain defined by \eqref{H2Dmod:eq:2D_1D_cylindrical_domains}. (b) The oblique domains $\genO_{\periodic, \cuti}$, $\genS_{\periodic, \cuti}$ defined by \eqref{H2Dmod:eq:oblique_cylindrical_domain}. The hatched domain is $\operatorname{T}_\cuti (\Qcuti \times \{0\})$.
      Any function in $H^1_\cutmat(\genO_{\periodic, \cuti})$ or $\Hgradper{\genO_\periodic}$ is $L^2$ in the $\eva_2$ direction, and $H^1$ in the $\eva_\xv$ and $\cuti_1\vts \eva_1 + \cuti_2\vts \eva_2$--directions indicated by the dashed arrows. (c) The point $\operatorname{T}_\cuti(\xv, s)$ defined by \eqref{H2Dmod:eq:oblique_mapping} for $(\xv, s) \in \Qcuti \times (0, 1)$ and for $|(\cuti_1, \cuti_2)| = 1$. We fix $I_\xvi = \R_-$ and $I_1 \subset \R_+$.\label{H2Dmod:fig:oblique_transformation}%
    }
  \end{minipage}
  \hfill
  \begin{minipage}[b]{0.45\textwidth}
    \raggedleft
    \begin{subfigure}[b]{1\textwidth}
      \vspace{2\baselineskip}
      \caption{\label{H2Dmod:fig:cylindrical_domain_c}}
      \begin{tikzpicture}[scale=1]
        \pgfmathsetmacro\xcoupemax{ 0.5+min(\BBymax/\coupevecY, \BBzmax/\coupevecZ)}
        \pgfmathsetmacro\xcoupemin{-0.5+max(\BBymin/\coupevecY, \BBzmin/\coupevecZ)}
        \def\spoint{0.675}
        \def\zpoint{2.}
        \fill[thick, fill=secondary!70, opacity=0.75] (\planymin, 0) -- (\planymin, 1) -- (\planymax, 1) -- (\planymax, 0);
        \fill[thick, fill=tertiary!70, opacity=0.75] (\planymin, 0) -- (\planymax, {\planymax*\coupevecZ/\coupevecY}) -- (\planymax, {1+(\planymax*\coupevecZ/\coupevecY)}) -- (\planymin, 1) -- cycle;
        \draw (\planymin, \spoint) node {$\bullet$} node[left] {$s$} -- ({\planymin + \zpoint*\coupevecY}, {\spoint + \zpoint*\coupevecZ}) node {$\bullet$};

        \draw[decorate,decoration={brace, mirror, amplitude=2mm}] (\planymin, \spoint) -- ({\planymin + \zpoint*\coupevecY}, {\spoint + \zpoint*\coupevecZ});
        \draw ({\planymin + 0.6*\zpoint*\coupevecY}, {\spoint - 0.45 + 0.6*\zpoint*\coupevecZ}) node {$\zvi$};
        \draw[dashed] ({\planymin + \zpoint*\coupevecY}, {\spoint + \zpoint*\coupevecZ}) -- ({\planymin + \zpoint*\coupevecY}, 0) node[below] {$\cuti_1\vts \zvi$};
        \draw[dashed] ({\planymin + \zpoint*\coupevecY}, {\spoint + \zpoint*\coupevecZ}) -- (0, {\spoint + \zpoint*\coupevecZ}) node[left] {$s + \cuti_2\vts \zvi$};
        %
        %
        \draw[-latex] (\BBymin-0.25, 0) -- (\BBymax+0.125, 0) node[right] {$\zvi_1$};
        \node (z2) at (0, {\BBzmax+0.25*(\BBzmax-\BBzmin)}) {$\zvi_2$};
        \draw[-latex] (0, {\BBzmin-0.25*(\BBzmax-\BBzmin)}) -- (z2);
        %
        \draw[latex-, secondary!80!black] (\planymax, 0.5) -- +(0.75, 0) node[right] {$\genS_\periodic$};
        \draw[latex-, tertiary!80!black] (\planymax, {0.5 + \planymax*\coupevecZ/\coupevecY}) -- +(0.75, 0) node[right] {$\genS_{\periodic, \cuti}$};
      \end{tikzpicture} 
    \end{subfigure}
  \end{minipage}
\end{figure}

\begin{defi}[Periodic extension]\label{H2Dmod:def:periodic_extension_e2}
  Let $p > 0$ and $\boldsymbol{V} \in [L^p (\genO_\periodic)]^d$. The \emph{periodic extension of $\boldsymbol{V}$ in the $\eva_2$--direction} is the function $\extper[2]{\boldsymbol{V}} \in [L^p_{\textit{loc}}(\genO)]^d$ defined by:
  \begin{equation}\label{H2Dmod:eq:periodic_extension_e2}
     \aeforall \xva = (\xvi, \zvi_1, \zvi_2) \in \genO_\periodic, \quad \spforall n \in \Z, \quad \extper[2]{\boldsymbol{V}} (\xva + n\vts \eva_2) := \boldsymbol{V}(\xva).
  \end{equation}
  For any $\Phi \in L^p (\genS_\periodic)$, we define $\extper[2]{\Phi} \in L^p_{\textit{loc}}(\genS)$ similarly, by replacing $\xva \in \genO_\periodic$ in \eqref{H2Dmod:eq:periodic_extension_e2} by $\xva \in \genS_\periodic$.
\end{defi}

\vspace{1\baselineskip} \noindent
Let $\mathscr{C}^\infty_{0, \periodic}(\overline{\genO_\periodic})$ denote the space of smooth functions in $\overline{\genO_\periodic}$ that are compactly supported in the $\eva_\xvi$ and $\eva_1$--directions, and $1$--periodic in the $\eva_2$--direction, that is,
\begin{equation}
  \label{H2Dmod:eq:def_Cinftyper}
  \mathscr{C}^\infty_{0, \periodic}(\overline{\genO_\periodic}) := \big\{V \in \mathscr{C}^\infty_0(\overline{\genO_\periodic}) \ /\ \extper[2]{V} \in \mathscr{C}^\infty(\genO)\big\}.
\end{equation}
Note that $\mathscr{C}^\infty_{0, \periodic}(\overline{\genO_\periodic})$ contains $\mathscr{C}^\infty_0(\genO_\periodic)$, and hence, is dense in $L^2(\genO_\periodic)$.

\vspace{1\baselineskip} \noindent
The change of variables \eqref{H2Dmod:eq:oblique_mapping} combined with the periodic extension along $\eva_2$ in Definition \ref{H2Dmod:def:periodic_extension_e2} allows us to introduce the so-called shear transform defined by
\begin{equation}
  \label{H2Dmod:eq:def_shear}
  \begin{array}{r@{\ } c@{\ }c@{\ }l}
    \shearmap: & [\mathscr{C}^\infty_{0, \periodic}(\overline{\genO_\periodic})]^d & \longrightarrow & [\mathscr{C}^\infty(\Qcuti \times \R)]^d
    \\[5pt]
    & \boldsymbol{V} & \longmapsto & \boldsymbol{V}_\cutmat,
  \end{array}
  \quad \boldsymbol{V}_\cutmat(\xv, s) := (\extper[2]{\boldsymbol{V}}) \circ \operatorname{T}_\cuti (\xv, s), \quad \spforall (\xv, s) \in \Qcuti \times \R.
\end{equation}
The definition of $\extper[2]{}$ implies that $\shearmap \boldsymbol{V} (\cdot, s + 1) = \shearmap \boldsymbol{V} (\cdot, s)$. For this reason, the study of $\shearmap$ will be restricted to $s \in (0, 1)$. As stated in the next proposition proved in Appendix \ref{H2Dmod:sec:proof_shear_properties}, $\shearmap$ extends by density to $L^2$--functions.
\begin{prop}\label{H2Dmod:prop:shear_properties}
  The mapping $\shearmap$ defined by \eqref{H2Dmod:eq:def_shear} extends by density to a mapping defined from $[L^2 (\genO_\periodic)]^d$ to $L^2(0, 1; [L^2(\Qcuti)]^d)$, with
  \begin{equation}
    \displaystyle
    \spforall \boldsymbol{U}, \boldsymbol{V} \in [L^2 (\genO_\periodic)]^d, \quad \frac{1}{\cuti_1}\, \int_0^1 \int_{\Qcuti} \shearmap \boldsymbol{U}(\xv, s)\cdot \overline{\shearmap \boldsymbol{V}(\xv, s)}\; d\xv ds = \int_{\genO_\periodic} \boldsymbol{U}\cdot \overline{\boldsymbol{V}}.
    \label{H2Dmod:eq:oblique_plancherel}
  \end{equation}
  Moreover, $\shearmap$ is an isomorphism from $[L^2 (\genO_\periodic)]^d$ to $L^2(0, 1; [L^2(\Qcuti)]^d)$, and its inverse is given for any $\boldsymbol{V}_\cutmat \in L^2(0, 1; [L^2(\Qcuti)]^d)$ by
  \begin{equation}
    \label{H2Dmod:eq:inverse_shear}
    \displaystyle
    \aeforall \xva = (\xvi, \zvi_1, \zvi_2) \in \genO_\periodic, \quad \invshearmap{\boldsymbol{V}_\cutmat} (\xva) = \sqrt{\cuti_1}\; \extper[s]{\boldsymbol{V}_\cutmat} \big(\xvi,\, \zvi_1 / \cuti_1,\, \zvi_2 - \zvi_1 (\cuti_2 /\cuti_1)\big),
  \end{equation}
  where $\extper[s]{\boldsymbol{V}_\cutmat} \in [L^2_{\textit{loc}}(\Qcuti \times \R)]^d$ denotes the periodic extension of $\boldsymbol{V}_\cutmat$ with respect to the variable $s$, defined for almost any $(\xv, s) \in \Qcuti \times (0, 1)$ and for any $n \in \Z$ by $\extper[s]{\boldsymbol{V}_\cutmat} (\xv, s + n) := \boldsymbol{V}_\cutmat (\xv, s)$.
\end{prop}

\vspace{1\baselineskip} %
\noindent
Now, by analogy with \eqref{H2Dmod:eq:fibered_characterization_H1cut_Hdivecut}, we define
\begin{equation}
  \displaystyle
  \label{H2Dmod:eq:def_H1Cper_HdivCper}
  \boxed{
    \begin{array}{r@{\ :=\ }l}
      \Hgradper{\genO_\periodic} & \big\{V \in L^2(\genO_\periodic)\ /\ \shearmap V \in L^2(0, 1; H^1(\Qcuti)) \big\},
      \\[15pt]
      \Hdiveper{\genO_\periodic} & \big\{\boldsymbol{W} \in [L^2(\genO_\periodic)]^2 \ /\ \shearmap \boldsymbol{W} \in L^2(0, 1; H(\dive; \Qcuti)) \big\}.
    \end{array}
  }
\end{equation}
It is not obvious from their definition that $\Hgradper{\genO_\periodic}$ and $\Hdiveper{\genO_\periodic}$ are respectively subspaces of $H^1_\cutmat(\genO_\periodic)$ and $H_\cutmat(\mathbf{div}; \genO_\periodic)$. This is however true from the chain rule, as highlighted in the next proposition proved in Appendix \ref{H2Dmod:sec:proof_link_H1Cper_H1C_and_HdivCper_HdivC}.
\begin{prop}\label{H2Dmod:prop:link_H1Cper_H1C_and_HdivCper_HdivC}
  Any $V \in \Hgradper{\genO_\periodic}$ belongs to $H^1_\cutmat(\genO_\periodic)$ and 
  \begin{equation}\label{H2Dmod:eq:shear_commutes_with_grad}
    \displaystyle
    \aeforall (\xv, s) \in \Qcuti \times (0, 1), \quad \shearmap (\transp{\cutmat} \bsnabla V)(\xv, s) = \nabla_\xv (\shearmap V)(\xv, s).
  \end{equation}
  Similarly, any $\boldsymbol{W} \in \Hdiveper{\genO_\periodic}$ belongs to $H_\cutmat(\mathbf{div}; \genO_\periodic)$ and 
  \begin{equation}\label{H2Dmod:eq:shear_commutes_with_dive}
    \displaystyle
    \aeforall (\xv, s) \in \Qcuti \times (0, 1), \quad \shearmap (\transp{\bsnabla} \cutmat \boldsymbol{W})(\xv, s) = \transp{\nabla}_\xv (\shearmap \boldsymbol{W})(\xv, s),
  \end{equation}
  where $\nabla_\xv$ is the gradient operator with respect to $\xv$.
\end{prop}

\vspace{1\baselineskip}\noindent
In what follows, we equip $\Hgradper{\genO_\periodic}$ with the scalar product $(\cdot, \cdot)_{H^1_{\cutmat}(\genO_\periodic)}$, and $\Hdiveper{\genO_\periodic}$ with the scalar product $(\cdot, \cdot)_{H_\cutmat(\mathbf{div}; \genO_\periodic)}$. Propositions \ref{H2Dmod:prop:shear_properties} and \ref{H2Dmod:prop:link_H1Cper_H1C_and_HdivCper_HdivC} then lead to the next result.
\begin{cor}\label{H2Dmod:cor:shear_properties}
  The shear map $\shearmap$ is an isomorphism from $\Hgradper{\genO_\periodic}$ to $L^2(0, 1; H^1(\Qcuti))$, and from $\Hdiveper{\genO_\periodic}$ to $L^2(0, 1; H(\dive; \Qcuti))$. Consequently, equipped with the $H^1_\cutmat$ (\emph{resp}. $H_\cutmat (\mathbf{div}))$) scalar product, $\Hgradper{\genO_\periodic}$ (\emph{resp}. $\Hdiveper{\genO_\periodic}$) is a Hilbert space.
\end{cor}

\vspace{1\baselineskip} \noindent
We finish this section with two density results which provide more insight on the nature of $\Hgradper{\genO_\periodic}$ and $\Hdiveper{\genO_\periodic}$. The proof is delayed to Appendix \ref{H2Dmod:sec:Cinftyper_dense_H1Cper}.
\begin{prop}\label{H2Dmod:prop:Cinftyper_dense_H1Cper}~
  \begin{enumerate}[label=$(\textit{\alph*})$., ref = \textit{\alph*}]
  \item The space $\mathscr{C}^\infty_{0, \periodic}(\overline{\genO_\periodic})$ is dense in $\Hgradper{\genO_\periodic}$.
  \item The space $[\mathscr{C}^\infty_{0, \periodic}(\overline{\genO_\periodic})]^2$ is dense in $\Hdiveper{\genO_\periodic}$.
  \end{enumerate}
\end{prop}

\begin{rmk}\label{H2Dmod:rmk:H1Cper_HdivCper_weak_periodicity}
  Proposition \ref{H2Dmod:prop:Cinftyper_dense_H1Cper} implies that for smooth functions, the definition of $\Hgradper{\genO_\periodic}$ corresponds to $1$--periodicity in the $\eva_2$--direction. More generally, one could show the characterization 
  \begin{equation*}
    \Hgradper{\genO_\periodic} = \{V \in H^1_\cutmat(\genO_\periodic)\ /\ V|_{\zvi_2 = 0} = V|_{\zvi_2 = 1}\},
  \end{equation*}
  even though the equality in this last characterization has to be understood not in the $L^2$ sense, but in the $L^2_{\textit{loc}}$ sense. In fact, by adapting the results in \cite[Section 3.2.1]{amenoagbadji2023wave}, it can be shown that the space of traces on $\{\zvi_2 = a\}$ involves a weighted $L^2$ space.
\end{rmk}

\subsection{Trace operator on transverse interfaces \texorpdfstring{$\genS^\tau_\periodic$}{Gamma}} \label{H2Dmod:sec:gamma0}
In order to prove Green's formulas, we need to define traces of functions in $\Hgradper{\genO_\periodic}$ on the transverse face $\genS^\tau_\periodic$ defined by \eqref{H2Dmod:eq:cylindrical_domain} for $\tau \in \overline{I_\xvi}$. Throughout the section, we assume for simplicity that $\tau = 0 \in \overline{I_\xvi}$, so that $\genS^\tau_\periodic = \genS_\periodic$. The extension to an arbitrary $\tau \in \overline{I_\xvi}$ is explained in Remark \ref{H2Dmod:rmk:gamma0_tau_neq_0}. 

\vspace{1\baselineskip} \noindent
Note that 
\begin{equation*}
  \spforall V \in \Hgradper{\genO_\periodic}, \quad V,\; \partial_\xvi V \in L^2(\genO_\periodic).
\end{equation*}
Therefore, the usual trace theorem applied to each $1$D function $\xvi \mapsto V(\xvi, \zvi_1, \zvi_2) \in H^1(I_\xvi)$ provides an estimate which can be integrated with respect to $(\zvi_1, \zvi_2) \in I_1 \times (0, 1)$, to obtain the following preliminary result.

\begin{prop}
  The trace application $\gamma_{0, \periodic}$, defined by $\gamma_{0, \periodic} V := V|_{\genS_\periodic}$ for $V \in \mathscr{C}^\infty_{0, \periodic}(\overline{\genO_\periodic})$, extends by continuity to a linear map still denoted by $\gamma_{0, \periodic}$, and defined from $\Hgradper{\genO_\periodic}$ to $L^2(\genS_\periodic)$.
\end{prop}

\noindent
In what follows, we will abusively write $V|_{\genS_\periodic}$ instead of $\gamma_{0, \periodic} V$ when referring to traces on $\genS_\periodic$.

\vspace{1\baselineskip} \noindent
We are next interested in characterizing the range of the trace map $\gamma_{0, \periodic}$. To this end, we first extend the definition of $\shearmap$ to functions defined on $\genS_\periodic$. For any $\Phi \in \mathscr{C}^\infty_0(\overline{\genS_\periodic})$ such that $\extper[2]{\Phi} \in \mathscr{C}^\infty(\genS)$, we can define $\shearmap \Phi = \Phi_\cutmat$ as in the volume case \eqref{H2Dmod:eq:def_shear} by choosing $(\xv, s) \in \Scuti \times \R$, where \surligner{$\Scuti = \Scuti^{\tau = 0}$} is the edge in \eqref{H2Dmod:eq:2D_1D_cylindrical_domains}. From
\[\displaystyle
  \spforall \xv = (0, \zvi) \in \Scuti, \quad \spforall s \in \R, \qquad \shearmap \Phi(\xv, s) := \extper[2]{\Phi}(0,\; \cuti_1\vts \zvi,\; \cuti_2\vts \zvi + s),
\]
we deduce that $\shearmap \Phi(\cdot, s)$ is simply the trace of $\extper[2]{\Phi}$ along the line $\{(\cuti_1\vts \zvi,\; \cuti_2\vts \zvi + s),\ \zvi \in \R\}$. Similarly to Proposition \ref{H2Dmod:prop:shear_properties}, $\shearmap$ extends by density as an isomorphism from $L^2(\genS_\periodic)$ to $L^2(0, 1; L^2(\Scuti))$. 

\vspace{1\baselineskip} \noindent
Consider the space
\begin{equation}\label{H2Dmod:eq:def_H12Cper}
  \boxed{\Honehalfper{\genS_\periodic} := \big\{\Phi \in L^2 (\genS_\periodic)\ /\ \shearmap \Phi \in L^2(0, 1; H^{1/2}(\Scuti)) \big\},}
\end{equation}
which is equipped with the norm $\Phi \mapsto \|\shearmap \Phi\|_{L^2(0, 1; H^{1/2}(\Scuti))}$. Then the next result holds.

\begin{figure}[ht!]
  \centering
  \begin{tikzpicture}
    \def\ecartHor{4}
    \def\ecartVer{-2}
    %
    \node (L2H1) at (0, 0) {$L^2(0, 1; H^1(\Qcuti))$};
    \node (H1cut) at (\ecartHor, 0) {$\Hgradper{\genO_\periodic}$};
    \node (L2H12) at (0, \ecartVer) {$L^2(0, 1; H^{1/2}(\Scuti))$};
    \node (H12cut) at (\ecartHor, \ecartVer) {$\Honehalfper{\genS_\periodic}$};
    %
    \draw[-latex] (H1cut) -- node[above, pos=0.5] {$\shearmap$} (L2H1);
    \draw[-latex] (H12cut) -- node[above, pos=0.5] {$\shearmap$} (L2H12);
    \draw[-latex] (H1cut) -- node[right, pos=0.5] {$\gamma_{0, \periodic}$} (H12cut);
    \draw[-latex] (L2H1) -- node[left, pos=0.5] {$\gamma_0$} (L2H12);
  \end{tikzpicture}
  \hfill
  \begin{tikzpicture}
    \def\ecartHor{4}
    \def\ecartVer{-2}
    %
    \node (L2H1) at (0, 0) {$L^2(0, 1; H(\dive; \Qcuti))$};
    \node (H1cut) at (\ecartHor, 0) {$\Hdiveper{\mathbb{\Omega}_\periodic}$};
    \node (L2Hminus12) at (0, \ecartVer) {$L^2(0, 1; H^{-1/2}(\Scuti))$};
    \node (Hminus12cut) at (\ecartHor, \ecartVer) {$\Hminusonehalfper{\mathbb{\Sigma}_\periodic}$};
    %
    \draw[-latex] (H1cut) -- node[above, pos=0.5] {$\shearmap$} (L2H1);
    \draw[-latex] (Hminus12cut) -- node[above, pos=0.5] {$\shearmap$} (L2Hminus12);
    \draw[-latex] (H1cut) -- node[right, pos=0.5] {$\gamma_{1, \periodic}$} (Hminus12cut);
    \draw[-latex] (L2H1) -- node[left, pos=0.5] {$\gamma_1$} (L2Hminus12);
  \end{tikzpicture}
  \caption{Left: In Proposition \ref{H2Dmod:prop:gamma0_continuous_surjective}, the range of $\gamma_{0, \periodic}$ is characterized using the shear map $\shearmap$ and the usual trace operator $\gamma_0$ defined by \eqref{H2Dmod:eq:usual_trace_operator}. More precisely, $\gamma_{0, \periodic} = \invshearmap\, \gamma_0 \shearmap$. Right: the same holds for the normal trace operator $\gamma_{1, \periodic}$ defined in Proposition \ref{H2Dmod:prop:Green_formula_strips} for $I_1 = \R$.\label{H2Dmod:fig:isomorphism_traces}}
\end{figure}

\begin{prop}\label{H2Dmod:prop:gamma0_continuous_surjective}
  The trace operator $\gamma_{0, \periodic}$ is continuous and surjective from $\Hgradper{\genO_\periodic}$ to $\Honehalfper{\genS_\periodic}$, and commutes with the shear map $\shearmap$:
  \begin{equation}\label{H2Dmod:eq:gamma0_commutes_with_S}
    \spforall V \in \Hgradper{\genO_\periodic}, \quad \spforall s \in (0, 1), \quad \shearmap (V|_{\genS_\periodic}) (\cdot, s) = \shearmap V (\cdot, s)|_{\Scuti}.
  \end{equation}
\end{prop}

\begin{dem}
  Equation \eqref{H2Dmod:eq:gamma0_commutes_with_S} is obtained easily for smooth functions $V \in \mathscr{C}^\infty_{0, \periodic}(\overline{\genO_\periodic})$ using the definition of $\shearmap$. It extends to $V \in \Hgradper{\genO_\periodic}$ using: $(i)$ the density of $\mathscr{C}^\infty_{0, \periodic}(\overline{\genO_\periodic})$ in $\Hgradper{\genO_\periodic}$, $(ii)$ the continuity of $\shearmap$ from $\Hgradper{\genO_\periodic}$ to $L^2(0, 1; H^1 (\Qcuti))$ (Corollary \ref{H2Dmod:cor:shear_properties}) and from $\Honehalfper{\genS_\periodic}$ to $L^2(0, 1; H^{1/2} (\Scuti))$ (by definition \eqref{H2Dmod:eq:def_H12Cper}), and $(iii)$ the continuity from $L^2(0, 1; H^1(\Qcuti))$ to $L^2(0, 1; H^{1/2}(\Scuti))$ of the trace map $\gamma_0$ given by: 
  \begin{equation}\label{H2Dmod:eq:usual_trace_operator}
    \displaystyle
    \spforall V_\cutmat \in L^2(0, 1; H^1(\Qcuti)), \quad \gamma_0 V_\cutmat (\cdot, s) := V_\cutmat (\cdot, s)|_{\Scuti},
  \end{equation}
  which is a direct consequence of the classical trace theorem in $H^1(\Qcuti)$. 
  
  \vspace{1\baselineskip} \noindent
  Since $\shearmap$ is an isomorphism, \eqref{H2Dmod:eq:gamma0_commutes_with_S} gives $\gamma_{0, \periodic} = \invshearmap\, \gamma_0 \shearmap$ \surligner{(see Figure \ref{H2Dmod:fig:isomorphism_traces} left)}. The continuity and surjectivity of $\gamma_{0, \periodic}$ then result from the continuity of $(\gamma_0, \shearmap, \invshearmap)$ and the surjectivity of $\gamma_0$.
\end{dem}

\begin{rmk}\label{H2Dmod:rmk:gamma0_tau_neq_0}
  Let $\tau \in \overline{I}_\xvi$. By applying the above arguments to $V(\cdot + \tau\, \eva_\xvi)$ where $V \in \Hgradper{\genO_\periodic}$, one can define the trace of $V$ on the face $\genS^\tau_\periodic$ for any $\tau \in \overline{I_\xvi}$.
\end{rmk}

\subsection{Normal trace operator and Green's formula for a strip}\label{H2Dmod:sec:gamma1_Greens_formula_strip}
In this section, we restrict ourselves to the case \eqref{H2Dmod:item:case_strip} of a strip\surligner{: more precisely, define} 
\[
  \displaystyle
  \mathbb{\Omega}_\periodic := I_\xvi \times \R \times (0, 1) \quad \textnormal{and} \quad \mathbb{\Sigma}^\tau_\periodic := \{\tau\} \times \R \times (0, 1).
\]
\surligner{Note that the strip $\mathbb{\Omega}_\periodic$ and the face $\mathbb{\Sigma}^\tau_\periodic$ correspond respectively to the domains $\genO_\periodic$ and $\genS^\tau_\periodic$ defined by \eqref{H2Dmod:eq:cylindrical_domain}, with
\[
  I_1 = \R.
\]
Moreover, the domain $Q$ and the edge $\Scuti^\tau$ defined by \eqref{H2Dmod:eq:2D_1D_cylindrical_domains} become
\begin{equation*}
  \Qcuti = I_\xvi \times \R \quad \textnormal{and} \quad \Scuti^\tau = \{\tau\} \times \R.
\end{equation*}
}
In what follows, we assume that $I_\xvi \neq \R$, so that $\partial I_\xvi \neq \varnothing$; typically $I_\xvi = \R_+$, $\R_-$, or $(a, b)$. We fix $\tau \in \partial I_\xvi$. Let $\nva = \pm \eva_\xvi$ (\emph{resp.} $\nv = \pm \ev_\xvi$) denote the unit normal vector on $\mathbb{\Sigma}^\tau_\periodic$ (\emph{resp.} $\Scuti^\tau$) which is outward with respect to $\mathbb{\Omega}_\periodic$ (\emph{resp.} $\Qcuti$). Our objective is to define a normal trace operator on $\mathbb{\Sigma}^\tau_\periodic$ and to prove the Green's formula given in Proposition \ref{H2Dmod:prop:Green_formula_strips}.

\vspace{1\baselineskip} \noindent
The topological dual space of $\Honehalfper{\mathbb{\Sigma}^\tau_\periodic}$ is denoted by $\Hminusonehalfper{\mathbb{\Sigma}^\tau_\periodic}$, and is equipped with the dual norm for now. Let $\langle \cdot,\, \cdot \rangle_{\mathbb{\Sigma}^\tau_\periodic}$ denote the duality product between $H^{-1/2}_{\cutmat\vts \periodic}(\mathbb{\Sigma}^\tau_\periodic)$ and $H^{1/2}_{\cutmat\vts \periodic}(\mathbb{\Sigma}^\tau_\periodic)$, defined as a natural extension of the $\vphantom{H^{-1/2}_{\cutmat\vts \periodic}(\mathbb{\Sigma}^\tau_\periodic)}L^2(\mathbb{\Sigma}^\tau_\periodic)$--scalar product. In order to define a normal trace operator, we first extend the shear map $\shearmap$ as an isomorphism from $H^{-1/2}_{\cutmat\vts \periodic}(\mathbb{\Sigma}^\tau_\periodic)$ to $[L^2(0, 1; H^{1/2}(\Scuti^\tau))]'$. This is achieved through duality, using the Parseval-like formula \eqref{H2Dmod:eq:oblique_plancherel}:
\begin{multline}
  \label{H2Dmod:eq:shearmap_duality}
  \displaystyle
  \spforall \Psi \in \Hminusonehalfper{\mathbb{\Sigma}^\tau_\periodic}, \qquad \frac{1}{\cuti_1} \big\langle \shearmap \Psi, \Phi_\cutmat \big\rangle_{[L^2(0, 1; H^{1/2}(\Scuti^\tau))]', L^2(0, 1; H^{1/2}(\Scuti^\tau))} 
  \\
  \quad := \big\langle \Psi, \invshearmap \Phi_\cutmat \big\rangle_{\mathbb{\Sigma}^\tau_\periodic}, \qquad \spforall \Phi_\cutmat \in L^2(0, 1; H^{1/2}(\Scuti^\tau)).
\end{multline}
The next lemma provides a more convenient characterization of $L^2(0, 1; H^{1/2}(\Scuti^\tau))'$.
\begin{lem}\label{H2Dmod:lem:dual_of_Bochner_space}
  Given $a < b \in \R$ and a Hilbert space $(\mathscr{X}, (\cdot, \cdot)_{\mathscr{X}})$, one has $[L^2(a, b; \mathscr{X})]' = L^2(a, b; \mathscr{X}')$ as well as the following identity: for any $\varPhi \in L^2(a, b; \mathscr{X})$ and $\varPsi \in L^2(a, b; \mathscr{X}')$:
  \begin{equation}\label{H2Dmod:eq:dual_of_Bochner_space}
    \displaystyle
    \big\langle \varPsi, \varPhi \big\rangle_{L^2(a, b; \mathscr{X})', L^2(a, b; \mathscr{X})} = \int_a^b \big\langle \varPsi(\cdot, s),\; \varPhi(\cdot, s) \big\rangle_{\mathscr{X}', \mathscr{X}} \; ds,
  \end{equation}
  where $\langle \cdot, \cdot \rangle_{\mathscr{X}', \mathscr{X}}$ denotes the duality product between $\mathscr{X}'$ and $\mathscr{X}$.
\end{lem}

\begin{dem}
  \surligner{This result holds for a more general class of Banach spaces (see e.g. \cite{diestel_vector_1977}), but we provide a simpler proof for Hilbert spaces, for ease of the reader.} One checks without any difficulty that the inclusion $[L^2(a, b; \mathscr{X})]' \supset L^2(a, b; \mathscr{X}')$ holds. Conversely, the inclusion $[L^2(a, b; \mathscr{X})]' \subset L^2(a, b; \mathscr{X}')$ and \eqref{H2Dmod:eq:dual_of_Bochner_space} result from the Riesz representation theorem applied to $\mathscr{X}$ and to the space $L^2(a, b; \mathscr{X})$ equipped with the scalar product 
  \begin{equation*}
    \displaystyle 
    (\varPhi, \varPsi) \mapsto \int_a^b (\varPhi(\cdot, s), \varPsi(\cdot, s))_{\mathscr{X}} \; ds.
  \end{equation*}

  \vspace{-1.5\baselineskip} \noindent
\end{dem}

\vspace{1\baselineskip} \noindent
By applying Lemma \ref{H2Dmod:lem:dual_of_Bochner_space} with $(a, b) = (0, 1)$ and $\mathscr{X} := H^{1/2}(\Scuti^\tau)$, we deduce that $\shearmap$ is an isomorphism from $H^{-1/2}_{\cutmat\vts \periodic}(\mathbb{\Sigma}^\tau_\periodic)$ to $L^2(0, 1; H^{-1/2}(\Scuti^\tau))$. Moreover, the bijectivity of $\shearmap$ implies that:
\begin{equation}\label{H2Dmod:eq:carac_Hminus12Cper}
  \boxed{
  \Hminusonehalfper{\mathbb{\Sigma}^\tau_\periodic} = \big\{\Psi\ /\ \shearmap \Psi \in L^2(0, 1; H^{-1/2}(\Scuti^\tau)) \big\}.}
\end{equation}
Finally, \eqref{H2Dmod:eq:shearmap_duality} combined with \eqref{H2Dmod:eq:dual_of_Bochner_space} leads to the following, which is analogous to \eqref{H2Dmod:eq:oblique_plancherel}:
\begin{equation}\label{H2Dmod:eq:shearmap_fibered_duality}
  \displaystyle
  \spforall (\Phi, \Psi) \in \Honehalfper{\mathbb{\Sigma}^\tau_\periodic} \times \Hminusonehalfper{\mathbb{\Sigma}^\tau_\periodic}, \quad \big\langle \Psi,\; \Phi \big\rangle_{\mathbb{\Sigma}^\tau_\periodic} = \frac{1}{\cuti_1} \int_0^1 \big\langle \shearmap \Psi(\cdot, s), \; \shearmap \Phi(\cdot, s) \big\rangle_{\Scuti^\tau} \; ds.
\end{equation} 
In what follows, $\Hminusonehalfper{\mathbb{\Sigma}^\tau_\periodic}$ is endowed with the norm $\Psi \mapsto \|\shearmap \Psi\|_{L^2(0, 1; H^{-1/2}(\Scuti^\tau))}$.

\vspace{1\baselineskip} \noindent
At last, we are able to define the normal trace operator thanks to the next result.
\begin{prop}\label{H2Dmod:prop:Green_formula_strips}
  Given $\tau \in \partial I_\xvi$, the normal trace application $\gamma^\tau_{1, \periodic}$ defined by $\gamma^\tau_{1, \periodic} \boldsymbol{W} := (\cutmat\, \boldsymbol{W} \cdot \nva)|_{\mathbb{\Sigma}^\tau_\periodic}$ for any $\boldsymbol{W} \in [\mathscr{C}^\infty_{0, \periodic}(\overline{\mathbb{\Omega}_\periodic})]^2$, extends by continuity to a linear and surjective map still denoted by $\gamma^\tau_{1, \periodic}$, and defined from $\Hdiveper{\mathbb{\Omega}_\periodic}$ to $\Hminusonehalfper{\mathbb{\Sigma}^\tau_\periodic}$.

  \vspace{1\baselineskip} \noindent
  Moreover, we have the Green's formula: for any $V \in \Hgradper{\mathbb{\Omega}_\periodic}$ and $\boldsymbol{W} \in \Hdiveper{\mathbb{\Omega}_\periodic}$,
  \begin{equation}\label{H2Dmod:eq:Green_formula_strips}
    \displaystyle
    \int_{\mathbb{\Omega}_\periodic} \big[(\transp{\bsnabla} \cutmat \boldsymbol{W}) \; \overline{V} + \boldsymbol{W} \cdot \overline{\transp{\cutmat} \bsnabla V} \big] = \sum_{\tau \in \partial I_\xvi} \big\langle \cutmat\, \boldsymbol{W} \cdot \nva,\; V \big\rangle_{\mathbb{\Sigma}^\tau_\periodic},
  \end{equation}
  where we have written abusively $\cutmat\, \boldsymbol{W} \cdot \nva \equiv (\cutmat\, \boldsymbol{W} \cdot \nva)|_{\mathbb{\Sigma}^\tau_\periodic}$ instead of $\gamma^\tau_{1, \periodic} \boldsymbol{W}$.
\end{prop}

\begin{rmk}
  Note that the integrals on $\partial \mathbb{\Omega}_\periodic \cap \{\zvi_2 \in \{0, 1\}\}$ do not appear in the Green's formula \eqref{H2Dmod:eq:Green_formula_strips} because the functions are \textquote{periodic} in the $\eva_2$--direction (see Remark \ref{H2Dmod:rmk:H1Cper_HdivCper_weak_periodicity}).
\end{rmk}

\begin{dem}[of Proposition \ref{H2Dmod:prop:Green_formula_strips}]
  Similarly to the proof of Proposition \ref{H2Dmod:prop:gamma0_continuous_surjective}, we begin by observing that $\shearmap$ commutes with $\gamma^\tau_{1, \periodic}$ in the sense that for any $\boldsymbol{W} \in [\mathscr{C}^\infty_{0, \periodic}(\overline{\mathbb{\Omega}_\periodic})]^2$ (with $\nva = \pm \eva_\xvi$ and $\nv = \pm \ev_\xvi$),
  \begin{equation}\label{H2Dmod:eq:gamma1_commutes_with_S}
    \displaystyle
    \spforall s \in (0, 1), \quad \shearmap \big((\cutmat\, \boldsymbol{W} \cdot \nva)|_{\mathbb{\Sigma}^\tau_\periodic}\big) (\cdot, s) = \big(\shearmap \boldsymbol{W} (\cdot, s) \cdot \nv\big) |_{\Scuti^\tau}.
  \end{equation}
  Equation \eqref{H2Dmod:eq:gamma1_commutes_with_S} extends to $\boldsymbol{W}\! \in\! \Hdiveper{\mathbb{\Omega}_\periodic}$ using $(i)$ the density of $[\mathscr{C}^\infty_{0, \periodic}(\overline{\mathbb{\Omega}_\periodic})]^2$ in $\Hdiveper{\mathbb{\Omega}_\periodic}$, $(ii)$ the continuity of the shear map $\shearmap$ from $\Hdiveper{\mathbb{\Omega}_\periodic}$ to $L^2(0, 1; H(\dive; \Qcuti))$ (Corollary \ref{H2Dmod:cor:shear_properties}) and from $\Hminusonehalfper{\mathbb{\Sigma}_\periodic}$ to the space $L^2(0, 1; H^{-1/2} (\Scuti^\tau))$ (by definition \eqref{H2Dmod:eq:def_H12Cper}), and $(iii)$ the continuity from $L^2(0, 1; H(\dive; \Qcuti))$ to $L^2(0, 1; H^{-1/2} (\Scuti^\tau))$ of the following normal trace application: 
  \[\displaystyle
    \spforall \boldsymbol{W}_\cutmat \in L^2(0, 1; H(\dive; \Qcuti)), \quad \gamma_1 \boldsymbol{W}_\cutmat (\cdot, s) := \big(\boldsymbol{W}_\cutmat (\cdot, s) \cdot \nv\big) |_{\Scuti^\tau},
  \]
  which follows directly from the classical trace theorem in $H(\dive; \Qcuti)$. 
  
  \vspace{1\baselineskip} \noindent
  Since $\shearmap$ is an isomorphism, \eqref{H2Dmod:eq:gamma0_commutes_with_S} implies that $\gamma^\tau_{1, \periodic} = \invshearmap\, \gamma_1 \shearmap$ \surligner{(see Figure \ref{H2Dmod:fig:isomorphism_traces} right)}. The continuity and surjectivity of $\gamma^\tau_{1, \periodic}$ then follows from the continuity of $(\gamma_1, \shearmap, \invshearmap)$ and from the surjectivity of $\gamma_1$. 
  
  \vspace{1\baselineskip} \noindent
  To finish we prove the Green's formula \eqref{H2Dmod:eq:Green_formula_strips}. For $V \in \Hgradper{\mathbb{\Omega}_\periodic}$ and $\boldsymbol{W} \in \Hdiveper{\mathbb{\Omega}_\periodic}$, the classical Green formula applied to $(\shearmap V (\cdot, s), \shearmap \boldsymbol{W}(\cdot, s)) \in H^1(\Qcuti) \times H(\dive; \Qcuti)$ for almost any $s \in (0, 1)$ and integrated with respect to $s$ leads to
  \begin{multline*}\displaystyle
    \int_0^1 \int_{\Qcuti} \big[\transp{\nabla}_\xv (\shearmap \boldsymbol{W})\; \overline{\shearmap V} + \shearmap \boldsymbol{W} \cdot \overline{\nabla_\xv \shearmap V}](\xv, s) \; d\xv ds 
    \\
    = \sum_{\tau \in \partial I_\xvi} \int_0^1 \big\langle (\shearmap \boldsymbol{W}(\cdot, s) \cdot \nv)|_{\Scuti^\tau},\; \shearmap V(\cdot, s)|_{\Scuti^\tau} \big\rangle_{\Scuti^\tau} \; ds,
  \end{multline*}
  where we recall that $\langle \cdot, \cdot \rangle_{\Scuti^\tau}$ is the duality product between $H^{-1/2}(\Scuti^\tau)$ and $H^{1/2}(\Scuti^\tau)$. To conclude, we use the properties (\ref{H2Dmod:eq:oblique_plancherel}, \ref{H2Dmod:eq:shear_commutes_with_grad}, \ref{H2Dmod:eq:shear_commutes_with_dive}) of $\shearmap$ for the left-hand side, as well as the identity \eqref{H2Dmod:eq:shearmap_fibered_duality} for the right-hand side.
\end{dem}

\subsection{Subspaces of periodic functions in a cylinder}\label{H2Dmod:sec:e1_e2_periodic_functions}
We now adress the case \eqref{H2Dmod:item:case_cylinder} of the cylinder. In addition to the strip $\mathbb{\Omega}_\periodic = I_\xvi \times \R \times (0, 1)$ and the interface $\mathbb{\Sigma}^\tau_\periodic = \{\tau\} \times \R \times (0, 1)$ defined in Section \ref{H2Dmod:sec:gamma1_Greens_formula_strip}, we also introduce the cylinder and the interface
\begin{equation}\label{H2Dmod:eq:Qperper_Sigmaperper}
  \displaystyle
  \mathbb{\Omega}_\perperiodic := I_\xvi \times (0, 1) \times (0, 1) \quad \textnormal{and} \quad \mathbb{\Sigma}^\tau_\perperiodic := \{\tau\} \times (0, 1) \times (0, 1),
\end{equation}
where the subscript "$\perperiodic$" refers to the boundedness of the domains in the $\eva_1$ and the $\eva_2$--directions. 

\vspace{1\baselineskip} \noindent
Note that $\mathbb{\Omega}_\perperiodic$ and $\mathbb{\Sigma}^\tau_\perperiodic$ correspond \surligner{respectively} to the domains $\genO_\periodic$ and $\genS^\tau_\periodic$ in \eqref{H2Dmod:eq:cylindrical_domain} with
\begin{equation*}
  \surligner{I_1 = (0, 1).}
\end{equation*}
Therefore, one can use the spaces $\Hgradper{\mathbb{\Omega}_\perperiodic}$, $\Hdiveper{\mathbb{\Omega}_\perperiodic}$ \surligner{defined by \eqref{H2Dmod:eq:def_H1Cper_HdivCper}, as well as the space $\Honehalfper{\mathbb{\Sigma}^\tau_\perperiodic}$ defined by \eqref{H2Dmod:eq:def_H12Cper}}. In the sequel, we assume that $I_\xvi \neq \R$, so that $\partial I_\xvi \neq \varnothing$.

\vspace{1\baselineskip} \noindent
For the purpose of Sections \ref{H2Dmod:sec:FB_transform} and \ref{H2Dmod:sec:DtN_waveguide}, we want to define spaces of $H^1_\cutmat$ or $H_\cutmat(\mathbf{div})$--functions in $\mathbb{\Omega}_\perperiodic$ which are $1$--periodic with respect to \surligner{both} $\zvi_1$ and $\zvi_2$. To do so, let us begin with the following definition.
\begin{defi}\label{H2Dmod:def:periodic_extension_e1}
  Let $\boldsymbol{V} \in [L^2 (\mathbb{\Omega}_\perperiodic)]^d$. The \emph{periodic extension $\extper[1]{\boldsymbol{V}} \in [L^2_{\textit{loc}}(\mathbb{\Omega}_\periodic)]^d$ of $\boldsymbol{V}$ in the $\eva_1$--direction} is defined by:
  \begin{equation}\label{H2Dmod:eq:periodic_extension_e1}
     \aeforall \xva = (\xvi, \zvi_1, \zvi_2) \in \mathbb{\Omega}_\perperiodic, \quad \spforall n \in \Z, \quad (\extper[1]{\boldsymbol{V}}) (\xva + n\vts \eva_1) := \boldsymbol{V}(\xva).
  \end{equation}
  For $\Phi \in L^2 (\mathbb{\Sigma}^\tau_\perperiodic)$, we define $\extper[1]{\Phi} \in L^2_{\textit{loc}}(\mathbb{\Sigma}^\tau_\periodic)$ similarly, with $\xva \in \mathbb{\Sigma}^\tau_\perperiodic$ instead of $\xva \in \mathbb{\Omega}_\perperiodic$.
\end{defi}

\vspace{1\baselineskip} \noindent
In addition, let $\mathscr{C}^\infty_{0, \perperiodic}(\overline{\mathbb{\Omega}_\perperiodic})$ be the space of smooth functions in $\overline{\mathbb{\Omega}_\perperiodic}$ that are compactly supported in the $\eva_\xvi$--direction and $1$--periodic in the $\eva_1$ and $\eva_2$--directions, namely
\begin{equation}
  \label{H2Dmod:eq:def_Cinftyperper}
  \mathscr{C}^\infty_{0, \perperiodic}(\overline{\mathbb{\Omega}_\perperiodic}) := \big\{V \in \mathscr{C}^\infty_0(\overline{\mathbb{\Omega}_\perperiodic}) \ /\ \extper[1]{V} \in \mathscr{C}^\infty_{0, \periodic}(\overline{\mathbb{\Omega}_\periodic})\big\},
\end{equation}
where $\mathscr{C}^\infty_{0, \periodic}(\overline{\mathbb{\Omega}_\periodic})$ is given by \eqref{H2Dmod:eq:def_Cinftyper}. Note that $\mathscr{C}^\infty_{0, \perperiodic}(\overline{\mathbb{\Omega}_\perperiodic})$ contains $\mathscr{C}^\infty_0(\mathbb{\Omega}_\perperiodic)$, and therefore, is dense in $L^2(\mathbb{\Omega}_\perperiodic)$.

\vspace{1\baselineskip} \noindent
Now, we introduce $\Hgradperper{\mathbb{\Omega}_\perperiodic}$ (\emph{resp.} $\Hdiveperper{\mathbb{\Omega}_\perperiodic}$), which is defined as the closure of $\mathscr{C}^\infty_{0, \perperiodic}(\overline{\mathbb{\Omega}_\perperiodic})$ (\emph{resp.} $[\mathscr{C}^\infty_{0, \perperiodic}(\overline{\mathbb{\Omega}_\perperiodic})]^2$) in $H^1_\cutmat(\mathbb{\Omega}_\perperiodic)$ (\emph{resp.} $H_\cutmat(\mathbf{div}; \mathbb{\Omega}_\perperiodic)$):
\begin{equation}
  \displaystyle
  \label{H2Dmod:eq:def_H1Cperper_HdivCperper}
  \boxed{
  \Hgradperper{\mathbb{\Omega}_\perperiodic} := \overline{\mathscr{C}^\infty_{0, \perperiodic}(\overline{\mathbb{\Omega}_\perperiodic})}^{{H^1_\cutmat(\mathbb{\Omega}_\perperiodic)}} %
  \quad \textnormal{and} \quad%
   \Hdiveperper{\mathbb{\Omega}_\perperiodic} := \overline{[\mathscr{C}^\infty_{0, \perperiodic}(\overline{\mathbb{\Omega}_\perperiodic})]^2}^{{H_\cutmat(\mathbf{div}; \mathbb{\Omega}_\perperiodic)}}. 
  }
\end{equation}
We note that $\Hgradperper{\mathbb{\Omega}_\perperiodic}$ and $\Hdiveperper{\mathbb{\Omega}_\perperiodic}$ are respectively closed subspaces of $\Hgradper{\mathbb{\Omega}_\perperiodic}$ and $\Hdiveper{\mathbb{\Omega}_\perperiodic}$, and hence, are Hilbert spaces when equipped with the respective norms of these last spaces. 

\begin{rmk}\label{H2Dmod:rmk:H1Cqp_HdivCqp_weak_qperiodicity}
  Similarly to Remark \ref{H2Dmod:rmk:H1Cper_HdivCper_weak_periodicity}, for smooth functions, the definition of $\Hgradperper{\mathbb{\Omega}_\perperiodic}$ corresponds to periodicity in the $\eva_1$ and $\eva_2$--directions. More generally, one could show the characterization 
  \begin{equation*}
    \Hgradperper{\mathbb{\Omega}_\perperiodic} = \{V \in H^1_\cutmat(\mathbb{\Omega}_\perperiodic)\ /\ V|_{\zvi_1 = 0} = V|_{\zvi_1 = 1} \ \textnormal{and}\ V|_{\zvi_2 = 0} = V|_{\zvi_2 = 1}\},
  \end{equation*}
  where the equalities have to be understood not in the $L^2$ sense, but in the $L^2_{\textit{loc}}$ sense. We refer to \cite[Section 3.2.1]{amenoagbadji2023wave} for similar considerations.
\end{rmk}

\noindent
\surligner{Section \ref{H2Dmod:sec:gamma0} and Proposition \ref{H2Dmod:prop:gamma0_continuous_surjective} applied for $I_1 = (0, 1)$ allow to define the trace operator on $\mathbb{\Sigma}^\tau_\perperiodic$ ($\tau \in \partial I_\xvi$) for functions in $\Hgradperper{\mathbb{\Omega}_\perperiodic}$. Introduce}
\begin{equation}\label{H2Dmod:eq:def_H12Cperper}
  \displaystyle
  \Honehalfperper{\mathbb{\Sigma}^\tau_\perperiodic} := \big\{V|_{\mathbb{\Sigma}^\tau_\perperiodic} \ /\ \ V \in \Hgradperper{\mathbb{\Omega}_\perperiodic} \big\},
\end{equation}
to which we associate the graph norm. Let $\Hminusonehalfperper{\mathbb{\Sigma}^\tau_\perperiodic}$ be the topological dual space of $\Honehalfperper{\mathbb{\Sigma}^\tau_\perperiodic}$, equipped with the dual norm. The duality product $\langle \cdot, \cdot\rangle_{\mathbb{\Sigma}^\tau_\perperiodic}$ between $\Hminusonehalfperper{\mathbb{\Sigma}^\tau_\perperiodic}$ and $\Honehalfperper{\mathbb{\Sigma}^\tau_\perperiodic}$ is defined as a natural extension of the $\smash{L^2(\mathbb{\Sigma}^\tau_\perperiodic)}$--scalar product.

\vspace{1\baselineskip} \noindent
%
%
The density of $\mathscr{C}^\infty_{0, \perperiodic}(\overline{\mathbb{\Omega}_\perperiodic})$ in $\Hgradperper{\mathbb{\Omega}_\perperiodic}$, the density of $[\mathscr{C}^\infty_{0, \perperiodic}(\overline{\mathbb{\Omega}_\perperiodic})]^2$ in $\Hdiveperper{\mathbb{\Omega}_\perperiodic}$, and the continuity of the trace operator on $\mathbb{\Sigma}^\tau_\perperiodic$ lead directly to the next result.
\begin{prop}
  For $\tau \in \partial I_\xvi$, the normal trace map $\gamma^\tau_{1, \perperiodic}$ defined by $\gamma^\tau_{1, \perperiodic} \boldsymbol{W} := (\cutmat\, \boldsymbol{W} \cdot \nva)|_{\mathbb{\Sigma}^\tau_\perperiodic}$ for any $\smash{\boldsymbol{W} \in [\mathscr{C}^\infty_{0, \perperiodic}(\overline{\mathbb{\Omega}_\perperiodic})]^2}$ extends by continuity to a linear and surjective map still denoted by $\gamma^\tau_{1, \perperiodic}$, and defined from $\Hdiveperper{\mathbb{\Omega}_\perperiodic}$ to $\Hminusonehalfperper{\mathbb{\Sigma}^\tau_\perperiodic}$.

  \vspace{1\baselineskip} \noindent
  Moreover, we have the Green's formula: for any $V \in \Hgradperper{\mathbb{\Omega}_\perperiodic}$ and $\boldsymbol{W} \in \Hdiveperper{\mathbb{\Omega}_\perperiodic}$,
  \begin{equation}\label{H2Dmod:eq:Green_cylinder}
    \displaystyle
    \int_{\mathbb{\Omega}_\perperiodic} \big[(\transp{\bsnabla} \cutmat \boldsymbol{W}) \; \overline{V} + \boldsymbol{W} \cdot \overline{\transp{\cutmat} \bsnabla V} \big] = \sum_{\tau \in \partial I_\xvi} \big\langle (\cutmat\, \boldsymbol{W} \cdot \nva),\; V \big\rangle_{\mathbb{\Sigma}^\tau_\perperiodic}, 
  \end{equation}
  where we have abusively written $(\cutmat\, \boldsymbol{W} \cdot \nva) \equiv (\cutmat\, \boldsymbol{W} \cdot \nva)|_{\mathbb{\Sigma}^\tau_\perperiodic}$ instead of $\gamma^\tau_{1, \perperiodic} \boldsymbol{W}$.
\end{prop}

%

\vspace{1\baselineskip}\noindent
We finish with the following result, which will be used in Section \ref{H2Dmod:sec:DtN_waveguide}. The proof is similar to the proof of the jump rule for isotropic Sobolev spaces, and relies on the Green's formula \eqref{H2Dmod:eq:Green_cylinder}.
\begin{prop}\label{H2Dmod:prop:formule_sauts}
  Assume that $\mathbb{\Omega}_\perperiodic = \mathbb{\Omega}^1_\perperiodic \cup \mathbb{\Omega}^2_\perperiodic$ where $\mathbb{\Omega}^1_\perperiodic$ and $\mathbb{\Omega}^2_\perperiodic$ are disjoint cylinders defined by
  \[
    \mathbb{\Omega}^\ell_\perperiodic := J^\ell_\xvi \times (0, 1) \times (0, 1) \subset \mathbb{\Omega}_\perperiodic, \quad \spforall \ell \in \{1, 2\}.
  \]
  where $J^1_\xvi \cup J^2_\xvi = I_\xvi$ and $\overline{J^1_\xvi} \cap \overline{J^2_\xvi} = \{\tau\}$ (see Figure \ref{H2Dmod:fig:formule_sauts}). Let $\mathbb{\Sigma}^\tau_\perperiodic := \overline{\mathbb{\Omega}^1_\perperiodic} \cap \overline{\mathbb{\Omega}^2_\perperiodic}$.
  
  \begin{enumerate}[label=$(\textit{\alph*})$., ref = \textit{\alph*}]
    \item Let $V_\ell \in \Hgradperper{\mathbb{\Omega}^\ell_\perperiodic}$, $\ell \in \{1, 2\}$. Then the function defined by
    \[
      \displaystyle
      \aeforall \xva \in \mathbb{\Omega}_\perperiodic, \quad V(\xva) := %
      \left\{%
        \begin{array}{cl}
          V_1(\xva) & \textnormal{if} \quad \xva \in \mathbb{\Omega}^1_\perperiodic
          \\[4pt]
          V_2(\xva) & \textnormal{if} \quad \xva \in \mathbb{\Omega}^2_\perperiodic
        \end{array}
      \right.
    \]
    belongs to $\Hgradperper{\mathbb{\Omega}_\perperiodic}$ if and only if $V_1|_{\mathbb{\Sigma}^\tau_\perperiodic} = V_2|_{\mathbb{\Sigma}^\tau_\perperiodic}$.
    \item Let $\boldsymbol{W}_\ell \in \Hdiveperper{\mathbb{\Omega}^\ell_\perperiodic}$, $\ell \in \{1, 2\}$. Then the function given by
    \[
      \displaystyle
      \aeforall \xva \in \mathbb{\Omega}_\perperiodic, \quad \boldsymbol{W}(\xva) := %
      \left\{%
        \begin{array}{cl}
          \boldsymbol{W}_1(\xva) & \textnormal{if} \quad \xva \in \mathbb{\Omega}^1_\perperiodic
          \\[4pt]
          \boldsymbol{W}_2(\xva) & \textnormal{if} \quad \xva \in \mathbb{\Omega}^2_\perperiodic
        \end{array}
      \right.
    \]
    belongs to $\Hdiveperper{\mathbb{\Omega}_\perperiodic}$ if and only if $(\cutmat\, \boldsymbol{W}_1 \cdot \eva_\xvi)|_{\mathbb{\Sigma}^\tau_\perperiodic} = (\cutmat\, \boldsymbol{W}_2 \cdot \eva_\xvi)|_{\mathbb{\Sigma}^\tau_\perperiodic}$.
  \end{enumerate}
\end{prop}

\begin{figure}[H]
  \noindent\makebox[\textwidth]{%
    \def\coteCyl{1.525}
    \begin{tikzpicture}[scale=0.675]
      \def\angleinterface{55}%
      \definecolor{coplancoupe}{RGB}{211, 227, 65}
      \tdplotsetmaincoords{70}{20}
      \def\BBxmin{-4}
      \def\BBxmax{+4}
      \def\BBymin{-1.25}
      \def\BBymax{+3.75}
      \def\BBzmin{0}
      \def\BBzmax{+2.75}
      \begin{scope}[tdplot_main_coords, scale=0.9]
        \node (xaxis) at ({\BBxmax+0.1*(\BBxmax-\BBxmin)}, 0, 0) {$\xvi$};%
        \node (yaxis) at (0, {\BBymax+0.75*(\BBymax-\BBymin)}, 0) {$\zvi_1$};%
        %
        %
        %
        \fill[secondary!80, opacity=0.5] (\BBxmin, 0, 0) -- +(-\BBxmin, 0, 0) -- +(-\BBxmin, \coteCyl, 0) -- +(0, \coteCyl, 0) -- cycle;
        \fill[secondary!90, opacity=0.625] (\BBxmin, \coteCyl, 0) -- +(-\BBxmin, 0, 0) -- +(-\BBxmin, 0, \coteCyl) -- +(0, 0, \coteCyl) -- cycle;
        \draw[dashed, white] (\BBxmin, \coteCyl, 0) -- (\BBxmax, \coteCyl, 0);
        \node[secondary!80!black] (Omega) at (\BBxmin-1.25, 0.5*\coteCyl, 1.5*\coteCyl) {$\mathbb{\Omega}^2_\perperiodic$};
        \draw[-latex, secondary!80!black] (Omega) to[bend right=40] (\BBxmin+1.5, 0.5*\coteCyl, 0.5*\coteCyl);
        \filldraw[draw=white, fill=secondary!90, opacity=0.625] (\BBxmin, 0, 0) -- +(-\BBxmin, 0, 0) -- +(-\BBxmin, 0, \coteCyl) -- +(0, 0, \coteCyl) -- cycle;
        \filldraw[draw=white, fill=secondary!80, opacity=0.5] (\BBxmin, 0, \coteCyl) -- +(-\BBxmin, 0, 0) -- +(-\BBxmin, \coteCyl, 0) -- +(0, \coteCyl, 0) -- cycle;
        \filldraw[draw=white, fill=secondary!95!black, opacity=0.875] (0, 0, 0) -- (0, \coteCyl, 0) -- (0, \coteCyl, \coteCyl) -- (0, 0, \coteCyl) -- cycle;
        %
        \draw[-latex] (0, {\BBymin-0.1*(\BBymax-\BBymin)}, 0) -- (yaxis);%
        \fill[tertiary!80, opacity=0.5] (\BBxmax, 0, 0) -- +(-\BBxmax, 0, 0) -- +(-\BBxmax, \coteCyl, 0) -- +(0, \coteCyl, 0) -- cycle;
        \fill[tertiary!90, opacity=0.625] (\BBxmax, \coteCyl, 0) -- +(-\BBxmax, 0, 0) -- +(-\BBxmax, 0, \coteCyl) -- +(0, 0, \coteCyl) -- cycle;
        \draw[dashed, white] (\BBxmax, \coteCyl, 0) -- (\BBxmax, \coteCyl, 0);
        \node[tertiary!80!black] (Omega) at (\BBxmax+1.25, 0.5*\coteCyl, 1.5*\coteCyl) {$\mathbb{\Omega}^1_\perperiodic$};
        \draw[-latex, tertiary!80!black] (Omega) to[bend left=40] (\BBxmax-1.5, 0.5*\coteCyl, 0.5*\coteCyl);
        \filldraw[draw=white, fill=tertiary!90, opacity=0.625] (\BBxmax, 0, 0) -- +(-\BBxmax, 0, 0) -- +(-\BBxmax, 0, \coteCyl) -- +(0, 0, \coteCyl) -- cycle;
        \filldraw[draw=white, fill=tertiary!80, opacity=0.5] (\BBxmax, 0, \coteCyl) -- +(-\BBxmax, 0, 0) -- +(-\BBxmax, \coteCyl, 0) -- +(0, \coteCyl, 0) -- cycle;
        %
        %
        %
        %
        \draw[-latex] (0, 0, {\BBzmin-0.25*(\BBzmax-\BBzmin)}) -- (0, 0, {\BBzmax+0.1*(\BBzmax-\BBzmin)}) node[above] {$\smash{\zvi_2}$};
        \draw[-latex] ({\BBxmin-0.1*(\BBxmax-\BBxmin)}, 0, 0) -- (xaxis);%
      \end{scope}%
    \end{tikzpicture} 
  }
  \caption{The domains $\mathbb{\Omega}^1_\perperiodic$ and $\mathbb{\Omega}^2_\perperiodic$ in Proposition \ref{H2Dmod:prop:formule_sauts} for $I_\xvi = \R$, $J^1_\xvi = \R_+$, and $J^2_\xvi = \R_-$.\label{H2Dmod:fig:formule_sauts}}
\end{figure}

}{%
}

\ifthenelse{\boolean{brouillon}}{%
  \section{The solution of the augmented periodic problem}
\label{H2Dmod:sec:DtN_approach_augmented_problem}
Using the function spaces introduced in Section \ref{H2Dmod:sec:functional_framework}, we propose in this section a rigorous \surligner{formulation of the augmented problem} \eqref{H2Dmod:eq:formal_augmented_PDE_recalled}. This problem is then studied and solved by exploiting the periodic nature of $(\aten_p, \rho_p)$, following the steps \eqref{H2Dmod:item:lifting_step_1}--\eqref{H2Dmod:item:lifting_step_3} described in Section \ref{H2Dmod:sec:formal_description_method}.

\subsection{The augmented strip problem and its quasi-2D structure}\label{H2Dmod:sec:analysis_augmented_problem}
We consider the $3$--dimensional strips $\Omegaper, \Omegaper^\pm$ and the 2--dimensional interface $\Sigmaper$ defined by
\begin{equation}
  \label{H2Dmod:eq:strips}
  \displaystyle
  \begin{array}{r@{\ :=\ }l}
    \Omegaper & \{(\xvi, \zvi_1, \zvi_2) \in \R^3\ /\ \zvi_2 \in (0, 1)\},
    \\[10pt]
    \Omegaper^\pm & \{(\xvi, \zvi_1, \zvi_2) \in \Omegaper\ /\  \pm \xvi > 0\},
    \\[10pt]
    \Sigmaper & \{(\xvi, \zvi_1, \zvi_2) \in \Omegaper\ /\ \xvi = 0\}.
  \end{array}
\end{equation}
These domains are represented in Figure \ref{H2Dmod:fig:strip}. We \surligner{shall use the indexing} $\Omegaper^\idom$, $\idom \in \{\varnothing, +, -\}$ where by convention, $\Omegaper^\idom = \Omegaper$ for $\idom = \varnothing$. Note that $\Omegaper^\idom$ and $\Sigmaper$ correspond respectively to the domains \surligner{$\mathbb{\Omega}_\periodic := I_\xvi \times \R \times (0, 1)$ and $\mathbb{\Sigma}^\tau_\periodic := \{\tau\} \times \R \times (0, 1)$ introduced in Section \ref{H2Dmod:sec:gamma1_Greens_formula_strip}}, with
\[
  I_\xvi := \R_\idom, \quad \Qcuti = \R_\idom \times \R, \quad \textnormal{and} \quad \Scuti = \{0\} \times \R.
\]
Therefore, Sections \ref{H2Dmod:sec:e2_periodic_functions}, \ref{H2Dmod:sec:gamma0}, and \ref{H2Dmod:sec:gamma1_Greens_formula_strip} enable us to use \surligner{the spaces $\Hgradper{\Omegaper^\idom}$, $\Hdiveper{\Omegaper^\idom}$ defined} by \eqref{H2Dmod:eq:def_H1Cper_HdivCper}, the space $\Honehalfper{\Sigmaper}$ given by \eqref{H2Dmod:eq:def_H12Cper} and its dual $\Hminusonehalfper{\Sigmaper}$ characterized by \eqref{H2Dmod:eq:carac_Hminus12Cper}, and the trace and normal trace applications on $\Sigmaper$.

\begin{figure}[H]
  \noindent\makebox[\textwidth]{%
    \begin{tikzpicture}[scale=0.675]
      \def\angleinterface{55}%
      \definecolor{coplancoupe}{RGB}{211, 227, 65}
      \tdplotsetmaincoords{80}{20}
      \def\coteCyl{1.525}
      \def\BBxmin{-4.5}
      \def\BBxmax{+4.5}
      \def\BBymin{-5.25}
      \def\BBymax{+5.5}
      \def\BBzmin{0}
      \def\BBzmax{+2.75}
      \def\coupevecY{1.414213562}%
      \def\coupevecZ{0.75}
      \pgfmathsetmacro\xcoupemax{min(\BBymax/\coupevecY, \BBzmax/\coupevecZ)}
      \pgfmathsetmacro\xcoupemin{max(\BBymin/\coupevecY, \BBzmin/\coupevecZ)}
      \begin{scope}[tdplot_main_coords, scale=0.9, shift={(-4.15cm, 1.15cm)}]%
        \node (xaxis) at ({\BBxmax+0.25*(\BBxmax-\BBxmin)}, 0, 0) {$\xvi$};%
        \draw[-latex] ({\BBxmin-0.1*(\BBxmax-\BBxmin)}, 0, 0) -- (xaxis);%
        \node (yaxis) at (0, {0.5+\BBymax+1*(\BBymax-\BBymin)}, 0) {$\zvi_1$};%
        %
        %
        %
        %
        %
        %
        \fill[canvas is zx plane at y=\BBymax, fill=secondary!70, opacity=0.5] (0, \BBxmin) rectangle (\coteCyl, 0);
        \fill[canvas is yz plane at x=\BBxmin, fill=secondary!60, opacity=0.675] (\BBymin, 0) rectangle (\BBymax, \coteCyl);
        \fill[secondary!80, opacity=0.675] (\BBxmin, \BBymin, 0) -- +(-\BBxmin, 0, 0) -- +(-\BBxmin, \BBymax-\BBymin, 0) -- +(0, \BBymax-\BBymin, 0) -- cycle;
        \node[secondary!80!black] (Omega) at (\BBxmin-2, 0, 2*\coteCyl) {$\Omegaper^-$};
        \draw[-latex, opacity=1, secondary!80!black] (Omega) to[bend right=40] (\BBxmin+1, 0, 0.5*\coteCyl);
        \fill[secondary!90!black, opacity=0.75] (\BBxmin, \BBymin, \coteCyl) -- +(-\BBxmin, 0, 0) -- +(-\BBxmin, \BBymax-\BBymin, 0) -- +(0, \BBymax-\BBymin, 0) -- cycle;
        \filldraw[thick, draw=white, fill=secondary!95, opacity=0.875] (0, \BBymin, 0) -- +(0, \BBymax-\BBymin, 0) -- +(0, \BBymax-\BBymin, \coteCyl) -- +(0, 0, \coteCyl) -- cycle;
        \fill[canvas is zx plane at y=\BBymin, fill=secondary!80, opacity=0.675] (0, \BBxmin) rectangle (\coteCyl, 0);
        %
        %
        \draw[-latex] (0, {\BBymin-0.1*(\BBymax-\BBymin)}, 0) -- (yaxis);%
        \draw[-latex] (0, 0, {\BBzmin-0.0*(\BBzmax-\BBzmin)}) -- (0, 0, {\BBzmax+0.25*(\BBzmax-\BBzmin)}) node[above] {$\smash{\zvi_2}$};
        \fill[canvas is xy plane at z=0, tertiary!60, opacity=0.675] (\BBxmax, \BBymin) -- +(-\BBxmax, 0) -- +(-\BBxmax, \BBymax-\BBymin) -- +(0, \BBymax-\BBymin) -- cycle;
        \draw[-latex, secondary!20!black] (\BBxmax, 0.75*\BBymax, 0.5*\coteCyl) -- (\BBxmax-1, 0.75*\BBymax, 0.5*\coteCyl);
        \fill[canvas is zx plane at y=\BBymax, fill=tertiary!70, opacity=0.25] (0, \BBxmax) rectangle (\coteCyl, 0);
        \fill[canvas is xy plane at z=\coteCyl, tertiary!95!black, opacity=0.75] (\BBxmax, \BBymin) -- +(-\BBxmax, 0) -- +(-\BBxmax, \BBymax-\BBymin) -- +(0, \BBymax-\BBymin) -- cycle;
        \fill[canvas is zx plane at y=\BBymin, fill=tertiary!80, opacity=0.675] (0, \BBxmax) rectangle (\coteCyl, 0);
        \fill[canvas is yz plane at x=\BBxmax, fill=tertiary!87, opacity=0.675] (\BBymin, 0) rectangle (\BBymax, \coteCyl);

        \node[tertiary!80!black] (Omega) at (\BBxmax+2, 0.75*\BBymax, 0.5*\coteCyl) {$\Omegaper^+$};
        \draw[tertiary!80!black] (Omega) -- (\BBxmax, 0.75*\BBymax, 0.5*\coteCyl);

      \end{scope}%
    \end{tikzpicture}
  }
  \caption{The half-strips $\Omegaper^+$ and $\Omegaper^-$ defined in \eqref{H2Dmod:eq:strips}\label{H2Dmod:fig:strip}}
\end{figure}

\noindent
Using the formal observations in the step \eqref{H2Dmod:item:lifting_step_1} of Section \ref{H2Dmod:sec:formal_description_method}, we introduce for any $G \in \Hminusonehalfper{\Sigmaper}$ the strip problem: \textit{find $U \in H^1_\cutmat (\Omegaper)$ such that}
\begin{equation}
  \label{H2Dmod:eq:augmented_transmission_problem}
  \left\{
    \begin{array}{r@{\ =\ }l}
      \displaystyle- \transp{\bsnabla}\, \cutmat\vts \aten_p \!  \transp{\cutmat} \bsnabla\, U - \rho_p\, \omega^2\, U & 0 \quad \textnormal{in}\ \  \Omegaper^+ \cup \Omegaper^-,
      \ret
      \multicolumn{2}{c}{\displaystyle U \in H^1_{\cutmat, \periodic} (\Omegaper), \ \  (\aten_p \! \transp{\cutmat} \bsnabla\, U)|_{\Omegaper^\pm} \in \Hdiveper{\Omegaper^\pm},}
      \ret
      \displaystyle \llbracket (\cutmat\vts \aten_p \! \transp{\cutmat} \bsnabla\, U) \cdot \eva_\xvi \rrbracket_{\Sigmaper} & G,
    \end{array}
  \right. \tag{$\mathscr{P}_\periodic$}
\end{equation}
{\setstretch{1.15}%
where $\aten_p, \rho_p$ are given by \eqref{H2Dmod:eq:Apm_rhopm_config_a} for Configuration \eqref{H2Dmod:item:config_a} and by (\ref{H2Dmod:eq:Ap_from_A0}, \ref{H2Dmod:eq:def_Am_rhom_config_b}) for Configuration \eqref{H2Dmod:item:config_b}. Note that from the first equation, one has $\smash{(\aten_p \! \transp{\cutmat} \bsnabla\, U) |_{\Omegaper^\pm} \in H_\cutmat(\mathbf{div}; \Omegaper^\pm)}$. Formally speaking, the second equation requires $U$ and $\smash{(\aten_p \! \transp{\cutmat} \bsnabla\, U)|_{\Omegaper^\pm}}$ to be $1$--periodic with respect to $\zvi_2$ (see Remark \ref{H2Dmod:rmk:H1Cper_HdivCper_weak_periodicity}). Finally, $\smash{\llbracket \cutmat\, \boldsymbol{W} \cdot \eva_\xvi \rrbracket_{\Sigmaper}}$ is the jump accross $\Sigmaper$ of the normal components of a function $\boldsymbol{W} \in [L^2(\Omegaper)]^2$ which satisfies $\boldsymbol{W}^\pm := \boldsymbol{W}|_{\Omegaper^\pm} \in \Hdiveper{\Omegaper^\pm}$:
\par}
\begin{equation}
  \llbracket \cutmat\, \boldsymbol{W} \cdot \eva_\xvi \rrbracket_{\Sigmaper} := (\cutmat\, \boldsymbol{W}^- \cdot \eva_\xvi)|_{\Sigmaper} - (\cutmat\, \boldsymbol{W}^+ \cdot \eva_\xvi)|_{\Sigmaper},
\end{equation}
and the third equation holds in $\Hminusonehalfper{\Sigmaper}$.

\vspace{1\baselineskip} \noindent
In what follows, we shall study Problem \eqref{H2Dmod:eq:augmented_transmission_problem} and emphasize its fibered structure as a "concatenation" of $2$--dimensional problems. For that latter purpose, it is useful to introduce the family of functions indexed by $s \in \R$:
\begin{equation}\label{H2Dmod:eq:def_As_rho_s}
  \displaystyle
  \spforall s \in \R, \quad \aten_s(\xv) := \aten_p(\cutmat\, \xv + s\, \eva_2) \quad \textnormal{and} \quad \rho_s(\xv) := \rho_p(\cutmat\, \xv + s\, \eva_2), \quad \spforall \xv \in \R^2.
\end{equation}
Note that the restrictions of $\aten_s$ and $\rho_s$ to $\R^2_\pm$ are well-defined and continuous since for both Configurations \eqref{H2Dmod:item:config_a} and \eqref{H2Dmod:item:config_b}, $\aten_p$ and $\rho_p$ depend in a simple manner on $\aten$ and $\rho$ which are continuous on $\R^2_\pm$. Furthermore, $s \mapsto \aten_s$ and $s \mapsto \rho_s$ are continuous from $\R$ to $L^\infty(\R^2; \R^{2\times 2})$ and $L^\infty(\R^2)$ respectively. Note also that for any $s \in \R$, $\aten_s$ and $\rho_s$ satisfy the ellipticity assumption \eqref{H2Dmod:eq:ellipticity_assumption}, and that $\aten_{s + 1} = \aten_s$ and $\rho_{s + 1} = \rho_s$. 

\vspace{1\baselineskip} \noindent
Given $\psi \in H^{-1/2}(\sigma)$, we consider the $2$--dimensional problem: \textit{Find $u_s \in H^1(\R^2)$ such that}
\begin{equation}
  \label{H2Dmod:eq:transmission_problems}
  \left\{
    \begin{array}{r@{\ =\ }l}
      \displaystyle- \transp{\nabla} \aten_s(\xv)\, \nabla u_s(\xv) - \rho_s(\xv)\, \omega^2\, u_s(\xv) & 0 \quad \textnormal{for}\ \ \xv \in \R^2_+ \cup \R^2_-,
      \\[8pt]
      \displaystyle \llbracket \aten_s\, \nabla u_s \cdot \ev_\xvi \rrbracket_{\sigma} & \psi.
    \end{array}
  \right.
  \tag{$\mathscr{P}_s$}
\end{equation}

\begin{rmk}\label{H2Dmod:rmk:transmission_problems_for_s_zero}
  One has $\aten_0 = \aten$ and $\rho_0 = \rho$, and \eqref{H2Dmod:eq:transmission_problem} corresponds to \eqref{H2Dmod:eq:transmission_problems} with $s = 0$ and $\psi = g$.
\end{rmk}

\vspace{1\baselineskip} \noindent
Because of the ellipticity assumption \eqref{H2Dmod:eq:ellipticity_assumption} and the presence of absorption \eqref{H2Dmod:eq:absorption}, Problem \eqref{H2Dmod:eq:transmission_problems} admits a unique solution $\linapp{u_s}{\psi} \in H^1(\R^2)$. Furthermore, the well-posedness of Problem \eqref{H2Dmod:eq:augmented_transmission_problem} and the link between its solution $U$ and the solutions $\linapp{u_s}{\psi}$ is given by the following result.
\begin{prop}\label{H2Dmod:prop:augmented_BVP_equiv_FV_and_Quasi2D_structure}
  For any $G \in \Hminusonehalfper{\Sigmaper}$, Problem \eqref{H2Dmod:eq:augmented_transmission_problem} is equivalent to the variational problem
  \begin{equation}
    \label{H2Dmod:eq:FV_augmented_transmission}
    \left|
    \begin{array}{l}
      \textit{Find $U \in \Hgradper{\Omegaper}$ such that $\spforall V \in \Hgradper{\Omegaper}$,}
      \ret
      \multicolumn{1}{c}{\displaystyle\int_{\Omegaper} \big[\,(\aten_p \! \transp{\cutmat} \bsnabla\, U) \cdot \overline{\vphantom{U^+}(\transp{\cutmat} \bsnabla\, V)} - \rho_p\, \omega^2\, U\, \overline{V}\,\big] = \big\langle G,\; V \big\rangle_{\Sigmaper}},
    \end{array}
    \right.\tag{$\mathrm{FV}_\periodic$}
  \end{equation}
  which is well-posed in $\Hgradper{\Omegaper}$. Furthermore, 
  \begin{equation}\label{H2Dmod:eq:fibered_structure_strip}
    \aeforall s \in (0, 1), \quad \extper[2] U (\xvi, \cuti_1\, \zvi, \cuti_2\, \zvi + s) = u_s (\xv), \quad \aeforall \xv = (\xvi, \zvi) \in \R^2,
  \end{equation}
  where $U = \linapp{U}{G}$ and $u_s = \linapp{u_s}{g_s}$ with $g_s := \shearmap {G}(\cdot, s)$ are the respective solutions of Problems \eqref{H2Dmod:eq:augmented_transmission_problem} and \eqref{H2Dmod:eq:transmission_problems}.
\end{prop}


\begin{dem}
  To obtain the variational formulation \eqref{H2Dmod:eq:FV_augmented_transmission}, one multiplies the volume equation in \eqref{H2Dmod:eq:augmented_transmission_problem} by $V \in \Hgradper{\Omegaper}$, and integrates over $\Omegaper^+$ and $\Omegaper^-$ separately. Since $\aten_p \! \transp{\cutmat} \bsnabla\, U \in \Hdiveper{\Omegaper^\pm}$, we then apply Green's formula \eqref{H2Dmod:eq:Green_formula_strips} on each domain $\Omegaper^+$ and $\Omegaper^-$, add the corresponding identities, and we finally use the transmission condition (that is, the third equation in \eqref{H2Dmod:eq:augmented_transmission_problem}) to conclude.

  \vspace{1\baselineskip} \noindent
  The well-posedness of the variational formulation \eqref{H2Dmod:eq:FV_augmented_transmission} is guaranteed by Lax-Milgram's theorem. In particular, the bilinear form $\mathscr{A}$ associated to \eqref{H2Dmod:eq:FV_augmented_transmission} satisfies
  \begin{equation}
    \spforall V \in \Hgradper{\Omegaper}, \quad \Imag \Big[\frac{\mathscr{A}(V, V)}{\omega}\Big] = - \Imag \omega \; \int_{\Omegaper} \Big[\, \frac{1}{|\omega|^2}\,\aten_p |\transp{\cutmat} \bsnabla\, V|^2 + \rho_p\, |V|^2\,\Big],
  \end{equation}
  and therefore is coercive in $\Hgradper{\Omegaper}$ due to the presence of absorption \eqref{H2Dmod:eq:absorption}, that is, $\Imag \omega > 0$.

  \vspace{1\baselineskip} \noindent
  {
  Now we shall prove \eqref{H2Dmod:eq:fibered_structure_strip}, which will enable us to deduce the equivalence between \eqref{H2Dmod:eq:FV_augmented_transmission} and \eqref{H2Dmod:eq:augmented_transmission_problem} at last. To this end, given $s \in (0, 1)$, we show that $\shearmap U(\cdot, s) \in H^1(\R^2)$ satisfies the same problem as $\linapp{u_s}{\shearmap G(\cdot, s)}$. For $V \in \Hgradper{\Omegaper}$, the properties (\ref{H2Dmod:eq:shear_commutes_with_grad}, \ref{H2Dmod:eq:oblique_plancherel}) of $\shearmap$ combined with \eqref{H2Dmod:eq:FV_augmented_transmission} and \eqref{H2Dmod:eq:shearmap_fibered_duality} imply that \par}
  \begin{align}
    \int_0^1 \int_{\R^2} \big[\,\aten_s \nabla_\xv \shearmap U(\cdot, s) \cdot \overline{\nabla_\xv \shearmap V(\cdot, s)}& - \rho_s\, \omega^2\, \shearmap U(\cdot, s)\, \overline{\shearmap V(\cdot, s)}\,\big] d\xv ds \nonumber
    \\
    &= \cuti_1 \int_{\Omegaper} \big[\,(\aten_p \! \transp{\cutmat} \bsnabla\, U) \cdot \overline{\vphantom{U^+}(\transp{\cutmat} \bsnabla\, V)} - \rho_p\, \omega^2\, U\, \overline{V}\,\big]\nonumber
    \\
    &= \cuti_1\, \big\langle G,\; V \big\rangle_{\Sigmaper} = \int_0^1 \big\langle \shearmap G(\cdot, s), \; \shearmap V(\cdot, s)|_{\sigma} \big\rangle_{\sigma} \; ds.\label{H2Dmod:eq:augmented_BVP_equiv_FV_and_Quasi2D_structure:dem0}
  \end{align}
  Now we choose $V$ such that $\shearmap V(\xv, s) = \varphi(s)\, v(\xv)$, with $\varphi \in L^2(0, 1)$ and $v \in H^1(\R^2)$. Then it follows that $\nabla_\xv \shearmap V(\xv, s) = \varphi(s)\, \nabla v(\xv)$, and \eqref{H2Dmod:eq:augmented_BVP_equiv_FV_and_Quasi2D_structure:dem0} being true for any $\varphi \in L^2(0, 1)$ leads for almost any $s \in (0, 1)$ to the equation
  \[
    \spforall v \in H^1(\R^2), \quad \int_{\R^2} \big[\,\aten_s \nabla_\xv \shearmap U(\cdot, s) \cdot \overline{\nabla v} - \rho_s\, \omega^2\, \shearmap U(\cdot, s)\, \overline{v}\,\big] \; d\xv = \big\langle \shearmap G(\cdot, s), \; v \big\rangle_{\sigma},
  \]
  which is equivalent to the transmission problem \eqref{H2Dmod:eq:transmission_problems} satisfied by $\linapp{u_s}{\shearmap G(\cdot, s)}$. From the uniqueness of the solution of \eqref{H2Dmod:eq:transmission_problems}, we then deduce $\shearmap U(\cdot, s) = \linapp{u_s}{\shearmap G(\cdot, s)}$ which is \eqref{H2Dmod:eq:fibered_structure_strip}.

  \vspace{1\baselineskip} \noindent
  It remains to prove the equivalence between Problem \eqref{H2Dmod:eq:augmented_transmission_problem} and its variational formulation \eqref{H2Dmod:eq:FV_augmented_transmission}. But we have shown in the last step that if $U$ satisfies \eqref{H2Dmod:eq:FV_augmented_transmission}, then $(\shearmap U) (\cdot, s) = \linapp{u_s}{\shearmap G(\cdot, s)}$. Combining this result with the properties (\ref{H2Dmod:eq:shear_commutes_with_grad}, \ref{H2Dmod:eq:shear_commutes_with_dive}) of $\shearmap$ thus leads to the first and second equations in \eqref{H2Dmod:eq:augmented_transmission_problem}, while the transmission condition (that is, the third equation in \eqref{H2Dmod:eq:augmented_transmission_problem}) follows from \eqref{H2Dmod:eq:shearmap_fibered_duality} and from the transmission condition satisfied by $\linapp{u_s}{\shearmap G(\cdot, s)}$ in \eqref{H2Dmod:eq:transmission_problems}.
\end{dem}

\vspace{1\baselineskip} \noindent
To conclude this section, we propose in Proposition \ref{H2Dmod:prop:rigorous_ansatz} a rigorous formulation of the ansatz \eqref{H2Dmod:eq:lifting_ansatz}. To begin, note that due to the uniform boundedness of $(\aten_s)_s$ and $(\rho_s)_s$ in $L^\infty(\R^2)$, $(\linapp{u_s}{\psi})_s$ is bounded uniformly with respect to $s$:
\begin{equation}\label{H2Dmod:eq:uniform_boundedness_us}
  \displaystyle
  \spexists c > 0, \quad \spforall s \in \R, \quad \big\|\linapp{u_s}{\psi} \big\|_{H^1(\R^2)} \leq c\; \|\psi\|_{H^{-1/2}(\sigma)}, \quad \spforall \psi \in H^{-1/2}(\sigma).
\end{equation}
By combining this property with the continuity of $\aten_s$ and $\rho_s$ in $L^\infty(\R^2)$ with respect to $s$, one can show that $s \mapsto \linapp{u_s}{\psi}$ is a continuous application. This is the object of the next proposition, whose proof is similar to the proof of \cite[Proposition 3.18]{amenoagbadji2023wave}.
\begin{prop}\label{H2Dmod:prop:uniform_continuity_us}
  For $\psi \in H^{-1/2}(\sigma)$, $s \mapsto \linapp{u_s}{\psi}$ is continuous from $\R$ to $H^1(\R^2)$ and is $1$--periodic.
\end{prop}

\vspace{1\baselineskip} \noindent
The next proposition shows how the solution $u$ of \eqref{H2Dmod:eq:transmission_problem} can be retrieved from the solution $\linapp{U}{G}$ of the augmented problem \eqref{H2Dmod:eq:augmented_transmission_problem} for a well-chosen boundary data $G$ (see Remark \ref{H2D:rmk:example_G}).
\begin{prop}\label{H2Dmod:prop:rigorous_ansatz}
  Consider $G \in \Hminusonehalfper{\Sigmaper}$ such that $s \mapsto \shearmap G (\cdot, s) \in H^{-1/2}(\sigma)$ is continuous at $0$, and that $\shearmap G (\cdot, 0) = g$, where $g \in H^{-1/2}(\sigma)$ is the jump data in \eqref{H2Dmod:eq:transmission_problem}. Then setting $U = \linapp{U}{G}$, the map $s \mapsto (\shearmap U) (\cdot, s) \in H^1(\R^2)$ is continuous as well at $0$, and
  \begin{equation}
    \label{H2Dmod:eq:rigorous_ansatz}
    \displaystyle
    \aeforall \xv = (\xvi, \zvi) \in \R^2, \quad u(\xv) = \extper{U} (\xvi, \cuti_1\, \zvi, \cuti_2\, \zvi).
  \end{equation}
\end{prop}

\begin{dem}
  Since $\psi \mapsto \linapp{u_s}{\psi}$ is continuous from $ H^{-1/2}(\sigma)$ to $H^1(\R^2)$ uniformly with respect to $s$ (according to \eqref{H2Dmod:eq:uniform_boundedness_us}) on one hand, and $s \mapsto \linapp{u_s}{\psi}$ is continuous (according to Proposition \ref{H2Dmod:prop:uniform_continuity_us}) on the other hand, it follows that $s \mapsto \linapp{u_s}{\shearmap G(\cdot, s)}$ is continuous at $0$ as a product of continuous maps. Consequently \eqref{H2Dmod:eq:fibered_structure_strip} becomes true for $s = 0$, and thus leads to \eqref{H2Dmod:eq:rigorous_ansatz}.
\end{dem}

\begin{rmk}\label{H2D:rmk:example_G}
  If $g \in L^2(\sigma)$, a class of augmented data $G \in \Hminusonehalfper{\Sigmaper}$ such that $\shearmap G (\cdot, 0) = g$ is given by
  \begin{equation*}
    \aeforall (0, \zvi_1, \zvi_2) \in \Sigmaper, \quad G(0, \zvi_1, \zvi_2) := f(\zvi_2 - \zvi_1 \cuti_2/\cuti_1)\, g(0, \zvi_1/\cuti_1),    
  \end{equation*}
  where $f \in \mathscr{C}^0(\R)$ is such that $f(s + 1) = f(s)$ and $f(0) = 1$ (e.g. $f(s) = \euler^{2\icplx \pi \ell s}$, $\ell \in \Z$). In fact, it can be seen that $\shearmap G (0, \zvi, 0) := \extper{G}(0, \cuti_1 \zvi, \cuti_2 \zvi) = f(0) g(0, \zvi) = g(0, \zvi)$. If $g \in H^{-1/2}(\sigma)$, then the above expression of $G$ can be extended by duality.
\end{rmk}

\subsection{Reduction to waveguide problems via the Floquet-Bloch transform}\label{H2Dmod:sec:FB_transform}
Now that we have shown that to solve the original problem \eqref{H2Dmod:eq:transmission_problem}, it is sufficient to study the augmented strip problem \eqref{H2Dmod:eq:augmented_transmission_problem}, we begin by taking advantage of the periodicity with respect to $\zvi_1$ by applying a Floquet-Bloch transform. %
Using the properties of the $1$D Floquet-Bloch transform (see \cite{kuchment2001mathematics} or \cite[Section 2.2]{kuchment1993floquet}), we introduce in Section \ref{H2Dmod:sec:partial_FB_anisotropic_spaces} the partial Floquet-Bloch transform with respect to $\zvi_1$, and present its properties in $\Hgradper{\Omegaper}$ or $\Hdiveper{\Omegaper}$. The transform is then applied to the solution of \eqref{H2Dmod:eq:augmented_transmission_problem} in Section \ref{H2Dmod:sec:application_FB_to_augmented_transmission_problem}, leading to a family of half-guide problems denoted \eqref{H2Dmod:eq:augmented_waveguide_problem} later on, and parameterized by the Floquet dual variable.

\subsubsection{The partial Floquet-Bloch transform with respect to \texorpdfstring{$\zvi_1$}{z1}}
\label{H2Dmod:sec:partial_FB_anisotropic_spaces}
We begin by introducing
\begin{equation}
  \label{H2Dmod:eq:cylinders}
  \displaystyle
  \begin{array}{l@{\quad}r@{\ :=\ \{(\xvi, \zvi_1, \zvi_2) \in \ }l@{\ /\ \zvi_1 \in (0, 1) \},}}
    \spforall \idom \in \{\varnothing, +, -\}, & \Omegaperper^\idom & \Omegaper^\idom
    \\[10pt]
    & \Sigmaperper & \Sigmaper
  \end{array}
\end{equation}
where the subscript "$\perperiodic$" indicates the boundedness of the domains in the $\eva_1$ and the $\eva_2$--directions. These domains are represented in Figure \ref{H2Dmod:fig:cylinder}. Note that $\Omegaperper^\idom$ and $\Sigmaperper$ correspond respectively to the domains \surligner{$\mathbb{\Omega}_\perperiodic := I_\xvi \times (0, 1) \times (0, 1)$ and $\mathbb{\Sigma}_\perperiodic := \{0\} \times (0, 1) \times (0, 1)$} defined by \eqref{H2Dmod:eq:Qperper_Sigmaperper} with
\[
  I_\xvi := \R_\idom, \quad \Qcuti = \R_\idom \times (0, 1/\cuti_1), \quad \textnormal{and} \quad \Scuti = \{0\} \times (0, 1/\cuti_1).
\]
Therefore, thanks to Section \ref{H2Dmod:sec:e1_e2_periodic_functions}, we can use the spaces $\Hgradperper{\Omegaperper}$, $\Hdiveperper{\Omegaperper}$ given by \eqref{H2Dmod:eq:def_H1Cperper_HdivCperper} (see also Remark \ref{H2Dmod:rmk:H1Cqp_HdivCqp_weak_qperiodicity} for the meaning of these spaces), as well as the space $\Honehalfperper{\Sigmaperper}$ defined by \eqref{H2Dmod:eq:def_H12Cperper}, its dual $\smash{\Hminusonehalfperper{\Sigmaperper}}$, and the trace and normal trace operators on $\Sigmaperper$.

\begin{figure}[H]
  \noindent\makebox[\textwidth]{%
    \def\coteCyl{1.525}
    \begin{tikzpicture}[scale=0.675]
      \def\angleinterface{55}%
      \definecolor{coplancoupe}{RGB}{211, 227, 65}
      \tdplotsetmaincoords{75}{20}
      \def\BBxmin{-4}
      \def\BBxmax{+4}
      \def\BBymin{-1.25}
      \def\BBymax{+4.75}
      \def\BBzmin{-0.75}
      \def\BBzmax{+2.75}
      \begin{scope}[tdplot_main_coords, scale=0.9]
        \node (xaxis) at ({\BBxmax+0.1*(\BBxmax-\BBxmin)}, 0, 0) {$\xvi$};%
        \node (yaxis) at (0, {0.5+\BBymax+0.75*(\BBymax-\BBymin)}, 0) {$\zvi_1$};%
        %
        %
        %
        \fill[secondary!80, opacity=0.5] (\BBxmin, 0, 0) -- +(-\BBxmin, 0, 0) -- +(-\BBxmin, \coteCyl, 0) -- +(0, \coteCyl, 0) -- cycle;
        \fill[secondary!90, opacity=0.625] (\BBxmin, \coteCyl, 0) -- +(-\BBxmin, 0, 0) -- +(-\BBxmin, 0, \coteCyl) -- +(0, 0, \coteCyl) -- cycle;
        \draw[dashed, white] (\BBxmin, \coteCyl, 0) -- (\BBxmax, \coteCyl, 0);
        \node[secondary!80!black] (Omega) at (\BBxmin-1.25, 0.5*\coteCyl, 1.5*\coteCyl) {$\Omegaperper^-$};
        \draw[-latex, secondary!80!black] (Omega) to[bend right=40] (\BBxmin+1.5, 0.5*\coteCyl, 0.5*\coteCyl);
        \filldraw[draw=white, fill=secondary!90, opacity=0.625] (\BBxmin, 0, 0) -- +(-\BBxmin, 0, 0) -- +(-\BBxmin, 0, \coteCyl) -- +(0, 0, \coteCyl) -- cycle;
        \filldraw[draw=white, fill=secondary!80, opacity=0.5] (\BBxmin, 0, \coteCyl) -- +(-\BBxmin, 0, 0) -- +(-\BBxmin, \coteCyl, 0) -- +(0, \coteCyl, 0) -- cycle;
        \filldraw[draw=white, fill=secondary!95!black, opacity=0.875] (0, 0, 0) -- (0, \coteCyl, 0) -- (0, \coteCyl, \coteCyl) -- (0, 0, \coteCyl) -- cycle;
        %
        \draw[-latex] (0, {\BBymin-0.0*(\BBymax-\BBymin)}, 0) -- (yaxis);%
        \fill[tertiary!80, opacity=0.5] (\BBxmax, 0, 0) -- +(-\BBxmax, 0, 0) -- +(-\BBxmax, \coteCyl, 0) -- +(0, \coteCyl, 0) -- cycle;
        \fill[tertiary!90, opacity=0.625] (\BBxmax, \coteCyl, 0) -- +(-\BBxmax, 0, 0) -- +(-\BBxmax, 0, \coteCyl) -- +(0, 0, \coteCyl) -- cycle;
        \draw[dashed, white] (\BBxmax, \coteCyl, 0) -- (\BBxmax, \coteCyl, 0);
        \node[tertiary!80!black] (Omega) at (\BBxmax+1.25, 0.5*\coteCyl, 1.5*\coteCyl) {$\Omegaperper^+$};
        \draw[-latex, tertiary!80!black] (Omega) to[bend left=40] (\BBxmax-1.5, 0.5*\coteCyl, 0.5*\coteCyl);
        \filldraw[draw=white, fill=tertiary!90, opacity=0.625] (\BBxmax, 0, 0) -- +(-\BBxmax, 0, 0) -- +(-\BBxmax, 0, \coteCyl) -- +(0, 0, \coteCyl) -- cycle;
        \filldraw[draw=white, fill=tertiary!80, opacity=0.5] (\BBxmax, 0, \coteCyl) -- +(-\BBxmax, 0, 0) -- +(-\BBxmax, \coteCyl, 0) -- +(0, \coteCyl, 0) -- cycle;
        \draw[-latex] (0, 0, {\BBzmin-0.0*(\BBzmax-\BBzmin)}) -- (0, 0, {\BBzmax+0.25*(\BBzmax-\BBzmin)}) node[above] {$\smash{\zvi_2}$};
        \draw[-latex] ({\BBxmin-0.1*(\BBxmax-\BBxmin)}, 0, 0) -- (xaxis);%
      \end{scope}%
    \end{tikzpicture} 
  }
  \caption{The half-cylinders $\Omegaperper^+$ and $\Omegaperper^-$ defined in \eqref{H2Dmod:eq:cylinders}\label{H2Dmod:fig:cylinder}}
\end{figure}

\noindent
Let $d \in \{1, 2\}$. The Floquet-Bloch transform in the $\eva_1$--direction is defined by
\begin{equation}
  \label{H2Dmod:eq:def_partial_TFB}
  \displaystyle
  \begin{array}{r@{\ }c@{\ }c@{\ }l}
    \fbtransform[]: & \mathscr{C}^\infty_0(\overline{\Omegaper})^d & \longrightarrow & [\mathscr{C}^\infty(\Omegaper \times \R)]^d
    \\[5pt]
    &\boldsymbol{V} & \longmapsto & \widehat{\boldsymbol{V}},
  \end{array}
  \qquad 
  \left\{
    \begin{array}{l}
      \displaystyle\widehat{\boldsymbol{V}}(\xva, \fbvar) := \frac{1}{\sqrt{2\pi}} \sum_{n \in \Z} \boldsymbol{V}(\xva + n\vts \eva_1)\, \euler^{-\icplx\vts \fbvar\vts (\zvi_1 + n)},
      \\[16pt]
      \quad \spforall \xva = (\xvi, \zvi_1, \zvi_2) \in \Omegaper, \  \fbvar \in \R.
    \end{array}
  \right.
\end{equation}
Note that $\xva \mapsto \widehat{\boldsymbol{V}}(\xva, \fbvar)$ is $\Z \eva_1$--periodic, and $\fbvar \mapsto \euler^{\icplx\vts \fbvar\vts \zvi_1}\, \widehat{\boldsymbol{V}}(\xva, \fbvar)$ is $2\pi$--periodic. Therefore, $\widehat{\boldsymbol{V}}$ can be fully constructed from its knowledge on the cell $\Omegaperper \times (-\pi, \pi)$. For this reason, the study of $\widehat{\boldsymbol{V}}$ will be restricted to $\Omegaperper \times (-\pi, \pi)$.

\vspace{1\baselineskip} \noindent
The transform $\fbtransform[]$ extends to an isometry from $[L^2(\Omegaper)]^d$ to $L^2(-\pi, \pi; [L^2(\Omegaperper)]^d)$, with
\begin{equation}
  \label{H2Dmod:eq:partial_TFB_plancherel}
  \displaystyle
  \spforall \boldsymbol{U}, \boldsymbol{V} \in [L^2 (\Omegaper)]^d, \quad \int_{-\pi}^\pi \int_{\Omegaperper} \fbtransform[] \boldsymbol{U}(\xva, \fbvar)\cdot \overline{\fbtransform[] \boldsymbol{V}(\xva, \fbvar)}\; d\xva d\fbvar = \int_{\Omegaper} \boldsymbol{U}\cdot \overline{\boldsymbol{V}}.
\end{equation}
Moreover, $\fbtransform[]$ is an isomorphism, with the inversion formula: for any $\widehat{\boldsymbol{U}} \in L^2(-\pi, \pi; [L^2(\Omegaperper)]^d)$,
\begin{equation}\label{H2Dmod:eq:partial_TFB_inversion}
  \displaystyle
  \aeforall \xva = (\xvi, \zvi_1, \zvi_2) \in \Omegaperper, \quad \spforall n \in \Z, \quad  (\invfbtransform[] \widehat{\boldsymbol{U}})(\xva + n\vts \eva_1) = \frac{1}{\sqrt{2\pi}} \, \int_{-\pi}^{\pi} \widehat{\boldsymbol{U}} (\xva, \fbvar)\;  \euler^{\icplx\vts \fbvar\vts (\zvi_1 + n)} \; d\fbvar.
\end{equation}
In addition, for $\Phi \in \mathscr{C}^\infty_0(\overline{\Sigmaper})$, we define $\widehat{\Phi} := \fbtransform[] \Phi \in \mathscr{C}^\infty_0(\Sigmaper \times \R)$ as in \eqref{H2Dmod:eq:def_partial_TFB} by choosing $\xva \in \Sigmaper$. Then $\fbtransform[]$ extends to an isometry between $L^2(\Sigmaper)$ and $L^2(-\pi, \pi; L^2(\Sigmaperper))$.

\vspace{1\baselineskip}\noindent
What makes the Floquet-Bloch transform a valuable tool in the context of this paper is on one hand the fact that it commutes with the multiplication by $\Z\eva_1$--periodic coefficients:
\begin{equation}\label{H2Dmod:eq:fb_mult}
  \displaystyle
  \spforall \mu \in L^\infty(\Omegaper) \ \textnormal{such that}\ \mu(\cdot + \eva_1) = \mu, \quad \spforall V \in L^2(\Omegaper), \quad \fbtransform[] (\mu V) = \mu \fbtransform[] V,
\end{equation}
and on the other hand its action on differential operators. More precisely, the next proposition gives the properties of $\fbtransform[]$ in $\Hgradper{\Omegaper}$, $\Hdiveper{\Omegaper}$, and $\Honehalfper{\Sigmaper}$. Its proof is delayed to Appendix \ref{H2Dmod:sec:proof_FB_properties}.
\begin{prop}\label{H2Dmod:prop:TFB_properties}~
  \begin{enumerate}[label=$(\textit{\alph*})$., ref = \theprop.\textit{\alph*}, wide, labelindent=0pt]
    \setlength\itemsep{8pt}
    \item $\fbtransform[]$ is an isomorphism from $\Hgradper{\Omegaper}$ to $L^2(-\pi, \pi; \Hgradperper{\Omegaperper})$, and
    \begin{equation} \label{H2Dmod:eq:fb_grad}
      \displaystyle
      \spforall V \in \Hgradper{\Omegaper}, \quad \aeforall \fbvar \in (-\pi, \pi), \quad \fbtransform[] (\gradcut V)\vts (\cdot, \fbvar) = \gradcutFB\, \fbtransform[] V (\cdot, \fbvar).
    \end{equation}
    \item $\fbtransform[]$ is an isomorphism from $\Hdiveper{\Omegaper}$ to $L^2(-\pi, \pi; \Hdiveperper{\Omegaperper})$, and
    \begin{equation} \label{H2Dmod:eq:fb_dive}
      \displaystyle
      \spforall \boldsymbol{W} \in \Hdiveper{\Omegaper}, \quad \spforall \fbvar \in (-\pi, \pi), \quad \fbtransform[] (\divecut\, \boldsymbol{W})\vts (\cdot, \fbvar) = \divecutFB\, \fbtransform[] \boldsymbol{W} (\cdot, \fbvar).
    \end{equation}

    \item $\fbtransform[]$ is an isomorphism from $\Honehalfper{\Sigmaper}$ to $L^2(-\pi, \pi; \Honehalfperper{\Sigmaperper})$.
  \end{enumerate}
\end{prop}

\vspace{1\baselineskip} \noindent
We highlight the following consequence of Proposition \ref{H2Dmod:prop:TFB_properties} for future reference.
\begin{cor}\label{H2Dmod:cor:TFB_properties}
  For any $U, V \in \Hgradper{\Omegaper}$,
  \begin{equation}
    \displaystyle
    \int_{-\pi}^\pi \int_{\Omegaperper} \gradcutFB\, \fbtransform[] U(\xva, \fbvar)\cdot \overline{\gradcutFB\, \fbtransform[] V(\xva, \fbvar)}\; d\xva d\fbvar = \int_{\Omegaper} \gradcut U \cdot \overline{\gradcut V}.
  \end{equation}
\end{cor}

\vspace{1\baselineskip} \noindent
Finally, $\fbtransform[]$ extends by duality as an isomorphism from $\Hminusonehalfper{\Sigmaper}$ to $\mathscr{X}' := [L^2(-\pi, \pi; \Honehalfperper{\Sigmaperper})]'$: 
\begin{equation*}
  \displaystyle
  \spforall \Psi \in \Hminusonehalfper{\Sigmaper}, \quad \big\langle \fbtransform[] \Psi,\; \widehat{\Phi} \big\rangle_{\mathscr{X}', \mathscr{X}} %
  := \big\langle \Psi,\; \invfbtransform[] \widehat{\Phi} \big\rangle_{\Sigmaper}, \quad \spforall \widehat{\Phi} \in \mathscr{X}.
\end{equation*}
Since $\Honehalfperper{\Sigmaperper}$ is a Hilbert space, Lemma \ref{H2Dmod:lem:dual_of_Bochner_space} implies $\mathscr{X}' = L^2(-\pi, \pi; \Hminusonehalfperper{\Sigmaperper})$ and 
\begin{equation}
  \displaystyle
  \spforall (\Phi, \Psi) \in \Honehalfper{\Sigmaper} \times \smash{\Hminusonehalfper{\Sigmaper}}, \quad \big\langle \Psi,\; \Phi \big\rangle_{\Sigmaper} = \int_{-\pi}^\pi \big\langle \fbtransform[] \Psi (\cdot, \fbvar),\; \fbtransform[] \Phi (\cdot, \fbvar) \big\rangle_{\Sigmaperper} \; d\fbvar.
\end{equation}

\subsubsection{Application to the augmented strip problem}
\label{H2Dmod:sec:application_FB_to_augmented_transmission_problem}
Thanks to the properties of the Floquet-Bloch transform $\fbtransform[]$ given in Section \ref{H2Dmod:sec:partial_FB_anisotropic_spaces}, we deduce directly the following proposition.
\begin{prop}\label{H2Dmod:prop:from_strip_to_waveguide}
Let $G \in \Hminusonehalfper{\Sigmaperper}$. Then the solution $\linapp{U}{G} \in \Hgradper{\Omegaper}$ of \eqref{H2Dmod:eq:augmented_transmission_problem} is given by
\begin{equation}\label{H2Dmod:eq:strip_from_waveguide}
  \displaystyle
  \aeforall \xva = (\xvi, \zvi_1, \zvi_2) \in \Omegaperper, \quad \spforall n \in \Z, \quad \linapp{U}{G}(\xva + n\vts \eva_1) = \frac{1}{\sqrt{2\pi}} \, \int_{-\pi}^{\pi} \linapp{\widehat{U}_\fbvar}{\widehat{G}_\fbvar}\vts (\xva)\;  \euler^{-\icplx\vts \fbvar\vts (\zvi_1 + n)} \; d\fbvar,
\end{equation}
where $\widehat{G}_\fbvar := \fbtransform[] G (\cdot, \fbvar) \in \Hminusonehalfperper{\Sigmaperper}$ $\aeforall \fbvar \in (-\pi, \pi)$, and where $\linapp{\widehat{U}_\fbvar}{\widehat{G}_\fbvar} := \fbtransform[] \linapp{U}{G}\vts (\cdot, \fbvar)$ is the unique solution of the well-posed waveguide problem: \textit{find $\widehat{U}_\fbvar \in H^1_\cutmat (\Omegaperper)$ such that}
\begin{equation}
  \label{H2Dmod:eq:augmented_waveguide_problem}
  \left\{
    \begin{array}{r@{\ =\ }l}
      \displaystyle- \divecutFB\, \aten_p \!  \gradcutFB\, \widehat{U}_\fbvar - \rho_p\, \omega^2\, \widehat{U}_\fbvar & 0 \quad \textnormal{in}\ \  \Omegaperper^+ \cup \Omegaperper^-,
      \ret
      \multicolumn{2}{c}{\displaystyle \widehat{U}_\fbvar \in \Hgradperper{\Omegaperper}, \ \  (\aten_p \! \gradcutFB\, \widehat{U}_\fbvar)|_{\Omegaperper^\pm} \in \Hdiveperper{\Omegaperper^\pm},}
      \ret
      \displaystyle \llbracket (\cutmat\vts \aten_p \! \gradcutFB\, \widehat{U}_\fbvar) \cdot \eva_\xvi \rrbracket_{\Sigmaperper} & \widehat{G}_\fbvar,
    \end{array}
  \right. \tag{$\mathscr{P}_\perperiodic$}
\end{equation}
whose variational formulation is given by
\begin{equation}
\label{H2Dmod:eq:FV_augmented_waveguide}
\left|
\begin{array}{l}
  \textit{Find $\widehat{U}_\fbvar \in \Hgradperper{\Omegaperper}$ such that $\spforall V \in \Hgradperper{\Omegaperper}$,}
  \\[10pt]
  \multicolumn{1}{c}{\displaystyle
  \int_{\Omegaperper} \big[\,(\aten_p \vts \gradcutFB\, \widehat{U}_\fbvar \cdot \overline{\vphantom{U^+}\gradcutFB\, V} - \rho_p\, \omega^2\, \widehat{U}_\fbvar\, \overline{V}\,\big] = \big\langle \widehat{G}_\fbvar,\; V \big\rangle_{\Sigmaperper}.}
\end{array}
\right.\tag{$\mathrm{FV}_\perperiodic$}
\end{equation}
\end{prop}

\begin{dem}
  On one hand, Proposition \ref{H2Dmod:prop:TFB_properties} applied to \eqref{H2Dmod:eq:augmented_transmission_problem} shows that $\linapp{\widehat{U}_\fbvar}{\widehat{G}_\fbvar} := \fbtransform[] \linapp{U}{G}\vts (\cdot, \fbvar)$ satisfies \eqref{H2Dmod:eq:augmented_waveguide_problem}. On the other hand, Corollary \ref{H2Dmod:cor:TFB_properties} applied to the variational formulation \eqref{H2Dmod:eq:FV_augmented_transmission} implies that $\widehat{U}_\fbvar$ satisfies \eqref{H2Dmod:eq:FV_augmented_waveguide}. The equivalence between \eqref{H2Dmod:eq:augmented_waveguide_problem} and \eqref{H2Dmod:eq:FV_augmented_waveguide} then follows from the equivalence between the strip problem \eqref{H2Dmod:eq:augmented_transmission_problem} and its variational formulation \eqref{H2Dmod:eq:FV_augmented_transmission} (Proposition \ref{H2Dmod:prop:augmented_BVP_equiv_FV_and_Quasi2D_structure}).
\end{dem}

\vspace{1\baselineskip} \noindent
Proposition \ref{H2Dmod:prop:from_strip_to_waveguide} shows that $\linapp{U}{G}$ can be reconstructed in the strip $\Omegaper$, provided that one knows how to solve the waveguide problem \eqref{H2Dmod:eq:augmented_waveguide_problem} for any $\fbvar \in (-\pi, \pi)$. The solution of this waveguide problem is the object of the next section.

\subsection{Reducing the waveguide problem to the interface \texorpdfstring{$\Sigmaperper$}{Sigmaperper}}\label{H2Dmod:sec:DtN_waveguide}
In this section, the Floquet variable $\fbvar \in (-\pi, \pi)$ is fixed. We investigate the waveguide problem \eqref{H2Dmod:eq:augmented_waveguide_problem}, whose solution is denoted by $\linapp{\widehat{U}_\fbvar}{\widehat{G}_\fbvar}$. To this end, we reformulate the problem as an equation on $\Sigmaperper$ involving DtN operators associated to half-guide problems defined in $\Omegaperper^\pm$.

\vspace{\baselineskip} \noindent
Given boundary data $\Phi \in \Honehalfperper{\Sigmaperper}$, we consider the half-guide problem: \textit{find $\widehat{U}^\pm_\fbvar \in H^1_\cutmat (\Omegaperper)$ such that}
\begin{equation}
  \label{H2Dmod:eq:augmented_halfguide_problem}
  \left\{
    \begin{array}{r@{\ =\ }l}
      \displaystyle- \divecutFB\, \aten^\pm_p \!  \gradcutFB\, \widehat{U}^\pm_\fbvar - \rho^\pm_p\, \omega^2\, \widehat{U}^\pm_\fbvar & 0 \quad \textnormal{in}\ \  \Omegaperper^\pm,
      \ret
      \multicolumn{2}{c}{\displaystyle \widehat{U}^\pm_\fbvar \in \Hgradperper{\Omegaperper^\pm}, \ \  \aten^\pm_p \! \gradcutFB\, \widehat{U}^\pm_\fbvar \in \Hdiveperper{\Omegaperper^\pm},}
      \ret
      \displaystyle \widehat{U}^\pm_\fbvar & \Phi \quad \textnormal{on}\ \  \Sigmaperper.
    \end{array}
  \right. \tag{$\mathscr{P}^\pm_\perperiodic$}
\end{equation}
Under Assumptions (\ref{H2Dmod:eq:absorption}, \ref{H2Dmod:eq:ellipticity_assumption}), Lax-Milgram's theorem combined with a lifting argument ensures that \eqref{H2Dmod:eq:augmented_halfguide_problem} admits a unique solution
\begin{equation*}
  \linapp{\widehat{U}^\pm_\fbvar}{\Phi} \in \Hgradperper{\Omegaperper^\pm}.
\end{equation*}
Let $\Lambdahat^\pm_\fbvar \in \mathscr{L}(\Honehalfperper{\Sigmaperper}, \Hminusonehalfperper{\Sigmaperper})$ be the DtN operator defined by
\begin{equation*}
  \displaystyle
  \spforall \Phi \in \Honehalfperper{\Sigmaperper}, \quad \Lambdahat^\pm_\fbvar\, \Phi := \big(\cutmat\vts \aten_p \! \gradcutFB \vts \, \linapp{\widehat{U}^\pm_\fbvar}{\Phi} \cdot \nva\big) |_{\Sigmaperper},
\end{equation*}
or equivalently (thanks to the Green's formula \eqref{H2Dmod:eq:Green_cylinder})
\begin{multline}\label{H2Dmod:eq:def_DtN_operator}
  \displaystyle
  \spforall \Phi, \Psi \in \Honehalfperper{\Sigmaperper}, \quad  
  \big\langle \Lambdahat^\pm_\fbvar\, \Phi,\; \Psi \big\rangle_{\Sigmaperper} = \mathscr{A}^\pm_\fbvar \big(\linapp{\widehat{U}^\pm_\fbvar}{\Phi},\, \linapp{\widehat{U}^\pm_\fbvar}{\Psi}\big),%
  \\[8pt]
  \textnormal{where} \quad \mathscr{A}^\pm_\fbvar (U, V) := \int_{\Omegaperper^\pm} \big[\,(\aten^\pm_p \vts \gradcutFB\, U \cdot \overline{\vphantom{U^+}\gradcutFB\, V} - \rho^\pm_p\, \omega^2\, U\, \overline{V}\,\big].
\end{multline}
%
%
The above weak form allows to show the following result.
\begin{prop}\label{H2Dmod:prop:properties_DtN}
  The operators $\Lambdahat^+_\fbvar$, $\Lambdahat^-_\fbvar$, and $\Lambdahat^+_\fbvar + \Lambdahat^-_\fbvar$ are coercive and thus invertible from $\Honehalfperper{\Sigmaperper}$ to $\Hminusonehalfperper{\Sigmaperper}$.
\end{prop}

\begin{dem}
  Let $\idom \in \{+, -\}$. From \eqref{H2Dmod:eq:def_DtN_operator}, we have for any $\Phi \in \Honehalfperper{\Sigmaperper}$ the equality
  \begin{equation*}
    \Imag \bigg[\frac{\big\langle \Lambdahat^\idom_\fbvar\, \Phi,\; \Phi \big\rangle_{\Sigmaperper}}{\omega}\bigg] = - \Imag \omega \; \int_{\Omegaper} \bigg[\, \frac{1}{|\omega|^2}\,\aten_p \big|\gradcutFB\, \linapp{\widehat{U}^\idom_\fbvar}{\Phi}\big|^2 + \rho_p\, |\linapp{\widehat{U}^\idom_\fbvar}{\Phi}|^2\,\bigg],
  \end{equation*}
  implying that $\Lambdahat^\idom_\fbvar$ is coercive, due to absorption ($\Imag \omega > 0$) and to the continuity of the trace operator on $\Sigmaperper$. The same holds for $\Lambdahat^+_\fbvar + \Lambdahat^-_\fbvar$ by summing the above equality over $\idom \in \{+, -\}$.
\end{dem}

\vspace{1\baselineskip} \noindent
By linearity and by uniqueness of the solutions \surligner{of} \eqref{H2Dmod:eq:augmented_waveguide_problem} and \eqref{H2Dmod:eq:augmented_halfguide_problem}, the waveguide solution $\linapp{\widehat{U}_\fbvar}{\widehat{G}_\fbvar}$ admits the expression
\begin{equation}\label{H2Dmod:eq:Uhat_from_Uplusminushat}
  \displaystyle
  \aeforall \xva \in \Omegaperper, \quad \linapp{\widehat{U}_\fbvar}{\widehat{G}_\fbvar}(\xva) = 
  \left\{%
  \begin{array}{cl}
    \linapp{\widehat{U}^+_\fbvar}{\Phi_\fbvar}\vts (\xva) & \textnormal{if} \quad \xva \in \Omegaperper^+,
    \\[4pt]
    \linapp{\widehat{U}^-_\fbvar}{\Phi_\fbvar}\vts (\xva) & \textnormal{if} \quad \xva \in \Omegaperper^-,
  \end{array}
  \right. 
\end{equation}
where $\Phi_\fbvar := \widehat{U}_\fbvar |_{\Sigmaperper}$ is \surligner{determined} by the jump condition \surligner{on} $\widehat{U}_\fbvar$ (i.e. the last equation in \eqref{H2Dmod:eq:augmented_waveguide_problem}):  
\begin{equation}\label{H2Dmod:eq:transmission_jump}
  (\Lambdahat^+_\fbvar + \Lambdahat^-_\fbvar)\; \Phi_\fbvar = \widehat{G}_\fbvar \quad \textnormal{in} \quad \Hminusonehalfperper{\Sigmaperper}.
\end{equation}
Note that this equation is well-posed since $\Lambdahat^+_\fbvar + \Lambdahat^-_\fbvar$ is coercive according to Proposition \ref{H2Dmod:prop:properties_DtN}. Figure \ref{H2Dmod:fig:Uhat_from_Uplusminushat} illustrates the link between $\linapp{\widehat{U}_\fbvar}{\widehat{G}_\fbvar}$ and $\linapp{\widehat{U}^\pm_\fbvar}{\Phi_\fbvar}$. 

\begin{figure}[H]
  \noindent\makebox[\textwidth]{%
    \def\coteCyl{1.525}
    \begin{tikzpicture}[scale=0.675]
      \def\angleinterface{55}%
      \definecolor{coplancoupe}{RGB}{211, 227, 65}
      \tdplotsetmaincoords{80}{25}
      \def\BBxmin{-4}
      \def\BBxmax{+4.5}
      \def\BBymin{-1.25}
      \def\BBymax{+2}
      \def\BBzmin{-0.25}
      \def\BBzmax{+2.5}
      \begin{scope}[tdplot_main_coords, scale=0.9]
        \node (xaxis) at ({\BBxmax+0.1*(\BBxmax-\BBxmin)}, 0, 0) {$\xvi$};%
        \node (yaxis) at (0, {1+\BBymax+0.75*(\BBymax-\BBymin)}, 0) {$\zvi_1$};%
        %
        %
        %
        \fill[secondary!80, opacity=0.5] (\BBxmin, 0, 0) -- +(-\BBxmin, 0, 0) -- +(-\BBxmin, \coteCyl, 0) -- +(0, \coteCyl, 0) -- cycle;
        \fill[secondary!90, opacity=0.625] (\BBxmin, \coteCyl, 0) -- +(-\BBxmin, 0, 0) -- +(-\BBxmin, 0, \coteCyl) -- +(0, 0, \coteCyl) -- cycle;
        \draw[dashed, white] (\BBxmin, \coteCyl, 0) -- (\BBxmax, \coteCyl, 0);
        \node[secondary!70!black] (Omega) at (\BBxmin-1.25, 0.5*\coteCyl, 1.5*\coteCyl) {$\Omegaperper^-$};
        \draw[-latex, secondary!80!black] (Omega) to[bend right=40] (\BBxmin+1, 0.5*\coteCyl, 0.5*\coteCyl);
        \filldraw[draw=white, fill=secondary!90, opacity=0.625] (\BBxmin, 0, 0) -- +(-\BBxmin, 0, 0) -- +(-\BBxmin, 0, \coteCyl) -- +(0, 0, \coteCyl) -- cycle;
        \filldraw[draw=white, fill=secondary!80, opacity=0.5] (\BBxmin, 0, \coteCyl) -- +(-\BBxmin, 0, 0) -- +(-\BBxmin, \coteCyl, 0) -- +(0, \coteCyl, 0) -- cycle;
        \filldraw[draw=white, fill=secondary!95!black, opacity=0.875] (0, 0, 0) -- (0, \coteCyl, 0) -- (0, \coteCyl, \coteCyl) -- (0, 0, \coteCyl) -- cycle;
        %
        \draw[-latex, opacity=0.8] (0, {\BBymin-0.1*(\BBymax-\BBymin)}, 0) -- (yaxis);%
        \fill[tertiary!80, opacity=0.5] (\BBxmax, 0, 0) -- +(-\BBxmax, 0, 0) -- +(-\BBxmax, \coteCyl, 0) -- +(0, \coteCyl, 0) -- cycle;
        \fill[tertiary!90, opacity=0.625] (\BBxmax, \coteCyl, 0) -- +(-\BBxmax, 0, 0) -- +(-\BBxmax, 0, \coteCyl) -- +(0, 0, \coteCyl) -- cycle;
        \draw[dashed, white] (\BBxmax, \coteCyl, 0) -- (\BBxmax, \coteCyl, 0);
        \node[tertiary!70!black] (Omega) at (\BBxmax+1.25, 0.5*\coteCyl, 1.5*\coteCyl) {$\Omegaperper^+$};
        \draw[-latex, tertiary!80!black] (Omega) to[bend left=40] (\BBxmax-1, 0.5*\coteCyl, 0.5*\coteCyl);
        \filldraw[draw=white, fill=tertiary!90, opacity=0.625] (\BBxmax, 0, 0) -- +(-\BBxmax, 0, 0) -- +(-\BBxmax, 0, \coteCyl) -- +(0, 0, \coteCyl) -- cycle;
        \filldraw[draw=white, fill=tertiary!80, opacity=0.5] (\BBxmax, 0, \coteCyl) -- +(-\BBxmax, 0, 0) -- +(-\BBxmax, \coteCyl, 0) -- +(0, \coteCyl, 0) -- cycle;
        \draw (1, 0.5, 2.5) node {$\linapp{\widehat{U}_\fbvar}{\widehat{G}_\fbvar}$};
        \draw[-latex] (0, 0, {\BBzmin-0.25*(\BBzmax-\BBzmin)}) -- (0, 0, {\BBzmax+0.25*(\BBzmax-\BBzmin)}) node[above] {$\smash{\zvi_2}$};
        \draw[-latex] ({\BBxmin-0.1*(\BBxmax-\BBxmin)}, 0, 0) -- (xaxis);%
      \end{scope}%
      \begin{scope}[shift={(11cm, 0.5cm)}, tdplot_main_coords, scale=0.9]
        %
        %
        \def\ecart{1}
        %
        \fill[secondary!80, opacity=0.5] (\BBxmin-\ecart, 0, 0) -- +(-\BBxmin, 0, 0) -- +(-\BBxmin, \coteCyl, 0) -- +(0, \coteCyl, 0) -- cycle;
        \fill[secondary!90, opacity=0.625] (\BBxmin-\ecart, \coteCyl, 0) -- +(-\BBxmin, 0, 0) -- +(-\BBxmin, 0, \coteCyl) -- +(0, 0, \coteCyl) -- cycle;
        \draw[dashed, white] (\BBxmin-\ecart, \coteCyl, 0) -- +(-\BBxmin, 0, 0);
        \draw[dashed, white] (\BBxmin+\ecart, \coteCyl, 0) -- +(-\BBxmin, 0, 0);
        \node[secondary!70!black] (Omega) at (\BBxmin-1.25-\ecart, 0.5*\coteCyl, 1.5*\coteCyl) {$\linapp{\widehat{U}_\fbvar^-}{\Phi_\fbvar}$};
        \draw[-latex, secondary!80!black] (Omega) to[bend right=40] (\BBxmin+1.5-\ecart, 0.5*\coteCyl, 0.5*\coteCyl);
        \filldraw[draw=white, fill=secondary!90, opacity=0.625] (\BBxmin-\ecart, 0, 0) -- +(-\BBxmin, 0, 0) -- +(-\BBxmin, 0, \coteCyl) -- +(0, 0, \coteCyl) -- cycle;
        \filldraw[draw=white, fill=secondary!80, opacity=0.5] (\BBxmin-\ecart, 0, \coteCyl) -- +(-\BBxmin, 0, 0) -- +(-\BBxmin, \coteCyl, 0) -- +(0, \coteCyl, 0) -- cycle;
        \filldraw[draw=white, fill=teal!95, opacity=0.875] (-\ecart, 0, 0) -- +(0, \coteCyl, 0) -- +(0, \coteCyl, \coteCyl) -- +(0, 0, \coteCyl) -- cycle;
        \node [teal!60!black] (Phi) at (0, 0.5*\coteCyl, 2.5*\coteCyl) {$\Phi_\fbvar$, solution of \eqref{H2Dmod:eq:transmission_jump}};
        \draw[-latex, teal!70!black] (Phi) to[bend left=15] (-\ecart, 0.5*\coteCyl, 0.5*\coteCyl);
        \draw[-latex, teal!70!black] (Phi) to[bend right=15] ( \ecart, 0.5*\coteCyl, 0.5*\coteCyl);
        \filldraw[draw=white, fill=teal!95, opacity=0.875] (\ecart, 0, 0) -- +(0, \coteCyl, 0) -- +(0, \coteCyl, \coteCyl) -- +(0, 0, \coteCyl) -- cycle;
        \fill[tertiary!80, opacity=0.5] (\BBxmax+\ecart, 0, 0) -- +(-\BBxmax, 0, 0) -- +(-\BBxmax, \coteCyl, 0) -- +(0, \coteCyl, 0) -- cycle;
        \fill[tertiary!90, opacity=0.625] (\BBxmax+\ecart, \coteCyl, 0) -- +(-\BBxmax, 0, 0) -- +(-\BBxmax, 0, \coteCyl) -- +(0, 0, \coteCyl) -- cycle;
        \draw[dashed, white] (\BBxmax+\ecart, \coteCyl, 0) -- +(-\BBxmax, 0, 0);
        \node[tertiary!70!black] (Omega) at (\BBxmax+1.25+\ecart, 0.5*\coteCyl, 1.5*\coteCyl) {$\linapp{\widehat{U}_\fbvar^+}{\Phi_\fbvar}$};
        \draw[-latex, tertiary!80!black] (Omega) to[bend left=40] (\BBxmax-1.5+\ecart, 0.5*\coteCyl, 0.5*\coteCyl);
        \filldraw[draw=white, fill=tertiary!90, opacity=0.625] (\BBxmax+\ecart, 0, 0) -- +(-\BBxmax, 0, 0) -- +(-\BBxmax, 0, \coteCyl) -- +(0, 0, \coteCyl) -- cycle;
        \filldraw[draw=white, fill=tertiary!80, opacity=0.5] (\BBxmax+\ecart, 0, \coteCyl) -- +(-\BBxmax, 0, 0) -- +(-\BBxmax, \coteCyl, 0) -- +(0, \coteCyl, 0) -- cycle;
        %
        %
        %
        %
      \end{scope}%
    \end{tikzpicture} 
  }
  \caption{Expression \eqref{H2Dmod:eq:Uhat_from_Uplusminushat} of the waveguide solution $\widehat{U}_\fbvar$ with respect to $\linapp{\widehat{U}^\pm_\fbvar}{\Phi_\fbvar}$.\label{H2Dmod:fig:Uhat_from_Uplusminushat}}
\end{figure}

\subsection{The DtN approach for the half-guide problems} \label{H2Dmod:sec:resolution_half_guide_problems}
Our objective is to solve the half-guide problems \eqref{H2Dmod:eq:augmented_halfguide_problem} defined in $\Omegaperper^+$ and $\Omegaperper^-$, and to compute the DtN operators $\Lambdahat^+_\fbvar$ and $\Lambdahat^-_\fbvar$. Since these problems are similar to each other, we will restrict ourselves to the half-guide problem \textnormal{(\hyperref[H2Dmod:eq:augmented_halfguide_problem]{$\mathscr{P}^+_\perperiodic$})} set in $\Omegaperper^+$, and the computation of $\smash{\Lambdahat^+_\fbvar}$. Exploiting the periodicity of $\aten^+_p$ and $\rho^+_p$ in the $\eva_\xvi$--direction, we resort to the method developed in \cite{fliss2009exact, joly2006exact} for the elliptic Helmholtz equation $-\transp{\bsnabla}\, \aten_p\, \bsnabla\; U - \rho\, \omega^2\; U = 0$, and in \cite{amenoagbadji2023wave} for the non-elliptic $2$D case.

\vspace{1\baselineskip} \noindent
We also introduce some additional notation:
\begin{equation}\label{H2Dmod:eq:cells}
  \displaystyle
  \cellperperi{0} := (0, 1)^3 \quad \textnormal{and} \quad \cellperperi{n} := \cellperperi{0} + n\vts \eva_\xvi \quad \spforall n \in \N, \quad \textnormal{so that} \quad \overline{\Omegaperper^+} = \bigcup_{n \in \N} \overline{\cellperperi{n}}.
\end{equation}
The interface between the cells $\cellperperi{n}$ and $\cellperperi{n + 1}$ is denoted by $\Sigmaperperi{n} := \Sigmaperper + n\vts \eva_\xvi$. These domains are represented in Figure \ref{H2Dmod:fig:cells}. By periodicity, one can identify each cell $\cellperperi{n}$ with the periodicity cell $\cellperperi{0}$ denoted $\cellperper$, and each interface $\Sigmaperperi{n}$ with the interface $\Sigmaperperi{0}$ denoted $\Sigmaperper$.

\begin{figure}[ht!]
  \noindent\makebox[\textwidth]{%
    \def\coteCyl{2.5}
    \begin{tikzpicture}[scale=0.675]
      \def\angleinterface{55}%
      \definecolor{coplancoupe}{RGB}{211, 227, 65}
      \tdplotsetmaincoords{75}{30}
      \def\BBxmin{0}
      \def\BBxmax{12}
      \def\BBymin{0}
      \def\BBymax{+4.75}
      \def\BBzmin{0}
      \def\BBzmax{+3}
      \def\numCellsminusone{4}
      \begin{scope}[tdplot_main_coords, scale=0.9]
        \node (xaxis) at ({\BBxmax+0.1*(\BBxmax-\BBxmin)}, 0, 0) {$\xvi$};%
        \node (yaxis) at (0, {\BBymax+0.75*(\BBymax-\BBymin)}, 0) {$\zvi_1$};%
        %
        %
        %
        \draw[-latex] (0, {\BBymin-0.1*(\BBymax-\BBymin)}, 0) -- (yaxis);%
        \fill[tertiary!80, opacity=0.5] (\BBxmax, 0, 0) -- +(-\BBxmax, 0, 0) -- +(-\BBxmax, \coteCyl, 0) -- +(0, \coteCyl, 0) -- cycle;
        \fill[tertiary!90, opacity=0.625] (\BBxmax, \coteCyl, 0) -- +(-\BBxmax, 0, 0) -- +(-\BBxmax, 0, \coteCyl) -- +(0, 0, \coteCyl) -- cycle;
        \draw[white, dashed, opacity=1] (0, \coteCyl, 0) -- +(\BBxmax, 0, 0);
        \foreach \xcell in {0, ..., \numCellsminusone} {
          \fill[draw=white, dashed, fill=teal!95!black, opacity=0.875] (\xcell*\coteCyl, 0, 0) -- +(0, \coteCyl, 0) -- +(0, \coteCyl, \coteCyl) -- +(0, 0, \coteCyl) -- cycle;
        }
        
        \foreach \xcell in {0, ..., 3} {
          \draw[-latex, dashed] (\xcell*\coteCyl+0.5*\coteCyl, 0.5*\coteCyl, \coteCyl) -- (\xcell*\coteCyl+0.5*\coteCyl, 0.5*\coteCyl, 0.5*\coteCyl);
        }
  
        %
        %
        %
        \fill[fill=tertiary!90, opacity=0.8] (\BBxmax, 0, 0) -- +(-\BBxmax, 0, 0) -- +(-\BBxmax, 0, \coteCyl) -- +(0, 0, \coteCyl) -- cycle;
        \fill[fill=tertiary!80, opacity=0.8] (\BBxmax, 0, \coteCyl) -- +(-\BBxmax, 0, 0) -- +(-\BBxmax, \coteCyl, 0) -- +(0, \coteCyl, 0) -- cycle;
        \foreach \xcell in {0, ..., 3} {
          \draw[white, opacity=0.5, thin] (\xcell*\coteCyl, 0, 0) -- +(\coteCyl, 0, 0) -- +(\coteCyl, 0, \coteCyl) -- +(0, 0, \coteCyl) -- cycle;
          \draw[white, opacity=0.5, thin] (\xcell*\coteCyl, 0, \coteCyl) -- +(0, \coteCyl, 0);
          \node (cellperper) at (\xcell*\coteCyl+0.5*\coteCyl, 0.5*\coteCyl, 4) {$\cellperperi{\xcell}$};
          \draw (cellperper) -- (\xcell*\coteCyl+0.5*\coteCyl, 0.5*\coteCyl, \coteCyl);
        }

        \foreach \xcell in {0, ..., \numCellsminusone} {
          \node (Sigmaperper) at (\xcell*\coteCyl, -5, 0) {$\Sigmaperperi{\xcell}$};
          \draw[-latex] (Sigmaperper) to[bend left=40] (\xcell*\coteCyl, 0, 0.5*\coteCyl);
        }
        \draw[white, opacity=0.5, thin] (\numCellsminusone*\coteCyl, \coteCyl, \coteCyl) -- +(0, -\coteCyl, 0) -- +(\BBxmax-\numCellsminusone*\coteCyl, -\coteCyl, 0);

        \draw[-latex] (0, 0, {\BBzmin-0.05*(\BBzmax-\BBzmin)}) -- (0, 0, {\BBzmax+0.25*(\BBzmax-\BBzmin)}) node[above] {$\smash{\zvi_2}$};
        \draw[-latex] ({\BBxmin-0.1*(\BBxmax-\BBxmin)}, 0, 0) -- (xaxis);%
      \end{scope}%
    \end{tikzpicture} 
  }
  \caption{The cells $\cellperperi{n}$ and the interfaces $\Sigmaperperi{n}$ in \eqref{H2Dmod:eq:cells}\label{H2Dmod:fig:cells}}
\end{figure}

\vspace{1\baselineskip} \noindent
Note that for any $n \in \N$, $\cellperperi{n}$ and $\Sigmaperperi{\tau}$, $\tau \in \{n, n+1\}$, correspond respectively to $\mathbb{\Omega}_\perperiodic := I_\xvi \times I_1 \times (0, 1)$ and $\mathbb{\Sigma}^\tau_\perperiodic := \{\tau\} \times I_1 \times (0, 1)$ defined by \eqref{H2Dmod:eq:Qperper_Sigmaperper} with
\[
  I_\xvi := (n, n + 1), \quad I_1 := (0, 1),
\]
{\setstretch{1.0}%
Thus Section \ref{H2Dmod:sec:e1_e2_periodic_functions} enables us to use the space $\Hgradperper{\cellperperi{n}}$ given by \eqref{H2Dmod:eq:def_H1Cperper_HdivCperper}, as well as $\Honehalfperper{\Sigmaperperi{n}}$, the space defined by \eqref{H2Dmod:eq:def_H12Cperper}, its dual $\smash{\Hminusonehalfperper{\Sigmaperperi{n}}}$, and the trace and the normal trace applications on $\Sigmaperperi{n}$.} In the sequel, we will systematically use the obvious identifications $\Hgradperper{\cellperperi{n}} \equiv \Hgradperper{\cellperper}$ and ${\Hplusminusonehalfperper{\Sigmaperperi{n}} \equiv \Hplusminusonehalfperper{\Sigmaperper}}$, even when not mentioned.

\subsubsection{Structure of the half-guide solution}\label{H2Dmod:sec:structure_half_guide_solution}
Consider the operator $\mathcal{P}_\fbvar \in \mathscr{L}(\Honehalfperper{\Sigmaperper})$:
\begin{equation}\label{H2Dmod:eq:def_propagation_operator}
  \displaystyle
  \spforall \Phi \in \Honehalfperper{\Sigmaperper}, \quad \mathcal{P}_\fbvar \vts \Phi := \linapp{\widehat{U}^+_\fbvar}{\Phi}|_{\Sigmaperperi{1}},
\end{equation}
where the boundedness property stems from the well-posedness of \textnormal{(\hyperref[H2Dmod:eq:augmented_halfguide_problem]{$\mathscr{P}^+_\perperiodic$})} and from the continuity of the trace map on $\Sigmaperperi{1}$ from $\Hgradperper{\Omegaperper^+}$ to $\smash{\Honehalfperper{\Sigmaperperi{1}} \equiv \Honehalfperper{\Sigmaperper}}$. The operator $\mathcal{P}_\fbvar$ is called the \emph{propagation operator}, because it determines how the half-guide solution propagates from one interface to the other, as the next result shows. 
\begin{prop}\label{H2Dmod:prop:structure_half_guide}
  For any $\Phi \in \Honehalfperper{\Sigmaperper}$, the solution $\linapp{\widehat{U}^+_\fbvar}{\Phi}$ of Problem \textnormal{(\hyperref[H2Dmod:eq:augmented_halfguide_problem]{$\mathscr{P}^+_\perperiodic$})} satisfies
  \begin{equation}\label{H2Dmod:eq:structure_half_guide}
    \displaystyle
    \spforall n \in \N, \quad \aeforall \xva \in \Omegaperper, \quad \linapp{\widehat{U}^+_\fbvar}{\Phi}\vts (\xva + n\vts \eva_\xvi) = \linapp{\widehat{U}^+_\fbvar}{(\mathcal{P}_\fbvar)^n \vts \Phi}\vts (\xva),
  \end{equation}
  so that 
  \begin{equation}
    \displaystyle
    \spforall n \in \N, \quad \linapp{\widehat{U}^+_\fbvar}{\Phi}|_{\Sigmaperperi{n}} = (\mathcal{P}_\fbvar)^n\, \Phi.
  \end{equation}
  Furthermore, $\mathcal{P}_\fbvar$ has a spectral radius denoted by $\rho(\mathcal{P}_\fbvar)$, which is strictly less than $1$.
\end{prop}

\begin{dem}
  The proof is a direct adaptation of \cite[Theorem 3.1]{joly2006exact}. To prove \eqref{H2Dmod:eq:structure_half_guide}, we begin with $n = 1$ and we define $\linapp{\widetilde{U}}{\Phi} := \linapp{\widehat{U}^+_\fbvar}{\Phi}\vts (\cdot + \eva_\xvi)$. Then one has $\linapp{\widetilde{U}}{\Phi}|_{\Sigmaperper} = \mathcal{P}_\fbvar \vts \Phi$. Moreover, by using the change of variables $\xva \mapsto \xva + \eva_1$ in \textnormal{(\hyperref[H2Dmod:eq:augmented_halfguide_problem]{$\mathscr{P}^+_\perperiodic$})} and the periodicity of $\aten^+_p$ and $\rho^+_p$ along $\eva_\xvi$, one obtains that $\linapp{\widetilde{U}}{\Phi}$ satisfies the same problem as $\smash{\linapp{\widehat{U}^+_\fbvar}{\mathcal{P}_\fbvar \vts \Phi}}$. Consequently, the uniqueness of \textnormal{(\hyperref[H2Dmod:eq:augmented_halfguide_problem]{$\mathscr{P}^+_\perperiodic$})} leads to \eqref{H2Dmod:eq:structure_half_guide} for $n = 1$. The extension to $n > 1$ follows by induction.

  \vspace{1\baselineskip} \noindent
  To prove that $\rho(\mathcal{P}_\fbvar) < 1$, we rely on the following estimate, \surligner{which describes the exponential decay of the half-guide solution $\widehat{U}^+_\fbvar$ in the $\eva_\xvi$--direction (see Remark \ref{H2D:rmk:proof_exponential_decay}):}
  \begin{equation}\label{H2Dmod:eq:exponential_decay_half_guide_solution}
    \spexists c,\, \alpha > 0, \quad \big\|\linapp{\widehat{U}^+_\fbvar}{\Phi} \, \exp(\alpha \Imag \omega\, |\xvi|)\big\|_{\Hgradperper{\Omegaperper^+}} \leq c\; \|\Phi\|_{\Honehalfperper{\Sigmaperper}}.
  \end{equation}
  Since $\smash{(\mathcal{P}_\fbvar)^n \vts \Phi = \linapp{\widehat{U}^+_\fbvar}{\Phi}|_{\Sigmaperperi{n}}}$ by \eqref{H2Dmod:eq:structure_half_guide}, the continuity of the trace \surligner{map on $\Sigmaperperi{n}$ and \eqref{H2Dmod:eq:exponential_decay_half_guide_solution} lead to}
  \[
    \displaystyle
    \spexists c', \alpha > 0, \quad \spforall n \in \N, \quad \|(\mathcal{P}_\fbvar)^n\| \leq c'\; \euler^{- \alpha \Imag \omega n}.
  \]
  One then concludes using Gelfand's formula $\displaystyle \rho(\mathcal{P}_\fbvar) = \lim_{n \to +\infty}\|(\mathcal{P}_\fbvar)^n\|^{1/n}$.
\end{dem}

\begin{rmk}\label{H2D:rmk:proof_exponential_decay}
  \surligner{The proof of \eqref{H2Dmod:eq:exponential_decay_half_guide_solution} is very similar to that of \eqref{H2Dmod:eq:exp_decay}: the idea is to examine the problem satisfied by the product $\widehat{U}^+_\fbvar\, \exp(\beta\, |\xvi|)$ for any $\beta \in \R$, and to identify the values of $\beta$ for which this problem is well-posed in $\Hgradperper{\Omegaperper^+}$ via coercivity.}
\end{rmk}

\vspace{1\baselineskip} \noindent
Proposition \ref{H2Dmod:prop:structure_half_guide} shows that the restrictions of $\linapp{\widehat{U}^+_\fbvar}{\Phi}$ to the interfaces $\Sigmaperperi{n}$ can be fully expressed with respect to the propagation operator $\mathcal{P}_\fbvar$. Therefore, knowing $\mathcal{P}_\fbvar$, one can construct $\smash{\linapp{\widehat{U}^+_\fbvar}{\Phi}}$ using solutions of problems defined in one periodicity cell, as shown in the next section.

\subsubsection{Local cell problems}
Given boundary data $\Phi \in \Honehalfperper{\Sigmaperper}$ and $j \in \{0, 1\}$, let us introduce the local Dirichlet cell problems: \textit{Find $E^j_\fbvar \in H^1_\cutmat (\cellperper)$ such that}
\begin{subequations}\label{H2Dmod:eq:local_cell_problems}
  \begin{equation}
    \label{H2Dmod:eq:local_cell_problem}
    \left\{
      \begin{array}{r@{\ =\ }l}
        \displaystyle- \divecutFB\, \aten_p \!  \gradcutFB\, E^j_\fbvar - \rho_p\, \omega^2\, E^j_\fbvar & 0 \quad \textnormal{in}\ \  \cellperper,
        \ret
        \multicolumn{2}{c}{\displaystyle E^j_\fbvar \in \Hgradperper{\cellperper}, \quad \aten_p \! \gradcutFB\, E^j_\fbvar \in \Hdiveperper{\cellperper},}
      \end{array}
    \right.
  \end{equation}
  combined with the Dirichlet boundary conditions (see Figure \ref{H2Dmod:fig:local_cell_Dirichlet})
  \begin{equation}\label{H2Dmod:eq:BC_local_cell_problems}
    \left\{
    \begin{array}{r@{\ =\ }l@{\quad \textnormal{and} \quad}r@{\ =\ }l}
      \linapp{E^0_\fbvar}{\Phi}|_{\Sigmaperperi{0}} & \Phi & \linapp{E^0_\fbvar}{\Phi}|_{\Sigmaperperi{1}} & 0
      \\[10pt]
      \linapp{E^1_\fbvar}{\Phi}|_{\Sigmaperperi{0}} & 0 & \linapp{E^1_\fbvar}{\Phi}|_{\Sigmaperperi{1}} & \Phi.
    \end{array}
    \right.
  \end{equation}
\end{subequations}
\begin{figure}[H]
  \noindent\makebox[\textwidth]{%
    \def\coteCyl{2}
    \begin{tikzpicture}[scale=0.825]
      \def\angleinterface{55}%
      \definecolor{coplancoupe}{RGB}{211, 227, 65}
      \tdplotsetmaincoords{80}{30}
      \def\axesrapport{3}
      \def\ecart{4.5}
      \begin{scope}[tdplot_main_coords, scale=0.9]
        \draw[latex-, canvas is zx plane at y=0.5*\coteCyl, dashed] (0.75*\coteCyl, -\ecart) arc[radius = 0.5*\coteCyl, start angle=180, end angle=270] node[above] {$\Phi$};
        \draw[-latex, canvas is zx plane at y=0.5*\coteCyl, tertiary!90!black] (0.25*\coteCyl, -\ecart) arc[radius = 0.5*\coteCyl, start angle=360, end angle=270] node[below] {$\mathcal{T}^{00}_\fbvar \Phi$};
        \fill[tertiary!80, opacity=0.5] (\coteCyl-\ecart, 0, 0) -- +(-\coteCyl, 0, 0) -- +(-\coteCyl, \coteCyl, 0) -- +(0, \coteCyl, 0) -- cycle;
        \fill[tertiary!90, opacity=0.625] (\coteCyl-\ecart, \coteCyl, 0) -- +(-\coteCyl, 0, 0) -- +(-\coteCyl, 0, \coteCyl) -- +(0, 0, \coteCyl) -- cycle;
        \draw[white, dashed, opacity=1] (-\ecart, \coteCyl, 0) -- +(\coteCyl, 0, 0);
        \fill[draw=white, dashed, fill=tertiary!95!black, opacity=0.75] (-\ecart, 0, 0) -- +(0, \coteCyl, 0) -- +(0, \coteCyl, \coteCyl) -- +(0, 0, \coteCyl) -- cycle;
        \fill[draw=white, fill=tertiary!95!black, opacity=0.875] (-\ecart+\coteCyl, 0, 0) -- +(0, \coteCyl, 0) -- +(0, \coteCyl, \coteCyl) -- +(0, 0, \coteCyl) -- cycle;
        \filldraw[draw=white, fill=tertiary!90, opacity=0.8] (\coteCyl-\ecart, 0, 0) -- +(-\coteCyl, 0, 0) -- +(-\coteCyl, 0, \coteCyl) -- +(0, 0, \coteCyl) -- cycle;
        \fill[fill=tertiary!80, opacity=0.8] (\coteCyl-\ecart, 0, \coteCyl) -- +(-\coteCyl, 0, 0) -- +(-\coteCyl, \coteCyl, 0) -- +(0, \coteCyl, 0) -- cycle;
        \draw[latex-, canvas is zx plane at y=0.5*\coteCyl, dashed] (0.75*\coteCyl, \coteCyl-\ecart) arc[radius = 0.5*\coteCyl, start angle=180, end angle=90] node[above] {$0$};
        \draw[-latex, canvas is zx plane at y=0.5*\coteCyl, tertiary!90!black] (0.25*\coteCyl, \coteCyl-\ecart) arc[radius = 0.5*\coteCyl, start angle=0, end angle=90] node[below] {$\mathcal{T}^{01}_\fbvar \Phi$};
        \draw (0.5*\coteCyl-\ecart, 0.5*\coteCyl, 1.5*\coteCyl) node {$\linapp{E^0_\fbvar}{\Phi}$};
        \node (xaxis) at (0.5*\axesrapport, \axesrapport, 0) {$\xvi$};
        \node (yaxis) at (0, 1.9*\axesrapport, 0) {$\zvi_1$};
        \node (zaxis) at (0, \axesrapport, 0.5*\axesrapport) {$\zvi_2$};
        \draw[-latex] (-0.1*\axesrapport, \axesrapport, 0) -- (xaxis);
        \draw[-latex] (0, 0.9*\axesrapport, 0) -- (yaxis);
        \draw[-latex] (0, \axesrapport, -0.1*\axesrapport) -- (zaxis);
        \draw[latex-, canvas is zx plane at y=0.5*\coteCyl, dashed] (0.75*\coteCyl, \ecart) arc[radius = 0.675*\coteCyl, start angle=180, end angle=270, dashed] node[above] {$0$};
        \draw[-latex, canvas is zx plane at y=0.5*\coteCyl, tertiary!90!black] (0.25*\coteCyl, \ecart) arc[radius = 0.5*\coteCyl, start angle=360, end angle=270] node[below] {$\mathcal{T}^{10}_\fbvar \Phi$};
        \fill[tertiary!80, opacity=0.5] (\coteCyl+\ecart, 0, 0) -- +(-\coteCyl, 0, 0) -- +(-\coteCyl, \coteCyl, 0) -- +(0, \coteCyl, 0) -- cycle;
        \fill[tertiary!90, opacity=0.625] (\coteCyl+\ecart, \coteCyl, 0) -- +(-\coteCyl, 0, 0) -- +(-\coteCyl, 0, \coteCyl) -- +(0, 0, \coteCyl) -- cycle;
        \draw[white, dashed, opacity=1] (\ecart, \coteCyl, 0) -- +(\coteCyl, 0, 0);
        \fill[draw=white, dashed, fill=tertiary!95!black, opacity=0.75] (\ecart, 0, 0) -- +(0, \coteCyl, 0) -- +(0, \coteCyl, \coteCyl) -- +(0, 0, \coteCyl) -- cycle;
        \fill[draw=white, fill=tertiary!95!black, opacity=0.875] (\ecart+\coteCyl, 0, 0) -- +(0, \coteCyl, 0) -- +(0, \coteCyl, \coteCyl) -- +(0, 0, \coteCyl) -- cycle;
        \filldraw[draw=white, fill=tertiary!90, opacity=0.8] (\coteCyl+\ecart, 0, 0) -- +(-\coteCyl, 0, 0) -- +(-\coteCyl, 0, \coteCyl) -- +(0, 0, \coteCyl) -- cycle;
        \fill[fill=tertiary!80, opacity=0.8] (\coteCyl+\ecart, 0, \coteCyl) -- +(-\coteCyl, 0, 0) -- +(-\coteCyl, \coteCyl, 0) -- +(0, \coteCyl, 0) -- cycle;
        \draw[latex-, canvas is zx plane at y=0.5*\coteCyl, dashed] (0.75*\coteCyl, \coteCyl+\ecart) arc[radius = 0.675*\coteCyl, start angle=180, end angle=90] node[above] {$\Phi$};
        \draw[-latex, canvas is zx plane at y=0.5*\coteCyl, tertiary!90!black] (0.25*\coteCyl, \coteCyl+\ecart) arc[radius = 0.5*\coteCyl, start angle=0, end angle=90] node[below] {$\mathcal{T}^{11}_\fbvar \Phi$};
        \draw (0.5*\coteCyl+\ecart, 0.5*\coteCyl, 1.5*\coteCyl) node {$\linapp{E^1_\fbvar}{\Phi}$};
      \end{scope}%
    \end{tikzpicture} 
  }
  \caption{Dirichlet conditions satisfied by $E^0_\fbvar$ and $E^1_\fbvar$, and the DtN operators $\mathcal{T}^{j \ell}_\fbvar$ given by \eqref{H2Dmod:eq:def_local_DtN_operator}. \label{H2Dmod:fig:local_cell_Dirichlet}}
\end{figure}

\noindent
These problems are well-posed thanks to Lax-Milgram's theorem combined with a lifting argument. Moreover, using the structure of $\linapp{\widehat{U}^+_\fbvar}{\Phi}$ given by Proposition \ref{H2Dmod:prop:structure_half_guide}, it follows by linearity that
\begin{equation}\label{H2Dmod:eq:halfguide_from_local_cells}
  \displaystyle
  \spforall n \in \N, \quad \linapp{\widehat{U}^+_\fbvar}{\Phi} (\cdot + n\vts \eva_\xvi) |_{\cellperper} = \linapp{E^0_\fbvar}{(\mathcal{P}_\fbvar)^n \, \Phi} + \linapp{E^1_\fbvar}{(\mathcal{P}_\fbvar)^{n + 1} \, \Phi}.
\end{equation}
Therefore, if $\mathcal{P}_\fbvar$ is known, $\linapp{\widehat{U}^+_\fbvar}{\Phi}$ can be constructed cell by cell using the solutions $\linapp{E^j_\fbvar}{\Phi}$ of the local cell problems \eqref{H2Dmod:eq:local_cell_problems}. On the other hand, these local cell problems are involved in the characterization of $\mathcal{P}_\fbvar$. This is the object of the next section.

\subsubsection{Characterization of the propagation operator via a Riccati equation}
The goal of this section is to characterize $\mathcal{P}_\fbvar$ using the \emph{local DtN operators} defined for any $j, \ell \in \{0, 1\}$ by
\begin{equation}\label{H2Dmod:eq:def_local_DtN_operator}
  \displaystyle
  \spforall \Phi \in \Honehalfperper{\Sigmaperper}, \quad \mathcal{T}^{j \ell}_\fbvar\, \Phi = (-1)^{\ell + 1}\, \big(\cutmat\vts \aten_p \! \gradcutFB \, \, \linapp{E^j_\fbvar}{\Phi} \cdot \eva_\xvi\big) |_{\Sigmaperperi{\ell}} \in \Hminusonehalfperper{\Sigmaperper},
\end{equation}
By applying the second item of Proposition \ref{H2Dmod:prop:formule_sauts} to $\cellperperi{0}$ and $\cellperperi{1}$ in the cell by cell expression \eqref{H2Dmod:eq:halfguide_from_local_cells} of $\widehat{U}^+_\fbvar$, it follows that $\aten_p \! \gradcutFB\, \widehat{U}^+_\fbvar \cdot \eva_\xvi$ is continuous accross the interface $\Sigmaperperi{1}$, that is,
\[
  \displaystyle
  \spforall \Phi \in \Honehalfperper{\Sigmaperper}, \quad (\aten_p \! \gradcutFB\, \linapp{\widehat{U}^+_\fbvar}{\Phi} \cdot \eva_\xvi)|_{\Sigmaperperi{1}} = (\aten_p \! \gradcutFB\, \linapp{\widehat{U}^+_\fbvar}{\mathcal{P}_\fbvar \Phi} \cdot \eva_\xvi )|_{\Sigmaperperi{0}},
\]
or equivalently, for any $\Phi \in \Honehalfperper{\Sigmaperper}$,
\begin{multline}\label{H2Dmod:eq:cty_normal_trace_U_halfguide}
  \displaystyle
  \big(\aten_p \! \gradcutFB\, \linapp{E^0_\fbvar}{\Phi} \cdot \eva_\xvi \big)|_{\Sigmaperperi{1}} + \big(\aten_p \! \gradcutFB\, \linapp{E^1_\fbvar}{\mathcal{P}_\fbvar\vts \Phi} \cdot \eva_\xvi \big)|_{\Sigmaperperi{1}} 
  \ret
  = \big(\aten_p \! \gradcutFB\, \linapp{E^0_\fbvar}{\mathcal{P}_\fbvar\vts \Phi} \cdot \eva_\xvi \big)|_{\Sigmaperperi{0}} + \big(\aten_p \! \gradcutFB\, \linapp{E^1_\fbvar}{(\mathcal{P}_\fbvar)^2\vts \Phi} \cdot \eva_\xvi \big)|_{\Sigmaperperi{0}}.
\end{multline}
This leads to a \emph{Riccati} equation, which characterizes uniquely $\mathcal{P}_\fbvar$ as stated by the next result.
\begin{prop}
  The propagation operator $\mathcal{P}_\fbvar$ defined by \eqref{H2Dmod:eq:def_propagation_operator} is the unique solution of the problem
  \begin{equation}\label{H2Dmod:eq:riccati}
    \left|
    \begin{array}{l}
      \textit{Find $P \in \mathscr{L}(\Honehalfperper{\Sigmaperper})$ such that $\rho(P) < 1$ and}
      \\[10pt]
      \multicolumn{1}{c}{\displaystyle \mathcal{T}^{10}_\fbvar\; P^2 + (\mathcal{T}^{00}_\fbvar + \mathcal{T}^{11}_\fbvar)\; P^{\vphantom{0}} + \mathcal{T}^{01}_\fbvar = 0.}
    \end{array}
    \right.
  \end{equation}
\end{prop}

\begin{dem}
  We only present the outline of the proof which is detailed in \cite[Theorem 4.1]{joly2006exact}. The existence is directly deduced from \eqref{H2Dmod:eq:cty_normal_trace_U_halfguide}. On the other hand, the uniqueness is achieved by considering an operator $\widetilde{\mathcal{P}}$ which satisfies \eqref{H2Dmod:eq:riccati}, and by showing that for any $\Phi \in \Honehalfperper{\Sigmaperper}$, the function defined cell by cell by $\linapp{\widetilde{U}}{\Phi} (\cdot + n\vts \eva_\xvi) |_{\cellperper} := \linapp{E^0_\fbvar}{(\widetilde{\mathcal{P}})^n \, \Phi} + \linapp{E^1_\fbvar}{(\widetilde{\mathcal{P}})^{n + 1} \, \Phi}$ satisfies the same half-guide problem \textnormal{(\hyperref[H2Dmod:eq:augmented_halfguide_problem]{$\mathscr{P}^+_\perperiodic$})} as $\linapp{\widehat{U}^+_\fbvar}{\Phi}$, so that both solutions coincide in the half-guide $\Omegaperper^+$, and in particular on the interface $\Sigmaperper^1$: $\widetilde{\mathcal{P}}\, \Phi = \mathcal{P}_\fbvar\, \Phi$.
\end{dem}

\vspace{1\baselineskip} \noindent
One important point to note from the practical point of view is that the Riccati equation \eqref{H2Dmod:eq:riccati} only involves the local DtN operators $\mathcal{T}^{\ell j}_\fbvar$ which are defined from the solutions of the local cell problems \eqref{H2Dmod:eq:local_cell_problems}. These problems are computable numerically since they are defined in $\cellperper$.

\vspace{1\baselineskip} \noindent
Finally, we deduce from \eqref{H2Dmod:eq:halfguide_from_local_cells} the following expression:
\begin{equation}\label{H2Dmod:eq:DtN_from_local_DtN}
  \displaystyle
  \Lambdahat^+_\fbvar = \mathcal{T}^{10}_\fbvar\; \mathcal{P}^{\vphantom{0}}_\fbvar + \mathcal{T}^{00}_\fbvar.
\end{equation}%
}{%
}

\ifthenelse{\boolean{brouillon}}{%
  \section{Algorithm and discretization} \label{H2Dmod:sec:algo_results}
The method developed in the previous sections can be summarized before discretization as follows:
\begin{enumerate}[label=$\arabic*$., ref = {$\arabic*$}]
  \item\label{H2Dmod:item:algo_step_1} For any $\fbvar \in [-\pi, \pi)$,
  \begin{enumerate}[label=$(\textit{\alph*})$., ref = \textit{\theenumi.(\alph*)}]
    \item\label{H2Dmod:item:algo_step_1_1} solve the local cell problems \eqref{H2Dmod:eq:local_cell_problems} and compute the local DtN operators $\mathcal{T}^{j \ell}_\fbvar$ given by \eqref{H2Dmod:eq:def_local_DtN_operator};
    \item\label{H2Dmod:item:algo_step_1_2} determine the propagation operator $\mathcal{P}_\fbvar$ by solving the Riccati equation \eqref{H2Dmod:eq:riccati};
    \item\label{H2Dmod:item:algo_step_1_3} deduce the DtN operator $\Lambdahat^+_\fbvar$ using \eqref{H2Dmod:eq:DtN_from_local_DtN};
    \item\label{H2Dmod:item:algo_step_1_4} adapt Steps \ref{H2Dmod:item:algo_step_1_1} -- \ref{H2Dmod:item:algo_step_1_3} for the half-guide $\Omegaperper^-$ in order to compute the DtN operator $\Lambdahat^-_\fbvar$ thanks to cell problems defined in $(-1, 0) \times (0, 1)^2$;
    \item\label{H2Dmod:item:algo_step_1_5} find the solution $\Phi_\fbvar$ of the interface equation \eqref{H2Dmod:eq:transmission_jump};
    \item\label{H2Dmod:item:algo_step_1_6} reconstruct the solution $\linapp{\widehat{U}_\fbvar}{\widehat{G}_\fbvar}$ of the waveguide problem \eqref{H2Dmod:eq:FV_augmented_waveguide} using \eqref{H2Dmod:eq:Uhat_from_Uplusminushat} and the cell by cell expression \eqref{H2Dmod:eq:halfguide_from_local_cells} of $\linapp{\widehat{U}^\pm_\fbvar}{\Phi}$;
  \end{enumerate}
  \item\label{H2Dmod:item:algo_step_2} apply the inverse Floquet-Bloch transform \eqref{H2Dmod:eq:strip_from_waveguide} to reconstruct the solution $\linapp{U}{G}$ of the augmented strip problem \eqref{H2Dmod:eq:augmented_transmission_problem};
  \item\label{H2Dmod:item:algo_step_3} deduce the solution $u$ of the $2$D transmission problem \eqref{H2Dmod:eq:transmission_problem} using \eqref{H2Dmod:eq:rigorous_ansatz}, namely $u(\xv) = U(\cutmat\, \xv)$.
\end{enumerate}

\noindent
\surligner{Since this algorithm operates at a continuous level, it has to be discretized with respect to the space and Floquet variables. In Section \ref{H2Dmod:sec:Floquet_disc}, we present the semi-discretization procedure with respect to the Floquet variable. In Section \ref{H2Dmod:sec:discretization_x}, we first present the spatial discretization, when one uses the $3$D Finite Element method for the cell problems.} \surligner{However, as this approach is computationally costly, we propose an alternative in Section \ref{H2Dmod:sec:Quasi2D_approach}, which turned out to be essential for obtaining our numerical results. More precisely, we show that, due to their anisotropic nature, each local cell problem essentially reduces to a family of $2$D problems. The advantage of this decomposition is that the $2$D problems can be solved independently from one another, and therefore in parallel. As a result, this so-called \emph{Quasi-2D approach} appears to be less costly, though it may seem more intricate at first reading.}

\subsection{Semi-discretization with respect to the Floquet variable}\label{H2Dmod:sec:Floquet_disc}
The solution $\linapp{U}{G}$ of the augmented strip problem \eqref{H2Dmod:eq:augmented_transmission_problem} is obtained using the inverse Floquet-Bloch transform \eqref{H2Dmod:eq:strip_from_waveguide}. In this expression, the integrand involves $\linapp{\widehat{U}_\fbvar}{\widehat{G}_\fbvar}$ which, in general, can only be computed numerically (modulo an approximation) for a finite number of values of $\fbvar \in [-\pi, \pi)$ (using Step \ref{H2Dmod:item:algo_step_1}). Thus, the integral with respect to $\fbvar$ in \eqref{H2Dmod:eq:strip_from_waveguide} has to be evaluated using a quadrature rule.

\vspace{1\baselineskip} \noindent
To this end, we consider a regular mesh of $[-\pi, \pi)$ consisting of $N_\fbvar$ intervals of equal size $\Delta \fbvar$, and of $N_\fbvar + 1$ equispaced points $(\fbvar_j)_{0 \leq j \leq N_\fbvar}$, with $N_\fbvar > 0$ and $\Delta \fbvar := 2\pi / N_\fbvar$. The integral in \eqref{H2Dmod:eq:strip_from_waveguide} is evaluated using a trapezoidal rule using $\fbvar_j$ as quadrature points, leading to an approximate solution $\linapp{U_{\Delta \fbvar}}{G}$ of the augmented transmission problem \eqref{H2Dmod:eq:FV_augmented_transmission}:
\begin{equation*}\label{H2Dmod:eq:strip_from_waveguide_discrete}
  \displaystyle
  \aeforall \xva = (\xvi, \zvi_1, \zvi_2) \in \Omegaperper, \quad \spforall n \in \Z, \quad \linapp{U_{\Delta \fbvar}}{G}(\xva + n\vts \eva_1) = \frac{\Delta \fbvar}{\sqrt{2\pi}} \, \sum_{j = 0}^{N_\fbvar} \linapp{\widehat{U}_{\fbvar_j}}{\widehat{G}_{\fbvar_j}}\vts (\xva)\;  \euler^{\icplx\vts {\fbvar_j} (\zvi_1 + n)}.
\end{equation*}
The choice of the trapezoidal rule is motivated in Remark \ref{H2Dmod:rmk:exponential_convergence_trapezoidal_rule}.

\begin{rmk}\label{H2Dmod:rmk:exponential_convergence_trapezoidal_rule}
  If $\Imag \omega > 0$ and $g \in L^2(\sigma)$ is compactly supported, then similarly to \eqref{H2Dmod:eq:exp_decay_all_directions}, the solution $\linapp{U}{G}$ of \eqref{H2Dmod:eq:augmented_transmission_problem} decays exponentially in the $\eva_1$--direction. This implies thanks to Paley-Wiener-type theorems (see e.g. \cite[Theorem 2.2.2]{kuchment1993floquet}) that the map $\fbvar \mapsto \linapp{\widehat{U}_\fbvar}{\widehat{G}_\fbvar} := \fbtransform[] \linapp{U}{G}\, (\cdot, \fbvar) \in \Hgradperper{\Omegaperper}$ is real analytic.
  
  Moreover, from the properties of the Floquet-Bloch transform, it is clear that $\fbvar \mapsto \euler^{\icplx\vts \fbvar\vts (\zvi_1 + n)} \linapp{\widehat{U}_\fbvar}{\widehat{G}_\fbvar}(\xva)$ is $2\pi$--periodic, $\aeforall \xva \in \Omegaperper$, $\spforall n \in \Z$. This formal observation justifies the use of the trapezoidal rule, which is known to converge exponentially for smooth periodic integrands (see for instance \cite{trefethen2014exponentially}).
  
  \vspace{1\baselineskip} \noindent
  One way to obtain a precise error estimate is to show that the function $\linapp{U_{\Delta \fbvar}}{G}$ obtained using the trapezoidal rule can be reinterpreted as the solution of a boundary value problem defined in a domain which is bounded in the $\eva_1$--direction. This approach has been studied in detail in \cite{coatleven2012helmholtz} for a different problem, but it can be extended to \eqref{H2Dmod:eq:augmented_transmission_problem}, thus leading to a quadrature error of the form
  \[
    \spexists \alpha > 0, \quad \big\|\linapp{U}{G} - \linapp{U_{\Delta \fbvar}}{G}\big\|_{H^1_\cutmat (\Omegaper)} = O(\euler^{-\alpha\, \Imag \omega\, N_\fbvar}).
  \]
\end{rmk}

\subsection{Discretization with respect to the space variable: a first approach}\label{H2Dmod:sec:discretization_x}
We fix $\fbvar \in [-\pi, \pi)$, and consider the discretization of Step \ref{H2Dmod:item:algo_step_1} \surligner{using the Finite Element method in $3$D}. The discretization is a direct adaptation of the procedure described extensively in \cite{fliss2009analyse, fliss2016solutions, joly2006exact}. For this reason, we do not intend to go into details. We begin with a tetrahedral mesh of $\cellperper := (0, 1)^3$ with step $h > 0$. This mesh is assumed to be \emph{periodic}, in the sense that the mesh nodes on the boundary $\xvi = 0$ (\emph{resp.} $\zvi_j = 0$) can be identified to those on $\xvi = 1$ (\emph{resp.} $\zvi_j = 1$) in a trivial manner. Then, using the classical $H^1$--conforming Lagrange Finite Element space of order $d > 0$ which we call $\mathcal{V}_h (\cellperper)$, an internal approximation of $\Hgradperper{\cellperper}$ and of $\Honehalfperper{\Sigmaperper^j}$ is provided by the respective subspaces
\[
  \displaystyle
  \begin{array}[pos]{r@{\ :=\ }l}
    \mathcal{V}_{h, \perperiodic} (\cellperper) & \big\{ V \in \mathcal{V}_h (\cellperper) \ /\ V_h|_{\zvi_j = 0} = V_h|_{\zvi_j = 1} \ \ \spforall j \in \{1, 2\}  \big\},
    \\[1em]
    \mathcal{V}_{h, \perperiodic} (\Sigmaperper^j) & \big\{ V|_{\Sigmaperper^j} \ /\ V_h \in \mathcal{V}_h (\cellperper) \big\}, \quad \spforall j \in \{0, 1\}.
  \end{array}
\]
The periodicity of the mesh allows the identification $\mathcal{V}_{h, \perperiodic} (\Sigmaperper^0) \equiv \mathcal{V}_{h, \perperiodic} (\Sigmaperper^1) \equiv \mathcal{V}_{h, \perperiodic} (\Sigmaperper)$, similarly to the continuous case. In what follows, let $N_h := \dim \mathcal{V}_{h, \perperiodic} (\Sigmaperper)$ \surligner{denote the number of degrees of freedom on the face $\Sigmaperper$}.

\vspace{1\baselineskip} \noindent
For any $\Phi_h \in \smash{\mathcal{V}_{h, \perperiodic} (\Sigmaperper^j)}$, we solve the Finite Element formulation of the local cell problems \eqref{H2Dmod:eq:local_cell_problems}, to obtain solutions $\linapp{E^j_{\fbvar, h}}{\Phi_h} \in \mathcal{V}_{h, \perperiodic} (\cellperper)$, $j \in \{0, 1\}$, and we deduce finite-dimensional DtN operators $\smash{\mathcal{T}^{j \ell}_{\fbvar, h}} \in \mathscr{L}(\smash{\mathcal{V}_{h, \perperiodic} (\Sigmaperper^j)})$ using the discrete analog of the weak formulation of \eqref{H2Dmod:eq:def_local_DtN_operator}. Note that these discrete DtN operators can be represented as $N_h \times N_h$ matrices.

\vspace{1\baselineskip} \noindent
To approximate the propagation operator $\mathcal{P}_\fbvar$, it is natural to introduce the solution $\mathcal{P}_{\fbvar, h}$ of the constrained Riccati equation \eqref{H2Dmod:eq:riccati} where $\mathcal{T}^{j \ell}_\fbvar$ is replaced by $\mathcal{T}^{j \ell}_{\fbvar, h}$. As described in \cite{joly2006exact}, one can solve the discrete Riccati equation using either $(1)$ a \emph{spectral} approach which consists in characterizing $\mathcal{P}_{\fbvar, h}$ by means of its eigenpairs which satisfy a quadratic eigenvalue problem, or $(2)$ a modified \emph{Newton} method, that is, a standard Newton method with an additional projection step in order to take the spectral radius constraint into account.

\vspace{1\baselineskip} \noindent
From the approximate operator $\mathcal{P}_{\fbvar, h}$, the discrete analog of \eqref{H2Dmod:eq:DtN_from_local_DtN} allows to define an approximate DtN operator $\Lambdahat^+_{\fbvar, h} \in \mathscr{L}(\smash{\mathcal{V}_{h, \perperiodic} (\Sigmaperper^j)})$ which can be represented as a $N_h \times N_h$ matrix. We compute $\smash{\Lambdahat^-_{\fbvar, h}}$ similarly and write the discrete version of the interface equation \eqref{H2Dmod:eq:transmission_jump} as a $N_h \times N_h$ linear system. The corresponding solution $\Phi_{\fbvar, h} \in \smash{\mathcal{V}_{h, \perperiodic} (\Sigmaperper^j)}$ is then used to deduce an approximation $\linapp{\widehat{U}_{\fbvar, h}}{\widehat{G}_\fbvar}$ of the waveguide solution of \eqref{H2Dmod:eq:FV_augmented_waveguide} thanks to \eqref{H2Dmod:eq:Uhat_from_Uplusminushat} and the cell by cell expression \eqref{H2Dmod:eq:halfguide_from_local_cells}.

\begin{rmk}\label{H2Dmod:rmk:error_analysis_spatial_disc}
  \begin{enumerate}[label=$(\textit{\alph*})$., ref = \theprop.\textit{\alph*}, wide=0pt]
  \item For $\fbvar \in [-\pi, \pi)$ fixed, an error analysis can be performed for the approximation of the waveguide solution ${\linapp{\widehat{U}_\fbvar}{\widehat{G}_\fbvar}}$, assuming that the Riccati equation is solved exactly. To do so, the idea is to show as in \cite[Section 2.3.1]{fliss2009analyse} that the discrete solution $\linapp{\widehat{U}_{\fbvar, h}}{\widehat{G}_\fbvar}$ satisfies a discrete half-guide problem on an infinite mesh of $\Omegaperper$. This observation then allows to use Céa's lemma and to derive classical Finite Element estimates. In particular, for Lagrange Finite Elements of order $1$, assuming that $\aten_p$ and $\rho_p$ are smooth enough, one shows that
  \[
     \big\|\linapp{\widehat{U}_\fbvar}{\widehat{G}_\fbvar} - \linapp{\widehat{U}_{\fbvar, h}}{\widehat{G}_\fbvar} \big\|_{L^2(\Omegaperper)} = O(h^2) \quad \textnormal{and} \quad  \big\|\linapp{\widehat{U}_\fbvar}{\widehat{G}_\fbvar} - \linapp{\widehat{U}_{\fbvar, h}}{\widehat{G}_\fbvar} \big\|_{H^1_\cutmat(\Omegaperper)} = O(h).
  \]
  \item For Configuration \eqref{H2Dmod:item:config_a}, we recall that $(\aten^+_p, \rho^+_p)$ are independent of $\zvi_2$. This property can be used to solve the local cell problems more efficiently. In fact, by expanding $\linapp{E^j_\fbvar}{\Phi}$ in a Fourier series with respect to $\zvi_2$:
  \begin{equation*}
    \displaystyle
    \aeforall \xva \in \cellperper, \quad \linapp{E^j_\fbvar}{\Phi} (\xva) = \sum_{n \in \Z} \widehat{E}^j_{\fbvar, n} (\xvi, \zvi_1)\, \euler^{2\icplx \pi n \zvi_2}, \quad \widehat{E}^j_{\fbvar, n} (\xvi, \zvi_1) := \int_0^1 \linapp{E^j_\fbvar}{\Phi} (\xvi, \zvi_1, t)\, \euler^{- 2\icplx \pi n t}\, dt,
  \end{equation*}
  we note that the $3$D local cell problem satisfied by $\linapp{E^j_\fbvar}{\Phi}$ reduces to a countable family (indexed by $n \in \Z$) of decoupled $2$D problems satisfied by the Fourier coefficients $\widehat{E}^j_{\fbvar, n}$. The same remark holds for the local cell problems defined on $(-1, 0) \times (0, 1)^2$, since $(\aten^-_p, \rho^-_p)$ are independent of $\zvi_1$.
  \end{enumerate}
\end{rmk}

\vspace{1\baselineskip} \noindent
From the discrete solution $\linapp{\widehat{U}_{\fbvar, h}}{\widehat{G}_\fbvar}$, we compute an approximate solution $\linapp{U_\hv}{G}$ of the strip problem using Section \ref{H2Dmod:sec:Floquet_disc}, with $\hv := (h, \Delta \fbvar)$. Finally, for Step \ref{H2Dmod:item:algo_step_3}, an approximate solution $u_\hv$ of the $2$--dimensional problem \eqref{H2Dmod:eq:transmission_problem} is given by $u_\hv (\xv) := \extper[2]{U_\hv}(\cutmat\vts \xv),\ \spforall \xv \in \R^2$, in the spirit of \eqref{H2Dmod:eq:rigorous_ansatz}. 

\vspace{1\baselineskip} \noindent
\surligner{The discretization of the local cell problems can be computationally costly, due to the use of $3$D Finite Elements. In the next section, we show that by exploiting the anisotropic nature of the differential operator $- \rho^{-1}_p \divecut \, \aten_p\! \gradcut$, the solution of the local cell problems reduces to a family of independent $2$D cell problems and an interface equation.
%
}

\subsection{A quasi--2D idea to solve the local cell problems}\label{H2Dmod:sec:Quasi2D_approach}
In this section, we focus on the step \ref{H2Dmod:item:algo_step_1_1} of the solution algorithm, that is, the computation of the local cell solutions $E^j_\fbvar$ and the local DtN operators $\mathcal{T}^{j \ell}_\fbvar$. Although the local cell problems \eqref{H2Dmod:eq:local_cell_problems} can be solved directly using $3$--dimensional Finite Elements as in Section \ref{H2Dmod:sec:discretization_x}, it is worth recalling that the differential operator $- \rho^{-1}_p \divecut \, \aten_p\! \gradcut$ is strongly linked to the family of operators $-\rho^{-1}_s \transp{\nabla} \aten_s \!\, \nabla$ ($s \in \R$) by means of the chain rule \eqref{H2Dmod:eq:chain_rule}. It is in particular this \textquote{\emph{fibered}} link that leads to the expression \eqref{H2Dmod:eq:fibered_structure_strip} recalled below:
\begin{equation}\label{H2Dmod:eq:fibered_structure_strip_recalled}
  \displaystyle
  \aeforall (\xvi, \zvi, s) \in \R^2 \times (0, 1), \quad \widetilde{U} (\xvi, \cuti_1\vts \zvi, \cuti_2\vts \zvi + s) = u_s (\xvi, \zvi),
\end{equation}
where $U = \linapp{U}{G}$ and $u_s = \linapp{u_s}{g_s}$ (with $g_s(0, \zvi) := \widetilde{G}(0, \cuti_1\vts \zvi, \cuti_2\vts \zvi + s)$) are the respective solutions of \eqref{H2Dmod:eq:augmented_transmission_problem} and \eqref{H2Dmod:eq:transmission_problems}, and where $\widetilde{U} := \extper[2]{U}$ and $\widetilde{G} := \extper[2]{G}$ are respectively the periodic extensions of $U$ and $G$ in the $\eva_2$--direction.

\vspace{1\baselineskip} \noindent
Up to now, in the context of the lifting approach, the relation \eqref{H2Dmod:eq:fibered_structure_strip_recalled} has been useful in practice to compute $\linapp{u_s}{g_s}$ (in particular for $s = 0$) from $\linapp{U}{G}$. The so-called \emph{quasi $2$--dimensional} (or \emph{quasi-$2$D}) approach developed in this section relies on the converse: in the spirit of \eqref{H2Dmod:eq:fibered_structure_strip_recalled}, our goal is to reduce the solution of the local cell problems to that of a family of $2$--dimensional decoupled cell problems. 

\vspace{1\baselineskip} \noindent
The principle of the quasi-2D method is very similar to the method developed in \cite[Section 5.2]{amenoagbadji2023wave} for the \surligner{$1$D} Helmholtz equation with quasiperiodic coefficients. In \cite{amenoagbadji2023wave}, \surligner{a quasi-1D approach was} used \surligner{to solve} a \surligner{$2$D} cell problem with Dirichlet boundary conditions on two opposite edges \surligner{and periodic conditions on the two other edges} (as shown in Figure \ref{H2Dmod:fig:quasi_1D_periodic} right). However, the extension to our 3D problems is more delicate, due to the periodicity conditions with respect to \textbf{both} $\zvi_1$ and $\zvi_2$. To illustrate this difficulty, it is useful to \surligner{extend} the quasi-1D approach for a $2$D cell problem with periodic conditions on \surligner{the two pairs of opposite edges}. This is the object of the next section.

\subsubsection{Illustration of the method in a 2D case}
The goal in this section is first to highlight the difficulty of the fibered approach in the case of a particular 2D cell problem with periodic conditions, and then propose a \surligner{quasi-1D method} which will be extended \surligner{into a quasi-2D method for the 3D cell problems}. Throughout this section, the subscript \textquote{$\petitcarre$} is used to emphasize that we are dealing with a 2D case. The generic point is denoted by $\zv = (\zvi_1, \zvi_2) \in \R^2$, and $\{\ev_1 = (1, 0), \ev_2 = (0, 1)\}$ denotes the canonical basis of $\R^2$.

\paragraph*{\surligner{The $2$D cell problem with periodic conditions}}~\\
We consider the problem defined in $C := \{\zv = (\zvi_1, \zvi_2) \in (0, 1)^2\}$: \textit{Find $E_\petitcarre$ such that}
\begin{equation}\label{H2Dmod:eq:2D_local_cell_problem_periodic}
  \left\{
    \begin{array}{r@{\ =\ }l@{\quad}l}
      \displaystyle\displaystyle- \Dt{} (\mu_\petitcarre\, \Dt{} E^{}_\petitcarre) - \rho_\petitcarre\, \omega^2\, E^{}_\petitcarre & f_\petitcarre & \textnormal{in}\ \ C,
      \ret
      \multicolumn{3}{c}{E^{}_\petitcarre|_{\zvi_j = 0} = E^{}_\petitcarre|_{\zvi_j = 1} \quad \textnormal{and} \quad (\mu_\petitcarre\, \Dt{} E^{}_\petitcarre)|_{\zvi_j = 0} = (\mu_\petitcarre\, \Dt{} E^{}_\petitcarre)|_{\zvi_j = 1} \quad \spforall j \in \{ 1, 2\},}
    \end{array}
    \right.
\end{equation}
where $\mu_\petitcarre, \rho_\petitcarre, f_\petitcarre \in \mathscr{C}^0(\R^2)$ are $\Z^2$--periodic, and with $\cutvec := (\cuti_1, \cuti_2)$ and $\Dt{} := \cutvec \cdot \nabla = \cuti_1\, \partial_{\zvi_1} + \cuti_2\, \partial_{\zvi_2}$. Let $\widetilde{E}_\petitcarre$ denote the periodic extension of $E_\petitcarre$ in the $\ev_2$--direction. Using the chain rule and the fact that $E_\petitcarre$ is $\Z\ev_2$--periodic, one shows that the $1$D function $\mathbf{e}_s: \zvi \mapsto \widetilde{E}_\petitcarre(\zvi\vts \cutvec + s\ev_2) \in H^1(0, 1/\cuti_1)$ is well-defined for any $s \in \R$, and satisfies
\begin{equation}\label{H2Dmod:eq:1D_local_cell_PDE}
  \displaystyle
  -\frac{d}{d\zvi} \Big(\mu_{\petitcarre, s}\, \frac{d \mathbf{e}_s}{d\zvi} \Big) - \rho_{\petitcarre, s}\, \omega^2\, \mathbf{e}_s = f_{\petitcarre, s} \quad \textnormal{in}\ \ (0, 1/\cuti_1),
\end{equation}
with $\mu_{\petitcarre, s} (\zvi) := \mu_\petitcarre(\zvi\vts \cutvec + s\ev_2)$, $\rho_{\petitcarre, s} (\zvi) := \rho_\petitcarre(\zvi\vts \cutvec + s\ev_2)$, and $f_{\petitcarre, s} (\zvi) := f_\petitcarre(\zvi\vts \cutvec + s\ev_2)$. Note also that $\mathbf{e}_{s + 1} = \mathbf{e}_s$, so that the study of $(\mathbf{e}_s)_s$ can be restricted to $s \in [0, 1]$.

\vspace{1\baselineskip} \noindent
By inverting the change of variables $(s, \zvi) \mapsto \zvi\vts \cutvec + s\ev_2$, we obtain
\begin{equation}\label{H2Dmod:eq:inversion_quasi_1D_structure}
  \aeforall \zv = (\zvi_1, \zvi_2) \in C, \quad E_\petitcarre(\zv) = \mathbf{e}_{\zvi_2 - \zvi_1\vts \surligner{\vartheta}}(\zvi_1 / \cuti_1) \quad \textnormal{with} \quad \surligner{\vartheta} := \cuti_2 / \cuti_1.
\end{equation}
The quasi-1D approach would consist in computing $E_\petitcarre$ from $\mathbf{e}_s$ using the above expression. In practice, this \surligner{expression would be interesting if the $1$D problems satisfied by the functions $\mathbf{e}_s$ (indexed by $s \in [0, 1]$) were decoupled from one another}. However, because $E_\petitcarre$ is also periodic with respect to $\zvi_1$, we deduce from \eqref{H2Dmod:eq:inversion_quasi_1D_structure} that:
\begin{equation}\label{H2Dmod:eq:coupling_relation}
  \displaystyle
  \spforall s \in \R, \quad \mathbf{e}_s(0) = \mathbf{e}_{s - \surligner{\vartheta}}(1 / \cuti_1),
\end{equation}
as shown in Figure \ref{H2Dmod:fig:quasi_1D_periodic} (left). Equation \eqref{H2Dmod:eq:coupling_relation} makes it impossible to compute the $\mathbf{e}_s$ independently from one another. The alternative we propose is to relax the coupling \eqref{H2Dmod:eq:coupling_relation} by first replacing the periodicity conditions with respect to $\zvi_1$ by Dirichlet conditions for instance \surligner{(in order to decouple the solutions of the cell problems)}, and then by imposing the periodicity afterwards. This is explored in the next paragraph.

\begin{figure}[ht!]
  \centering
  \def\pentecoupe{0.33}
  \def\sorigin{0.1}
  \begin{tikzpicture}[scale=3]
    \begin{scope}
      \draw (0, 0) rectangle (1, 1);
      \draw[tertiary, very thick] (0, 0) node[above left]  {\emph{per}} -- +(0, 1);
      \draw[tertiary, very thick] (1, 0) node[above right] {\emph{per}} -- +(0, 1);
      \draw (0.5, 0) node[below] {\emph{per}};
      \draw (0.5, 1) node[above] {\emph{per}};
      \draw (0, \sorigin) -- +(1, \pentecoupe);
      \draw (0, \sorigin+\pentecoupe) -- +(1, \pentecoupe);
      \draw[dashed, tertiary] (0, \sorigin+\pentecoupe) -- +(1, 0) node[tertiary] {$\bullet$} node[right] {$\mathbf{e}_{s-\surligner{\vartheta}}(1/\cuti_1)$};
      \draw (0, \sorigin+\pentecoupe) node[tertiary] {$\bullet$} node[left, tertiary] {$\mathbf{e}_s(0)$};
      \draw (0.5, {0.5*\sorigin+0.5*\pentecoupe}) node[below] {$\mathbf{e}_{s-\surligner{\vartheta}}$};
      \draw (0.5, {\sorigin+1.5*\pentecoupe}) node[above] {$\mathbf{e}_s$};
    \end{scope}

    \begin{scope}[shift={(3.25, 0)}]
      \draw (0, 0) rectangle (1, 1);
      \draw[tertiary, very thick] (0, 0) node[above left]  {\emph{Dirichlet}} -- +(0, 1);
      \draw[tertiary, very thick] (1, 0) node[above right] {\emph{Dirichlet}} -- +(0, 1);
      \draw (0.5, 0) node[below] {\emph{per}};
      \draw (0.5, 1) node[above] {\emph{per}};
      \draw (0, \sorigin+\pentecoupe) -- +(1, \pentecoupe);
      \draw (0, \sorigin+\pentecoupe) node[tertiary] {$\bullet$} node[left, tertiary] {$\mathbf{f}^j_s (0) = \delta_{j0}$};
      \draw (1, \sorigin+2*\pentecoupe) node[tertiary] {$\bullet$} node[right, tertiary] {$\mathbf{f}^j_s (1/\cuti_1) = \delta_{j1}$};
      \draw (0.5, {\sorigin+1.5*\pentecoupe}) node[above] {$\surligner{\mathbf{f}^j_s}$};
    \end{scope}
    \begin{scope}[shift={(2.125, 0)}]
      \draw[-latex] (-0.15, 0) -- +(0.3, 0) node[right] {$\zvi_1$};
      \draw[-latex] (-0.15, 0) -- +(0, 0.3) node[above] {$\zvi_2$};
    \end{scope}
  \end{tikzpicture}
  \caption{Left: Since $E_\petitcarre$ is periodic with respect to both $\zvi_1$ and $\zvi_2$, the 1D traces $\mathbf{e}_s$ and $\mathbf{e}_{s - \surligner{\vartheta}}$ are coupled by sharing a common boundary value \eqref{H2Dmod:eq:coupling_relation}. Right: Since $F_\petitcarre$ satisfies Dirichlet conditions on $\zvi_1 \in \{0, 1\}$, each trace $\surligner{\mathbf{f}^j_s}$ satisfies independently a 1D problem with Dirichlet conditions.\label{H2Dmod:fig:quasi_1D_periodic}}
\end{figure}

\paragraph*{Auxiliary local cell problems}~\\
To overcome the difficulty induced by the periodic coupling, one can reformulate $E_\petitcarre$ in terms of the solutions of \emph{auxiliary} local Dirichlet cell problems for which the quasi-1D method can be applied. For this purpose, we introduce an auxiliary unknown $\varphi \in L^2(0, 1)$, namely
\begin{equation}\label{H2Dmod:eq:auxiliary_unknown_quasi_1D}
  \varphi := E_\petitcarre|_{\zvi_1 = 0}\ \surligner{\textnormal{(}}= E_\petitcarre|_{\zvi_1 = 1} \ \ \surligner{\textnormal{by periodicity of $\zvi_1 \mapsto E_\petitcarre (\zvi_1, \zvi_2)$)}}.
\end{equation}
With this choice, by introducing the well-posed cell problems: \textit{Find $F_\petitcarre \equiv \linapp{F_\petitcarre}{\varphi}$ such that}
\begin{equation*}\label{H2Dmod:eq:2D_local_cell_problem_Dirichlet_1}
  \left\{
    \begin{array}{r@{\ =\ }l@{\quad}l}
      \displaystyle\displaystyle- \Dt{} (\mu_\petitcarre\, \Dt{} F_\petitcarre) - \rho_\petitcarre\, \omega^2\, F_\petitcarre & 0 & \textnormal{in}\ \ C,
      \\[8pt]
      \multicolumn{3}{c}{F_\petitcarre|_{\zvi_2 = 0} = F_\petitcarre|_{\zvi_2 = 1} \quad \textnormal{and} \quad (\mu_\petitcarre\, \Dt{} F_\petitcarre)|_{\zvi_2 = 0} = (\mu_\petitcarre\, \Dt{} F_\petitcarre)|_{\zvi_2 = 1},}
      \\[8pt]
      \multicolumn{3}{c}{F_\petitcarre|_{\zvi_1 = 0} = \varphi \quad \textnormal{and} \quad F_\petitcarre|_{\zvi_1 = 1} = \varphi}
    \end{array}
    \right.
\end{equation*}
and: \textit{Find $\mathcal{G}_{\petitcarre} \equiv \linapp{\mathcal{G}_{\petitcarre}}{f_\petitcarre}$ such that}
\begin{equation*}\label{H2Dmod:eq:2D_local_cell_problem_Dirichlet_2}
  \left\{
    \begin{array}{r@{\ =\ }l@{\quad}l}
      \displaystyle\displaystyle- \Dt{} (\mu_\petitcarre\, \Dt{} \mathcal{G}^{}_{\petitcarre}) - \rho_\petitcarre\, \omega^2\, \mathcal{G}^{}_{\petitcarre} & f_\petitcarre & \textnormal{in}\ \ C,
      \\[8pt]
      \multicolumn{3}{c}{\mathcal{G}^{}_{\petitcarre}|_{\zvi_2 = 0} = \mathcal{G}^{}_{\petitcarre}|_{\zvi_2 = 1} \quad \textnormal{and} \quad (\mu_\petitcarre\, \Dt{} \mathcal{G}^{}_{\petitcarre})|_{\zvi_2 = 0} = (\mu_\petitcarre\, \Dt{} \mathcal{G}^{}_{\petitcarre})|_{\zvi_2 = 1},}
      \\[8pt]
      \multicolumn{3}{c}{\mathcal{G}^{}_{\petitcarre}|_{\zvi_1 = 0} = 0 \quad \textnormal{and} \quad \mathcal{G}^{}_{\petitcarre}|_{\zvi_1 = 1} = 0,}
    \end{array}
    \right.
\end{equation*}
we obtain by linearity that
\begin{equation}\label{H2Dmod:eq:expression_per_dir_quasi1D}
  E_\petitcarre = \linapp{F_{\petitcarre}}{\varphi} + \linapp{\mathcal{G}_{\petitcarre}}{f_\petitcarre}.
\end{equation}
The important difference between the problems satisfied by $F_\petitcarre$, $\mathcal{G}_{\petitcarre}$, and the one satisfied by $E_\petitcarre$, is that the periodicity conditions on $\{\zvi_1 = 0\}$ and $\{\zvi_1 = 1\}$ have been replaced by Dirichlet boundary conditions. Hence, $F_\petitcarre$ and $\mathcal{G}_{\petitcarre}$ can be solved using the quasi-1D approach described in \cite{amenoagbadji2023wave}, as we will \surligner{describe} later, in Proposition \ref{H2Dmod:eq:quasi_1D_structure_Dirichlet}.

\paragraph*{The auxiliary unknown $\varphi$}~\\
We now explain how to characterize the unknown $\varphi$. By imposing the periodicity of $\mu_\petitcarre\, \Dt{} E^{}_\petitcarre$ with respect to $\zvi_1$ in \eqref{H2Dmod:eq:expression_per_dir_quasi1D}, we deduce that
\begin{equation}\label{H2Dmod:eq:controle_equation}
  \displaystyle
  (\Xi^{0}_\petitcarre + \Xi^{1}_\petitcarre)\, \varphi = -(\Upsilon^0_{\petitcarre} + \Upsilon^1_{\petitcarre})\, f_\petitcarre,
\end{equation}
where the auxiliary local DtN operators $\Xi^{j}_\petitcarre$ and the right-hand sides $\Upsilon^j_{\petitcarre}\, f$, $j \in \{0, 1\}$ are given by
\begin{equation*}
  \left\{
    \begin{array}{r@{\ :=\ }l}
      \Xi^{j}_\petitcarre\, \varphi & (-1)^{j+1}\, (\mu_\petitcarre\, \Dt{} \linapp{F_{\petitcarre}}{\varphi})|_{\zvi_1 = j},
      \ret
      \Upsilon^j_{\petitcarre}\, f & (-1)^{j+1}\, (\mu_\petitcarre\, \Dt{} \linapp{\mathcal{G}_{\petitcarre}}{f})|_{\zvi_1 = j}.
    \end{array}
  \right.
\end{equation*}
Conversely, it can be shown using the presence of absorption that $\Xi^{0}_\petitcarre + \Xi^{1}_\petitcarre$ is invertible, so that \eqref{H2Dmod:eq:controle_equation} admits a unique solution. 

\paragraph{The fibered structure of the auxiliary local cell problems}~\\
In order to emphasize the fibered structure of $F_\petitcarre, \mathcal{G}_\petitcarre$, consider the well-posed $1$D problems given for any $s \in \R$ and $j \in \{0, 1\}$ by: \textit{Find $\surligner{\mathbf{f}^j_s} \in H^1(0, 1/\cuti_1)$ such that}
\begin{equation*}
  \left\{
    \begin{array}{r@{\ =\ }l}
      \displaystyle -\frac{d}{d\zvi} \Big(\mu_{\petitcarre, s}\, \frac{d \surligner{\mathbf{f}^j_s}}{d\zvi} \Big) - \rho_{\petitcarre, s}\, \omega^2\, \surligner{\mathbf{f}^j_s} & 0 \quad \textnormal{in}\ \ (0, 1/\cuti_1),
      \\[4pt]
      \surligner{\mathbf{f}^0_s (0) = 1,\quad \mathbf{f}^0_s (1/\cuti_1)} & \surligner{0},
      \\[4pt]
      \surligner{\mathbf{f}^1_s (0) = 0,\quad \mathbf{f}^1_s (1/\cuti_1)} & \surligner{1},
    \end{array}
  \right.
\end{equation*}
and: \textit{Find $\mathbf{g}_s = \linapp{\mathbf{g}_s}{f_{\petitcarre, s}} \in H^1(0, 1/\cuti_1)$ such that}
\begin{equation*}
  \left\{
    \begin{array}{r@{\ =\ }l}
      \displaystyle -\frac{d}{d\zvi} \Big(\mu_{\petitcarre, s}\, \frac{d \mathbf{g}_s}{d\zvi} \Big) - \rho_{\petitcarre, s}\, \omega^2\, \mathbf{g}_s & f_{\petitcarre, s} \quad \textnormal{in}\ \ (0, 1/\cuti_1),
      \\[4pt]
      \mathbf{g}_s (0) = \mathbf{g}_s (1/\cuti_1) & 0.
    \end{array}
  \right.
\end{equation*}
Note that $\surligner{\mathbf{f}^j_{s + 1} = \mathbf{f}^j_s}$ and $\mathbf{g}_{s + 1} = \mathbf{g}_s$. The advantage of resorting to the auxiliary solutions $F_\petitcarre, \mathcal{G}_\petitcarre$ lies in the following result, which highlight their fibered or so-called \emph{quasi-1D} structure.
\begin{prop}[\textnormal{\cite[Proposition 5.1]{amenoagbadji2023wave}}]\label{H2Dmod:eq:quasi_1D_structure_Dirichlet}
For any $(\varphi, f_\petitcarre) \in L^2(0, 1) \times L^2(C)$, we have
\begin{equation}
  \aeforall (\zvi, s) \in (0, 1/\cuti_1) \times (0, 1), \quad %
  \left\{%
  \begin{array}{r@{\ =\ }l}
    \widetilde{F}_\petitcarre (\cuti_1 \zvi, \cuti_2 \zvi + s) & \surligner{\varphi(s)\, \mathbf{f}^0_s (\zvi) + \varphi(s + \surligner{\vartheta})\, \mathbf{f}^1_s (\zvi)}, \quad \surligner{\vartheta} := \cuti_2/\cuti_1,
    \ret
    \widetilde{\mathcal{G}}_\petitcarre (\cuti_1 \zvi, \cuti_2 \zvi + s) & \mathbf{g}_s (\zvi),
  \end{array}
  \right.
\end{equation}
\surligner{with $\mathbf{g}_s = \linapp{\mathbf{g}_s}{f_{\petitcarre, s}}$, and where $\widetilde{F}_\petitcarre$ (\emph{resp.} $\widetilde{\mathcal{G}}_\petitcarre$) denotes the periodic extension of $F_\petitcarre$ (\emph{resp.} $\mathcal{G}_\petitcarre$) in the $\ev_2$--direction, with $F_\petitcarre = \linapp{F_\petitcarre}{\varphi}$ and $\mathcal{G}_\petitcarre = \linapp{\mathcal{G}_\petitcarre}{f_\petitcarre}$.} 
\end{prop}
\noindent
Proposition \ref{H2Dmod:eq:quasi_1D_structure_Dirichlet} allows to compute $\linapp{F_\petitcarre}{\varphi}$ (\emph{resp.} $\linapp{\mathcal{G}_\petitcarre}{f_\petitcarre}$) from the $1$D functions $\surligner{\mathbf{f}^j_s}$ (\emph{resp.} $\linapp{\mathbf{g}_s}{f_{\petitcarre, s}}$), $s \in [0, 1]$. The difference with $\mathbf{e}_s$ is that $\surligner{\mathbf{f}^j_s}$ and $\linapp{\mathbf{g}_s}{f_{\petitcarre, s}}$ satisfy $1$D problems which can be solved independently from one another in practice, and hence in parallel (see Figure \ref{H2Dmod:fig:quasi_1D_periodic} right). Note also that the operators $\Xi^j_\petitcarre, \Upsilon^j_\petitcarre$ can be computed using $\surligner{\mathbf{f}^j_s}$ and $\linapp{\mathbf{g}_s}{f_{\petitcarre, s}}$. We refer to \cite[Proposition 5.2]{amenoagbadji2023wave}, where a similar situation is presented. 

\paragraph*{Algorithm for the quasi-1D method}~\\
The approach presented above to compute the solution $E_\petitcarre$ of \eqref{H2Dmod:eq:2D_local_cell_problem_periodic} can be summarized as follows:
\begin{itemize}
  \item Compute the solutions $F_\petitcarre, \mathcal{G}_\petitcarre$ of the auxiliary cell problems using the solutions $\surligner{\mathbf{f}^j_s}, \mathbf{g}_s$ of the $1$D problems indexed by $s \in [0, 1]$, and which can be solved in parallel;
  \item Construct $\Xi^j_\petitcarre, \Upsilon^j_\petitcarre$, and deduce the auxiliary unknown $\varphi$ by solving \eqref{H2Dmod:eq:controle_equation}.
\end{itemize}
The discretization of the first item is described in \cite[Section 5.2.2]{amenoagbadji2023wave}. The difference with the quasi-1D method in \cite{amenoagbadji2023wave} is the second item above. \surligner{In practice, \eqref{H2Dmod:eq:controle_equation} is a dense system of size $N \times N$, where $N$ is the number of degrees of freedom on the edge $\zvi_1 = 0$.} \surligner{Therefore, instead of solving a $2$D problem of size $N^2$, one can solve $N$ independent $1$D problems of size $N$ in parallel, and invert a $N$--sized dense matrix, which encodes the coupling between the $1$D solutions.}

\subsubsection{Extension to the 3D local cell problems}
In what follows, the Floquet variable $\fbvar \in (-\pi, \pi)$ is fixed, and is omitted in some of the notations. We now extend the idea developed in the previous section to the case of the 3D local cell problems \eqref{H2Dmod:eq:local_cell_problems} which are recalled under the formal form:
\begin{equation}
  \label{H2Dmod:eq:recalled_local_cell_problems}
  \left\{
    \begin{array}{r@{\ =\ }l}
      \displaystyle- \divecutFB\, \aten_p \!  \gradcutFB\, E^j_\fbvar - \rho_p\, \omega^2\, E^j_\fbvar & 0 \quad \textnormal{in}\ \  \cellperper,
      \\[8pt]
      \multicolumn{2}{c}{\textnormal{$E^j_\fbvar$ is $1$--periodic with respect to $\zvi_1$ and $\zvi_2$,}}
      \\[8pt]
      \multicolumn{2}{c}{\linapp{E^0_\fbvar}{\Phi}|_{\Sigmaperperi{0}} = \Phi \quad \textnormal{and} \quad \linapp{E^0_\fbvar}{\Phi}|_{\Sigmaperperi{1}} = 0,}
      \\[8pt]
      \multicolumn{2}{c}{\linapp{E^1_\fbvar}{\Phi}|_{\Sigmaperperi{0}} = 0 \quad \textnormal{and} \quad \linapp{E^1_\fbvar}{\Phi}|_{\Sigmaperperi{1}} = \Phi,}
    \end{array}
  \right.
\end{equation}
where $\Phi$ \surligner{is a fixed surface data, which plays a similar role than the volume} source term $f_\petitcarre$ in the previous section (see \eqref{H2Dmod:eq:1D_local_cell_PDE}). By analogy with the previous section, we want to take advantage of the quasi-2D nature of these \surligner{3D} problems, to deduce $E^j_\fbvar$ from the solution of $2$D cell problems set in $(\cutmat\, \R^2 + s\, \eva_2) \cap \cellperper$ \surligner{and parameterized by $s \in [0, 1]$}. However, \surligner{due to} the periodicity conditions with respect to $\zvi_1$ and $\zvi_2$, the traces of $E^j_\fbvar$ on the hyperplane $\cutmat\, \R^2 + s\, \eva_2$, \surligner{$s \in [0, 1]$,} satisfy a coupling relation that \surligner{prevents} them from being computed independently from one another. For this reason, we introduce auxiliary Dirichlet cell problems that have the desired fibered structure. The link between $E^j_\fbvar$ and the solutions of these auxiliary problems is related to an auxiliary unknown, \surligner{chosen as} the trace of $E^j_\fbvar$ on the interface $\{\zvi_1 = \surligner{0}\}$, \surligner{and} which satisfies an equation similar to \eqref{H2Dmod:eq:controle_equation}.

\paragraph*{The auxiliary local cell problems}~\\
Consider $E^0_\fbvar$ for simplicity. We introduce a Dirichlet-to-Dirichlet operator defined by
\begin{equation} \label{H2Dmod:eq:def_DtD}
  \displaystyle
  \spforall \Phi \in \Honehalfperper{\Sigmaperper}, \quad \mathcal{R}^0 \Phi := \linapp{E^0_\fbvar}{\Phi}|_{\zvi_1 = 0} \ \ (=  \linapp{E^0_\fbvar}{\Phi}|_{\zvi_1 = 1}\ \ \textnormal{by periodicity of $\zvi_1 \mapsto \linapp{E^0_\fbvar}{\Phi} (\xva)$}).
\end{equation}
Then for any $\Phi, \Psi$, using the respective solutions $\linapp{F}{\Psi}$ and $\linapp{\mathcal{G}^0}{\Phi}$ of the well-posed \surligner{Dirichlet} cell problems \surligner{(see Figure \ref{H2Dmod:fig:quasi2D_auxiliary_Dirichlet})}:
\begin{equation}\label{H2Dmod:eq:auxiliary_local_cell_problem_F}
  \left\{
    \begin{array}{r@{\ =\ }l@{\quad}l}
      \displaystyle\displaystyle- \divecutFB\, \aten_p\!  \gradcutFB\, F - \rho_p\, \omega^2\, F & 0 & \textnormal{in}\ \ \cellperper,
      \\[8pt]
      \multicolumn{2}{c}{\textnormal{$F = \linapp{F}{\Psi}$ is $1$--periodic with respect to $\zvi_2$,}}
      \\[8pt]
      \multicolumn{3}{c}{\linapp{F}{\Psi}|_{\xvi = 0} = 0 \quad \textnormal{and} \quad \linapp{F}{\Psi}|_{\xvi = 1} = 0,}
      \\[8pt]
      \multicolumn{3}{c}{\linapp{F}{\Psi}|_{\zvi_1 = 0} = \Psi \quad \textnormal{and} \quad \linapp{F}{\Psi}|_{\zvi_1 = 1} = \Psi,}
    \end{array}
    \right.
\end{equation} 
and
\begin{equation}\label{H2Dmod:eq:auxiliary_local_cell_problem_G}
  \left\{
    \begin{array}{r@{\ =\ }l@{\quad}l}
      \displaystyle\displaystyle- \divecutFB\, \aten_p\!  \gradcutFB\, \mathcal{G}^0 - \rho_p\, \omega^2\, \mathcal{G}^0 & 0 & \textnormal{in}\ \ \cellperper,
      \\[8pt]
      \multicolumn{2}{c}{\textnormal{$\mathcal{G}^0 = \linapp{\mathcal{G}^0}{\Phi}$ is $1$--periodic with respect to $\zvi_2$,}}
      \\[8pt]
      \multicolumn{3}{c}{\linapp{\mathcal{G}^0}{\Phi}|_{\xvi = 0} = \Phi \quad \textnormal{and} \quad \linapp{\mathcal{G}^0}{\Phi}|_{\xvi = 1} = 0,}
      \\[8pt]
      \multicolumn{3}{c}{\linapp{\mathcal{G}^0}{\Phi}|_{\zvi_1 = 0} = 0 \quad \textnormal{and} \quad \linapp{\mathcal{G}^0}{\Phi}|_{\zvi_1 = 1} = 0,}
    \end{array}
    \right.
\end{equation} 
it follows by linearity that 
\begin{equation}\label{H2Dmod:eq:from_auxiliary_local_cell_to_local_cell}
  \displaystyle
  \linapp{E^0_\fbvar}{\Phi} = \linapp{F}{\mathcal{R}^0 \Phi} + \linapp{\mathcal{G}^0}{\Phi}.
\end{equation}
\surligner{The main interest of the problems (\ref{H2Dmod:eq:auxiliary_local_cell_problem_F}, \ref{H2Dmod:eq:auxiliary_local_cell_problem_G}) satisfied by $\linapp{F}{\Psi}$ and $\linapp{\mathcal{G}^0}{\Phi}$ is that although they remain periodic in the $\eva_2$--direction, they involve Dirichlet conditions on the faces $\{\xvi= \ell\}$ and $\{\zvi_1 = \ell\}$. For this reason, these problems can be interpreted as a concatenation of $2$D cell problems, as we shall describe later, in Proposition \ref{H2Dmod:prop:fibered_structure_auxiliary_problems}.}

\begin{rmk}
  For a fixed $\Phi$, $\mathcal{R}^0 \Phi$ plays the same role as the auxiliary unknown $\varphi$ defined by \eqref{H2Dmod:eq:auxiliary_unknown_quasi_1D}.
\end{rmk}

\begin{figure}[ht!]
  \noindent\makebox[\textwidth]{%
    \def\coteCyl{1.5}
    \begin{tikzpicture}[scale=0.825]
      \def\angleinterface{55}%
      \definecolor{coplancoupe}{RGB}{211, 227, 65}
      \tdplotsetmaincoords{75}{40}
      \def\axesrapport{3}
      \def\ecart{4.5}
      \begin{scope}[tdplot_main_coords, scale=0.9]
        \draw[latex-, canvas is zx plane at y=0.5*\coteCyl, tertiary!80!black] (0.5*\coteCyl, 0) -- (0.5*\coteCyl, -1*\coteCyl) node[left] {$\Phi$};
        \draw[latex-, canvas is yz plane at x=0.5*\coteCyl, rougeTerre!80!black] (\coteCyl, 0.5*\coteCyl) -- (2.25*\coteCyl, 0.5*\coteCyl) node[pos=1.5] {\emph{per}};
        %
        %
        \fill[gray!40, opacity=0.5] (\coteCyl, 0, 0) -- +(-\coteCyl, 0, 0) -- +(-\coteCyl, \coteCyl, 0) -- +(0, \coteCyl, 0) -- cycle;
        \fill[rougeTerre!90, opacity=0.625] (\coteCyl, \coteCyl, 0) -- +(-\coteCyl, 0, 0) -- +(-\coteCyl, 0, \coteCyl) -- +(0, 0, \coteCyl) -- cycle;
        \draw[white, dashed, opacity=1] (0, \coteCyl, 0) -- +(\coteCyl, 0, 0);
        \fill[draw=white, dashed, fill=tertiary!95!black, opacity=0.75] (0, 0, 0) -- +(0, \coteCyl, 0) -- +(0, \coteCyl, \coteCyl) -- +(0, 0, \coteCyl) -- cycle;
        \fill[draw=white, fill=tertiary!95!black, opacity=0.75] (\coteCyl, 0, 0) -- +(0, \coteCyl, 0) -- +(0, \coteCyl, \coteCyl) -- +(0, 0, \coteCyl) -- cycle;
        \filldraw[draw=white, fill=rougeTerre!90, opacity=0.8] (\coteCyl, 0, 0) -- +(-\coteCyl, 0, 0) -- +(-\coteCyl, 0, \coteCyl) -- +(0, 0, \coteCyl) -- cycle;
        \fill[fill=gray!40, opacity=0.8] (\coteCyl, 0, \coteCyl) -- +(-\coteCyl, 0, 0) -- +(-\coteCyl, \coteCyl, 0) -- +(0, \coteCyl, 0) -- cycle;
        \draw[latex-, canvas is zx plane at y=0.5*\coteCyl, tertiary!80!black] (0.5*\coteCyl, \coteCyl) -- (0.5*\coteCyl, 2*\coteCyl) node[right] {$0$};
        \draw[latex-, canvas is yz plane at x=0.5*\coteCyl, rougeTerre!80!black] (0, 0.5*\coteCyl) -- (-1.25*\coteCyl, 0.5*\coteCyl) node[pos=1.4] {\emph{per}};
        \draw (0.5*\coteCyl, 0.5*\coteCyl, 1.5*\coteCyl) node {$\linapp{E^0_\fbvar}{\Phi}$};
      \end{scope}%
      \begin{scope}[shift={(5.75, 0)}]
        \begin{scope}[tdplot_main_coords, scale=0.9]
          \draw[latex-, canvas is zx plane at y=0.5*\coteCyl] (0.5*\coteCyl, 0) -- (0.5*\coteCyl, -1*\coteCyl) node[left] {$0$};
          \draw[latex-, canvas is yz plane at x=0.5*\coteCyl, rougeTerre!80!black] (\coteCyl, 0.5*\coteCyl) -- (2.25*\coteCyl, 0.5*\coteCyl) node[pos=1.25] {$\Psi$};
          %
          %
          \fill[gray!40, opacity=0.5] (\coteCyl, 0, 0) -- +(-\coteCyl, 0, 0) -- +(-\coteCyl, \coteCyl, 0) -- +(0, \coteCyl, 0) -- cycle;
          \fill[rougeTerre!90, opacity=0.625] (\coteCyl, \coteCyl, 0) -- +(-\coteCyl, 0, 0) -- +(-\coteCyl, 0, \coteCyl) -- +(0, 0, \coteCyl) -- cycle;
          \draw[white, dashed, opacity=1] (0, \coteCyl, 0) -- +(\coteCyl, 0, 0);
          \fill[draw=white, dashed, fill=gray!55, opacity=0.75] (0, 0, 0) -- +(0, \coteCyl, 0) -- +(0, \coteCyl, \coteCyl) -- +(0, 0, \coteCyl) -- cycle;
          \fill[draw=white, fill=gray!95, opacity=0.75] (\coteCyl, 0, 0) -- +(0, \coteCyl, 0) -- +(0, \coteCyl, \coteCyl) -- +(0, 0, \coteCyl) -- cycle;
          \filldraw[draw=white, fill=rougeTerre!90, opacity=0.8] (\coteCyl, 0, 0) -- +(-\coteCyl, 0, 0) -- +(-\coteCyl, 0, \coteCyl) -- +(0, 0, \coteCyl) -- cycle;
          \fill[fill=gray!40, opacity=0.8] (\coteCyl, 0, \coteCyl) -- +(-\coteCyl, 0, 0) -- +(-\coteCyl, \coteCyl, 0) -- +(0, \coteCyl, 0) -- cycle;
          \draw[latex-, canvas is zx plane at y=0.5*\coteCyl] (0.5*\coteCyl, \coteCyl) -- (0.5*\coteCyl, 2*\coteCyl) node[right] {$0$};
          \draw[latex-, canvas is yz plane at x=0.5*\coteCyl, rougeTerre!80!black] (0, 0.5*\coteCyl) -- (-1.25*\coteCyl, 0.5*\coteCyl) node[pos=1.25] {$\Psi$};
          \draw (0.5*\coteCyl, 0.5*\coteCyl, 1.5*\coteCyl) node {$\linapp{F}{\Psi}$};
        \end{scope}%
      \end{scope}
      \begin{scope}[shift={(11.5, 0)}]
        \begin{scope}[tdplot_main_coords, scale=0.9]
          \draw[latex-, canvas is zx plane at y=0.5*\coteCyl, tertiary!80!black] (0.5*\coteCyl, 0) -- (0.5*\coteCyl, -1*\coteCyl) node[left] {$\Phi$};
          \draw[latex-, canvas is yz plane at x=0.5*\coteCyl] (\coteCyl, 0.5*\coteCyl) -- (2.25*\coteCyl, 0.5*\coteCyl) node[pos=1.25] {$0$};
          %
          %
          \fill[gray!40, opacity=0.5] (\coteCyl, 0, 0) -- +(-\coteCyl, 0, 0) -- +(-\coteCyl, \coteCyl, 0) -- +(0, \coteCyl, 0) -- cycle;
          \fill[gray!90, opacity=0.625] (\coteCyl, \coteCyl, 0) -- +(-\coteCyl, 0, 0) -- +(-\coteCyl, 0, \coteCyl) -- +(0, 0, \coteCyl) -- cycle;
          \draw[white, dashed, opacity=1] (0, \coteCyl, 0) -- +(\coteCyl, 0, 0);
          \fill[draw=white, dashed, fill=tertiary!95!black, opacity=0.75] (0, 0, 0) -- +(0, \coteCyl, 0) -- +(0, \coteCyl, \coteCyl) -- +(0, 0, \coteCyl) -- cycle;
          \fill[draw=white, fill=tertiary!95!black, opacity=0.75] (\coteCyl, 0, 0) -- +(0, \coteCyl, 0) -- +(0, \coteCyl, \coteCyl) -- +(0, 0, \coteCyl) -- cycle;
          \filldraw[draw=white, fill=gray!90, opacity=0.8] (\coteCyl, 0, 0) -- +(-\coteCyl, 0, 0) -- +(-\coteCyl, 0, \coteCyl) -- +(0, 0, \coteCyl) -- cycle;
          \fill[fill=gray!40, opacity=0.8] (\coteCyl, 0, \coteCyl) -- +(-\coteCyl, 0, 0) -- +(-\coteCyl, \coteCyl, 0) -- +(0, \coteCyl, 0) -- cycle;
          \draw[latex-, canvas is zx plane at y=0.5*\coteCyl, tertiary!80!black] (0.5*\coteCyl, \coteCyl) -- (0.5*\coteCyl, 2*\coteCyl) node[right] {$0$};
          \draw[latex-, canvas is yz plane at x=0.5*\coteCyl] (0, 0.5*\coteCyl) -- (-1.25*\coteCyl, 0.5*\coteCyl) node[pos=1.25] {$0$};
          \draw (0.5*\coteCyl, 0.5*\coteCyl, 1.5*\coteCyl) node {$\linapp{\mathcal{G}^0}{\Phi}$};
        \end{scope}%
      \end{scope}
      \begin{scope}[shift={(16, 0)}]
        \begin{scope}[tdplot_main_coords, scale=0.9]
          \node (xaxis) at (\coteCyl, 0, 0) {$\xvi$};
          \node (yaxis) at (0, \coteCyl, 0) {$\zvi_1$};
          \node (zaxis) at (0, 0, \coteCyl) {$\zvi_2$};
          \draw[-latex] (0, 0, 0) -- (xaxis);
          \draw[-latex] (0, 0, 0) -- (yaxis);
          \draw[-latex] (0, 0, 0) -- (zaxis);
        \end{scope}
      \end{scope}
    \end{tikzpicture} 
  }
  \caption{Dirichlet conditions satisfied by $\linapp{E^0_\fbvar}{\Phi}$, $\linapp{F}{\Psi}$, and $\linapp{\mathcal{G}^0}{\Phi}$. \label{H2Dmod:fig:quasi2D_auxiliary_Dirichlet}}
\end{figure}

\paragraph*{Characterization of the Dirichlet-to-Dirichlet operator $\mathcal{R}^0$}~\\
Before highlighting the advantage of introducing $F$ and $\mathcal{G}^0$, let us derive an equation to characterize $\mathcal{R}^0$. The periodicity of $E^0_\fbvar$ in $\zvi_1$ leads to the following equality
\begin{equation*}
  \displaystyle
  \big(\cutmat\vts \aten_p \! \gradcutFB\, \linapp{E^0_\fbvar}{\Phi}\big) \cdot \eva_1|_{\zvi_1 = 0} = \big(\cutmat\vts \aten_p \! \gradcutFB\, \linapp{E^0_\fbvar}{\Phi}\big) \cdot \eva_1|_{\zvi_1 = 1},
\end{equation*}
which reformulates as
\begin{equation}\label{H2Dmod:eq:DtD_from_auxiliary_DtN}
  \spforall \Phi \in \Honehalfperper{\Sigmaperper}, \quad (\Xi^{0} + \Xi^{1})\, \mathcal{R}^0 \Phi = -(\Upsilon^{00} + \Upsilon^{01})\, \Phi,
\end{equation}
with $\Xi^{\ell}$ and $\Upsilon^{0\ell}$, $\ell \in \{0, 1\}$, being the local auxiliary DtN operators defined by
\begin{equation}\label{H2Dmod:eq:def_auxiliary_DtN}
  \displaystyle
  \spforall \Phi, \Psi, \quad \left\{
    \begin{array}{r@{\ :=\ }l}
      \Xi^{\ell}\, \Psi & (-1)^{\ell+1}\, \big(\cutmat\vts \aten_p \! \gradcutFB\, \linapp{F}{\Psi}\big) \cdot \eva_1 |_{\zvi_1 = \ell},
      \ret
      \Upsilon^{0\ell}\, \Phi & (-1)^{\ell+1}\, \big(\cutmat\vts \aten_p \! \gradcutFB\, \linapp{\mathcal{G}^0}{\Phi}\big) \cdot \eva_1 |_{\zvi_1 = \ell}.
    \end{array}
  \right.
\end{equation}
Conversely, it can be shown using the presence of absorption that $\Xi^{0} + \Xi^{1}$ is invertible, so that \eqref{H2Dmod:eq:DtD_from_auxiliary_DtN} is well-posed.

\vspace{1\baselineskip} \noindent
The above arguments extend naturally to $E^1_\fbvar$, to which we can associate $\mathcal{R}^1$, $\mathcal{G}^1$, and $\Upsilon^{1\ell}$ by adapting respectively \eqref{H2Dmod:eq:def_DtD}, \eqref{H2Dmod:eq:auxiliary_local_cell_problem_G}, and \eqref{H2Dmod:eq:def_auxiliary_DtN}. Then one has
\begin{equation}\label{H2Dmod:eq:DtD_from_auxiliary_DtN_1}
  \spforall \Phi \in \Honehalfperper{\Sigmaperper}, \quad \linapp{E^1_\fbvar}{\Phi} = \linapp{F}{\mathcal{R}^1 \Phi} + \linapp{\mathcal{G}^1}{\Phi}, \quad \textnormal{where} \quad (\Xi^{0} + \Xi^{1})\, \mathcal{R}^1 = -(\Upsilon^{10} + \Upsilon^{11}).  
\end{equation}
From $\mathcal{R}^0$ and $\mathcal{R}^1$, one can also deduce the local DtN operators $\mathcal{T}^{j \ell}_\fbvar$ defined by \eqref{H2Dmod:eq:def_local_DtN_operator}. More precisely, \surligner{using (\ref{H2Dmod:eq:from_auxiliary_local_cell_to_local_cell}, \ref{H2Dmod:eq:DtD_from_auxiliary_DtN}) and \eqref{H2Dmod:eq:DtD_from_auxiliary_DtN_1}}, it follows by linearity that
\begin{equation}\label{H2Dmod:eq:DtN_from_auxiliary_DtN}
  \spforall j, \ell \in \{0, 1\}, \quad \mathcal{T}^{j\ell}_\fbvar = \widetilde{\Upsilon}^{j\ell} + \widetilde{\Xi}^\ell\, \mathcal{R}^j,
\end{equation}
where $\widetilde{\Xi}^\ell$ and $\widetilde{\Upsilon}^{j \ell}$ are given by
\begin{equation}
  \displaystyle
  \spforall \Phi, \Psi, \quad \left\{
    \begin{array}{r@{\ :=\ }l}
      \widetilde{\Xi}^{\ell}\, \Psi & (-1)^{\ell+1}\, \big(\cutmat\vts \aten_p \! \gradcutFB\, \linapp{F}{\Psi}\big) \cdot \eva_\xvi |_{\xvi = \ell},
      \ret
      \widetilde{\Upsilon}^{j \ell}\, \Phi & (-1)^{\ell+1}\, \big(\cutmat\vts \aten_p \! \gradcutFB\, \linapp{\mathcal{G}^j}{\Phi}\big) \cdot \eva_1 |_{\xvi = \ell}.
    \end{array}
  \right.
\end{equation}
It is worth noting that the operators $\Xi^\ell$, $\widetilde{\Xi}^\ell$, $\Upsilon^{j \ell}$, $\widetilde{\Upsilon}^{j \ell}$ can all be obtained by computing $F$ and $\mathcal{G}^j$, which satisfy cell problems with Dirichlet conditions on \emph{both} the boundaries $\{\xvi = \ell\}$ and $\{\zvi_1 = \ell\}$. We now highlight the structure of $F$ and $\mathcal{G}^j$.

\paragraph*{The fibered structure of the auxiliary cell problems}~\\
Given $s \in \R$, using the definition \eqref{H2Dmod:eq:def_As_rho_s} of $(\aten_s, \rho_s)$, let us introduce the 2D cell problems defined in \surligner{the cell} $\Qcuti := \{\xv = (\xvi, \zvi) \in (0, 1) \times (0, 1/\cuti_1)\}$ \surligner{and for any $j \in \{0, 1\}$ as}
\begin{equation}\label{H2Dmod:eq:2D_auxiliary_local_cell_problem_F}
  \left\{
    \begin{array}{r@{\ =\ }l@{\quad}l}
      \displaystyle\displaystyle- \divecutFBtwoD\, \aten_s\, \gradcutFBtwoD\, \surligner{\mathbf{f}^j_s} - \rho_s\, \omega^2\, \surligner{\mathbf{f}^j_s} & 0 & \textnormal{in}\ \ \Qcuti,
      \\[8pt]
      \multicolumn{3}{c}{\surligner{\linapp{\mathbf{f}^j_s}{\psi}}|_{\xvi = 0} = 0 \quad \textnormal{and} \quad \surligner{\linapp{\mathbf{f}^j_s}{\psi}}|_{\xvi = 1} = 0,}
      \\[8pt]
      \multicolumn{3}{c}{\surligner{\linapp{\mathbf{f}^j_s}{\psi}|_{\zvi = \ell} = \delta_{j \ell}\, \psi,}}
    \end{array}
    \right.
\end{equation}
and
\begin{equation}\label{H2Dmod:eq:2D_auxiliary_local_cell_problem_G}
  \left\{
    \begin{array}{r@{\ =\ }l@{\quad}l}
      \displaystyle\displaystyle- \divecutFBtwoD\, \aten_s\, \gradcutFBtwoD\, \mathbf{g}^j_s - \rho_s\, \omega^2\, \mathbf{g}^j_s & 0 & \textnormal{in}\ \ \Qcuti,
      \\[8pt]
      \multicolumn{3}{c}{\linapp{\mathbf{g}^j_s}{\varphi}|_{\xvi = \ell} = \delta_{j \ell}\, \varphi,}
      \\[8pt]
      \multicolumn{3}{c}{\linapp{\mathbf{g}^j_s}{\varphi}|_{\zvi = 0} = 0 \quad \textnormal{and} \quad \linapp{\mathbf{g}^j_s}{\varphi}|_{\zvi = 1/\cuti_1} = 0.}
    \end{array}
    \right.
\end{equation}
We also introduce for $\ell \in \{0, 1\}$, $\Gamma_\perperiodic^\ell := (0, 1) \times \{\ell\} \times (0, 1)$ and the edges (see Figure \ref{H2Dmod:fig:coupled_traces_faces})
\begin{equation*}
  \Scuti^\mathtt{X}_\ell := \{\ell\} \times (0, 1/\cuti_1), \quad \textnormal{and} \quad \Scuti^\mathtt{Z}_\ell := (0, 1) \times \{\ell/\cuti_1\}.
\end{equation*}
Then, the link between the solution $F$ (\textit{resp} $\mathcal{G}^j$) of \eqref{H2Dmod:eq:auxiliary_local_cell_problem_F} (\textit{resp} \eqref{H2Dmod:eq:auxiliary_local_cell_problem_G} for $j = 0$) and $\mathbf{f}_s$ (\textit{resp} $\mathbf{g}^j_s$) is given by the next result, which mirrors Proposition \ref{H2Dmod:eq:quasi_1D_structure_Dirichlet}.
\begin{prop}\label{H2Dmod:prop:fibered_structure_auxiliary_problems}
  Let $\Psi \in L^2(\Gamma_\perperiodic^j)$, for $j \in \{0, 1\}$, be such that $\psi_s := \shearmap \Psi (\cdot, s): \surligner{\xvi \mapsto \extper{\Psi}(\xvi, 0, s)}$ belongs to $H^{1/2}(\Scuti^\mathtt{Z}_j)$ for almost any $s \in \R$. \surligner{Recall that $\vartheta := \cuti_2 / \cuti_1$}. Then one has
  \begin{equation}\label{H2Dmod:eq:link_auxiliary_2D_3D_solutions_F}
    \displaystyle
    \surligner{\aeforall (\xvi, \zvi, s) \in \Qcuti \times (0, 1), \quad \extper{F} (\xvi, \cuti_1\, \zvi, \cuti_2\, \zvi + s) = \linapp{\mathbf{f}^0_s}{\psi_s} (\xvi, \zvi) + \linapp{\mathbf{f}^1_s}{\psi_{s + \vartheta}} (\xvi, \zvi)},
  \end{equation}
  where $F = \linapp{F}{\Psi}$ and \surligner{$(\mathbf{f}^0_s, \mathbf{f}^1_s)$} denote the respective solutions of \eqref{H2Dmod:eq:auxiliary_local_cell_problem_F} and \eqref{H2Dmod:eq:2D_auxiliary_local_cell_problem_F}. Similarly, for any $\Phi \in \Honehalfperper{\Sigmaperper^j}$, \surligner{by setting} $\varphi_s := \shearmap \Phi (\cdot, s)\surligner{: \zvi \mapsto \extper{\Phi}(0, \cuti_1\, \zvi, \cuti_2\, \zvi + s)}$ for almost any $s \in \R$, one has
  \begin{equation}\label{H2Dmod:eq:link_auxiliary_2D_3D_solutions_G}
    \displaystyle
    \aeforall (\xvi, \zvi, s) \in \Qcuti \times (0, 1), \quad \extper{\mathcal{G}^j} (\xvi, \cuti_1\, \zvi, \cuti_2\, \zvi + s) = \mathbf{g}^j_s (\xvi, \zvi),
  \end{equation}
  where $\mathcal{G}^j = \linapp{\mathcal{G}^j}{\Phi}$ and $\mathbf{g}^j_s = \linapp{\mathbf{g}^j_s}{\varphi_s}$ are the respective solutions of \eqref{H2Dmod:eq:auxiliary_local_cell_problem_G} and \eqref{H2Dmod:eq:2D_auxiliary_local_cell_problem_G}.
\end{prop}

\begin{figure}[ht!]
  \noindent\makebox[\textwidth]{%
    \def\coteCyl{4}
    \def\arrowlen{1.5mm}
    \begin{tikzpicture}[scale=0.7]
      \def\angleinterface{55}%
      \definecolor{coplancoupe}{RGB}{211, 227, 65}
      \tdplotsetmaincoords{70}{30}
      \def\BBxmax{+7}
      \def\BBymax{+13}
      \def\BBzmax{+6}
      \def\cslope{0.35355339}
      \def\szero{0.5}
      \begin{scope}[tdplot_main_coords, scale=0.9]
        %
        \node (xaxis) at (\BBxmax, 0, 0) {$\xvi$};%
        \draw[-latex] (0, 0, 0) -- (xaxis);%
        \node (yaxis) at (0, \BBymax, 0) {$\zvi_1$};%
        \draw[-latex] (0, 0, 0) -- (0, 0, \BBzmax) node[above] {$\smash{\zvi_2}$};
        
        %
        \node[tertiary] (facept) at (-0.75*\coteCyl, 0.5*\coteCyl, 0.5*\coteCyl) {$\Sigmaperper^0$};
        \draw[-{Latex[length=\arrowlen, width=\arrowlen]}, tertiary] (facept) -- (0, 0.5*\coteCyl, 0.5*\coteCyl);
        \fill[canvas is xy plane at z=0, gray!80, opacity=0.5] (0, 0) rectangle (\coteCyl, \coteCyl);
        \fill[canvas is yz plane at x=0, fill=tertiary!95, opacity=0.75] (0, 0) rectangle (\coteCyl, \coteCyl);
        \draw (0, 0, 0) -- (0, \coteCyl, 0); \draw[-latex] (0, \coteCyl, 0) -- (yaxis);%
        \draw[dashed] (0, \coteCyl, \coteCyl) -- +(0, 0, -\coteCyl) -- +(\coteCyl, 0, -\coteCyl);
        \draw[canvas is yz plane at x=0.5*\coteCyl, {Latex[length=\arrowlen, width=\arrowlen]}-, rougeTerre!80!black] (\coteCyl, 0.5*\coteCyl) arc[radius =0.675*\coteCyl, start angle=-90, end angle=0] node[above] {$\Gamma_\perperiodic^1$};
        \fill[canvas is zx plane at y=\coteCyl, rougeTerre!90, opacity=0.625] (0, 0) rectangle (\coteCyl, \coteCyl);
        \fill[canvas is yz plane at x=\coteCyl, fill=tertiary!85, opacity=0.875] (0, 0) rectangle (\coteCyl, \coteCyl);
        \fill[canvas is zx plane at y=0, rougeTerre!90, opacity=0.75] (0, 0) rectangle (\coteCyl, \coteCyl);
        \fill[canvas is xy plane at z=\coteCyl, fill=gray!50, opacity=0.8] (0, 0) rectangle (\coteCyl, \coteCyl);
        \node[rougeTerre!80!black] (facept) at (0.5*\coteCyl, -1.5*\coteCyl, 0.5*\coteCyl) {$\Gamma_\perperiodic^0$};
        \draw[-{Latex[length=\arrowlen, width=\arrowlen]}, rougeTerre!80!black] (facept) -- (0.5*\coteCyl, 0, 0.5*\coteCyl);
        \node[tertiary] (facept) at (1.75*\coteCyl, 0.5*\coteCyl, 0.5*\coteCyl) {$\Sigmaperper^1$};
        \draw[-{Latex[length=\arrowlen, width=\arrowlen]}, tertiary] (facept) -- (\coteCyl, 0.5*\coteCyl, 0.5*\coteCyl);
      \end{scope}%
    \end{tikzpicture} 
    \hfill
    \begin{tikzpicture}[scale=0.7]
      \def\angleinterface{55}%
      \definecolor{coplancoupe}{RGB}{211, 227, 65}
      \tdplotsetmaincoords{70}{30}
      \def\BBxmax{+7}
      \def\BBymax{+14}
      \def\BBzmax{+6}
      \def\cslope{0.675}
      \def\szero{0}
      \begin{scope}[tdplot_main_coords, scale=0.9]
        %
        \node (xaxis) at (\BBxmax, 0, 0) {$\xvi$};%
        \draw[-latex] (0, 0, 0) -- (xaxis);%
        \node[opacity=0.2] (yaxis) at (0, \BBymax, 0) {$\zvi_1$};%
        \draw[-latex, opacity=0.2] (0, 0, 0) -- (0, 0, \BBzmax) node[above] {$\smash{\zvi_2}$};
        
        \node[tertiary] (S0X) at (-0.75*\coteCyl, 0.5*\coteCyl, {\szero + 0.5*\cslope*\coteCyl}) {$\Scuti^\texttt{X}_0$};
        \draw[-{Latex[length=\arrowlen, width=\arrowlen]}, tertiary] (S0X) -- (0, 0.5*\coteCyl, {\szero + 0.5*\cslope*\coteCyl});

        %
        \draw[line width=0.5mm, tertiary] (0, 0, \szero) -- +(0, \coteCyl, \cslope*\coteCyl);
        \fill[canvas is yz plane at x=0, tertiary!20, opacity=0.2] (0, 0) rectangle (\coteCyl, \coteCyl);
        \draw[dashed, opacity=0.2] (0, 0, 0) -- (0, \coteCyl, 0); 
        \draw[-latex, opacity=0.2] (0, \coteCyl, 0) -- (yaxis);%
        \draw[dashed, opacity=0.2] (0, \coteCyl, \coteCyl) -- +(0, 0, -\coteCyl) -- +(\coteCyl, 0, -\coteCyl);
        \fill[canvas is zx plane at y=\coteCyl, rougeTerre!20, opacity=0.2] (0, 0) rectangle (\coteCyl, \coteCyl);

        \node (zvi) at (0, 1.75*\coteCyl, {\szero+1.75*\cslope*\coteCyl}) {$\zvi$};
        \draw[-latex] (0, 0, \szero) -- (zvi);
        \filldraw[draw=black, fill=teal!40, opacity=0.8] (0, 0, \szero) -- +(0, \coteCyl, \cslope*\coteCyl) -- +(\coteCyl, \coteCyl, \cslope*\coteCyl) -- +(\coteCyl, 0, 0) -- cycle;
        \draw (0.5*\coteCyl, 0.5*\coteCyl, \szero+0.5*\cslope*\coteCyl) node[scale=0.875] {$\Qcuti$};
        %
        %
        %
        \fill[canvas is yz plane at x=\coteCyl, tertiary!80, opacity=0.1] (0, 0) rectangle (\coteCyl, \coteCyl);
        \fill[canvas is zx plane at y=0, rougeTerre!80, opacity=0.1] (0, 0) rectangle (\coteCyl, \coteCyl);
        \draw[canvas is yz plane at x=\coteCyl, dotted] (0, 0) rectangle (\coteCyl, \coteCyl);
        \draw[canvas is xy plane at z=\coteCyl, dotted] (0, 0) rectangle (\coteCyl, \coteCyl);
        \draw[line width=0.5mm, tertiary] (\coteCyl, 0, \szero) -- +(0, \coteCyl, \cslope*\coteCyl);
        \node[tertiary] (S1X) at (1.75*\coteCyl, 0.5*\coteCyl, {\szero + 0.5*\cslope*\coteCyl}) {$\Scuti^\texttt{X}_1$};
        \draw[-{Latex[length=\arrowlen, width=\arrowlen]}, tertiary] (S1X) -- (\coteCyl, 0.5*\coteCyl, {\szero + 0.5*\cslope*\coteCyl});
        \draw[line width=0.5mm, rougeTerre!80!black] (0, 0, 0) -- (\coteCyl, 0, 0);
        \draw[line width=0.5mm, rougeTerre!80!black] (0, \coteCyl, {\szero + \cslope*\coteCyl}) -- (\coteCyl, \coteCyl, {\szero + \cslope*\coteCyl});
        \node[rougeTerre!80!black] (S0Z) at (0.5*\coteCyl, -0.5*\coteCyl, {\szero-0.5*\cslope*\coteCyl}) {$\Scuti^\texttt{Z}_0$};
        \draw[-{Latex[length=\arrowlen, width=\arrowlen]}, rougeTerre!80!black] (S0Z) -- (0.5*\coteCyl, 0, 0);
        \node[rougeTerre!80!black] (S1Z) at (0.5*\coteCyl, 1.5*\coteCyl, {\szero + 1.5*\cslope*\coteCyl}) {$\Scuti^\texttt{Z}_1$};
        \draw[-{Latex[length=\arrowlen, width=\arrowlen]}, rougeTerre!80!black] (S1Z) -- (0.5*\coteCyl, \coteCyl, {\szero + \cslope*\coteCyl});
      \end{scope}%
    \end{tikzpicture} 
  }
  \caption{Left: The cell $\cellperper$ with its faces $\Sigmaperper^\ell, \Gamma_\perperiodic^\ell$. Right: $\Qcuti$ and its edges $\Scuti^\texttt{X}_\ell, \Scuti^\texttt{Z}_\ell, \ell \in \{0, 1\}$.\label{H2Dmod:fig:coupled_traces_faces}}
\end{figure}

\vspace{0\baselineskip} \noindent
Proposition \ref{H2Dmod:prop:fibered_structure_auxiliary_problems} shows that computing $(F, \mathcal{G}^j)$ reduces to finding $(\mathbf{f}_s, \mathbf{g}^j_s)$ for any $s \in (0, 1)$. The advantage in solving the problems satisfied by $(\mathbf{f}_s, \mathbf{g}^j_s)$ is that they are 2-dimensional, and can be solved independently from one another (with respect to $s$), and therefore in parallel.

\vspace{1\baselineskip} \noindent
Finally, let $\surligner{\mathbf{t}^{j \ell}_s}$, $\surligner{\widetilde{\mathbf{t}}^{j \ell}_s}$, and $\boldsymbol{\upsilon}^{j \ell}$, $\widetilde{\boldsymbol{\upsilon}}^{j \ell}$, $j, \ell \in \{0, 1\}$ be the auxiliary edge DtN operators defined by
\begin{equation}\label{H2Dmod:eq:auxiliary_edge_DtN_operators}
  \left\{
    \begin{array}{r@{\ :=\ }l}
      \surligner{\mathbf{t}^{j \ell}_s \psi} & \surligner{(-1)^{\ell+1} \big(\aten_s\, \gradcutFBtwoD\, \linapp{\mathbf{f}^j_s}{\psi}\big) \cdot \ev_\zvi|_{\Scuti^\mathtt{Z}_\ell}},
      \ret
      \surligner{\widetilde{\mathbf{t}}^{j \ell}_s \psi} & \surligner{(-1)^{\ell+1} \big(\aten_s\, \gradcutFBtwoD\, \linapp{\mathbf{f}^j_s}{\psi}\big) \cdot \ev_\xvi|_{\Scuti^\mathtt{X}_\ell}},
      \ret
      \boldsymbol{\upsilon}^{j \ell}_s \varphi & (-1)^{\ell+1} \big(\aten_s\, \gradcutFBtwoD\, \linapp{\mathbf{g}^j_s}{\varphi}\big) \cdot \ev_\zvi|_{\Scuti^\mathtt{Z}_\ell},
      \ret
      \widetilde{\boldsymbol{\upsilon}}^{j \ell}_s \varphi & (-1)^{\ell+1} \big(\aten_s\, \gradcutFBtwoD\, \linapp{\mathbf{g}^j_s}{\varphi}\big) \cdot \ev_\zvi|_{\Scuti^\mathtt{X}_\ell}.
    \end{array}
  \right.
\end{equation}
Then the DtN operators $\Xi^\ell$, $\widetilde{\Xi}^\ell$, $\Upsilon^{j \ell}$, $\widetilde{\Upsilon}^{j \ell}$ can be derived from the edge DtN operators $\surligner{\mathbf{t}^{j \ell}_s}$, $\surligner{\widetilde{\mathbf{t}}^{j \ell}_s}$, and $\boldsymbol{\upsilon}^{j \ell}$, $\widetilde{\boldsymbol{\upsilon}}^{j \ell}$, as highlighted in the next result, which follows directly from the weak forms of these operators and from the duality property \eqref{H2Dmod:eq:shearmap_fibered_duality} which can also be extended to the faces $\Gamma_\perperiodic^\ell$.
\begin{prop}\label{H2Dmod:prop:fibered_structure_auxiliary_DtN}
  \surligner{Recall that $\vartheta := \cuti_2 / \cuti_1$}. For $\ell \in \{0, 1\}$, let $(\Phi^\ell, \Psi^\ell) \in L^2(\Sigmaperper^\ell) \times L^2(\Gamma^\ell_\perperiodic)$ be such that for almost any $s \in \R$, $\varphi^\ell_s := \shearmap \Phi^\ell(\cdot, s) \in H^{1/2}(\Scuti^\mathtt{X}_\ell)$ and $\psi^\ell_s := \shearmap \Psi(\cdot, s) \in H^{1/2}(\Scuti^\mathtt{Z}_\ell)$. Then
  \begin{equation}\label{H2Dmod:eq:fibered_structure_auxiliary_DtN}
    \spforall j, \ell \in \{0, 1\}, \quad %
    \left\{
    \begin{array}{r@{\ =\ }l}
      \big\langle \Xi^\ell\, \Psi^j,\; \Psi^\ell \big\rangle_{\Gamma^\ell_\perperiodic} &\displaystyle \frac{1}{\cuti_1} \int_0^1 \big\langle \surligner{\mathbf{t}^{0 \ell}_s\, \psi^j_s + \mathbf{t}^{1 \ell}_s\, \psi^j_{s + \vartheta}},\; \psi^\ell_s \big\rangle_{\Scuti^\mathtt{Z}_\ell} \; ds,
      \ret
      \big\langle \widetilde{\Xi}^\ell\, \Psi^j,\; \Phi^\ell \big\rangle_{\Sigmaperper^\ell} &\displaystyle \frac{1}{\cuti_1} \int_0^1 \big\langle \surligner{\widetilde{\mathbf{t}}^{0 \ell}_s\, \psi^j_s + \widetilde{\mathbf{t}}^{1 \ell}_s\, \psi^j_{s + \vartheta}},\; \varphi^\ell_s \big\rangle_{\Scuti^\mathtt{X}_\ell} \; ds,
      \ret
      \big\langle \Upsilon^{j \ell}\, \Phi^j,\; \Psi^\ell \big\rangle_{\Gamma^\ell_\perperiodic} &\displaystyle \frac{1}{\cuti_1} \int_0^1 \big\langle \boldsymbol{\upsilon}^{j \ell}_s\, \varphi^j_s,\; \psi^\ell_s \big\rangle_{\Scuti^\mathtt{Z}_\ell} \; ds,
      \ret
      \big\langle \widetilde{\Upsilon}^{j \ell}\, \Phi^j,\; \Phi^\ell \big\rangle_{\Sigmaperper^\ell} &\displaystyle \frac{1}{\cuti_1} \int_0^1 \big\langle \widetilde{\boldsymbol{\upsilon}}^{j \ell}_s\, \varphi^j_s,\; \varphi^\ell_s \big\rangle_{\Scuti^\mathtt{X}_\ell} \; ds.
    \end{array}
  \right.
  \end{equation}
\end{prop}

\paragraph*{The algorithm for the quasi-2D method}~\\
Let us summarize the quasi-2D method in the following algorithm.
\begin{itemize}
  \item solve the 2D local cell problems (\ref{H2Dmod:eq:2D_auxiliary_local_cell_problem_F}, \ref{H2Dmod:eq:2D_auxiliary_local_cell_problem_G}) and compute the auxiliary edge DtN operators $\surligner{\mathbf{t}^{j \ell}_s}$, $\surligner{\widetilde{\mathbf{t}}^{j \ell}_s}$, and $\boldsymbol{\upsilon}^{j \ell}$, $\widetilde{\boldsymbol{\upsilon}}^{j \ell}$ given by \eqref{H2Dmod:eq:auxiliary_edge_DtN_operators} for any $s \in [0, 1]$;
  \item Compute the auxiliary DtN operators $\Xi^j$, $\widetilde{\Xi}^\ell$, $\Upsilon^{j \ell}$, $\widetilde{\Upsilon}^{j \ell}$ using \eqref{H2Dmod:eq:fibered_structure_auxiliary_DtN};
  \item Determine the DtD operator $\mathcal{R}^j$, $j \in \{0, 1\}$, by solving the linear equation \eqref{H2Dmod:eq:DtD_from_auxiliary_DtN};
  \item Deduce $E^j_\fbvar$, and compute the local DtN operators $\mathcal{T}^{j \ell}_\fbvar$ using their expression \eqref{H2Dmod:eq:DtN_from_auxiliary_DtN} with respect to $\mathcal{R}^j$ and $\Xi^j$, $\widetilde{\Xi}^\ell$, $\Upsilon^{j \ell}$, $\widetilde{\Upsilon}^{j \ell}$.
\end{itemize}

\noindent
We do not describe the discretization of this algorithm, since it is an extension of the quasi-1D method. We note that \surligner{in practice, the discrete counterpart of the equation \eqref{H2Dmod:eq:DtD_from_auxiliary_DtN} used to compute $\mathcal{R}^j$ is a dense system of size $N^2 \times N^2$, where $N$ is the number of degrees of freedom on the edge $\{\xvi \in [0, 1],\ \zvi_1 = \zvi_2 = 0\}$}. \surligner{In other words, instead of solving $3$D problems of size $N^3$, the quasi-2D method allows to solve $N$ independent $2$D cell problems of size $N^2$ in parallel, and invert a $N^2$--sized dense matrix which encodes the coupling between the $2$D solutions.} Once the DtN operators $\mathcal{T}^{j \ell}_\fbvar$ have been approximated, they can be used to derive the propagation operator $\mathcal{P}_\fbvar$ (using the Riccati equation \eqref{H2Dmod:eq:riccati}), and the DtN operator $\Lambdahat^+_\fbvar$ (using \eqref{H2Dmod:eq:DtN_from_local_DtN}). \surligner{Moreover, the local cell solutions $E^j_\fbvar$ can be used to reconstruct the half-guide solution $\widehat{U}^+_\fbvar$ cell by cell.}%

}{%
}

\ifthenelse{\boolean{brouillon}}{%
  \section{Numerical results}\label{H2Dmod:sec:numerical_results}
This section provides a series of numerical results with the goal to validate the method in various situations. For the sole sake of simplicity, simulations are performed with the tensor $\aten = \mathbb{I}_2$. Unless otherwise specified, the jump data $g$ is a cut-off function, and the augmented data $G$ is constant with respect to $\zvi_2$:
\begin{equation}\label{H2Dmod:eq:jump_data_and_augmented_jump_data}
  \displaystyle
  \left\{
    \begin{array}{r@{\quad}r@{\ :=\ }l}
      \spforall \zvi \in \R, & g(0, \zvi) & 100\, \phi(2 \zvi), \quad \textnormal{with} \quad \phi(\zvi) := \exp\big( 1 - 1 / (1 - \zvi^2)\big) \, \mathds{1}_{[-1, 1]} (\zvi),
      \ret
      \spforall \zvi_1 \in \R, & G(0, \zvi_1, \zvi_2) & g(0, \zvi_1 / \cuti_1).
    \end{array}
  \right.
\end{equation}
Simulations are carried out using Lagrange finite elements of order $1$ on a regular triangular mesh, and we use the quasi-2D method to solve the local cell problems.

\subsection{Validation in the homogeneous setting}
In the first example, we consider the case where $\rho$ is piecewise constant:
\[
  \spforall \xv \in \R^2, \quad \rho^+(\xv) := 1 \quad \textnormal{and} \quad \rho^-(\xv) := 2.
\]
This coefficient falls within the scope of both Configurations \eqref{H2Dmod:item:config_a} and \eqref{H2Dmod:item:config_b} with any periodicity parameter. For Configuration \eqref{H2Dmod:item:config_a}, we choose $p^+_\zvi := 1$ and $p^-_\zvi := \sqrt{2}$ to be the periods on both sides of the interface. For Configuration \eqref{H2Dmod:item:config_b}, we choose $\pv^+ := (\sqrt{2}, 1)$ as the periodicity vector in $\R^2_+$. The jump data $g$ and its extension $G$ are given by \eqref{H2Dmod:eq:jump_data_and_augmented_jump_data}.

\vspace{1\baselineskip} \noindent
The reference solution $u_{\textit{ref}}$ used in this context is obtained by applying the partial Fourier transform with respect to $\zvi$, and by solving analytically a family of transmission problems defined on $\R$ and parameterized by the Fourier dual variable. In fact, it can be computed that
\[\displaystyle 
  \aeforall (\xvi, \zvi) \in \R^2_\pm, \quad u_{\textit{ref}}\vts (\xvi, \zvi) = \frac{1}{2\pi} \int_{\R} \frac{\widehat{g}_\zeta}{r^+_\zeta + r^-_\zeta} \; \euler^{\mp r^\pm_\zeta\vts \xvi + \icplx \vts \zeta\vts \zvi} \; d\zeta, \quad \textnormal{with} \quad \widehat{g}_\zeta := \int_{\R} g(\zvi)\; \euler^{-\icplx\vts \zeta\vts \zvi}\; d\zvi,
\]
and where $r^\pm_\zeta$ are defined for any $\zeta \in \R$ by $(r^\pm_\zeta)^2 = \zeta^2 - \rho^\pm\, \omega^2$, $\Real r^\pm_\zeta \geq 0$.

\vspace{1\baselineskip} \noindent
The solutions that follow from applying the lifting approach to Configurations \eqref{H2Dmod:item:config_a} and \eqref{H2Dmod:item:config_b} are compared in Figure \ref{H2Dmod:fig:solution_homogeneous_setting} to $u_\textit{ref}$. We first choose $\omega = 8 + 0.25\,\icplx$ and a mesh step $h = 0.025$, which corresponds approximately to $31$ points per wavelength in $\R^2_+$ and $22$ points per wavelength in $\R^2_-$. The number of Floquet points is set to $N_\fbvar = 64$. The visual similarity between the results validates qualitatively the method in the homogeneous setting.

\begin{figure}[ht!]
  \centering
  \makebox[0pt][c]{%
    \centering
    \def\Ldomain{5}
    \begin{tikzpicture}
      \begin{groupplot}[
        group style={
          group size=3 by 1,
          horizontal sep=0.75cm,
          vertical sep=1.5cm,
        },
        enlargelimits=false,
        axis on top,
        width=0.3\textwidth,
        height=0.3\textwidth,
        scale only axis,
        xmin=-\Ldomain, xmax=\Ldomain, ymin=-\Ldomain, ymax=\Ldomain,
        xtick = {-4, -2, 0, 2, 4}, ytick=\empty,
        disabledatascaling,
        axis equal,
        clip mode=individual, 
      ]
        \nextgroupplot[title={$u_{\textit{ref}}$}]
        \addplot graphics [xmin=-\Ldomain, xmax=\Ldomain, ymin=-\Ldomain, ymax=\Ldomain] {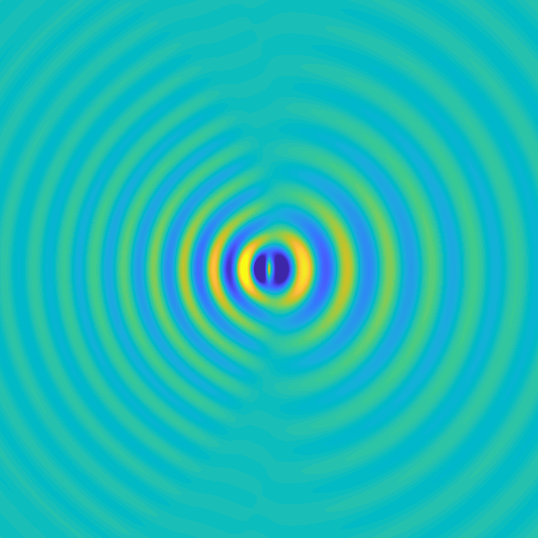};
        \draw[thick, dashed] (axis cs:0, -\Ldomain) -- (axis cs:0, \Ldomain);
        \nextgroupplot[%
          title={$u_\hv$, Configuration \eqref{H2Dmod:item:config_a}},
          colorbar,
          point meta min=-2.5, point meta max=2.5,
          colorbar horizontal,
          colorbar style={
            xticklabel style={
              /pgf/number format/precision=2,
              /pgf/number format/fixed,
            },
            xtick={-2, 0, 2},
            axis equal=false,
            major tick length=0.025*\pgfkeysvalueof{/pgfplots/parent axis width},
            scaled x ticks=false,
            at = {(0, -0.75cm)},
            anchor=south west,
            width = \pgfkeysvalueof{/pgfplots/parent axis  width},
            height = 0.025*\pgfkeysvalueof{/pgfplots/parent axis height},
          },
        ]
        \addplot graphics [xmin=-6, xmax=6, ymin=-6, ymax=6]{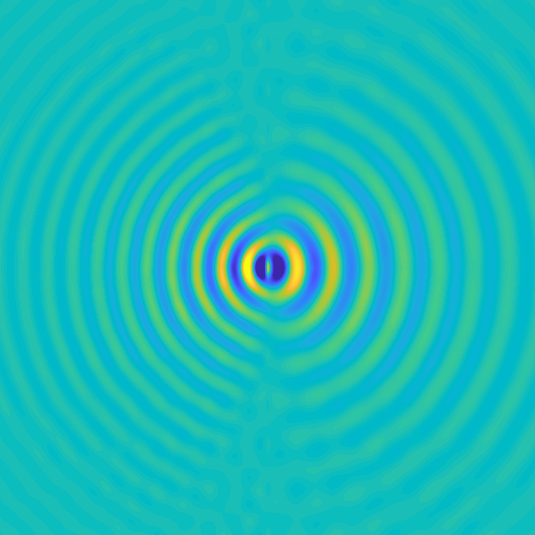};
        \draw[thick, dashed] (axis cs:0, -\Ldomain) -- (axis cs:0, \Ldomain);
        \nextgroupplot[title={$u_\hv$, Configuration \eqref{H2Dmod:item:config_b}}]
        \addplot graphics [xmin=-6, xmax=6, ymin=-6, ymax=6]{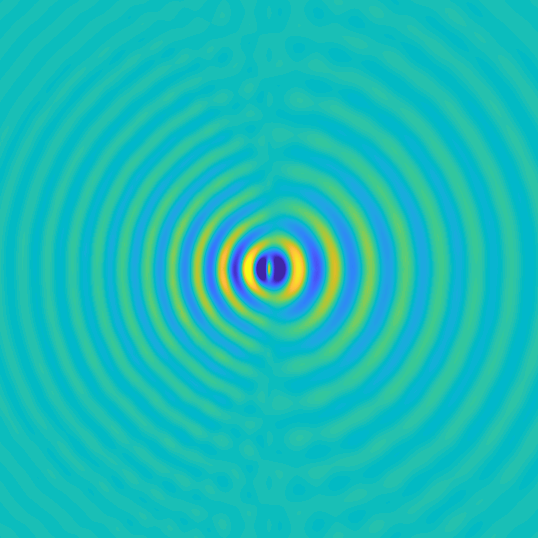};
        \draw[thick, dashed] (axis cs:0, -\Ldomain) -- (axis cs:0, \Ldomain);
      \end{groupplot}
    \end{tikzpicture}
  }
  \caption{Approximate solution $u_\hv$ (real part) computed for Configurations \eqref{H2Dmod:item:config_a} and \eqref{H2Dmod:item:config_b}, and compared to a reference solution $\textit{u}_\textit{ref}$, for $\omega = 8 + 0.25\vts \icplx$. We use order $1$ Lagrange finite elements with $h = 0.025$ and $N_\fbvar = 64$.\label{H2Dmod:fig:solution_homogeneous_setting}}
\end{figure}

\vspace{1\baselineskip}\noindent
Next, we fix $\omega = 1 + 0.25\,\icplx$, $N_\fbvar = 64$ and we restrict ourselves to Configuration \eqref{H2Dmod:item:config_a}. Figure \ref{H2Dmod:fig:relative_errors_homogeneous_setting} (right) shows the relative errors
\begin{equation}\label{H2Dmod:eq:def_relative_errors}
  \displaystyle
  \veps^{(0)}_h := \frac{\|u_\hv - u_\textit{ref}\|_{L^2(\Omega_0)}}{\|u_\textit{ref}\|_{L^2(\Omega_0)}} \quad \textnormal{and} \quad \veps^{(1)}_h := \frac{\|\nabla (u_\hv - u_\textit{ref})\|_{[L^2(\Omega_0)]^2}}{\|\nabla u_\textit{ref}\|_{[L^2(\Omega_0)]^2}}, \quad \Omega_0 := (-1, 1)^2,
\end{equation}
which decay as the mesh step $h$ tends to $0$. As expected for order $1$ Lagrange finite elements, $\veps^{(0)}_h$ tends to $0$ as $h^2$, whereas $\veps^{(1)}_h$ tends to $0$ as $h^2$ (instead of the expected $O(h)$ decay). \surligner{We believe that this superconvergence phenomenon is linked to our use of a regular mesh, although its precise cause remains unclear.}

\begin{figure}[ht!]
	\centering
  \def\Ldomain{5}
	\makebox[\textwidth][c]{
    \centering
    \begin{subfigure}[t]{0.48\textwidth}
      \begin{tikzpicture}
        \begin{axis}[%
          enlargelimits=false,
          axis on top,
          width=0.625\textwidth, 
          height=0.625\textwidth, 
          title={$u_\hv - u_\textit{ref}$,\quad $h = 0.025$},
          scale only axis,
          xmin=-\Ldomain, xmax=\Ldomain, ymin=-\Ldomain, ymax=\Ldomain,
          xtick = {-4, -2, 0, 2, 4}, ytick=\empty,
          disabledatascaling,
          axis equal,
          clip mode=individual, 
          colorbar,
          point meta min=-0.05, point meta max=0.05,
          colorbar style={
            ytick={-0.05, 0.05},
            axis equal=false,
            major tick length=0.025*\pgfkeysvalueof{/pgfplots/parent axis width},
            scaled y ticks=false,
            at = {(0.25cm+\pgfkeysvalueof{/pgfplots/parent axis  width}, 0)},
            anchor=south west,
            width  = 0.025*\pgfkeysvalueof{/pgfplots/parent axis  width},
            height = \pgfkeysvalueof{/pgfplots/parent axis height},
          },
        ]
          \addplot graphics [xmin=-5, xmax=5, ymin=-5, ymax=5]{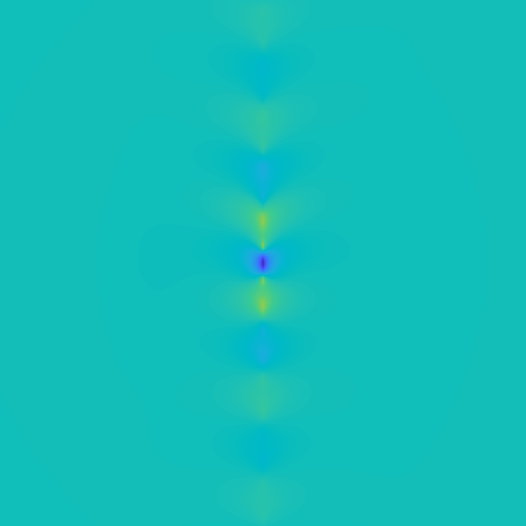};
          \draw[thick, dashed] (axis cs:0, -\Ldomain) -- (axis cs:0, \Ldomain);
        \end{axis}
      \end{tikzpicture}
    \end{subfigure}
    \hfill
    \begin{subfigure}[t]{0.5\textwidth}    
      \includegraphics[page=7]{H2Dmod/tikzpicture_H2Dmod.pdf}
    \end{subfigure}
  }
	\caption{$\omega = 1 + 0.25\vts \icplx$ and $N_k = 64$ are fixed. Left: real part of $u_\hv - u_{\textit{ref}}$ for $h = 0.025$. Right: relative errors $\veps^{(0)}_h$ and $\veps^{(1)}_h$ (see \eqref{H2Dmod:eq:def_relative_errors}) with respect to the mesh step $h$, for Configuration \eqref{H2Dmod:item:config_a}. \label{H2Dmod:fig:relative_errors_homogeneous_setting}}
\end{figure}

\vspace{1\baselineskip} \noindent
In the next sections, numerical experiments are performed with a variable coefficient $\rho$. Using the cut-off function $\phi \in \mathscr{C}^\infty_0(\R)$ defined in \eqref{H2Dmod:eq:jump_data_and_augmented_jump_data}, we start from the $\Z^2$--periodic functions defined in one periodicity cell by
\begin{equation*}
  \displaystyle
  \spforall \accentset{\circ}{\xv} = (\accentset{\circ}{\xvi}, \accentset{\circ}{\zvi}) \in (0, 1)^2, \quad \accentset{\circ}{\rho}^-(\accentset{\circ}{\xv}) := 0.5 + \phi(4\accentset{\circ}{\xvi})\, \phi(4\accentset{\circ}{\zvi}) \quad \textnormal{and} \quad \accentset{\circ}{\rho}^+(\accentset{\circ}{\xv}) := 0.5 + \phi(2.5 |\accentset{\circ}{\xv}|).
\end{equation*}
Then, for Configuration \eqref{H2Dmod:item:config_a}, we use
\begin{equation}\label{H2Dmod:eq:def_rho_0_config_a}
  \spforall \xv = (\xvi, \zvi) \in \R^2, \quad \rho^\pm (\xv) := \accentset{\circ}{\rho}^\pm (\xvi, \zvi / p^\pm_\zvi),
\end{equation}
which is $p^\pm_\zvi$--periodic with respect to $\zvi$. For Configuration \eqref{H2Dmod:item:config_b}, we set
\begin{equation}\label{H2Dmod:eq:def_rho_0_config_b}
  \displaystyle
  \spforall \xv = (\xvi, \zvi) \in \R^2, \quad \rho^-(\xv) = 1 \quad \textnormal{and} \quad \rho^+ (\xv) := \accentset{\circ}{\rho}(\xvi - (p^+_\xvi / p^+_\zvi)\vts \zvi, \zvi / p^+_\zvi),
\end{equation}
so that $\rho^+$ is $\Z\vts \ev_\xvi + \Z\, \pv^+$--periodic with $\pv^+ = (p^+_\xvi, p^+_\zvi)$. These coefficients are shown in Figure \ref{H2Dmod:fig:coefficients}. 
\begin{figure}[ht!]
  \hspace{-11pt}
  \centering
  \makebox[0pt][c]{%
    \def\echelle{0.75} 
    \def\Ldomain{5}
    \begin{tikzpicture}
      \begin{groupplot}[
        group style={
          group size=2 by 1,
          horizontal sep=3.5cm,
          vertical sep=1.5cm,
        },
        enlargelimits=false,
        axis on top,
        width=0.3\textwidth,
        height=0.3\textwidth,
        scale only axis,
        xmin=-\Ldomain, xmax=\Ldomain, ymin=-\Ldomain, ymax=\Ldomain,
        xtick = {-4, -2, 0, 2, 4}, ytick = {-4, -2, 0, 2, 4},
        disabledatascaling,
        axis equal,
        clip mode=individual, 
      ]
        \nextgroupplot[title={$\rho$, Configuration \eqref{H2Dmod:item:config_a}},
          colorbar,
          point meta min=0.5, point meta max=1.5,
          colorbar style={
            yticklabel style={
              /pgf/number format/precision=2,
              /pgf/number format/fixed,
            },
            ytick={0.6, 1, 1.4},
            axis equal=false,
            major tick length=0.025*\pgfkeysvalueof{/pgfplots/parent axis width},
            scaled y ticks=false,
            at = {(0.25cm+\pgfkeysvalueof{/pgfplots/parent axis  width}, 0)},
            anchor=south west,
            width  = 0.025*\pgfkeysvalueof{/pgfplots/parent axis  width},
            height = \pgfkeysvalueof{/pgfplots/parent axis height},
          },]
        \addplot graphics [xmin=-6, xmax=6, ymin=-6, ymax=6]{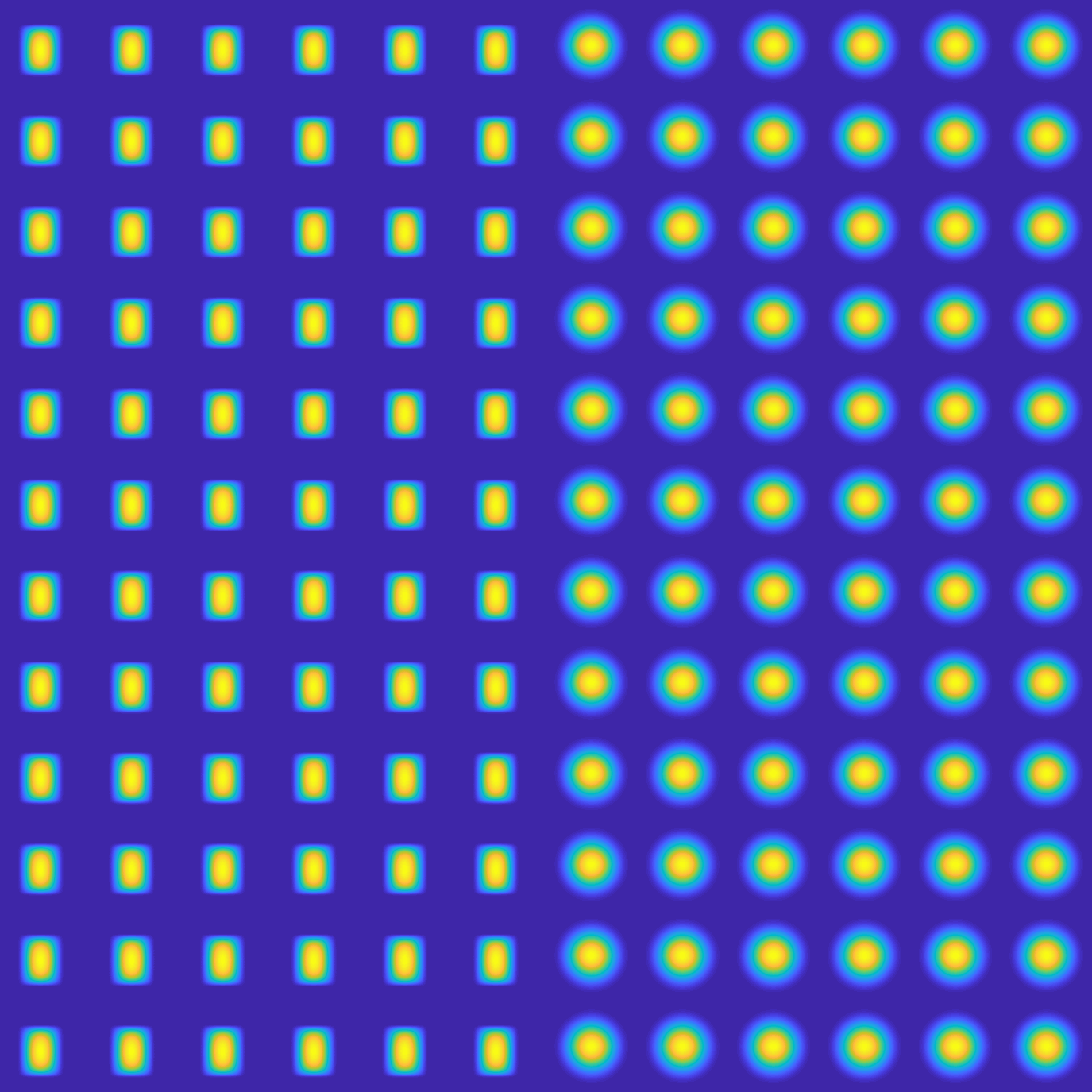};
        \pgfplotsinvokeforeach{1,2,...,5}{
          \draw[white, opacity=0.4] (axis cs:#1, -\Ldomain) -- (axis cs:#1, \Ldomain);
          \draw[white, opacity=0.4] (axis cs:-#1, -\Ldomain) -- (axis cs:-#1, \Ldomain);
        }
        \pgfplotsinvokeforeach{-5,-4,...,5}{
          \draw[white, opacity=0.4] (axis cs:-\Ldomain, #1) -- (axis cs:\Ldomain, #1);
        }
        \draw[thick, white, dashed] (axis cs:0, -\Ldomain) -- (axis cs:0, \Ldomain);
        \nextgroupplot[
          title={$\rho$, Configuration \eqref{H2Dmod:item:config_b}},
          colorbar,
          point meta min=0.5, point meta max=1.5,
          colorbar style={
            yticklabel style={
              /pgf/number format/precision=2,
              /pgf/number format/fixed,
            },
            ytick={0.6, 1, 1.4},
            axis equal=false,
            major tick length=0.025*\pgfkeysvalueof{/pgfplots/parent axis width},
            scaled y ticks=false,
            at = {(0.25cm+\pgfkeysvalueof{/pgfplots/parent axis  width}, 0)},
            anchor=south west,
            width  = 0.025*\pgfkeysvalueof{/pgfplots/parent axis  width},
            height = \pgfkeysvalueof{/pgfplots/parent axis height},
          },
        ]
        \addplot graphics [xmin=-6, xmax=6, ymin=-6, ymax=6]{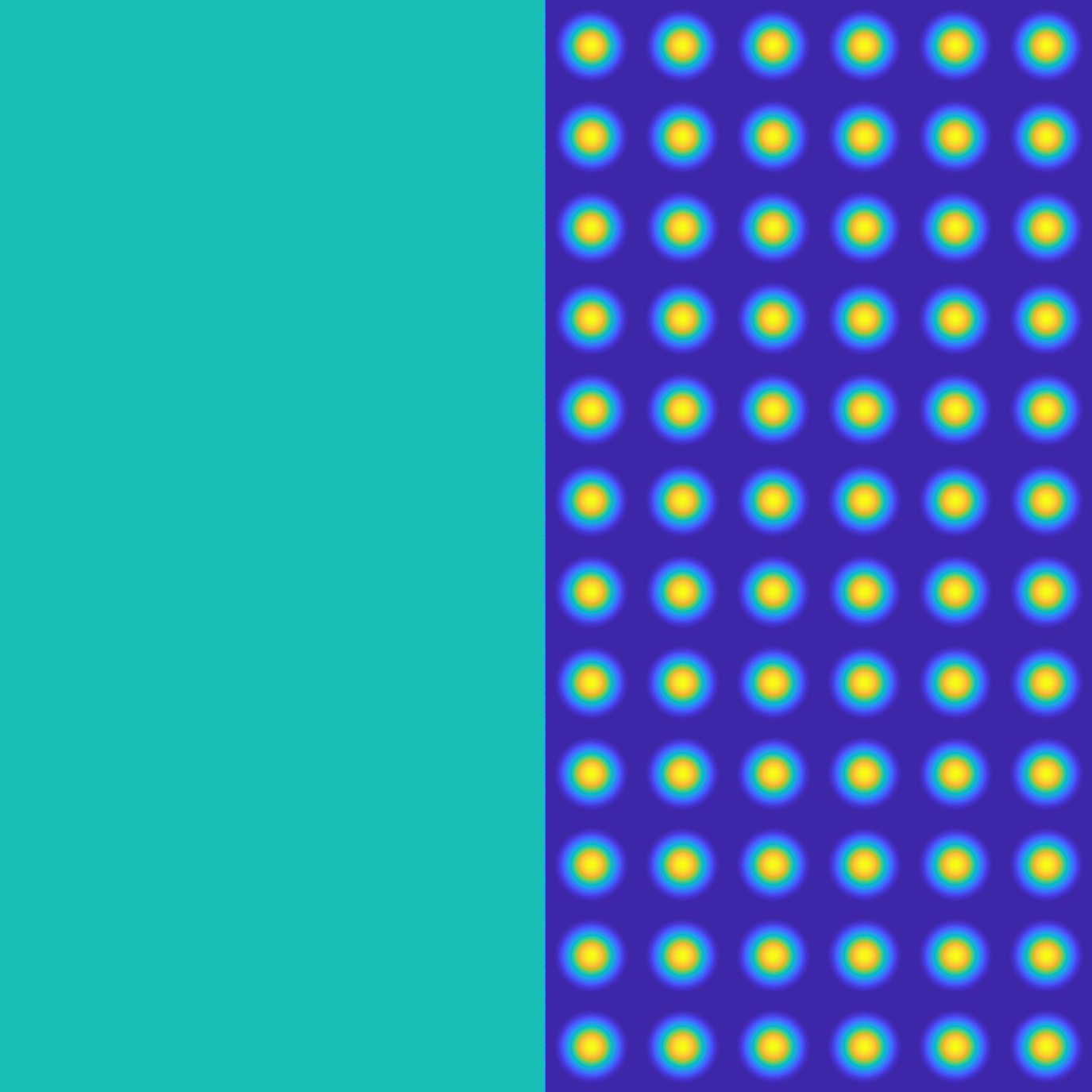};
        \pgfplotsinvokeforeach{1,2,...,5}{
          \draw[white, opacity=0.4] (axis cs:#1, -\Ldomain) -- (axis cs:#1, \Ldomain);
        }
        \pgfplotsinvokeforeach{-5,-4,...,5}{
          \draw[white, opacity=0.4] (axis cs:0, #1) -- (axis cs:\Ldomain, #1);
        }
        \draw[thick, white, dashed] (axis cs:0, -\Ldomain) -- (axis cs:0, \Ldomain);
      \end{groupplot}
    \end{tikzpicture}
  }
  \caption{The coefficient $\rho$ defined for Configuration \eqref{H2Dmod:item:config_a} by \eqref{H2Dmod:eq:def_rho_0_config_a}, with $p^+_\zvi = p^-_\zvi = 1$, and for Configuration \eqref{H2Dmod:item:config_b} by \eqref{H2Dmod:eq:def_rho_0_config_b}, with $\pv^+ = \ev_\zvi$.\label{H2Dmod:fig:coefficients}}
\end{figure}

\subsection{Validation in the rational setting}
We consider a coefficient $\rho$ which is $1$--periodic in the direction of the interface (see Figures \ref{H2Dmod:fig:solution_rational_setting_config_a} and \ref{H2Dmod:fig:solution_rational_setting_config_b} left). In this case, as done in \cite{flisscassanbernier2010}, one can directly apply a Floquet-Bloch transform in the direction of the interface, leading to a family of 2D transmission problems defined in $\R \times (0, 1)$ and parameterized by the Floquet dual variable. Each of the waveguide problems can then be reduced to an interface equation featuring 2D DtN operators, which we obtain by computing the solution of half-guide problems defined in $\R_\pm \times (0, 1)$. Solving these half-guide problems involves 2D local cell problems and a propagation operator, similarly to Section \ref{H2Dmod:sec:resolution_half_guide_problems}. We use this approach to construct a reference solution $u_{\textit{ref}}$, to which we compare the solution $u$ obtained using the lifting approach. We fix $\omega = 8 + 0.25\vts \icplx$.

\vspace{1\baselineskip} \noindent
For Configuration \eqref{H2Dmod:item:config_a}, we define $\rho$ using \eqref{H2Dmod:eq:def_rho_0_config_a}, with $p^+_\zvi = p^-_\zvi = 1$. Figure \ref{H2Dmod:fig:solution_rational_setting_config_a} shows the solution $u_\hv$ computed numerically using the lifting approach for a mesh step $h = 0.025$ and $N_\fbvar = 64$. As expected, this solution is close to the reference solution $u_{\textit{ref}}$.
\begin{figure}[ht!]
  \hspace{-11pt}
  \centering
  \makebox[0pt][c]{%
    \def\echelle{0.75}
    \def\Ldomain{5}
    \begin{tikzpicture}
      \begin{groupplot}[
        group style={
          group size=3 by 1,
          horizontal sep=0.75cm,
          vertical sep=1.5cm,
        },
        enlargelimits=false,
        axis on top,
        width=0.3\textwidth,
        height=0.3\textwidth,
        scale only axis,
        xmin=-\Ldomain, xmax=\Ldomain, ymin=-\Ldomain, ymax=\Ldomain,
        xtick = {-4, -2, 0, 2, 4}, ytick=\empty,
        disabledatascaling,
        axis equal,
        clip mode=individual, 
      ]
        \nextgroupplot[%
          title={$\rho$}, 
          colorbar,
          point meta min=0.5, point meta max=1.5,
          colorbar horizontal,
          colorbar style={
            xticklabel style={
              /pgf/number format/precision=2,
              /pgf/number format/fixed,
            },
            xtick={0.6, 1, 1.4},
            axis equal=false,
            major tick length=0.025*\pgfkeysvalueof{/pgfplots/parent axis width},
            scaled x ticks=false,
            at = {(0, -0.75cm)},
            anchor=south west,
            width = \pgfkeysvalueof{/pgfplots/parent axis  width},
            height = 0.025*\pgfkeysvalueof{/pgfplots/parent axis height},
          },
        ]
        \addplot graphics [xmin=-6, xmax=6, ymin=-6, ymax=6] {rho0_A_6_6.png};
        \pgfplotsinvokeforeach{1,2,...,5}{
          \draw[white, opacity=0.4] (axis cs:#1, -\Ldomain) -- (axis cs:#1, \Ldomain);
          \draw[white, opacity=0.4] (axis cs:-#1, -\Ldomain) -- (axis cs:-#1, \Ldomain);
        }
        \pgfplotsinvokeforeach{-5,-4,...,5}{
          \draw[white, opacity=0.4] (axis cs:-\Ldomain, #1) -- (axis cs:\Ldomain, #1);
        }
        \draw[thick, white, dashed] (axis cs:0, -\Ldomain) -- (axis cs:0, \Ldomain);
        \nextgroupplot[%
          title={$u_{\textit{ref}}$},
          colorbar,
          point meta min=-3, point meta max=3,
          colorbar horizontal,
          colorbar style={
            xticklabel style={
              /pgf/number format/precision=2,
              /pgf/number format/fixed,
            },
            xtick={-2, 0, 2},
            axis equal=false,
            major tick length=0.025*\pgfkeysvalueof{/pgfplots/parent axis width},
            scaled x ticks=false,
            at = {(0, -0.75cm)},
            anchor=south west,
            width = \pgfkeysvalueof{/pgfplots/parent axis  width},
            height = 0.025*\pgfkeysvalueof{/pgfplots/parent axis height},
          },
        ]
        \addplot graphics [xmin=-6, xmax=6, ymin=-6, ymax=6]{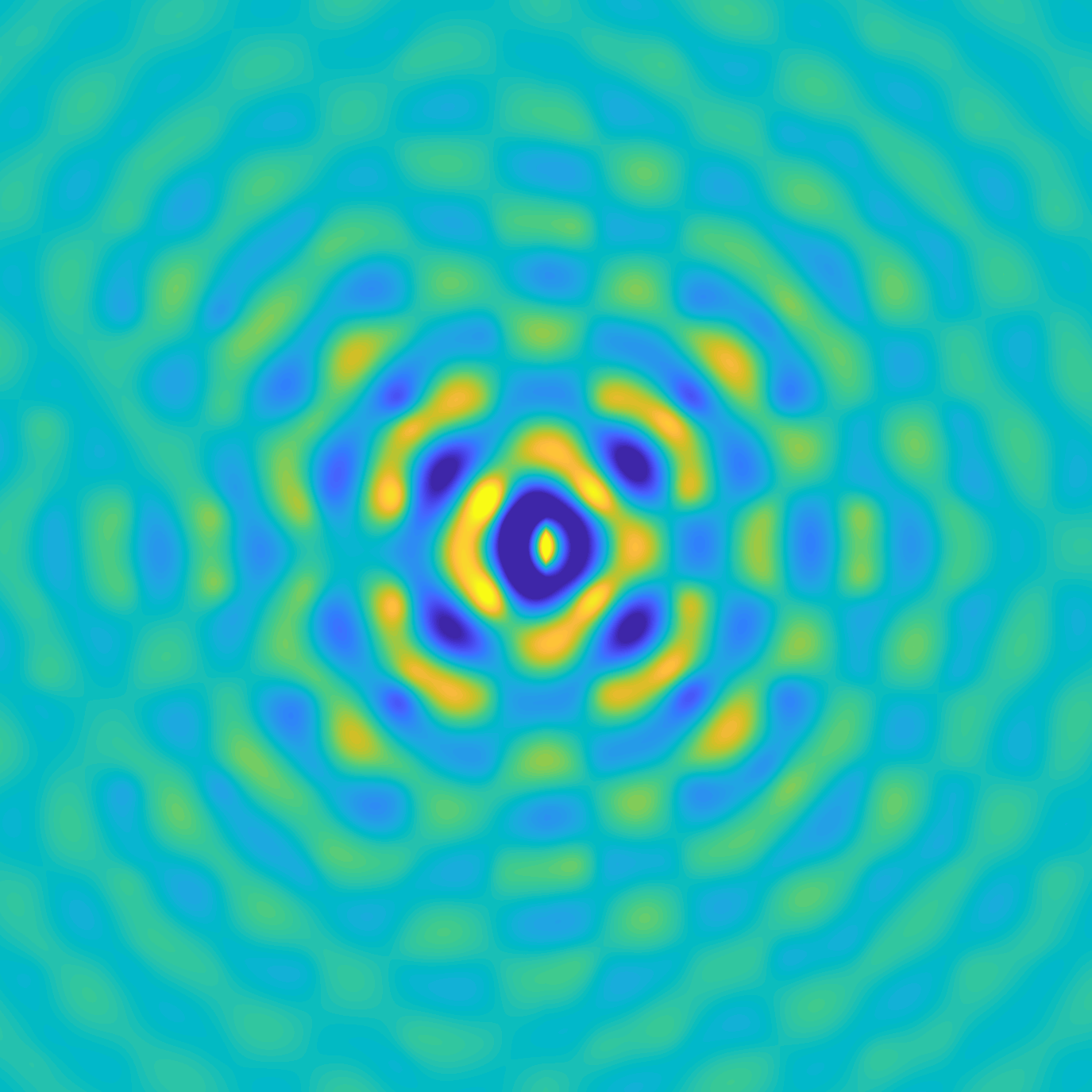};
        \draw[thick, dashed] (axis cs:0, -\Ldomain) -- (axis cs:0, \Ldomain);
        \nextgroupplot[%
          title={$u_\hv$},
          colorbar,
          point meta min=-3, point meta max=3,
          colorbar horizontal,
          colorbar style={
            xticklabel style={
              /pgf/number format/precision=2,
              /pgf/number format/fixed,
            },
            xtick={-2, 0, 2},
            axis equal=false,
            major tick length=0.025*\pgfkeysvalueof{/pgfplots/parent axis width},
            scaled x ticks=false,
            at = {(0, -0.75cm)},
            anchor=south west,
            width = \pgfkeysvalueof{/pgfplots/parent axis  width},
            height = 0.025*\pgfkeysvalueof{/pgfplots/parent axis height},
          },
        ]
        \addplot graphics [xmin=-6, xmax=6, ymin=-6, ymax=6]{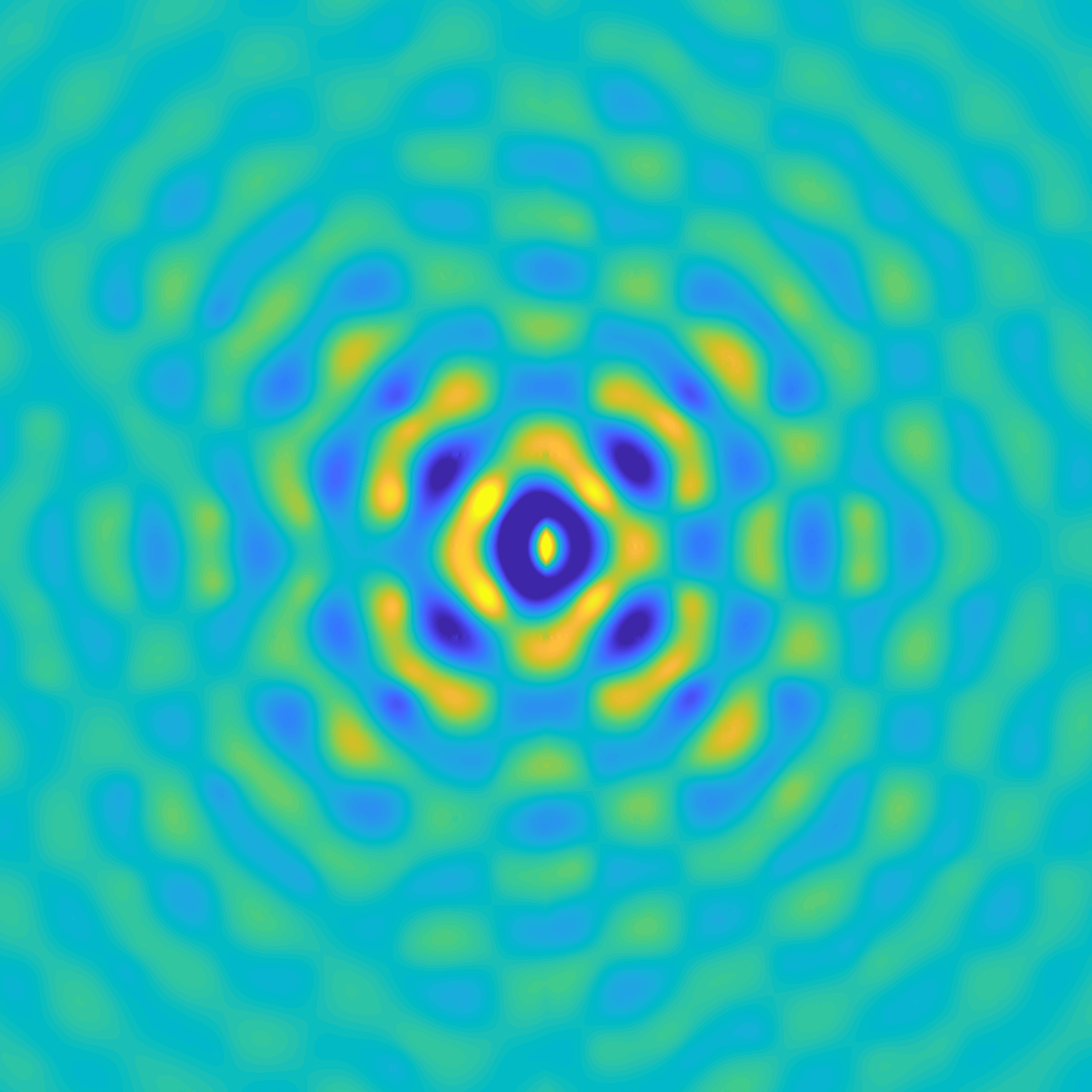};
        \draw[thick, dashed] (axis cs:0, -\Ldomain) -- (axis cs:0, \Ldomain);
      \end{groupplot}
    \end{tikzpicture}
  }
  \caption{Approximate solution $u_\hv$ (real part) computed for Configuration \eqref{H2Dmod:item:config_a} in the rational case ($p^+_\zvi = p^-_\zvi = 1$). Here, $\omega = 8 + 0.25\vts \icplx$, $h = 0.025$, and $N_\fbvar = 64$.\label{H2Dmod:fig:solution_rational_setting_config_a}}
\end{figure}

\vspace{0\baselineskip} \noindent
For Configuration \eqref{H2Dmod:item:config_b}, we define $\rho$ using \eqref{H2Dmod:eq:def_rho_0_config_b}, with $\pv^+ = (1, 1)$. Figure \ref{H2Dmod:fig:solution_rational_setting_config_b} shows the solution $u_\hv$ computed numerically using the lifting approach for a mesh step $h = 0.025$ and $N_\fbvar = 64$. As expected, this solution is close to the reference solution $u_{\textit{ref}}$.

\begin{figure}[ht!]
  \hspace{-11pt}
  \centering
  \makebox[0pt][c]{%
    \def\echelle{0.75}
    \def\Ldomain{5}
    \begin{tikzpicture}
      \begin{groupplot}[
        group style={
          group size=3 by 1,
          horizontal sep=0.75cm,
          vertical sep=1.5cm,
        },
        enlargelimits=false,
        axis on top,
        width=0.3\textwidth,
        height=0.3\textwidth,
        scale only axis,
        xmin=-\Ldomain, xmax=\Ldomain, ymin=-\Ldomain, ymax=\Ldomain,
        xtick = {-4, -2, 0, 2, 4}, ytick=\empty,
        disabledatascaling,
        axis equal,
        clip mode=individual, 
      ]
        \nextgroupplot[%
          title={$\rho$}, 
          colorbar,
          point meta min=0.5, point meta max=1.5,
          colorbar horizontal,
          colorbar style={
            xticklabel style={
              /pgf/number format/precision=2,
              /pgf/number format/fixed,
            },
            xtick={0.6, 1, 1.4},
            axis equal=false,
            major tick length=0.025*\pgfkeysvalueof{/pgfplots/parent axis width},
            scaled x ticks=false,
            at = {(0, -0.75cm)},
            anchor=south west,
            width = \pgfkeysvalueof{/pgfplots/parent axis  width},
            height = 0.025*\pgfkeysvalueof{/pgfplots/parent axis height},
          },
        ]
        \addplot graphics [xmin=-6, xmax=6, ymin=-6, ymax=6] {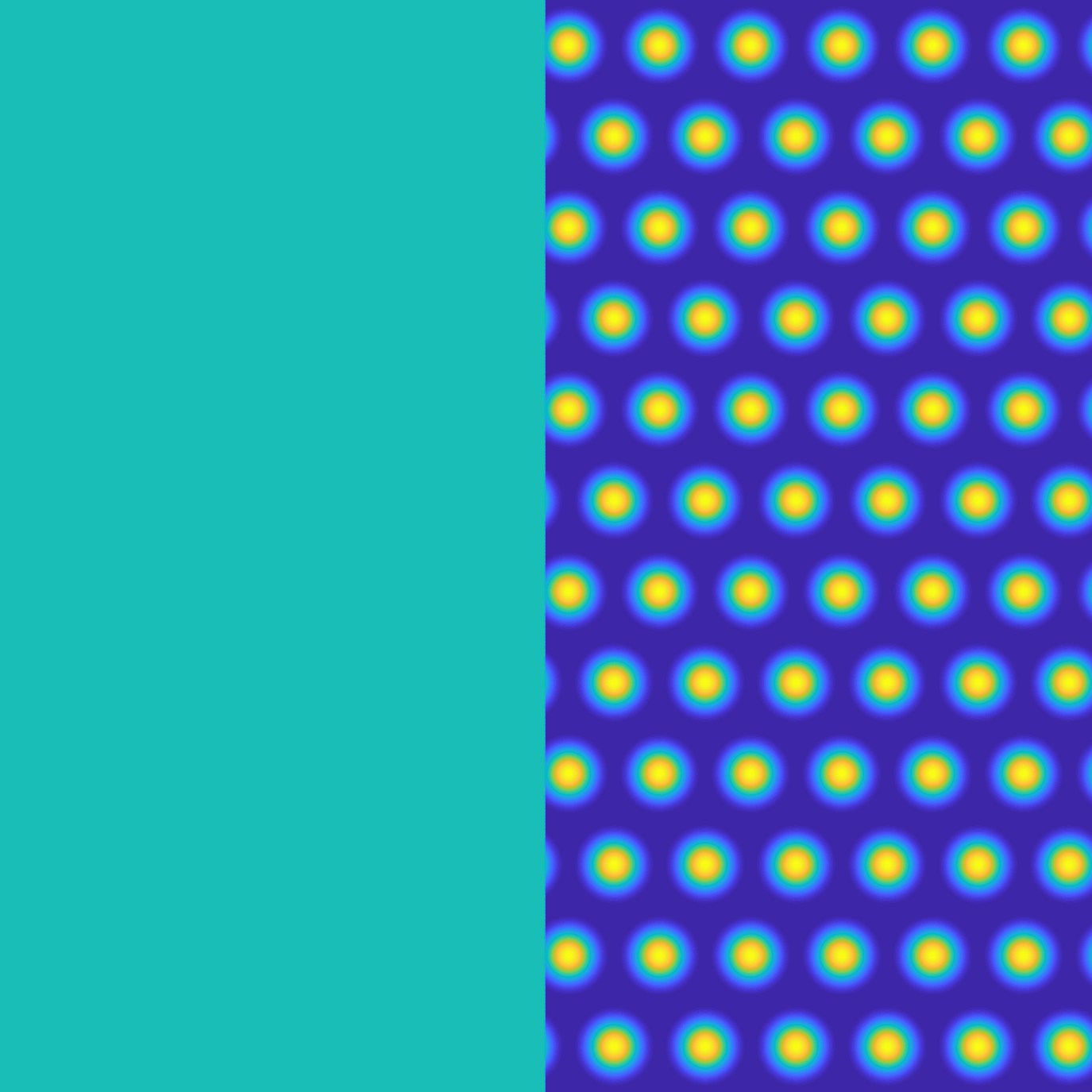};
        \draw[thick, white, dashed] (axis cs:0, -\Ldomain) -- (axis cs:0, \Ldomain);
        
        \def\vecperX{1}
        \def\vecperZ{1}
        \draw[thick, white] (axis cs:1.75, 2) -- +(1, 0) -- +({1+\vecperX}, {\vecperZ}) -- +(\vecperX, \vecperZ) -- cycle;
        \draw[thick, white, dashed] (axis cs:2.25, -2) -- +(1, 0) -- +(1, 1) -- +(0, 1) -- cycle;
        \draw[thick, white, dashed] (axis cs:0, -\Ldomain) -- (axis cs:0, \Ldomain);
        \nextgroupplot[%
          title={$u_{\textit{ref}}$},
          colorbar,
          point meta min=-3, point meta max=3,
          colorbar horizontal,
          colorbar style={
            xticklabel style={
              /pgf/number format/precision=2,
              /pgf/number format/fixed,
            },
            xtick={-2, 0, 2},
            axis equal=false,
            major tick length=0.025*\pgfkeysvalueof{/pgfplots/parent axis width},
            scaled x ticks=false,
            at = {(0, -0.75cm)},
            anchor=south west,
            width = \pgfkeysvalueof{/pgfplots/parent axis  width},
            height = 0.025*\pgfkeysvalueof{/pgfplots/parent axis height},
          },
        ]
        \addplot graphics [xmin=-6, xmax=6, ymin=-6, ymax=6]{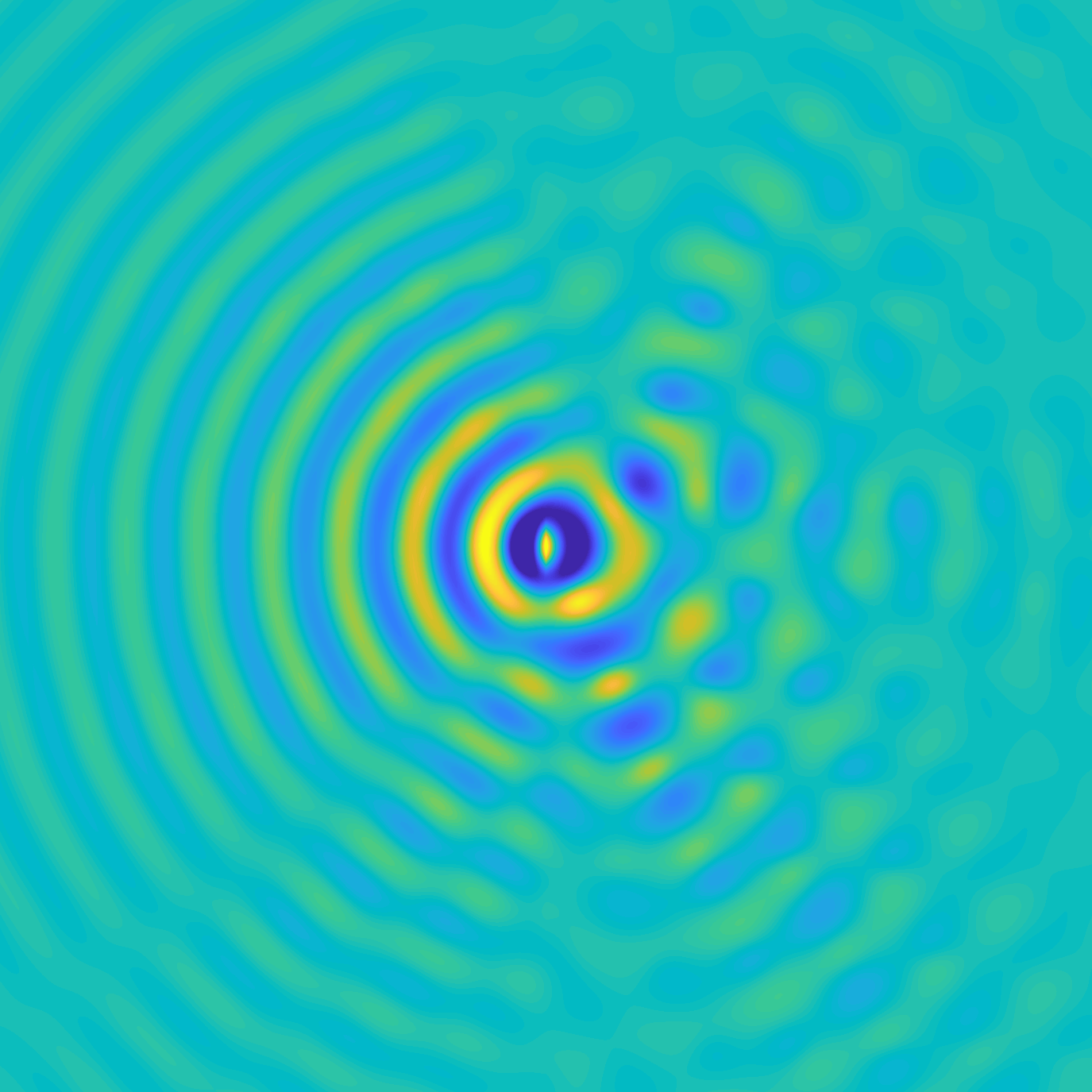};
        \draw[thick, white, dashed] (axis cs:2.25, 0) -- +(1, 0) -- +(1, 1) -- +(0, 1) -- cycle;
        \draw[thick, dashed] (axis cs:0, -\Ldomain) -- (axis cs:0, \Ldomain);
        \nextgroupplot[%
          title={$u_\hv$},
          colorbar,
          point meta min=-2.5, point meta max=2.5,
          colorbar horizontal,
          colorbar style={
            xticklabel style={
              /pgf/number format/precision=2,
              /pgf/number format/fixed,
            },
            xtick={-2, 0, 2},
            axis equal=false,
            major tick length=0.025*\pgfkeysvalueof{/pgfplots/parent axis width},
            scaled x ticks=false,
            at = {(0, -0.75cm)},
            anchor=south west,
            width = \pgfkeysvalueof{/pgfplots/parent axis  width},
            height = 0.025*\pgfkeysvalueof{/pgfplots/parent axis height},
          },
        ]
        \addplot graphics [xmin=-6, xmax=6, ymin=-6, ymax=6]{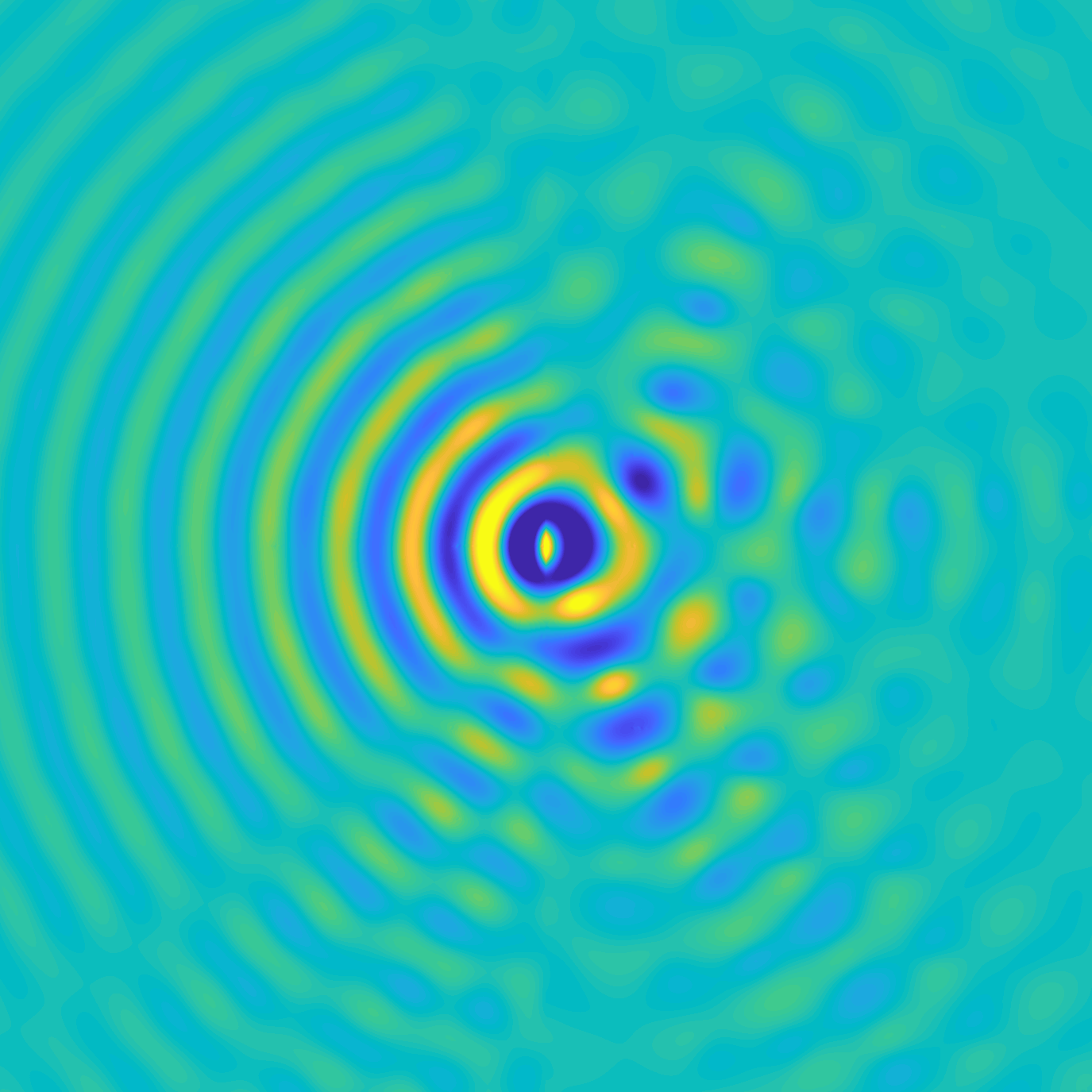};
        \def\vecperX{1}
        \def\vecperZ{1}
        \draw[thick, white] (axis cs:1.75, 0) -- +(1, 0) -- +({1+\vecperX}, {\vecperZ}) -- +(\vecperX, \vecperZ) -- cycle;
        \draw[thick, dashed] (axis cs:0, -\Ldomain) -- (axis cs:0, \Ldomain);
      \end{groupplot}
    \end{tikzpicture}
  }
  \caption{Approximate solution $u_\hv$ (real part) computed for Configuration \eqref{H2Dmod:item:config_b} in the rational case ($\pv^+ = (1, 1)$). Here, $\omega = 8 + 0.25\vts \icplx$, $h = 0.025$, and $N_\fbvar = 64$. The periodicity cells used to compute $u_{\textit{ref}}$ and $u_\hv$ are the square $(0, 1)^2$ (in dashed line) and the parallelogram $\{y_1 \pv^+ + y_2 \ev_\xvi,\ y_1, y_2 \in (0, 1)\}$ (in solid line) respectively. \label{H2Dmod:fig:solution_rational_setting_config_b}}
\end{figure}

\vspace{1\baselineskip} \noindent
\surligner{%
Finally, in Figure \ref{H2Dmod:fig:relative_errors_rational_setting}, we represent the $L^2$ and $H^1$--relative errors $\veps^{(0)}_h$ and $\veps^{(1)}_h$ (see \eqref{H2Dmod:eq:def_relative_errors}) with respect to $1/h$, for Configuration \eqref{H2Dmod:item:config_b}, $\omega = 8 + 0.25\vts \icplx$, and $N_\fbvar = 100$. The $L^2$--error decays as $h^2$, as expected for order $1$ Lagrange finite elements. However we observe that the $H^1$--error also decays as $h^2$ (instead of the order $h$ decay predicted by theory). Similarly to the homogeneous case, we believe that this superconvergence phenomenon is caused by our use of a regular mesh, although its precise cause remains unclear.

\begin{figure}[ht!]
	\centering
  \def\Ldomain{5}
	\makebox[\textwidth][c]{
    \centering
    \begin{subfigure}[t]{0.48\textwidth}
      \begin{tikzpicture}
        \begin{axis}[%
          enlargelimits=false,
          axis on top,
          width=0.625\textwidth, 
          height=0.625\textwidth, 
          title={$u_\hv - u_\textit{ref}$,\quad $h = 0.025$},
          scale only axis,
          xmin=-\Ldomain, xmax=\Ldomain, ymin=-\Ldomain, ymax=\Ldomain,
          xtick = {-4, -2, 0, 2, 4}, ytick=\empty,
          disabledatascaling,
          axis equal,
          clip mode=individual, 
          colorbar,
          point meta min=-0.002, point meta max=0.002,
          colorbar style={
            ytick={-0.002, 0.002},
            axis equal=false,
            major tick length=0.025*\pgfkeysvalueof{/pgfplots/parent axis width},
            scaled y ticks=false,
            at = {(0.25cm+\pgfkeysvalueof{/pgfplots/parent axis  width}, 0)},
            anchor=south west,
            width  = 0.025*\pgfkeysvalueof{/pgfplots/parent axis  width},
            height = \pgfkeysvalueof{/pgfplots/parent axis height},
          },
        ]
          \addplot graphics [xmin=-6, xmax=6, ymin=-6, ymax=6]{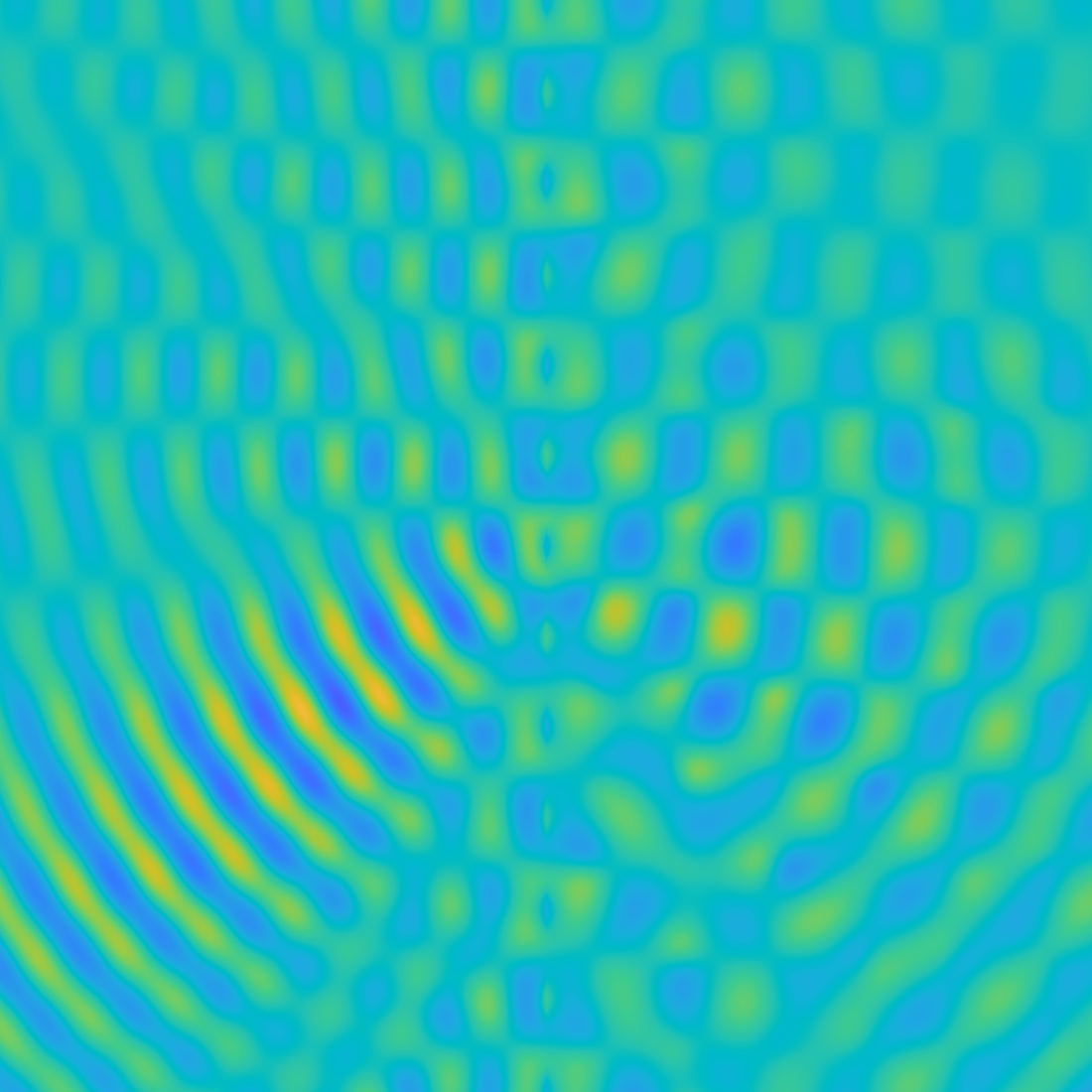};
          \draw[thick, dashed] (axis cs:0, -\Ldomain) -- (axis cs:0, \Ldomain);
        \end{axis}
      \end{tikzpicture}
    \end{subfigure}
    \hfill
    \begin{subfigure}[t]{0.52\textwidth}    
      \includegraphics[page=8]{H2Dmod/tikzpicture_H2Dmod.pdf}
    \end{subfigure}
  }
	\caption{$\omega = 8 + 0.25\vts \icplx$ and $N_k = 100$ are fixed. Left: real part of $u_\hv - u_{\textit{ref}}$ for $h = 0.025$. Right: relative errors $\veps^{(0)}_h$ and $\veps^{(1)}_h$ (see \eqref{H2Dmod:eq:def_relative_errors}) with respect to the mesh step $h$, for Configuration \eqref{H2Dmod:item:config_b}. \label{H2Dmod:fig:relative_errors_rational_setting}}
\end{figure}
}

\subsection{Qualitative validation and results in the irrational setting}\label{H2Dmod:sec:validation_irrational_setting}
This section is devoted to the case where $\rho$ is not periodic in the direction of the interface. For Configuration \eqref{H2Dmod:item:config_a}, we use the definition \eqref{H2Dmod:eq:def_rho_0_config_a} with $(p^+_\zvi, p^-_\zvi) = (1, \sqrt{2})$, and for Configuration \eqref{H2Dmod:item:config_b}, we use the definition \eqref{H2Dmod:eq:def_rho_0_config_b} with $\pv^+ = (\cos \alpha, \sin \alpha)$, $\alpha = 3\pi/5$. The corresponding coefficient $\rho$ is represented in Figure \ref{H2Dmod:fig:coefficients_irrational}.

\begin{figure}[ht!]
  \hspace{-11pt}
  \centering
  \makebox[0pt][c]{%
    \def\echelle{0.75} 
    \def\Ldomain{5}
    \begin{tikzpicture}
      \begin{groupplot}[
        group style={
          group size=2 by 1,
          horizontal sep=3.5cm,
          vertical sep=1.5cm,
        },
        enlargelimits=false,
        axis on top,
        width=0.3\textwidth,
        height=0.3\textwidth,
        scale only axis,
        xmin=-\Ldomain, xmax=\Ldomain, ymin=-\Ldomain, ymax=\Ldomain,
        xtick = {-4, -2, 0, 2, 4}, ytick = {-4, -2, 0, 2, 4},
        disabledatascaling,
        axis equal,
        clip mode=individual, 
      ]
        \nextgroupplot[title={$\rho$, Configuration \eqref{H2Dmod:item:config_a}},
          colorbar,
          point meta min=0.5, point meta max=1.5,
          colorbar style={
            yticklabel style={
              /pgf/number format/precision=2,
              /pgf/number format/fixed,
            },
            ytick={0.6, 1, 1.4},
            axis equal=false,
            major tick length=0.025*\pgfkeysvalueof{/pgfplots/parent axis width},
            scaled y ticks=false,
            at = {(0.25cm+\pgfkeysvalueof{/pgfplots/parent axis  width}, 0)},
            anchor=south west,
            width  = 0.025*\pgfkeysvalueof{/pgfplots/parent axis  width},
            height = \pgfkeysvalueof{/pgfplots/parent axis height},
          },]
        \addplot graphics [xmin=-6, xmax=6, ymin=-6, ymax=6]{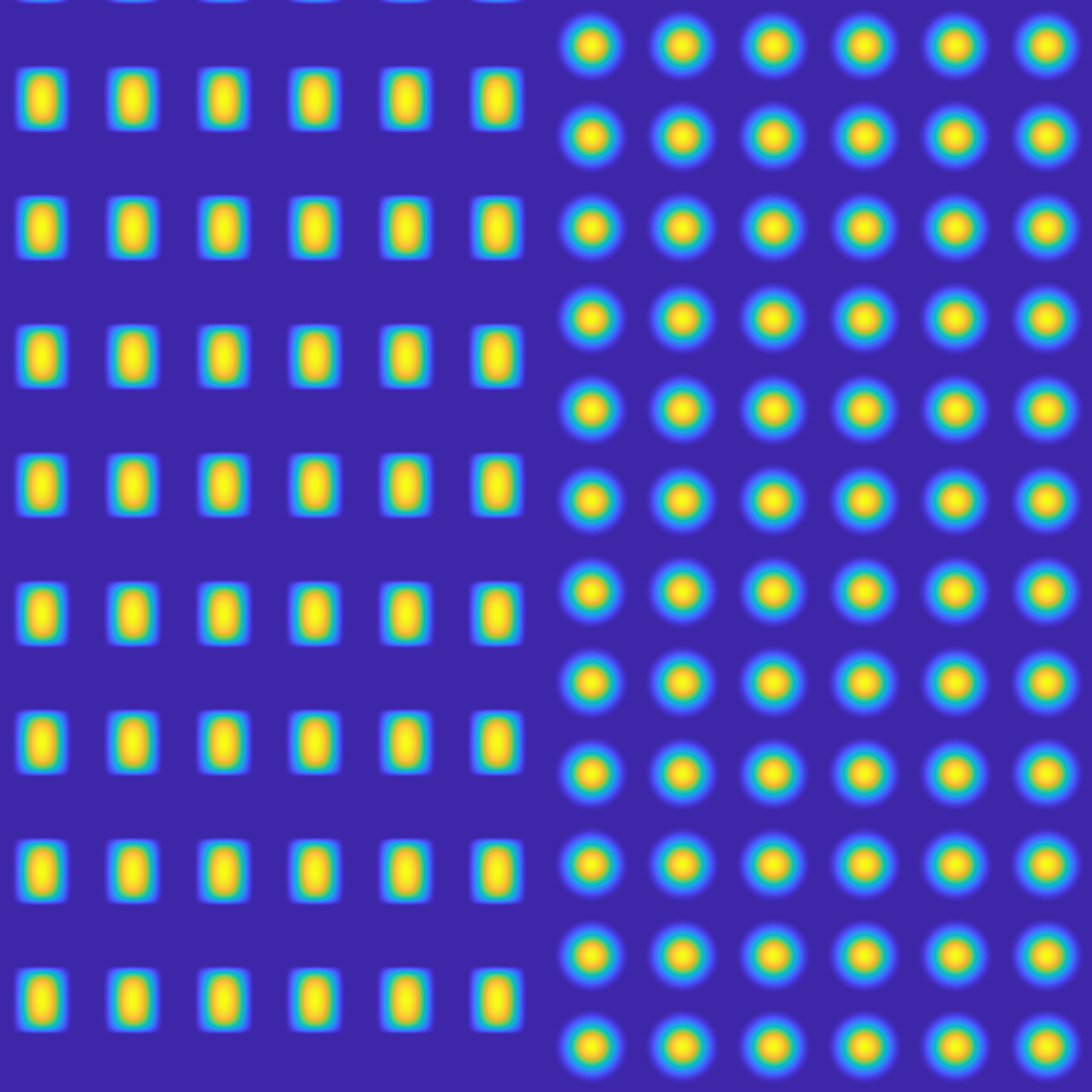};
        \pgfplotsinvokeforeach{1,2,...,5}{
          \draw[white, opacity=0.4] (axis cs:#1, -\Ldomain) -- (axis cs:#1, \Ldomain);
          \draw[white, opacity=0.4] (axis cs:-#1, -\Ldomain) -- (axis cs:-#1, \Ldomain);
        }
        \pgfplotsinvokeforeach{-5,-4,...,5}{
          \draw[white, opacity=0.4] (axis cs:0, #1) -- (axis cs:\Ldomain, #1);
          \draw[white, opacity=0.4] (axis cs:-\Ldomain, {#1*sqrt(2)}) -- (axis cs:0, {#1*sqrt(2)});
        }
        \draw[thick, white, dashed] (axis cs:0, -\Ldomain) -- (axis cs:0, \Ldomain);
        \nextgroupplot[
          title={$\rho$, Configuration \eqref{H2Dmod:item:config_b}},
          colorbar,
          point meta min=0.5, point meta max=1.5,
          colorbar style={
            yticklabel style={
              /pgf/number format/precision=2,
              /pgf/number format/fixed,
            },
            ytick={0.6, 1, 1.4},
            axis equal=false,
            major tick length=0.025*\pgfkeysvalueof{/pgfplots/parent axis width},
            scaled y ticks=false,
            at = {(0.25cm+\pgfkeysvalueof{/pgfplots/parent axis  width}, 0)},
            anchor=south west,
            width  = 0.025*\pgfkeysvalueof{/pgfplots/parent axis  width},
            height = \pgfkeysvalueof{/pgfplots/parent axis height},
          },
        ]
        \addplot graphics [xmin=-6, xmax=6, ymin=-6, ymax=6]{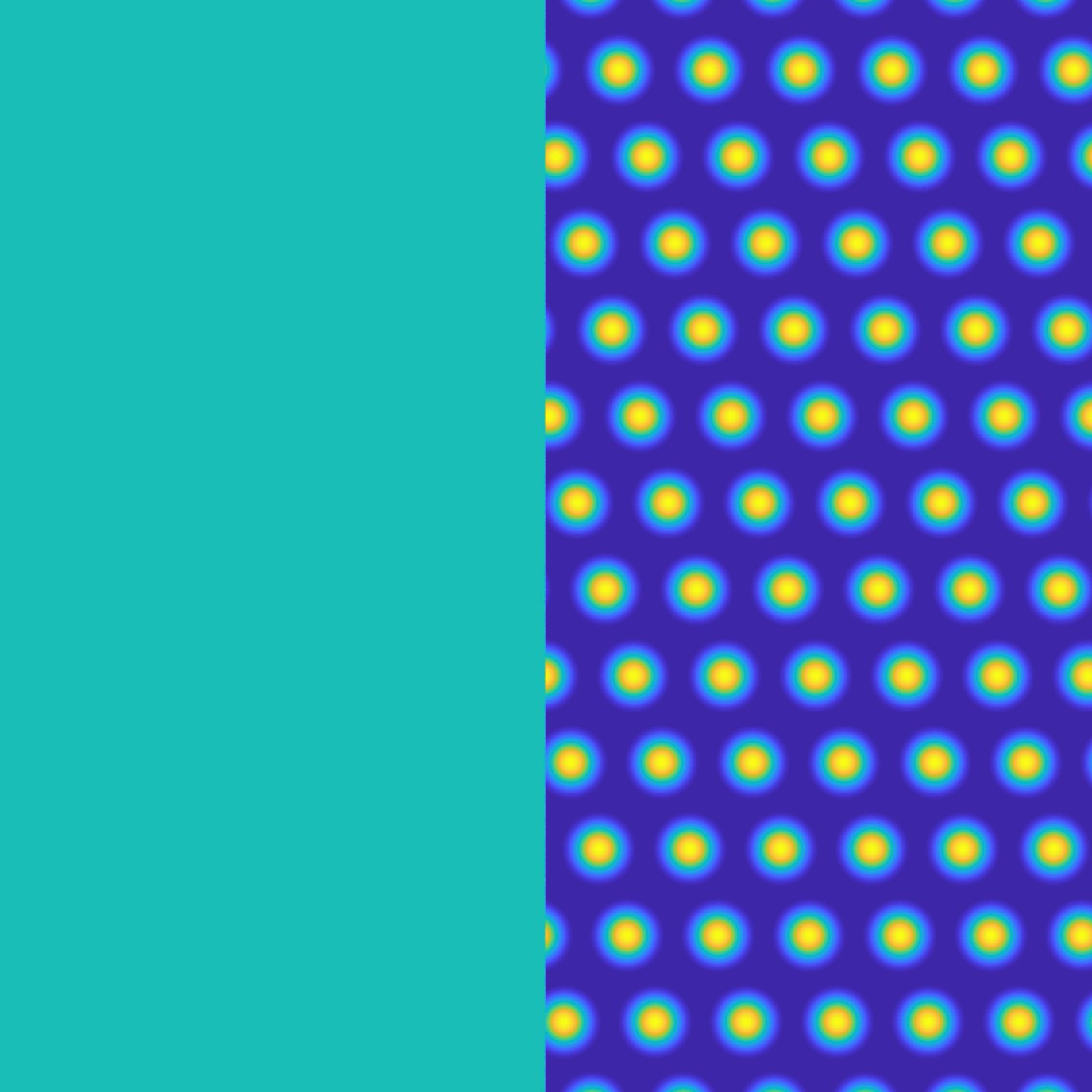};
        \draw[thick, white, dashed] (axis cs:0, -\Ldomain) -- (axis cs:0, \Ldomain);

        \clip (axis cs:0, -\Ldomain) rectangle (axis cs:\Ldomain, \Ldomain);
        \pgfmathsetmacro\vecperX{cos(108)}
        \pgfmathsetmacro\vecperZ{sin(108)}
        \pgfplotsinvokeforeach{-2,...,7}{
          \draw[white, opacity=0.4] (axis cs:0, {-#1*\vecperZ/\vecperX}) -- (axis cs:\Ldomain, {(\Ldomain-#1)*\vecperZ/\vecperX});
        }
        \pgfplotsinvokeforeach{-5,-4,...,5}{
          \draw[white, opacity=0.4] (axis cs:0, {#1*\vecperZ}) -- (axis cs:\Ldomain, {#1*\vecperZ});
        }
        \draw[thick, white, dashed] (axis cs:0, -\Ldomain) -- (axis cs:0, \Ldomain);
      \end{groupplot}
    \end{tikzpicture}
  }
  \caption{The coefficient $\rho$ for Configurations \eqref{H2Dmod:item:config_a} and \eqref{H2Dmod:item:config_b}.\label{H2Dmod:fig:coefficients_irrational}}
\end{figure}

\paragraph*{Invariance with respect to the period} We only consider Configuration \eqref{H2Dmod:item:config_a} for simplicity. As the coefficient $\rho^\pm$ is $p^\pm_\zvi$--periodic with respect to $\zvi$, it is also \surligner{$\ell^\pm\, p^\pm_\zvi$--periodic} with respect to $\zvi$ for any $\ell^\pm \in \Z^*$. Choosing $\ell^\pm\, p^\pm_\zvi$ to be the period along the interface will leave the $2$D solution $u$ of \eqref{H2Dmod:eq:transmission_problem} unchanged since $(\aten, \rho)$ remain the same, whereas the $3$D solution $\linapp{U}{G}$ of \eqref{H2Dmod:eq:augmented_transmission_problem} will be modified, due to the different expressions of the augmented coefficients $(\aten_p, \rho_p)$, the cut matrix $\cutmat$, and the augmented jump data $G$. To see if the approximate solution $u_\hv$ has the same invariances as $u$ with respect to the periods, we compute $u_\hv$ and the $3$D approximate solution $U_\hv$ in Figure \ref{H2Dmod:fig:invariance_period_G_config_A} for $(p^+_\zvi, p^-_\zvi) = (1, \sqrt{2})$ (middle column) and $(p^+_\zvi, p^-_\zvi) = (3, 2\sqrt{2})$ (left column), with $\omega = 8 + 0.25\,\icplx$, $h = 0.025$, and $N_\fbvar = 64$. For both values of the periods, the trace of $u_\hv$ on the interface $\sigma = \{\xvi = 0\}$ is shown in Figure \ref{H2Dmod:fig:invariance_period_G_config_A_trace}. One sees that $U_\hv$ changes with respect to $(p^+_\zvi, p^-_\zvi)$, whereas $u_\hv$ remains approximately the same.

\begin{figure}[ht!]
  \hspace{-11pt}
  \centering
  \makebox[0pt][c]{%
    \centering
    \def\echelle{0.75} 
    \def\Ldomain{5}
    \begin{tikzpicture}
      \begin{groupplot}[
        group style={
          group size=3 by 2,
          horizontal sep=0.75cm,
          vertical sep=1.5cm,
        },
        enlargelimits=false,
        axis on top,
        scale only axis,
        disabledatascaling,
        axis equal,
        clip mode=individual, 
      ]
        \nextgroupplot[
          title={$(a)$-- $(p^+_\zvi, p^-_\zvi) = (3, 2\sqrt{2})$, $G_1$},
          width=0.3\textwidth,
          height=0.3\textwidth,
          xmin=-\Ldomain, xmax=\Ldomain, ymin=-\Ldomain, ymax=\Ldomain,
          xtick = {-4, -2, 0, 2, 4}, ytick=\empty,
        ]
        \addplot graphics [xmin=-6, xmax=6, ymin=-6, ymax=6] {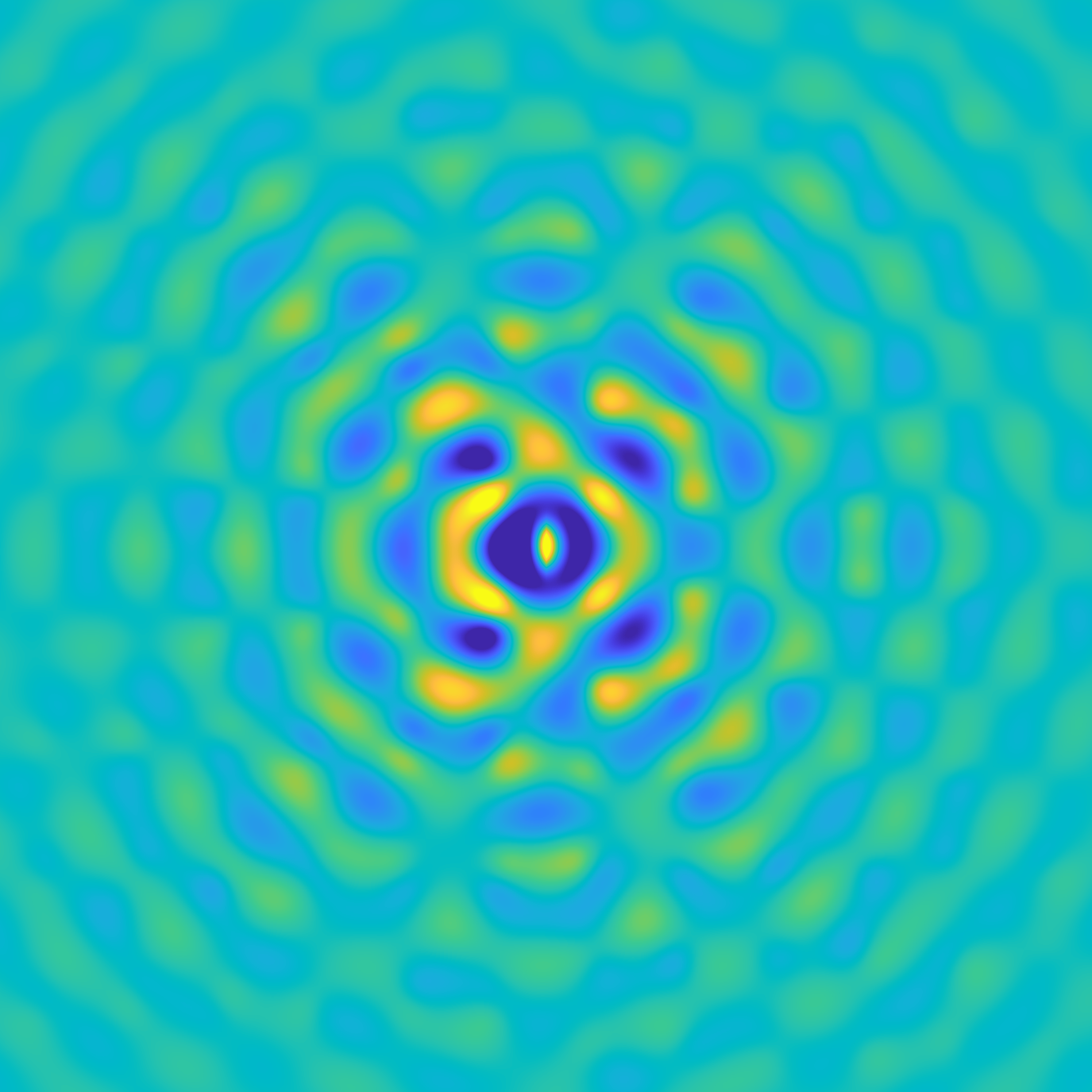};
        \draw[thick, dashed] (axis cs:0, -\Ldomain) -- (axis cs:0, \Ldomain);
        \nextgroupplot[%
          title={$(b)$-- $(p^+_\zvi, p^-_\zvi) = (1, \sqrt{2})$, $G_1$},
          width=0.3\textwidth,
          height=0.3\textwidth,
          xmin=-\Ldomain, xmax=\Ldomain, ymin=-\Ldomain, ymax=\Ldomain,
          xtick = {-4, -2, 0, 2, 4}, ytick=\empty,
          colorbar,
          point meta min=-3, point meta max=3,
          colorbar horizontal,
          colorbar style={
            xticklabel style={
              /pgf/number format/precision=2,
              /pgf/number format/fixed,
            },
            xtick={-2, 0, 2},
            axis equal=false,
            major tick length=0.025*\pgfkeysvalueof{/pgfplots/parent axis width},
            scaled x ticks=false,
            at = {(0, -0.875cm)},
            anchor=south west,
            width = \pgfkeysvalueof{/pgfplots/parent axis  width},
            height = 0.025*\pgfkeysvalueof{/pgfplots/parent axis height},
          },
        ]
        \addplot graphics [xmin=-6, xmax=6, ymin=-6, ymax=6]{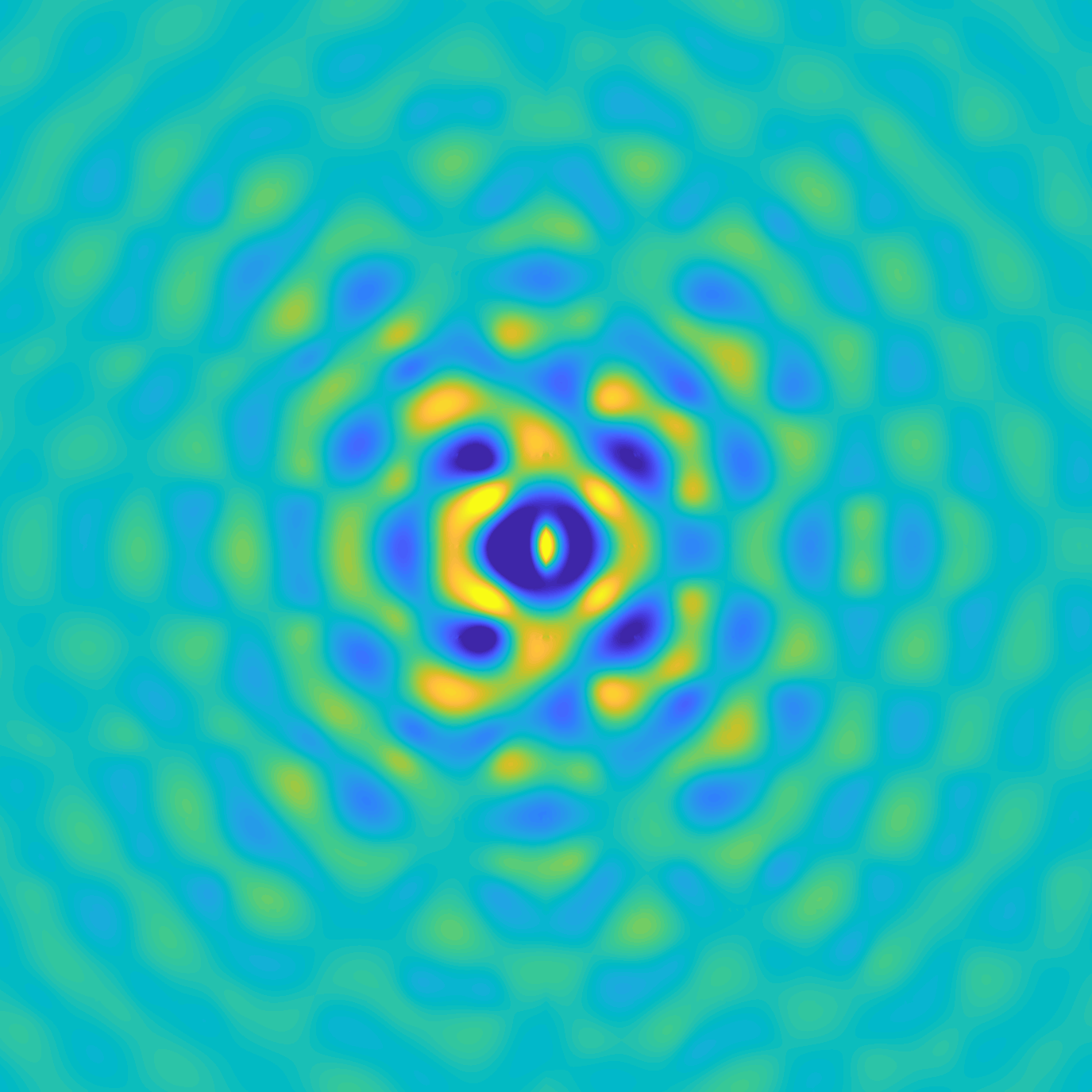};
        \draw[thick, dashed] (axis cs:0, -\Ldomain) -- (axis cs:0, \Ldomain);
        \nextgroupplot[
          title={$(c)$-- $(p^+_\zvi, p^-_\zvi) = (1, \sqrt{2})$, $G_2$},
          width=0.3\textwidth,
          height=0.3\textwidth,
          xmin=-\Ldomain, xmax=\Ldomain, ymin=-\Ldomain, ymax=\Ldomain,
          xtick = {-4, -2, 0, 2, 4}, ytick=\empty,
        ]
        \addplot graphics [xmin=-6, xmax=6, ymin=-6, ymax=6]{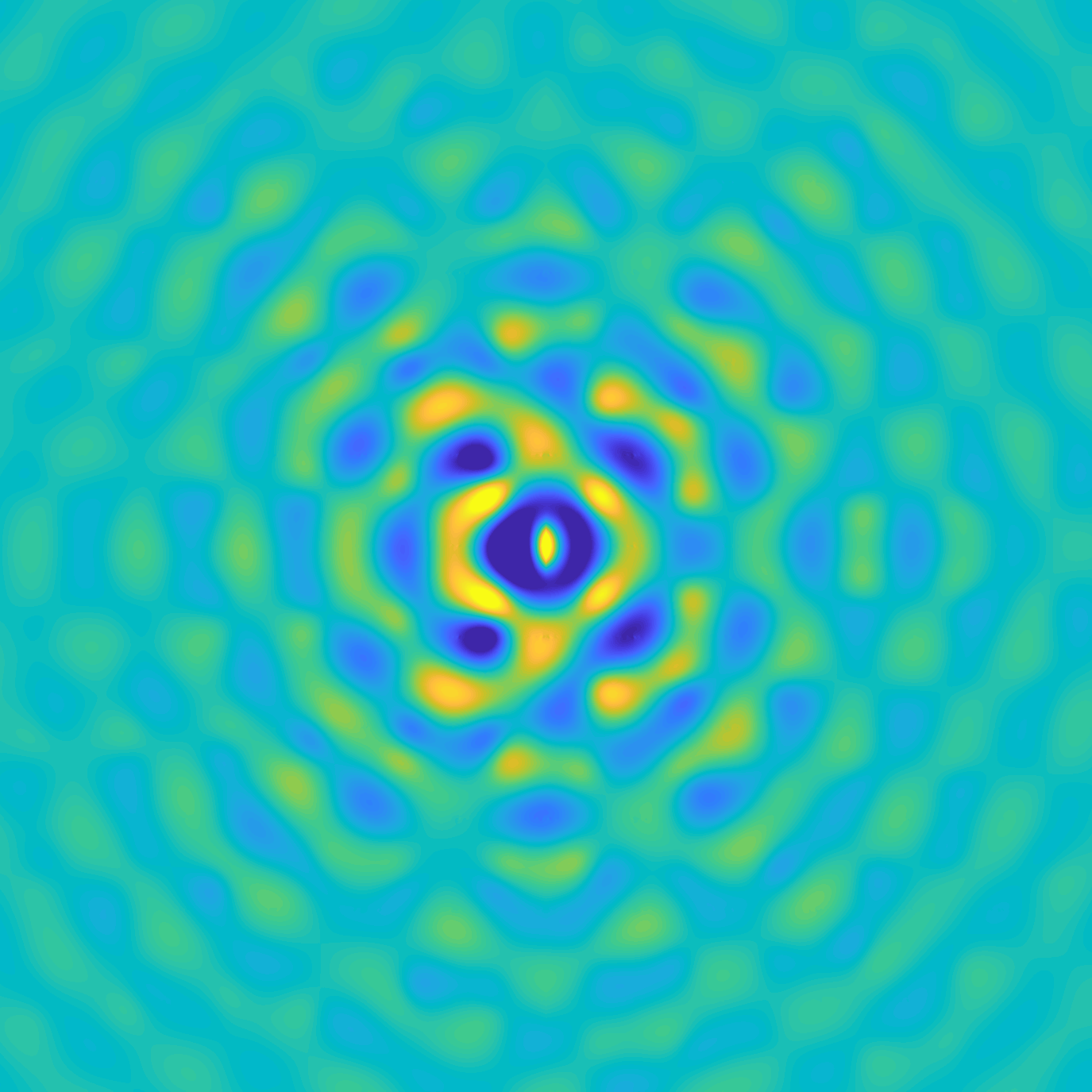};
        \draw[thick, dashed] (axis cs:0, -\Ldomain) -- (axis cs:0, \Ldomain);
        \nextgroupplot[%
          width=0.3\textwidth,
          height=0.2\textwidth,
          xmin=-2.0, xmax=2.0, ymin=-2.0, ymax=2.0, zmin=0, zmax=1.0, %
          xtick = {-2, 0, 2}, ytick=\empty, ztick={0, 1},
        ]
        \path (axis cs:0,0,0) coordinate (O1) (axis cs:1,0,0) coordinate (X1)  
              (axis cs:0,1,0) coordinate (Y1) (axis cs:0,0,1) coordinate (Z1);
        \nextgroupplot[%
          width=0.3\textwidth,
          height=0.2\textwidth,
          xmin=-2.0, xmax=2.0, ymin=-2.0, ymax=2.0, zmin=0, zmax=1.0, %
          xtick = {-2, 0, 2}, ytick=\empty, ztick={0, 1},
        ]
        \path (axis cs:0,0,0) coordinate (O2) (axis cs:1,0,0) coordinate (X2)  
              (axis cs:0,1,0) coordinate (Y2) (axis cs:0,0,1) coordinate (Z2);
        \nextgroupplot[%
          width=0.3\textwidth,
          height=0.2\textwidth,
          xmin=-2.0, xmax=2.0, ymin=-2.0, ymax=2.0, zmin=0, zmax=1.0, %
          xtick = {-2, 0, 2}, ytick={-2, 0, 2}, ztick={0, 1},
        ]
        \path (axis cs:0,0,0) coordinate (O3) (axis cs:1,0,0) coordinate (X3)  
              (axis cs:0,1,0) coordinate (Y3) (axis cs:0,0,1) coordinate (Z3);
      \end{groupplot}

      \begin{scope}[shift=(O1),x={($(X1)-(O1)$)},y={($(Y1)-(O1)$)},z={($(Z1)-(O1)$)}, transform shape]
        \draw[canvas is xy plane at z=1] (0, 0) node{\includegraphics[width=4cm,height=4cm]{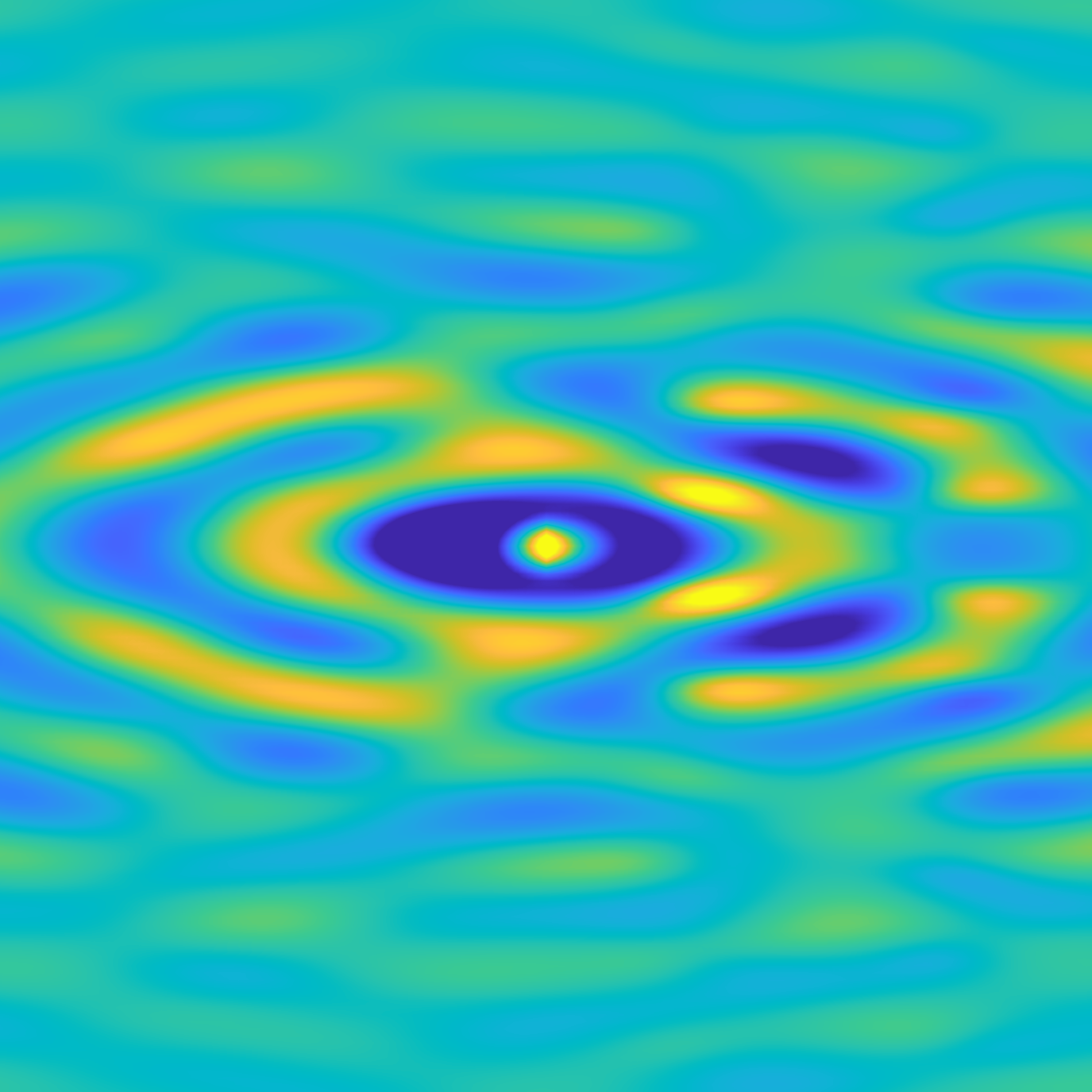}};
        \draw[canvas is yz plane at x=2] (0, 0.5) node{\includegraphics[width=4cm,height=1cm]{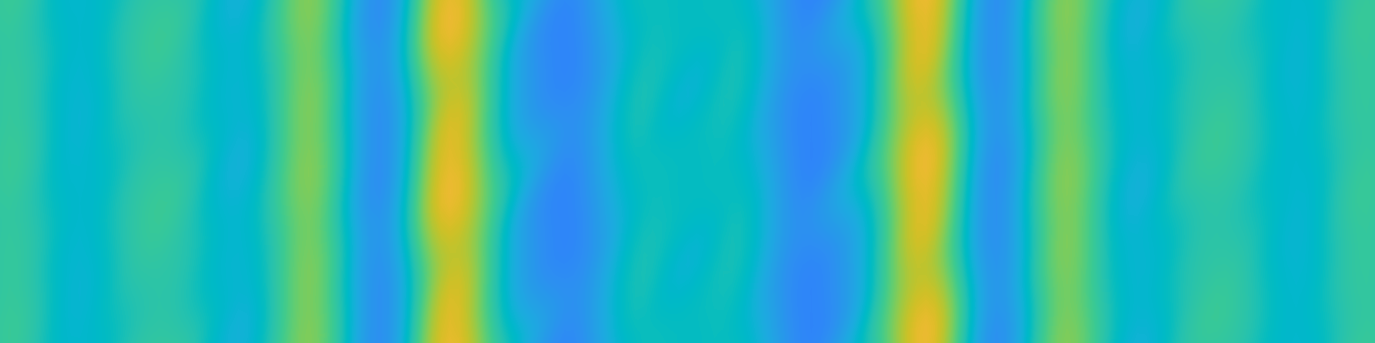}};
        \draw[canvas is zx plane at y=-2] (0.5, 0) node {\includegraphics[width=1cm,height=4cm]{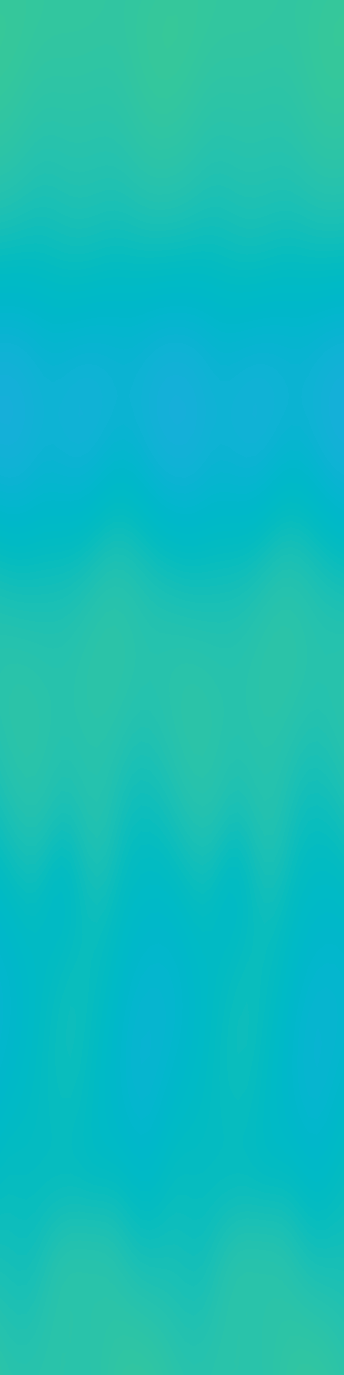}};
        
        \draw[canvas is zx plane at y=-2, black!80] (0, -2) rectangle (1, 2);
        \draw[canvas is zx plane at y=-2, black!80, dashed, thick] (0, 0) -- (1, 0);
        \draw[canvas is yz plane at x=2, black!80] (-2, 0) rectangle (2, 1);
        \draw[canvas is xy plane at z=1, black!80] (-2, -2) rectangle (2, 2);
        \draw[canvas is xy plane at z=1, black!80, dashed, thick] (0, -2) -- (0, 2);
      \end{scope}
      \begin{scope}[shift=(O2),x={($(X2)-(O2)$)},y={($(Y2)-(O2)$)},z={($(Z2)-(O2)$)}, transform shape]
        \draw[canvas is xy plane at z=1] (0, 0) node{\includegraphics[width=4cm,height=4cm]{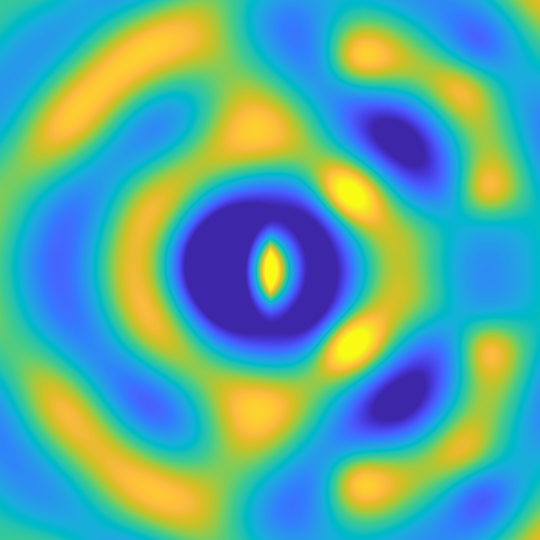}};
        \draw[canvas is yz plane at x=2] (0, 0.5) node{\includegraphics[width=4cm,height=1cm]{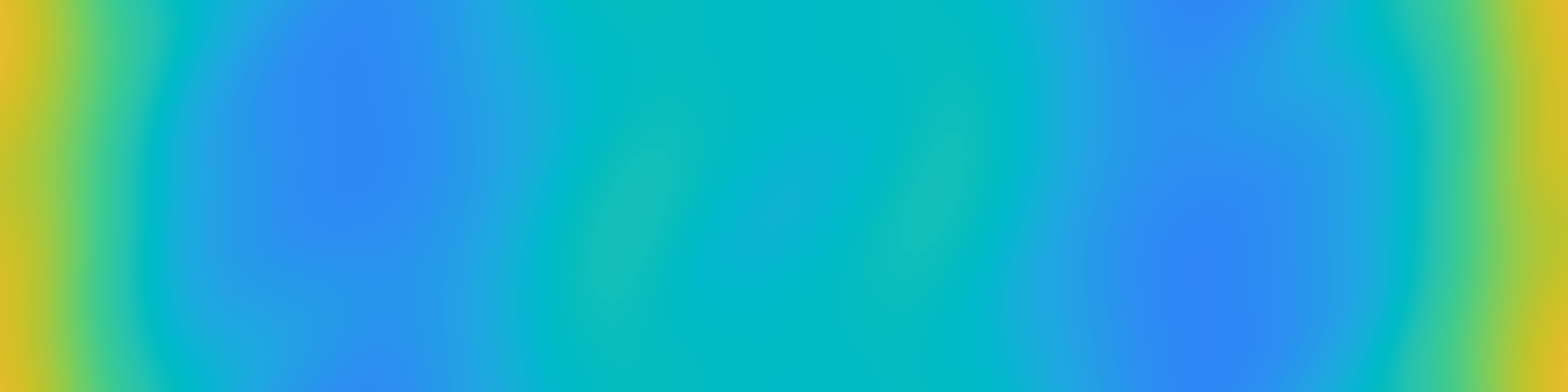}};
        \draw[canvas is zx plane at y=-2] (0.5, 0) node {\includegraphics[width=1cm,height=4cm]{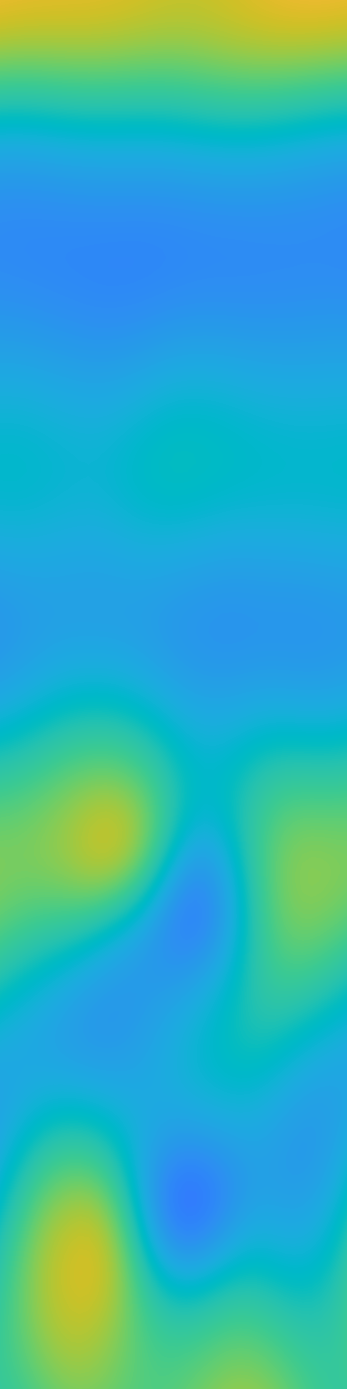}};
        
        \draw[canvas is zx plane at y=-2, black!80] (0, -2) rectangle (1, 2);
        \draw[canvas is zx plane at y=-2, black!80, dashed, thick] (0, 0) -- (1, 0);
        \draw[canvas is yz plane at x=2, black!80] (-2, 0) rectangle (2, 1);
        \draw[canvas is xy plane at z=1, black!80] (-2, -2) rectangle (2, 2);
        \draw[canvas is xy plane at z=1, black!80, dashed, thick] (0, -2) -- (0, 2);
      \end{scope}
      \begin{scope}[shift=(O3),x={($(X3)-(O3)$)},y={($(Y3)-(O3)$)},z={($(Z3)-(O3)$)}, transform shape]
        \draw[canvas is xy plane at z=1] (0, 0) node{\includegraphics[width=4cm,height=4cm]{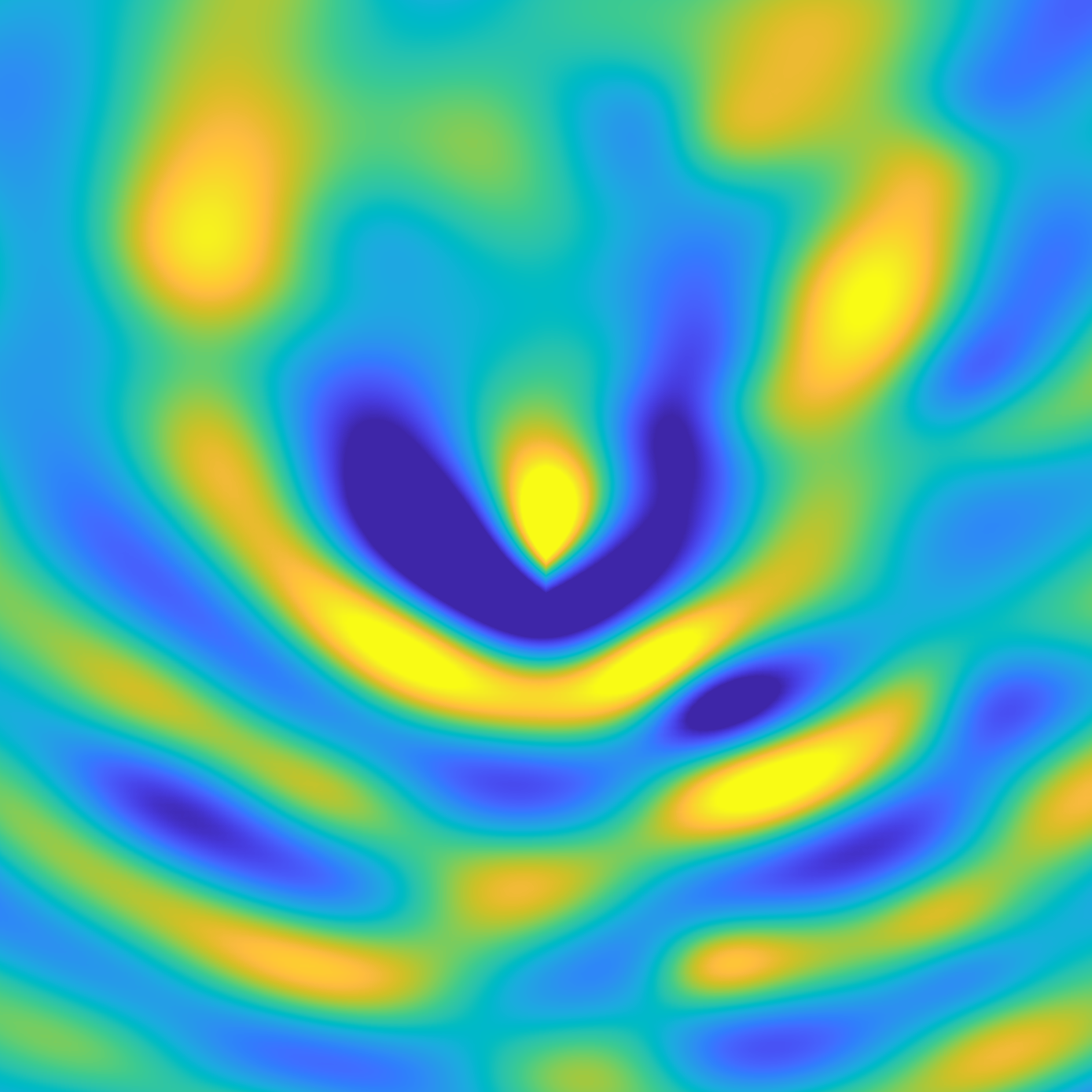}};
        \draw[canvas is yz plane at x=2] (0, 0.5) node{\includegraphics[width=4cm,height=1cm]{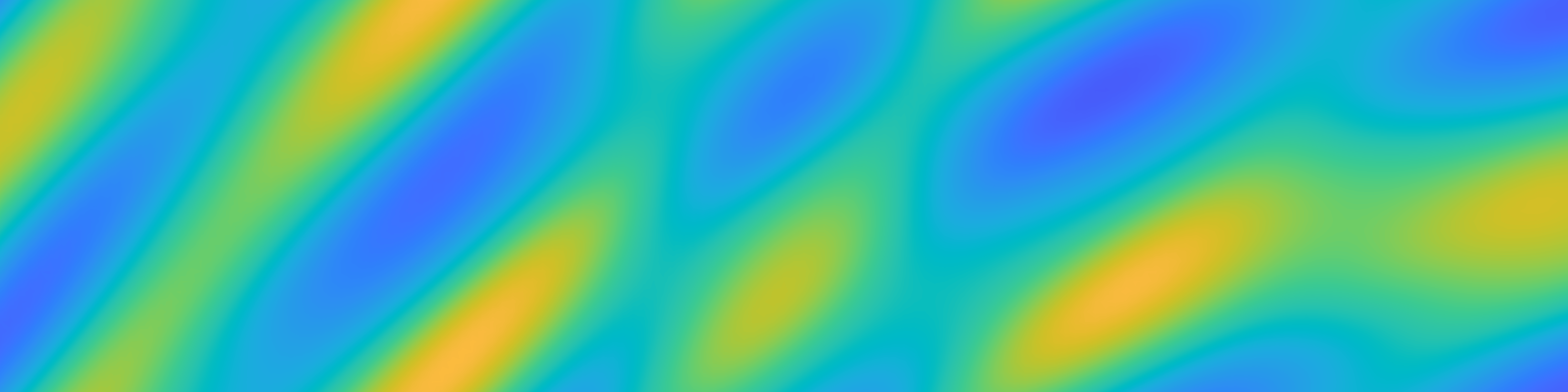}};
        \draw[canvas is zx plane at y=-2] (0.5, 0) node {\includegraphics[width=1cm,height=4cm]{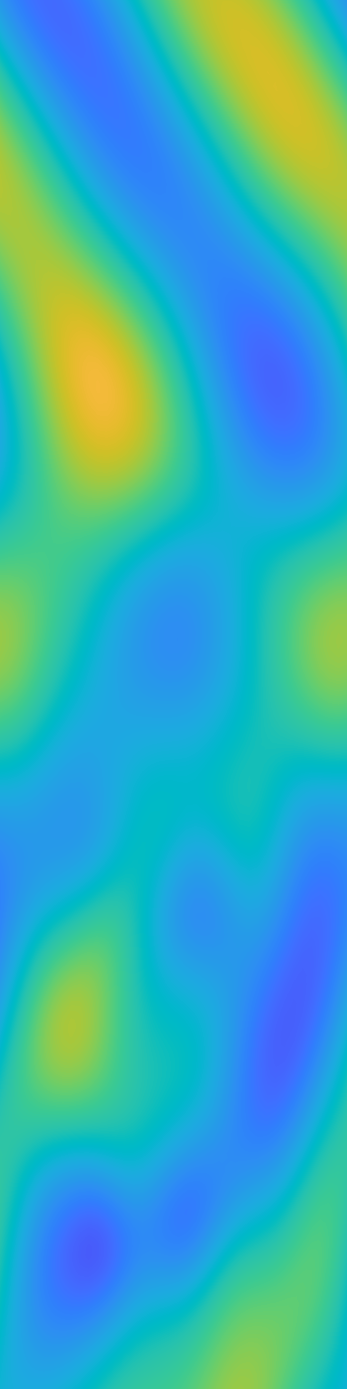}};
        
        \draw[canvas is zx plane at y=-2, black!80] (0, -2) rectangle (1, 2);
        \draw[canvas is zx plane at y=-2, black!80, dashed, thick] (0, 0) -- (1, 0);
        \draw[canvas is yz plane at x=2, black!80] (-2, 0) rectangle (2, 1);
        \draw[canvas is xy plane at z=1, black!80] (-2, -2) rectangle (2, 2);
        \draw[canvas is xy plane at z=1, black!80, dashed, thick] (0, -2) -- (0, 2);
      \end{scope}
    \end{tikzpicture}
  }
  \caption{The $2$D solution $u_\hv$ (first row, real part) and the augmented $3$D solution $U_\hv$ (second row, real part), for Configuration \eqref{H2Dmod:item:config_a} (with the coefficient $\rho$ represented in Figure \ref{H2Dmod:fig:coefficients_irrational}), and for different values of the periods along the interface and the augmented jump data. We use $\omega = 8 + 0.25\vts \icplx$, $h = 0.025$, and $N_\fbvar = 64$.\label{H2Dmod:fig:invariance_period_G_config_A}}
\end{figure}

\begin{figure}[ht!]
  \centering
  \makebox[0pt][c]{%
    \centering
    \includegraphics[page=6]{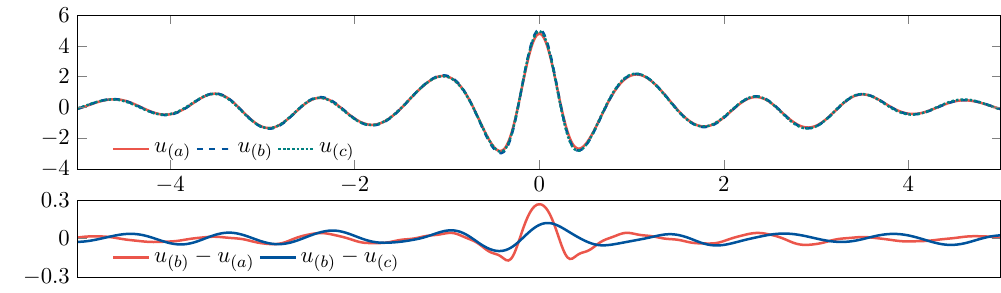}
  }
  \caption{Top: trace of the approximate solution $u_\hv$ on the interface $\sigma = \{\xvi = 0\}$ for Configuration \eqref{H2Dmod:item:config_a}, with the cases $(a)$, $(b)$, and $(c)$ outlined in Figure \ref{H2Dmod:fig:invariance_period_G_config_A}. Bottom: differences between the traces. $\omega = 8 + 0.25\vts \icplx$, $h = 0.025$, and $N_\fbvar = 64$ are fixed.\label{H2Dmod:fig:invariance_period_G_config_A_trace}}
\end{figure}

\paragraph*{Invariance with respect to augmented jump data} We still consider Configuration \eqref{H2Dmod:item:config_a} for simplicity. Now, let us look at the invariance of the solution $u_\hv$ with respect to the augmented jump data $G$. Since the data $g$ used for our experiments is smooth, we recall that the assumptions on the augmented jump data reduce to $G$ satisfying $G(0, \cuti_1\, \zvi, \cuti_2\, \zvi) = g(0, \zvi)$ and $G(\cdot + \eva_2) = G$. These requirements are satisfied by the following functions (see Remark \ref{H2D:rmk:example_G}):
\[\displaystyle
  G_1 (0, \zvi_1, \zvi_2) := g(0, \zvi_1 / \cuti_1) \quad \textnormal{and} \quad G_2 (0, \zvi_1, \zvi_2) := \exp\big(2\icplx \pi \vts (\zvi_2 - \zvi_1 \cuti_2 / \cuti_1)\big)\; g(0, \zvi_1 / \cuti_1). 
\]
Figure \ref{H2Dmod:fig:invariance_period_G_config_A} (middle and right columns) and Figure \ref{H2Dmod:fig:invariance_period_G_config_A_trace} show the solution $u_\hv$ and the augmented solution $U_\hv$ obtained for $G_1$ and $G_2$, with $p^+_\zvi = 1$, $p^-_\zvi = \sqrt{2}$, $h = 0.025$, $N_\fbvar = 64$, and $\omega = 8 + 0.25\,\icplx$. As expected, $U_\hv$ depends on the choice of the augmented data and changes accordingly, whereas $u_\hv$ is visually invariant.

\paragraph*{Dependence with respect to the frequency}
We finish by solving \eqref{H2Dmod:eq:transmission_problem} for different values of $\omega$. As expected for the Helmholtz equation, the solution, represented in Figure \ref{H2Dmod:fig:dependence_frequency}, oscillates more as $\Real \omega$ increases, and decays less at infinity as $\Imag \omega$ decreases.

\begin{figure}[ht!]
  \hspace{-11pt}
  \centering
  \makebox[0pt][c]{%
    \centering
    \def\echelle{0.75} 
    \def\Ldomain{5}
    \begin{tikzpicture}
      \begin{groupplot}[
        group style={
          group size=3 by 2,
          horizontal sep=0.75cm,
          vertical sep=1.5cm,
        },
        enlargelimits=false,
        axis on top,
        width=0.3\textwidth,
        height=0.3\textwidth,
        scale only axis,
        xmin=-\Ldomain, xmax=\Ldomain, ymin=-\Ldomain, ymax=\Ldomain,
        xtick = {-4, -2, 0, 2, 4}, ytick=\empty,
        disabledatascaling,
        axis equal,
        clip mode=individual, 
      ]
        \nextgroupplot[title={$\omega = 8 + 0.25\vts \icplx$}, ylabel={Configuration \eqref{H2Dmod:item:config_a}}]
        \addplot graphics [xmin=-6, xmax=6, ymin=-6, ymax=6] {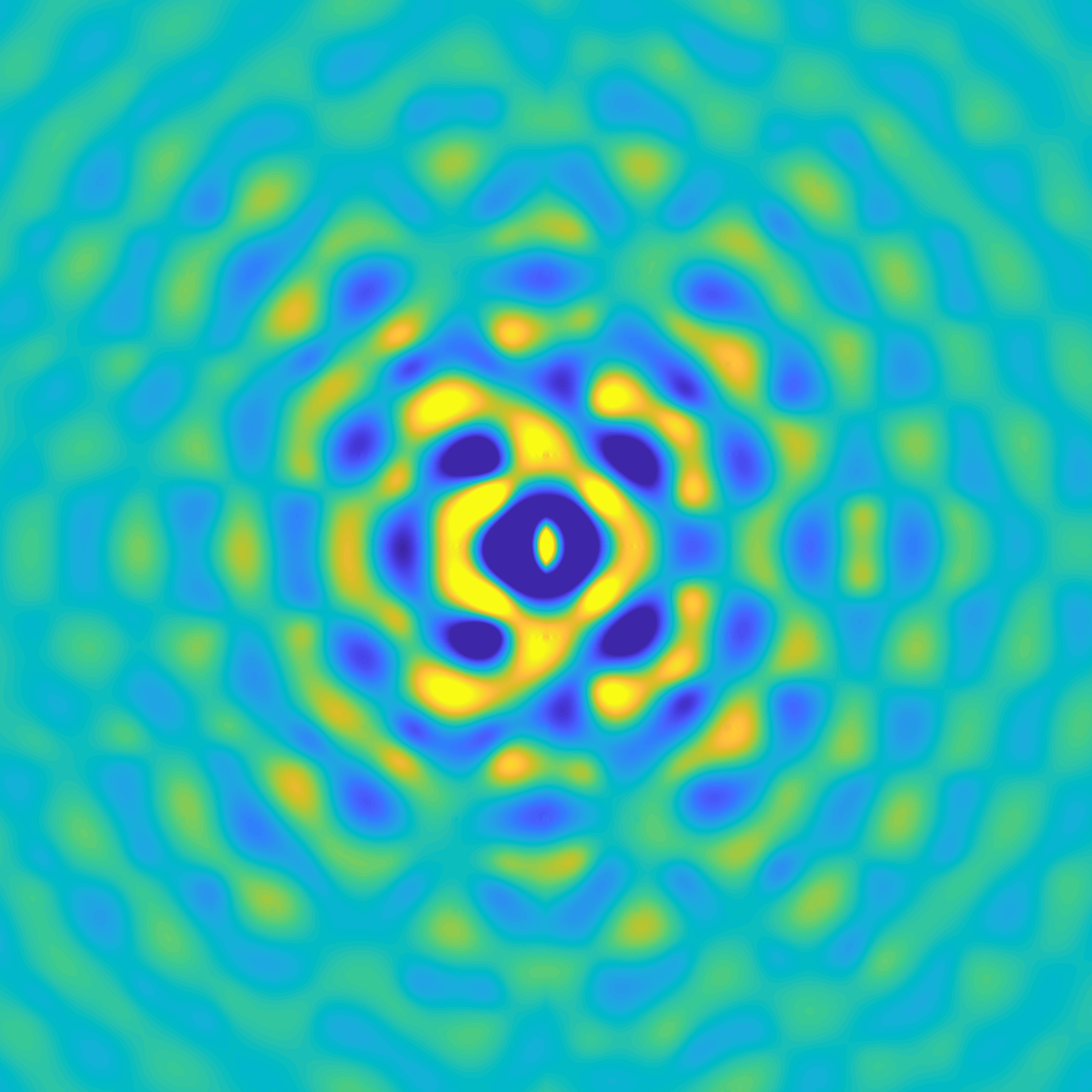};
        \draw[thick, dashed] (axis cs:0, -\Ldomain) -- (axis cs:0, \Ldomain);
        \nextgroupplot[%
          title={$\omega = 20 + 0.25\vts \icplx$},
          colorbar,
          point meta min=-2, point meta max=2,
          colorbar horizontal,
          colorbar style={
            xticklabel style={
              /pgf/number format/precision=2,
              /pgf/number format/fixed,
            },
            xtick={-1, 0, 1},
            axis equal=false,
            major tick length=0.025*\pgfkeysvalueof{/pgfplots/parent axis width},
            scaled x ticks=false,
            at = {(0, -0.75cm)},
            anchor=south west,
            width = \pgfkeysvalueof{/pgfplots/parent axis  width},
            height = 0.025*\pgfkeysvalueof{/pgfplots/parent axis height},
          },
        ]
        \addplot graphics [xmin=-6, xmax=6, ymin=-6, ymax=6]{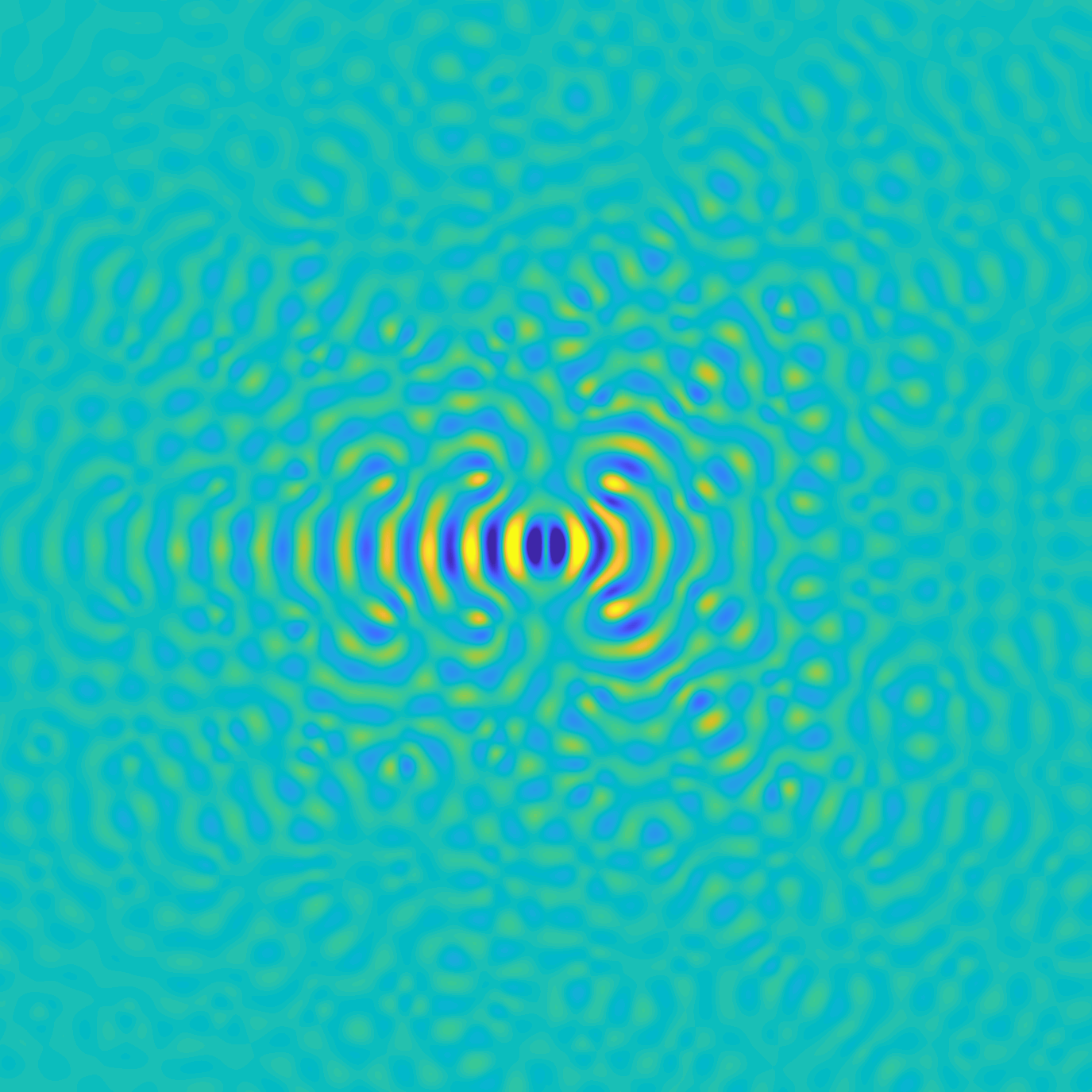};
        \draw[thick, dashed] (axis cs:0, -\Ldomain) -- (axis cs:0, \Ldomain);
        \nextgroupplot[title={$\omega = 20 + 0.05\, \icplx$}]
        \addplot graphics [xmin=-6, xmax=6, ymin=-6, ymax=6]{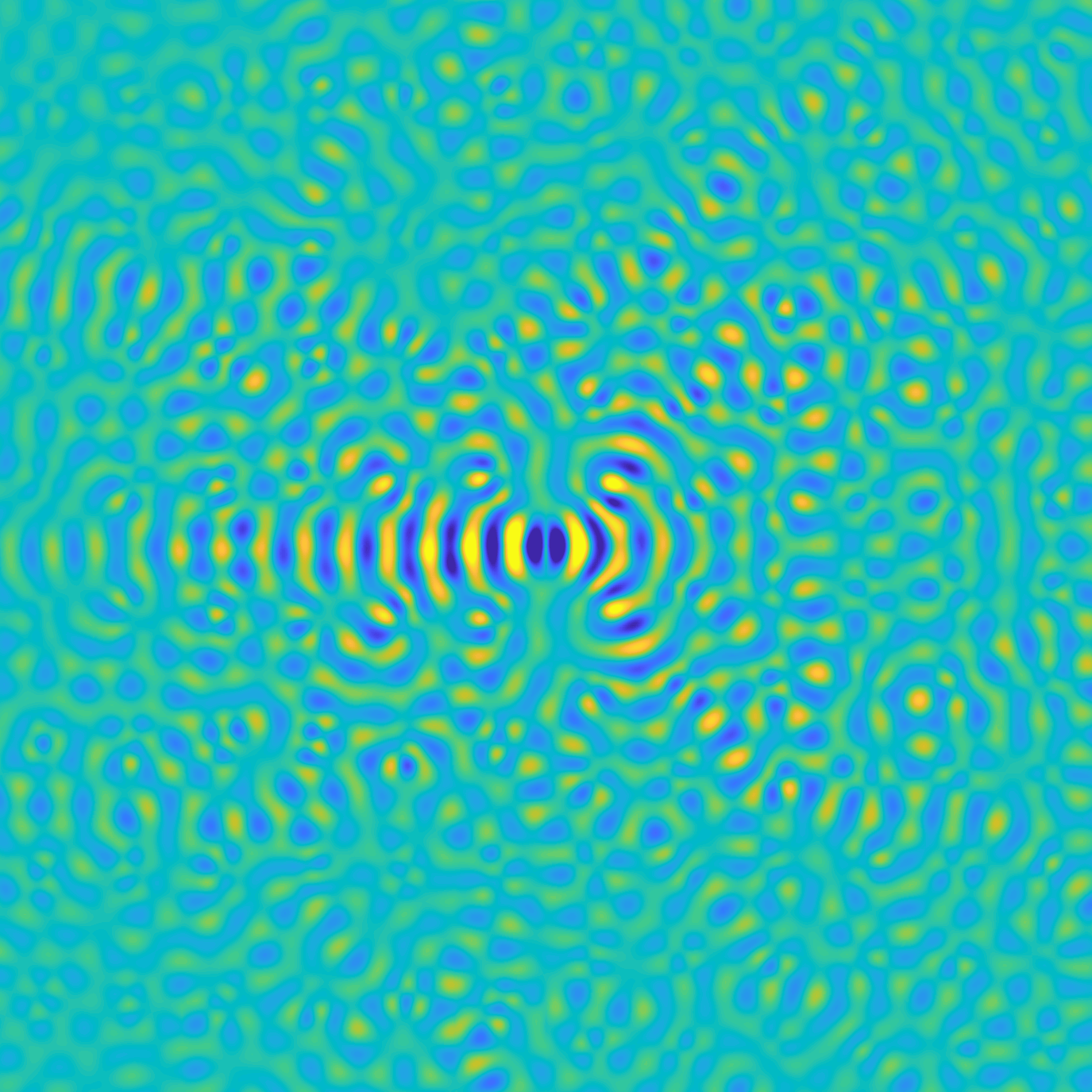};
        \draw[thick, dashed] (axis cs:0, -\Ldomain) -- (axis cs:0, \Ldomain);
        \nextgroupplot[ylabel={Configuration \eqref{H2Dmod:item:config_b}}]
        \addplot graphics [xmin=-6, xmax=6, ymin=-6, ymax=6] {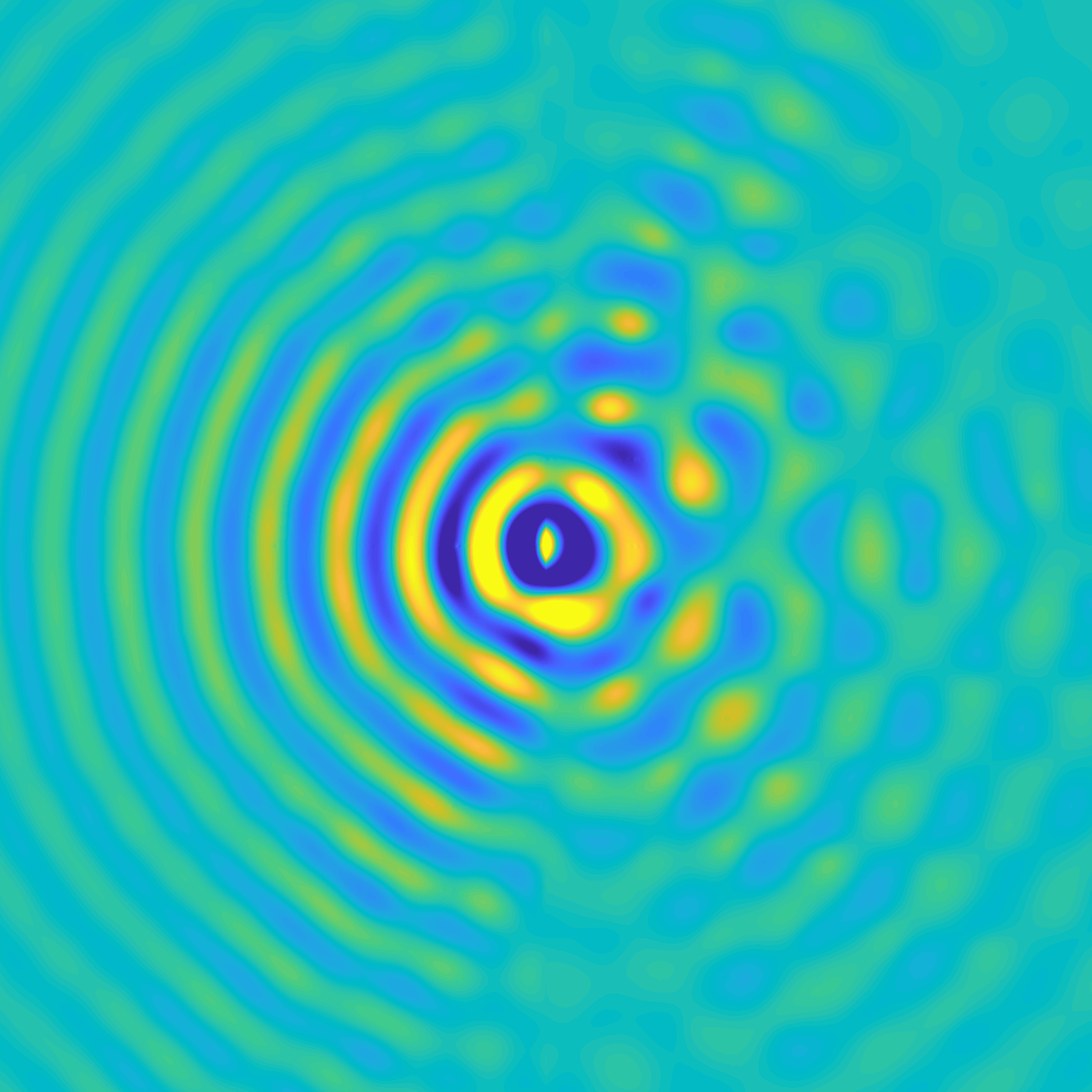};
        \draw[thick, dashed] (axis cs:0, -\Ldomain) -- (axis cs:0, \Ldomain);
        \nextgroupplot[]
        \addplot graphics [xmin=-6, xmax=6, ymin=-6, ymax=6]{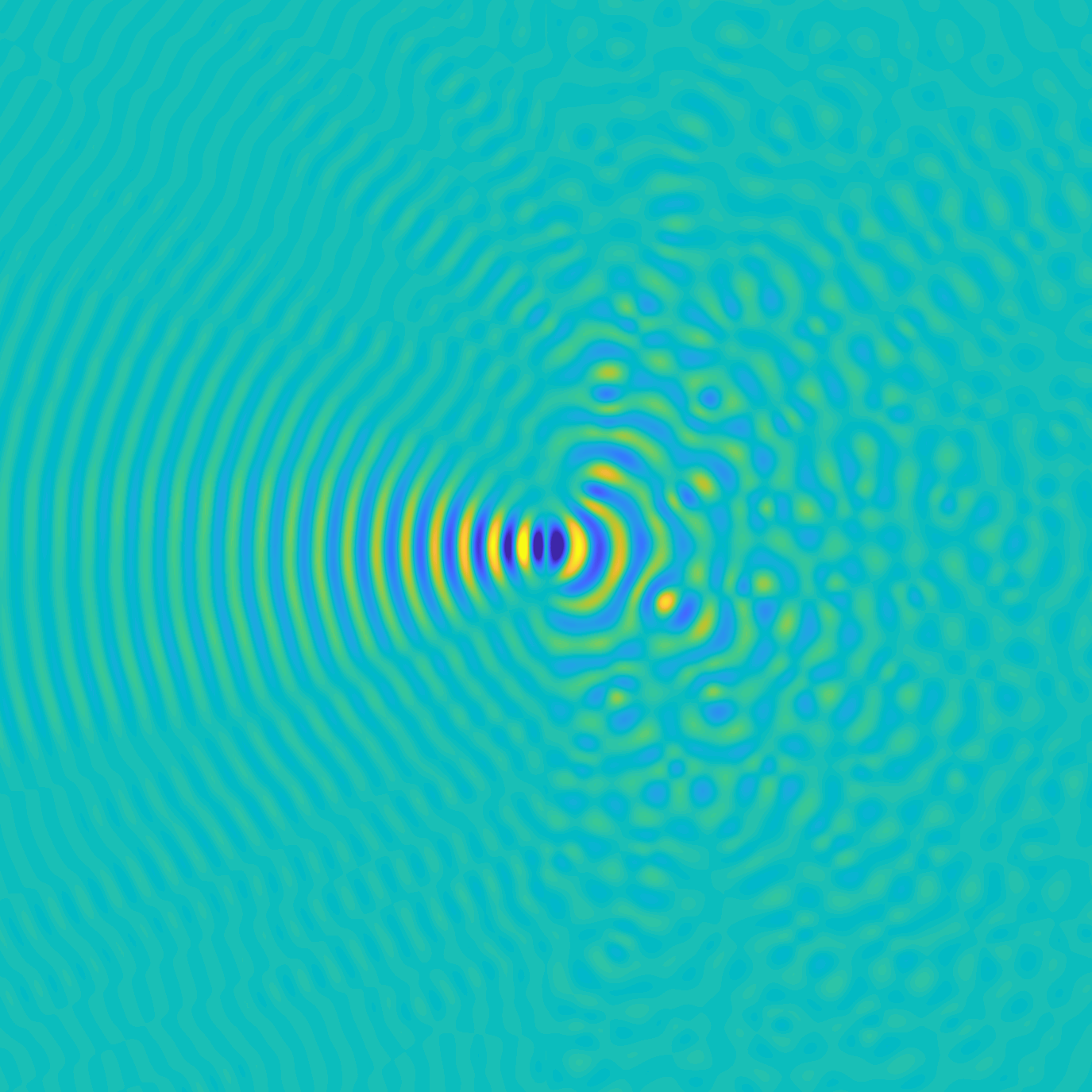};
        \draw[thick, dashed] (axis cs:0, -\Ldomain) -- (axis cs:0, \Ldomain);
        \nextgroupplot[]
        \addplot graphics [xmin=-6, xmax=6, ymin=-6, ymax=6]{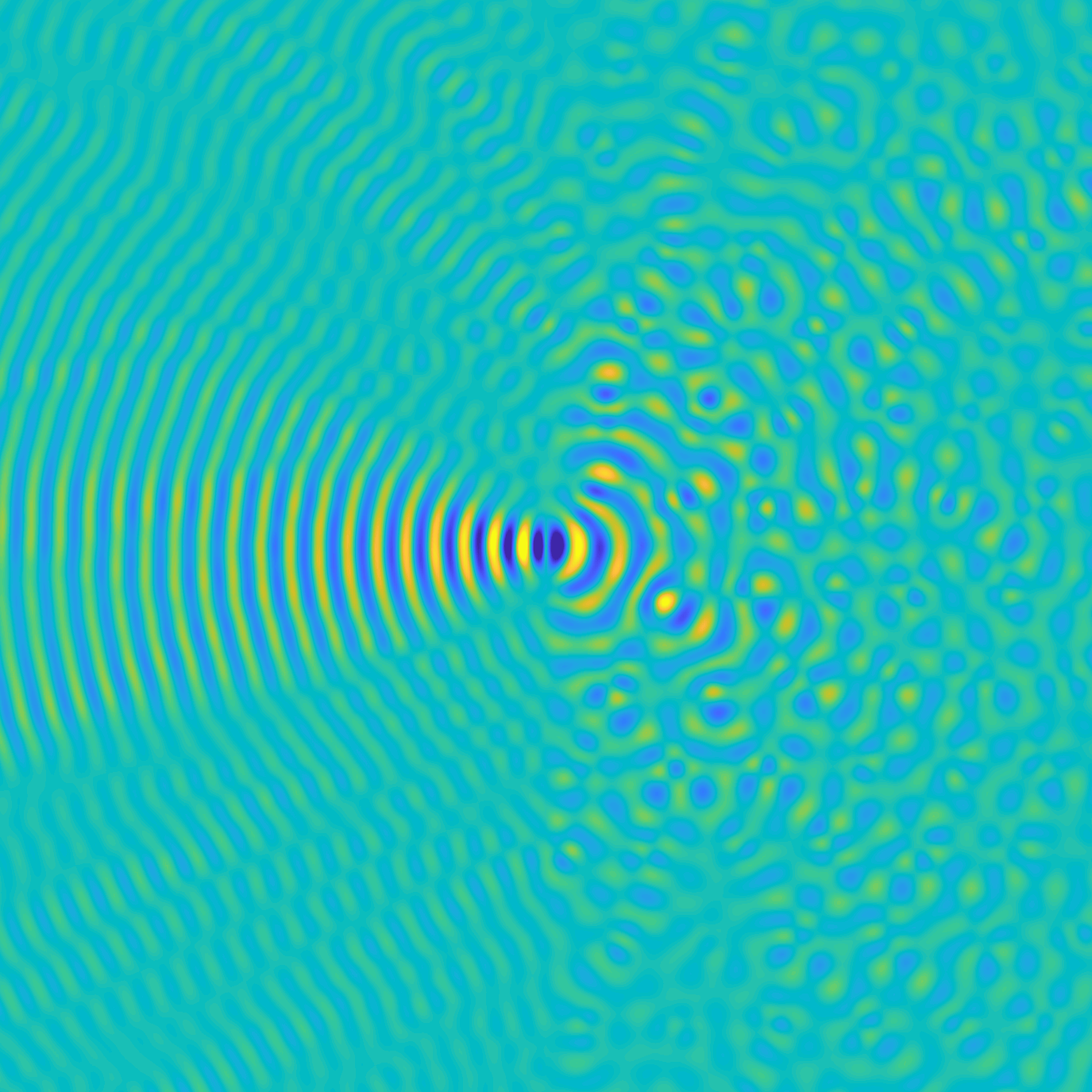};
        \draw[thick, dashed] (axis cs:0, -\Ldomain) -- (axis cs:0, \Ldomain);
      \end{groupplot}
    \end{tikzpicture}
  }
  \caption{The approximate solution $u_\hv$ (real part) in the irrational setting for Configuration \eqref{H2Dmod:item:config_a} ($p^+_\zvi = 1$, $p^-_\zvi = \sqrt{2}$) and Configuration \eqref{H2Dmod:item:config_b} ($\pv^+ = (\cos \alpha, \sin \alpha)$, $\alpha = 3\pi/5$). The discretization parameters $h = 0.025$ and $N_\fbvar = 64$ are fixed, and different values of $\omega$ are considered.\label{H2Dmod:fig:dependence_frequency}}
\end{figure}

%
}{%
}

\ifthenelse{\boolean{brouillon}}{%
  \section{Extensions and perspectives}\label{H2Dmod:sec:extensions_perspectives}
In this section, we explore potential extensions and perspectives arising from the present work.

\subsection{Direct extension to a source term}\label{H2Dmod:sec:extension_source_term}
The approach developed in this paper could be extended to the classical Helmholtz equation
\begin{equation*}
  -\transp{\nabla}\, \aten\, \nabla u^0 - \rho\, \omega^2\, u^0 = f, \quad \textnormal{in}\ \ \R^2,
\end{equation*}
where the source term $f \in L^2(\R^2)$ has a compact support in the $\ev_\xvi$--direction. In particular, using the fact that $\aten(\xv) = \aten_p (\cutmat\vts \xv)$ and $\rho(\xv) = \rho_p (\cutmat\vts \xv)$, one may seek $u^0$ under the form $u^0 (\xv) = U^0 (\cutmat\vts \xv)$, where $U^0$ satisfies 
\begin{equation*}
  \displaystyle- \transp{\bsnabla} \cutmat\vts \aten_p\!  \transp{\cutmat} \bsnabla\, U^0 - \rho_p\, \omega^2\, U^0 = F, \quad \textnormal{in}\ \  \R^3,
\end{equation*}
with $F \in L^2_{\textit{loc}}(\R^3)$ being compactly supported in the $\eva_\xvi$--direction, and satisfying $F(\cutmat\vts \xv) = f(\xv)$ for almost any $\xv \in \R^2$. Similarly to the procedure summarized in Section \ref{H2Dmod:sec:formal_description_method} and developed throughout the paper, the augmented equation satisfied by $U^0$ may be solved as follows.
\begin{itemize}
  \item We first choose an augmented source term $F$ such that $F(\cdot + \eva_2) = F$, so that $U^0(\cdot + \eva_2) = U^0$ satisfies \surligner{a strip problem with periodicity conditions with respect to $\zvi_2$:}
  \begin{equation*}
    \left\{
      \begin{array}{r@{\ =\ }l}
        \displaystyle- \transp{\bsnabla}\, \cutmat\vts \aten_p \!  \transp{\cutmat} \bsnabla\, U^0 - \rho_p\, \omega^2\, U^0 & F \quad \textnormal{in}\ \  \Omegaper,
        \ret
        \multicolumn{2}{c}{\displaystyle U^0 \in H^1_{\cutmat, \periodic} (\Omegaper), \ \  \aten_p \! \transp{\cutmat} \bsnabla\, U^0 \in \Hdiveper{\Omegaper}.}
      \end{array}
    \right.
  \end{equation*}
  \item We apply a partial Floquet-Bloch transform with respect to $\zvi_1$, to obtain a family of waveguide problems \surligner{parameterized by $\fbvar \in [-\pi, \pi]$, and} defined in the cylinder \surligner{$\Omegaperper := \R \times (0, 1)^2$:
  \begin{equation*}
    \left\{
      \begin{array}{r@{\ =\ }l}
        \displaystyle- \divecutFB\, \aten_p \!  \gradcutFB\, \widehat{U}^0_\fbvar - \rho_p\, \omega^2\, \widehat{U}^0_\fbvar & \widehat{F}_\fbvar \quad \textnormal{in}\ \  \Omegaperper,
        \ret
        \multicolumn{2}{c}{\displaystyle \widehat{U}^0_\fbvar \in \Hgradperper{\Omegaperper}, \ \  \aten_p \! \gradcutFB\, \widehat{U}^0_\fbvar \in \Hdiveperper{\Omegaperper},}
      \end{array}
    \right.
  \end{equation*}
  where $\widehat{F}_\fbvar := \fbtransform[] F (\cdot, \fbvar)$ is compactly supported in the $\eva_\xvi$--direction, because $F$ is compactly supported in the $\eva_\xvi$--direction.}
  \item Finally, we use the DtN method developed in \cite{joly2006exact,fliss2009analyse,fliss2009exact,fliss2012wave} \surligner{to compute the waveguide solution $\widehat{U}^0_\fbvar$ for any $\fbvar \in [-\pi, \pi]$. More precisely, assuming that $(a_-, a_+) \subset \R$ is the support of $\widehat{F}_\fbvar$ in the $\eva_\xvi$--direction, by analogy with Section \ref{H2Dmod:sec:DtN_waveguide} and the interface equation \eqref{H2Dmod:eq:transmission_jump}, we characterize the restriction of $\widehat{U}^0_\fbvar$ to the interior domain $\Omegaperper^{\emph{int}} := (a_-, a_+) \times (0, 1)^2$ as the solution of:
  \begin{equation*}
    \left\{
      \begin{array}{r@{\ =\ }l}
        \displaystyle- \divecutFB\, \aten_p \!  \gradcutFB\, \widehat{U}^0_\fbvar - \rho_p\, \omega^2\, \widehat{U}^0_\fbvar & \widehat{F}_\fbvar \quad \textnormal{in}\ \  \Omegaperper^{\emph{int}},
        \ret
        \pm \big(\cutmat\, \aten_p \! \gradcutFB\, \widehat{U}^0_\fbvar \cdot \eva_\xvi\big) + \Lambdahat^\pm_\fbvar\, \widehat{U}^0_\fbvar & 0, \quad \textnormal{on}\ \ \{\xvi = a^\pm\},
        \ret
        \multicolumn{2}{c}{\displaystyle \widehat{U}^0_\fbvar \in \Hgradperper{\Omegaperper^{\emph{int}}}, \ \  \aten_p \! \gradcutFB\, \widehat{U}^0_\fbvar \in \Hdiveperper{\Omegaperper^{\emph{int}}};}
      \end{array}
    \right.
  \end{equation*}
  see Figure \ref{H2Dmod:fig:interior_domain}. The conditions on $\{\xvi = a^\pm\}$ feature DtN operators $\Lambdahat^\pm_\fbvar$ which can be computed from solutions of half-guide problems set in the half-cylinders $\Omegaperper \cap \R^3_\pm$. Since the source term $\widehat{F}_\fbvar$ vanishes in these half-cylinders, the half-guide problems can be solved using the local cell problems and the propagation operator introduced in Section \ref{H2Dmod:sec:resolution_half_guide_problems}.}

  \begin{figure}[H]
    \noindent\makebox[\textwidth]{%
      \def\coteCyl{1.525}
      \begin{tikzpicture}[scale=0.675]
        \def\angleinterface{55}%
        \definecolor{coplancoupe}{RGB}{211, 227, 65}
        \tdplotsetmaincoords{78}{20}
        \def\BBxmin{-4.5}
        \def\BBxmax{+4.5}
        \def\BBymin{-1.25}
        \def\BBymax{+4.75}
        \def\BBzmin{-0.75}
        \def\BBzmax{+2.75}
        \def\amin{-1.5}
        \def\amax{+1.5}
        \begin{scope}[tdplot_main_coords, scale=0.9]
          \node (xaxis) at ({\BBxmax+\amax+0.2*(\BBxmax-\BBxmin)}, 0, 0) {$\xvi$};%
          \node (yaxis) at (0, {\BBymax+1.2*(\BBymax-\BBymin)}, 0) {$\zvi_1$};%
          %
          %
          %
          \fill[secondary!80, opacity=0.5] (\BBxmin+\amin, 0, 0) -- +(-\BBxmin, 0, 0) -- +(-\BBxmin, \coteCyl, 0) -- +(0, \coteCyl, 0) -- cycle;
          \fill[secondary!90, opacity=0.625] (\BBxmin+\amin, \coteCyl, 0) -- +(-\BBxmin, 0, 0) -- +(-\BBxmin, 0, \coteCyl) -- +(0, 0, \coteCyl) -- cycle;
          \draw[dashed, white] (\BBxmin+\amin, \coteCyl, 0) -- (\BBxmax, \coteCyl, 0);
          \filldraw[draw=white, fill=secondary!90, opacity=0.625] (\BBxmin+\amin, 0, 0) -- +(-\BBxmin, 0, 0) -- +(-\BBxmin, 0, \coteCyl) -- +(0, 0, \coteCyl) -- cycle;
          \filldraw[draw=white, fill=secondary!80, opacity=0.5] (\BBxmin+\amin, 0, \coteCyl) -- +(-\BBxmin, 0, 0) -- +(-\BBxmin, \coteCyl, 0) -- +(0, \coteCyl, 0) -- cycle;
          \filldraw[draw=white, fill=teal!90, opacity=0.875] (\amin, 0, 0) -- (\amin, \coteCyl, 0) -- (\amin, \coteCyl, \coteCyl) -- (\amin, 0, \coteCyl) -- cycle;
          \fill[teal!80, opacity=0.5] (\amin, 0, 0) -- +({\amax-\amin}, 0, 0) -- +({\amax-\amin}, \coteCyl, 0) -- +(0, \coteCyl, 0) -- cycle;
          \draw[-latex] (0, {\BBymin-0.0*(\BBymax-\BBymin)}, 0) -- (yaxis);%
          \fill[teal!90, opacity=0.625] (\amin, \coteCyl, 0) -- +({\amax-\amin}, 0, 0) -- +({\amax-\amin}, 0, \coteCyl) -- +(0, 0, \coteCyl) -- cycle;
          \node[teal!80!black] (Omega) at (0, 4.5*\coteCyl, 1.5*\coteCyl) {$\Omegaperper^{\emph{int}}$};
          \draw[dashed, white] (\amin, \coteCyl, 0) -- (\amax, \coteCyl, 0);
          \filldraw[draw=white, fill=teal!90, opacity=0.625] (\amin, 0, 0) -- +({\amax-\amin}, 0, 0) -- +({\amax-\amin}, 0, \coteCyl) -- +(0, 0, \coteCyl) -- cycle;
          \filldraw[draw=white, fill=teal!80, opacity=0.5] (\amin, 0, \coteCyl) -- +({\amax-\amin}, 0, 0) -- +({\amax-\amin}, \coteCyl, 0) -- +(0, \coteCyl, 0) -- cycle;
          \draw[-latex, teal!80!black] (Omega) to[bend right=30] (0, 0.5*\coteCyl, \coteCyl);
          \filldraw[draw=white, fill=tertiary!90, opacity=0.875] (\amax, 0, 0) -- (\amax, \coteCyl, 0) -- (\amax, \coteCyl, \coteCyl) -- (\amax, 0, \coteCyl) -- cycle;
          \fill[tertiary!80, opacity=0.5] (\BBxmax+\amax, 0, 0) -- +(-\BBxmax, 0, 0) -- +(-\BBxmax, \coteCyl, 0) -- +(0, \coteCyl, 0) -- cycle;
          \fill[tertiary!90, opacity=0.625] (\BBxmax+\amax, \coteCyl, 0) -- +(-\BBxmax, 0, 0) -- +(-\BBxmax, 0, \coteCyl) -- +(0, 0, \coteCyl) -- cycle;
          \draw[dashed, white] (\BBxmax+\amax, \coteCyl, 0) -- +(-\BBxmax, 0, 0);
          \filldraw[draw=white, fill=tertiary!90, opacity=0.625] (\BBxmax+\amax, 0, 0) -- +(-\BBxmax, 0, 0) -- +(-\BBxmax, 0, \coteCyl) -- +(0, 0, \coteCyl) -- cycle;
          \filldraw[draw=white, fill=tertiary!80, opacity=0.5] (\BBxmax+\amax, 0, \coteCyl) -- +(-\BBxmax, 0, 0) -- +(-\BBxmax, \coteCyl, 0) -- +(0, \coteCyl, 0) -- cycle;
          \draw[-latex] (0, 0, {\BBzmin-0.0*(\BBzmax-\BBzmin)}) -- (0, 0, {\BBzmax+0.25*(\BBzmax-\BBzmin)}) node[above] {$\smash{\zvi_2}$};
          \draw[-latex] ({\BBxmin+\amin-0.1*(\BBxmax-\BBxmin)}, 0, 0) -- (xaxis);%
          \draw (\amin, 0, -0.4) node {$a^-$};
          \draw (\amax, 0, -0.4) node {$a^+$};
        \end{scope}%
      \end{tikzpicture} 
    }
    \caption{\surligner{The waveguide solution $\widehat{U}^0_\fbvar$ can be constructed by solving a problem set in the interior domain $\Omegaperper^{\emph{int}}$, and half-guide problems similar to \eqref{H2Dmod:eq:augmented_halfguide_problem}.}\label{H2Dmod:fig:interior_domain}}
  \end{figure}
\end{itemize}

\subsection{Short-term perspective on junctions of general periodic half-spaces}\label{H2Dmod:sec:extension_general_config}
\surligner{We recall that the transmission configurations studied throughout the paper are characterized by a pair $(\aten, \rho)$, where $\aten: \R^2 \to \R^{2 \times 2}$ (\emph{resp.} $\rho: \R^2 \to \R^2$) coincides in $\R^2_\pm$ with a function $\aten^\pm$ (\emph{resp.} $\rho^\pm$) with periodicity vectors $\pv^\pm_1, \pv^\pm_2 \in \R^2$, as illustrated in Figure \ref{H2Dmod:fig:general_configuration} (see also Figure \ref{H2Dmod:fig:all_general_configurations} top). Configuration \eqref{H2Dmod:item:config_a} and Configuration \eqref{H2Dmod:item:config_b} were particular cases for which the associated transmission problem \eqref{H2Dmod:eq:transmission_problem} could be lifted into a $3$D augmented problem.

\vspace{1\baselineskip} \noindent
In this section, we describe how one can also derive an augmented formulation for arbitrary vectors $\pv^\pm_1$ and $\pv^\pm_2$. This will be the focus of Sections \ref{H2Dmod:sec:general_particular_A_B} and \ref{H2Dmod:sec:other_configurations}. We will rely on a notion of equivalence between two transmission configurations, which is introduced in Section \ref{H2Dmod:sec:equivalence_transmission_problems}.

\subsubsection{A notion of equivalence between two transmission problems}\label{H2Dmod:sec:equivalence_transmission_problems}
In this section, we work with pairs $(\mathbb{T}_+, \mathbb{T}_-)$ of matrices such that 
\begin{equation}\label{H2Dmod:eq:pties_T_pm}
  \textnormal{$\mathbb{T}_\pm \in \R^{2 \times 2}$ is invertible}, \quad \mathbb{T}_\pm\, \R^2_\pm = \R^2_\pm \quad \textnormal{and} \quad \mathbb{T}_+\, \ev_\zvi = \mathbb{T}_-\, \ev_\zvi = \ev_\zvi.
\end{equation}
From each pair, we can construct the piecewise-defined map $\mathbb{T} : \R^2 \to \R^2$:
\begin{equation}\label{H2Dmod:eq:def_T}
  \displaystyle
  \spforall \xv = (\xvi, \zvi) \in \R^2, \quad \mathbb{T}(\xv) := %
  \left\{
    \begin{array}{l@{\quad \textnormal{if}\ \ }l}
      \mathbb{T}_+\vts \xv & x > 0
      \\[5pt]
      \mathbb{T}_-\vts \xv & x < 0.
    \end{array}
  \right.
\end{equation}
Note that by \eqref{H2Dmod:eq:pties_T_pm}, $\mathbb{T}$ is a continuous bijection from $\R^2$ to $\R^2$, which is piecewise differentiable, with Jacobian matrix
\[\displaystyle
  \spforall \xv = (\xvi, \zvi) \in \R^2, 
  \quad \mathbb{J}_{\mathbb{T}} (\xv) = %
  \left\{
    \begin{array}{l@{\quad \textnormal{if}\ \ }l}
      \mathbb{T}_+ & x > 0
      \\[5pt]
      \mathbb{T}_- & x < 0,
    \end{array}
  \right.  
\]
so that it defines a continuous and piecewise smooth change of variables $\xv \mapsto \mathbb{T} \xv$. We shall say that two transmission configurations associated to pairs $(\aten, \rho)$ and $(\widetilde{\aten}, \widetilde{\rho})$ are \emph{equivalent} if and only if there exists a transformation $\mathbb{T}$ defined by (\ref{H2Dmod:eq:pties_T_pm}, \ref{H2Dmod:eq:def_T}), such that
\begin{equation}
  \spforall \xv \in \R^2, \qquad \widetilde{\aten}(\xv) = \mathbb{J}_{\mathbb{T}} (\xv)\; \aten\big(\mathbb{T}^{-1} (\xv)\big) \transp{[\mathbb{J}_{\mathbb{T}} (\xv)]} \quad \textnormal{and} \quad \widetilde{\rho}(\xv) = \rho\big(\mathbb{T}^{-1} (\xv)\big).
  \label{H2Dmod:eq:A_rho_transfo}
\end{equation}
The first interest of the notion of equivalence introduced above is that the solutions $u$ and $\widetilde{u}$ of the transmission problems \eqref{H2Dmod:eq:transmission_problem} associated to $(\aten, \rho)$ and $(\widetilde{\aten}, \widetilde{\rho})$, respectively, are simply related through the change of variables $\xv \mapsto \mathbb{T} \xv$. This is the object of the next proposition, which is easily proved using the chain rule.

\begin{prop}\label{H2Dmod:prop:alternative_configuration}
  Let $u$ and $\widetilde{u}$ be the respective solutions of the transmission problems \eqref{H2Dmod:eq:transmission_problem} associated to the pairs $(\aten, \rho)$ and $(\widetilde{\aten}, \widetilde{\rho})$ related by \eqref{H2Dmod:eq:A_rho_transfo}. Then
  \begin{equation}
    \displaystyle
    \widetilde{u}(\xv) = u\big(\mathbb{T}^{-1} (\xv)\big).
    \label{H2Dmod:eq:change_var_unknowns}
  \end{equation}
  %
  %
\end{prop}

\begin{figure}[ht!]
  \centering  
  \begin{subfigure}[t]{\textwidth}
    \makebox[\textwidth][c]{
      \begin{tikzpicture}
        \def\xbb{7.5}
        \def\ybb{3.5}
        \def\vectxneg{1}
        \def\vectyneg{0.25}
        \def\vecnxneg{-0.4}
        \def\vecnyneg{1.13137085}
        \def\vectxpos{0.8}
        \def\vectypos{-0.4}
        \def\vecnxpos{0.8}
        \def\vecnypos{0.848528137}
        \def\nbCellsXneg{6}
        \def\nbCellsYneg{3}
        \def\nbCellsXpos{7}
        \def\nbCellsYpos{4}
        \def\xinterneg{0}
        \def\xinterpos{0}

        \def\mgcg{0.2}
        \def\anglecoupe{60}
        \def\dimCells{1.25}

        \begin{scope}
          \clip (\xinterneg, 0) -- (\xinterpos, \ybb) -- (\xbb, \ybb) -- (\xbb, 0) -- cycle;
          \begin{scope}[scale=\dimCells, shift={(-0.32*\xbb/\dimCells, 0.3*\ybb/\dimCells)}]
            \foreach \x in {0, 1, ..., \nbCellsXpos} {
              \foreach \y in {0, 1, ..., \nbCellsYpos} {
                \filldraw[fill=gray!28, draw=gray!10, thick] ({\x*\vectxpos+\y*\vecnxpos}, {\x*\vectypos+\y*\vecnypos}) -- +(\vectxpos, \vectypos) -- +({\vectxpos+\vecnxpos}, {\vectypos+\vecnypos}) -- +(\vecnxpos, \vecnypos) -- cycle;
                \filldraw[fill=gray, draw=black!85, thick]  ({(\x+0.5)*\vectxpos+(\y+0.5)*\vecnxpos}, {(\x+0.5)*\vectypos+(\y+0.5)*\vecnypos}) circle[radius=0.22];
                %
              }
            }
            \draw[dashed, black!85] ({3*\vectxpos+2*\vecnxpos}, {3*\vectypos+2*\vecnypos}) -- +(\vectxpos, \vectypos) -- +({\vectxpos+\vecnxpos}, {\vectypos+\vecnypos}) -- +(\vecnxpos, \vecnypos) -- cycle;
            \draw[-latex, thick] ({3*\vectxpos+2*\vecnxpos}, {3*\vectypos+2*\vecnypos}) -- +(\vectxpos, \vectypos) node[pos=1.3] {$\qv^+_1$}; 
            \draw[-latex, thick] ({3*\vectxpos+2*\vecnxpos}, {3*\vectypos+2*\vecnypos}) -- +(\vecnxpos, \vecnypos) node[pos=1.3] {$\qv^+_2$}; 
          \end{scope}
        \end{scope}
        \begin{scope}
          \clip (\xinterneg, 0) -- (\xinterpos, \ybb) -- (-\xbb, \ybb) -- (-\xbb, 0) -- cycle;
          \begin{scope}[scale=\dimCells, shift={(-0.95*\xbb/\dimCells, -0.7*\ybb/\dimCells)}]
            \foreach \x in {0, 1, ..., \nbCellsXneg} {
              \foreach \y in {0, 1, ..., \nbCellsYneg} {
                \filldraw[fill=gray!28, draw=gray!10, thick] ({\x*\vectxneg+\y*\vecnxneg}, {\x*\vectyneg+\y*\vecnyneg}) -- +(\vectxneg, \vectyneg) -- +({\vectxneg+\vecnxneg}, {\vectyneg+\vecnyneg}) -- +(\vecnxneg, \vecnyneg) -- cycle;
                \filldraw[fill=gray, draw=black!85, thick] ({(\x+0.25)*\vectxneg+(\y+0.25)*\vecnxneg}, {(\x+0.25)*\vectyneg+(\y+0.25)*\vecnyneg}) -- +(0.475*\vectxneg, 0.475*\vectyneg) -- +({0.475*\vectxneg+0.475*\vecnxneg}, {0.475*\vectyneg+0.475*\vecnyneg}) -- +(0.475*\vecnxneg, 0.475*\vecnyneg)-- cycle;
                %
              }
            }
            \draw[dashed, black!85] ({3*\vectxneg+2*\vecnxneg}, {3*\vectyneg+2*\vecnyneg}) -- +(\vectxneg, \vectyneg) -- +({\vectxneg+\vecnxneg}, {\vectyneg+\vecnyneg}) -- +(\vecnxneg, \vecnyneg) -- cycle;
            \draw[-latex, thick] ({3*\vectxneg+2*\vecnxneg}, {3*\vectyneg+2*\vecnyneg}) -- +(\vectxneg, \vectyneg) node[pos=1.3] {$\qv_1^-$};
            \draw[-latex, thick] ({3*\vectxneg+2*\vecnxneg}, {3*\vectyneg+2*\vecnyneg}) -- +(\vecnxneg, \vecnyneg) node[pos=1.3] {$\qv_2^-$};
          \end{scope}
        \end{scope}
        \draw[line width=1mm, black!85] (\xinterneg, 0) -- (\xinterpos, \ybb);
        \draw[-latex] (-\xbb-0.25, 0) -- (\xbb+0.25, 0) node[right] {$x$};
        \draw[-latex] (0, -0.25) -- (0, \ybb+0.5) node[left] {$z$};

      \end{tikzpicture}
    }
  \end{subfigure}

  \vspace{1\baselineskip}
  \hfill
  \begin{subfigure}[t]{\textwidth}
    \makebox[\textwidth][c]{
      \begin{tikzpicture}
        \def\xbb{7.5}
        \def\ybb{3.5}
        \def\vectxneg{1.2}
        \def\vectyneg{0}
        \def\vecnxneg{-0.4}
        \def\vecnyneg{1.13137085}
        \def\vectxpos{1.2}
        \def\vectypos{0}
        \def\vecnxpos{0.8}
        \def\vecnypos{0.848528137}
        \def\nbCellsXneg{6}
        \def\nbCellsYneg{3}
        \def\nbCellsXpos{7}
        \def\nbCellsYpos{4}
        \def\xinterneg{0}
        \def\xinterpos{0}

        \def\mgcg{0.2}
        \def\anglecoupe{60}
        \def\dimCells{1.25}

        \begin{scope}
          \clip (\xinterneg, 0) -- (\xinterpos, \ybb) -- (\xbb, \ybb) -- (\xbb, 0) -- cycle;
          \begin{scope}[scale=\dimCells, shift={(-0.5*\xbb/\dimCells, 0*\ybb/\dimCells)}]
            \foreach \x in {0, 1, ..., \nbCellsXpos} {
              \foreach \y in {0, 1, ..., \nbCellsYpos} {
                \filldraw[fill=gray!28, draw=gray!10, thick] ({\x*\vectxpos+\y*\vecnxpos}, {\x*\vectypos+\y*\vecnypos}) -- +(\vectxpos, \vectypos) -- +({\vectxpos+\vecnxpos}, {\vectypos+\vecnypos}) -- +(\vecnxpos, \vecnypos) -- cycle;
                \filldraw[fill=gray, draw=black!85, thick]  ({(\x+0.5)*\vectxpos+(\y+0.5)*\vecnxpos}, {(\x+0.5)*\vectypos+(\y+0.5)*\vecnypos}) circle[radius=0.22];
                %
              }
            }
            \draw[dashed, black!85] ({4*\vectxpos+1*\vecnxpos}, {4*\vectypos+1*\vecnypos}) -- +(\vectxpos, \vectypos) -- +({\vectxpos+\vecnxpos}, {\vectypos+\vecnypos}) -- +(\vecnxpos, \vecnypos) -- cycle;
            \draw[-latex, thick] ({4*\vectxpos+1*\vecnxpos}, {4*\vectypos+1*\vecnypos}) -- +(\vectxpos, \vectypos) node[pos=1.2] {$\ev_\xvi$}; 
            \draw[-latex, thick] ({4*\vectxpos+1*\vecnxpos}, {4*\vectypos+1*\vecnypos}) -- +(\vecnxpos, \vecnypos) node[pos=1.2] {$\pv^+$}; 
          \end{scope}
        \end{scope}
        \begin{scope}
          \clip (\xinterneg, 0) -- (\xinterpos, \ybb) -- (-\xbb, \ybb) -- (-\xbb, 0) -- cycle;
          \begin{scope}[scale=\dimCells, shift={(-1.25*\xbb/\dimCells, 0.0*\ybb/\dimCells)}]
            \foreach \x in {0, 1, ..., \nbCellsXneg} {
              \foreach \y in {0, 1, ..., \nbCellsYneg} {
                \filldraw[fill=gray!28, draw=gray!10, thick] ({\x*\vectxneg+\y*\vecnxneg}, {\x*\vectyneg+\y*\vecnyneg}) -- +(\vectxneg, \vectyneg) -- +({\vectxneg+\vecnxneg}, {\vectyneg+\vecnyneg}) -- +(\vecnxneg, \vecnyneg) -- cycle;
                \filldraw[fill=gray, draw=black!85, thick] ({(\x+0.25)*\vectxneg+(\y+0.25)*\vecnxneg}, {(\x+0.25)*\vectyneg+(\y+0.25)*\vecnyneg}) -- +(0.475*\vectxneg, 0.475*\vectyneg) -- +({0.475*\vectxneg+0.475*\vecnxneg}, {0.475*\vectyneg+0.475*\vecnyneg}) -- +(0.475*\vecnxneg, 0.475*\vecnyneg)-- cycle;
                %
              }
            }
            \draw[dashed, black!85] ({4*\vectxneg+1*\vecnxneg}, {4*\vectyneg+1*\vecnyneg}) -- +(\vectxneg, \vectyneg) -- +({\vectxneg+\vecnxneg}, {\vectyneg+\vecnyneg}) -- +(\vecnxneg, \vecnyneg) -- cycle;
            \draw[-latex, thick] ({4*\vectxneg+1*\vecnxneg}, {4*\vectyneg+1*\vecnyneg}) -- +(\vectxneg, \vectyneg) node[pos=1.2] {$\ev_\xvi$};
            \draw[-latex, thick] ({4*\vectxneg+1*\vecnxneg}, {4*\vectyneg+1*\vecnyneg}) -- +(\vecnxneg, \vecnyneg) node[pos=1.2] {$\pv^-$};
          \end{scope}
        \end{scope}
        \draw[line width=1mm, black!85] (\xinterneg, 0) -- (\xinterpos, \ybb);
        \draw[-latex] (-\xbb-0.25, 0) -- (\xbb+0.25, 0) node[right] {$\xvi$};
        \draw[-latex] (0, -0.25) -- (0, \ybb+0.5) node[left] {$\zvi$};

      \end{tikzpicture}
    }
  \end{subfigure}
  \caption{\surligner{Junction of arbitrary periodic half-spaces (top). Thanks to Propositions \ref{H2Dmod:prop:alternative_configuration} and \ref{H2Dmod:prop:alternative_configuration_2}, this configuration can always be reduced to a junction of periodic half-spaces with one periodicity direction being $\ev_\xvi$ (bottom).}\label{H2Dmod:fig:all_general_configurations}}
\end{figure}

\vspace{1\baselineskip} \noindent
A useful property is that any transmission configuration $(\aten, \rho)$, with periodicity vectors $(\pv^\pm_1, \pv^\pm_2)$, is equivalent, in the sense of \eqref{H2Dmod:eq:A_rho_transfo}, to a transmission configuration $(\widetilde{\aten}, \widetilde{\rho})$ where one of the periodicity vectors in each half-space is $\ev_\xvi$ (see in Figure \ref{H2Dmod:fig:all_general_configurations} bottom).

\begin{prop}\label{H2Dmod:prop:alternative_configuration_2}
  Given a transmission configuration $(\aten, \rho)$ with periodicity vectors $(\pv^\pm_1, \pv^\pm_2)$, there exists a pair $(\widetilde{\aten}, \widetilde{\rho})$ equivalent to $(\aten, \rho)$ in the sense of \eqref{H2Dmod:eq:A_rho_transfo}, such that $\widetilde{\aten}$ (\emph{resp.} $\widetilde{\rho}$) concides in $\R^2_\pm$ with a $\Z \ev_\xvi + \Z \widetilde{\pv}^{\, \pm}$--periodic function $\widetilde{\aten}^{\, \pm}$ (\emph{resp.} $\widetilde{\rho}^{\, \pm}$), with $\widetilde{\pv}^{\, \pm} \in \R^2$.
\end{prop}

\begin{dem}
  Consider the matrices $\mathbb{T}_+, \mathbb{T}_- \in \R^{2 \times 2}$ defined by
  \begin{equation} \label{H2Dmod:eq:def_Tpm}
    \mathbb{T}_\pm\, \ev_\xvi = \pv^\pm_1 \quad \textnormal{and} \quad \mathbb{T}_\pm\, \ev_\zvi = \ev_\zvi.
  \end{equation}
  Such matrices are invertible if and only if $\smash{\pv^\pm_1 = (p^\pm_{1, \xvi}, p^\pm_{1, \zvi})}$ satisfies $p^\pm_{1, \xvi} \neq 0$. This assumption can be made without any loss of generality: in fact, since $\smash{\pv^\pm_1}$ and $\smash{\pv^\pm_2}$ are non-collinear, at least one of them is not collinear with $\ev_\zvi$. Moreover, because $(|p^\pm_{1, \xvi}|, p^\pm_{1, \zvi})$ is also a periodicity vector of $(\aten^\pm, \rho^\pm)$, we can assume that $p^\pm_{1, \xvi} > 0$, so that $\mathbb{T}_\pm\, \R^2_\pm = \R^2_\pm$. Finally, note that $\mathbb{T}_+\, \ev_\zvi = \mathbb{T}_-\, \ev_\zvi = \ev_\zvi$, by definition. Consequently, the pair $(\mathbb{T}_+, \mathbb{T}_-)$ satisfies the properties \eqref{H2Dmod:eq:pties_T_pm}. 

  \vspace{1\baselineskip}\noindent
  Finally, let $(\widetilde{\aten}, \widetilde{\rho})$ be defined by \eqref{H2Dmod:eq:A_rho_transfo}. Then using \eqref{H2Dmod:eq:def_Tpm}, one easily checks that $\widetilde{\aten}$ (\emph{resp.} $\widetilde{\rho}$) concides in $\R^2_\pm$ with a $\Z \ev_\xvi + \Z \widetilde{\pv}^{\, \pm}$--periodic function $\widetilde{\aten}^{\, \pm}$ (\emph{resp.} $\widetilde{\rho}^{\, \pm}$), with $\widetilde{\pv}^{\, \pm} := \mathbb{T}^{-1}_\pm\, \pv^\pm_2$.
\end{dem}
}

\subsubsection{General configurations that are equivalent to Configurations \texorpdfstring{\eqref{H2Dmod:item:config_a}}{(A)} and \texorpdfstring{\eqref{H2Dmod:item:config_b}}{(B)}}\label{H2Dmod:sec:general_particular_A_B}
\surligner{%
Proposition \ref{H2Dmod:prop:alternative_configuration_2} implies that any arbitrary junction of periodic half-spaces can be reduced to the case of Figure \ref{H2Dmod:fig:all_general_configurations} bottom. In particular,
\begin{enumerate}
  \item Any junction of $\Z \qv^\pm + \Z (p^\pm_\zvi \ev_\zvi)$--periodic half-spaces can be reduced to Configuration \eqref{H2Dmod:item:config_a};
  \item Any junction of a homogeneous half-space and a  $\Z \qv + \Z \pv$--periodic half-space can be reduced to Configuration \eqref{H2Dmod:item:config_b}.
\end{enumerate}

\noindent
\surligner{Therefore, by Proposition \ref{H2Dmod:prop:alternative_configuration}, the solution of \eqref{H2Dmod:eq:transmission_problem} for a junction of $\Z \qv^\pm + \Z (p^\pm_\zvi \ev_\zvi)$--periodic half-spaces (\emph{resp.} the junction of a homogeneous half-space and a  $\Z \qv + \Z \pv$--periodic half-space) can be obtained using the method developed in this paper for Configuration \eqref{H2Dmod:item:config_a} (\emph{resp.} \eqref{H2Dmod:item:config_b}).}

\subsubsection{Other configurations}\label{H2Dmod:sec:other_configurations}
\surligner{We now address the general case of a transmission configuration $(\aten, \rho)$ with arbitrary periodicity vectors.} By Proposition \ref{H2Dmod:prop:alternative_configuration_2}, one can assume without any loss of generality that the restriction $\aten^\pm$ of $\aten$ to $\R^2_\pm$ is $\Z \ev_\xvi + \Z \pv^\pm$--periodic, with $\pv^\pm \in \R^2$ (and similarly for $\rho$). By using the arguments developed in Section \ref{H2Dmod:sec:augmented_structure_config_b} for Configuration \eqref{H2Dmod:item:config_b}, we have the existence of $\cuti^\pm_1, \cuti^\pm_2 \in \R$ and $\Z^3$--periodic functions $\aten^\pm_p$ such that
\begin{equation*}\label{H2Dmod:eq:each_halfspace_lift}
  \displaystyle
  \spforall \xv \in \R^2, \quad 
  \left\{%
    \begin{array}{r@{\ =\ }l}
      \aten^+ (\xvi, \zvi) & \aten^+_p (\xvi, \cuti^+_1\vts \zvi, \cuti^+_2\vts \zvi) = \aten^+_p (\cutmat^+\, \xv),
      \ret
      \aten^- (\xvi, \zvi) & \aten^-_p (\xvi, \cuti^-_1\vts \zvi, \cuti^-_2\vts \zvi) = \aten^-_p (\cutmat^-\, \xv), %
    \end{array}
  \right. %
  \qquad %
  \cutmat^\pm := %
    \begin{pmatrix}
      1 & 0
      \\
      0 & \cuti^\pm_1
      \\
      0 & \cuti^\pm_2
    \end{pmatrix}.
\end{equation*}
%
\surligner{It is then natural to ask if one can lift the overall coefficient $\aten$ in a similar manner.} This corresponds to finding $(i)$ an integer $\ord \geq 2$, $(ii)$ a $\ord$D coefficient $\widetilde{\aten}_p$ which coincides in $\R^\ord_\pm$ with a continuous function $\widetilde{\aten}^\pm_p$ which is $\Z^\ord$--periodic, and $(iii)$ a matrix $\widetilde{\cutmat} \in \R^{\ord \times 2}$, such that
\begin{equation*}
  \displaystyle
  \spforall \xv \in \R^2, \quad \aten (\xv) = \widetilde{\aten}_p (\widetilde{\cutmat}\, \xv).
\end{equation*}
\surligner{Note that in general, one cannot choose $\widetilde{\aten}^\pm_p = \aten^\pm_p$, since the matrices $\cutmat^+$ and $\cutmat^-$ do not coincide. This is the reason why we need to introduce additional variables to construct $(n, \widetilde{\aten}_p, \widetilde{\cutmat})$.}

\begin{enumerate}
  \item \textbf{\surligner{General case: reduction to an augmented $5$D problem}}. One observes that in general, $\aten$ admits an extension $\widetilde{\aten}_p$ which is periodic, but defined in $\R^5$. In fact, it can be easily checked that $\aten (\xv) = \widetilde{\aten}_p (\widetilde{\cutmat}\, \xv)$ for any $\xv \in \R^2$, where
  \begin{equation*}
    \displaystyle
    \widetilde{\cutmat} := %
    \begin{pmatrix}
      1 & 0
      \\
      0 & \cuti^+_1
      \\
      0 & \cuti^+_2
      \\
      0 & \cuti^-_1
      \\
      0 & \cuti^-_2
    \end{pmatrix}
    \quad \textnormal{and} \quad 
    \left\{
      \begin{array}{r@{\ :=\ }l}
        \widetilde{\aten}^+_p (\xvi, \zvi^+_1, \zvi^+_2, \zvi^-_1, \zvi^-_2) & \aten^-_p (\xvi, \zvi^+_1, \zvi^+_2)
        \ret
        \widetilde{\aten}^-_p (\xvi, \zvi^+_1, \zvi^+_2, \zvi^-_1, \zvi^-_2) & \aten^-_p (\xvi, \zvi^-_1, \zvi^-_2).
      \end{array}
    \right.
  \end{equation*}
  \surligner{The process for constructing $\widetilde{\aten}$ and $\widetilde{\cutmat}$ will be presented in a subsequent paper.}

  \item \textbf{\surligner{Further reduction to a $4$D augmented problem in the rationally dependent case}}. We say that the family $\{\cuti^+_1, \cuti^+_2, \cuti^-_1, \cuti^-_2\}$ is \emph{rationally dependent} if there exists a non-zero vector of integers $(m^+_1, m^+_2, m^-_1, m^-_2) \in \Z^4 \setminus \{0\}$ such that $m^+_1\, \cuti^+_1 + m^+_2\, \cuti^+_2 + m^-_1\, \cuti^-_1 + m^-_2\, \cuti^-_2 = 0$. In this case, if we assume that $m^+_1 \neq 0$ for simplicity, then one can avoid introducing an additional variable, for instance $\zvi^+_1$, by observing that 
  $\aten (\xv) = \widetilde{\aten}_p (\widetilde{\cutmat}\, \xv)$ with $\ord = 4$,
  \begin{equation*}
    \widetilde{\cutmat} := %
    \begin{pmatrix}
      1 & 0
      \\
      0 & \cuti^+_2 / m^+_1
      \\
      0 & \cuti^-_1 / m^+_1
      \\
      0 & \cuti^-_2 / m^+_1
    \end{pmatrix}
    \ \ \textnormal{and} \ \
    \left\{
      \begin{array}{r@{\ :=\ }l}
        \widetilde{\aten}^+_p (\xvi, \zvi^+_2, \zvi^-_1, \zvi^-_2) & \aten^+_p (x, -m^+_2\vts \zvi^+_2 - m^-_1 \zvi^-_1 - m^-_2 \zvi^-_2,  m^+_1 \zvi^+_2)
        \ret
        \widetilde{\aten}^-_p (\xvi, \zvi^+_2, \zvi^-_1, \zvi^-_2) & \aten^-_p (x, m^+_1 \zvi^-_1, m^+_1 \zvi^-_2).
      \end{array}
    \right.
  \end{equation*}
  Note that one could find an augmented structure with dimension $\ord < 4$ if there exists another non-trivial integer combination $\widehat{m}^+_1\, \cuti^+_1 + \widehat{m}^+_2\, \cuti^+_2 + \widehat{m}^-_1\, \cuti^-_1 + \widehat{m}^-_2\, \cuti^-_2 = 0$. \surligner{This would correspond to the transmission configurations in Section \ref{H2Dmod:sec:general_particular_A_B}.}
\end{enumerate}

\noindent
\surligner{If the family $\{\cuti^+_1, \cuti^+_2, \cuti^-_1, \cuti^-_2\}$ is \emph{rationally independent}, that is $m^+_1 \cuti^+_1 + m^+_2 \cuti^+_2 + m^-_1 \cuti^-_1 + m^-_2 \cuti^-_2 \neq 0$ for any $(m^+_1, m^+_2, m^-_1, m^-_2) \in \Z^4 \setminus \{0\}$, then we believe that one cannot find an augmented structure with dimension $\ord < 5$, unless if $\aten$ is constant in one or multiple directions. This will be addressed in a future work.}
}

\subsection{The case without absorption}\label{H2Dmod:sec:extension_non_absorbing}
In the case where there is no absorption, the application of the method in this paper raises several challenges. First, some waves guided by the interface may exist (depending on the coefficients), so that the physical solution may propagate along the interface without attenuation. The use of the Floquet-Bloch transform along the 3D lifted interface is then not possible. If the solution is $L^2$ along the interface (in the absence of these interface guided modes), the Floquet-Bloch transform can be applied and the numerical method can be formally extended (although one has to add conditions to select the outgoing waves, see \cite{fliss2016solutions,fliss2021dirichlet}). A corresponding formal approach has already been used for the rational case \cite{flisscassanbernier2010}. However, even in the rational case, the rigorous justification of the numerical method is still a difficult challenge.

}{%
}

\appendix

\ifthenelse{\boolean{brouillon}}{%
  \section{Anisotropic Sobolev spaces}\label{H2Dmod:sec:proof_functional_framework}
\subsection{Proof of Proposition \ref{H2Dmod:prop:shear_properties}: properties of the shear map \texorpdfstring{$\shearmap$}{Stheta}}\label{H2Dmod:sec:proof_shear_properties}
To prove Proposition \ref{H2Dmod:prop:shear_properties}, we rely on the following easy but useful lemma, which links the integral of $\boldsymbol{V}$ on $\mathbb{\Omega}_\periodic$ and of $\extper[2]{\boldsymbol{V}}$ on $\mathbb{\Omega}_{\periodic, \cuti}$ for $\boldsymbol{V} \in [L^1(\mathbb{\Omega}_\periodic)]^d$. We refer to \cite[Lemma 3.12]{amenoagbadji2023wave} for its proof.
\begin{lem}\label{H2Dmod:lem:intg_Qtheta_Q_are_equal}
  For any $\boldsymbol{V} \in [L^1(\mathbb{\Omega}_\periodic)]^d$, we have
  \begin{equation}
    \displaystyle
    \int_{\mathbb{\Omega}_\periodic} \boldsymbol{V} = \int_{\mathbb{\Omega}_{\periodic, \cuti}} \extper[2]{\boldsymbol{V}},
  \end{equation}
  where $\extper[2]{\boldsymbol{V}} \in L^1_{\textit{loc}}(\mathbb{\Omega})^d$ is defined by \eqref{H2Dmod:eq:periodic_extension_e2}.
\end{lem}

\begin{dem}[of Proposition \ref{H2Dmod:prop:shear_properties}]
  The mapping $\operatorname{T}_\cuti : (\xv, s) \in \Qcuti \times (0, 1) \mapsto \cutmat\vts \xv + s\vts \eva_2 \in \mathbb{\Omega}_{\periodic, \cuti}$ is a $\mathscr{C}^1$--diffeomorphism with a non-vanishing \surligner{Jacobian determinant $\cuti_1 \neq 0$. Hence}, the associated change of variables leads to
  \[
    \displaystyle
    \spforall \boldsymbol{U}, \boldsymbol{V} \in [\mathscr{C}^\infty_{0, \periodic}(\overline{\mathbb{\Omega}_\periodic})]^d, \quad \frac{1}{\cuti_1} \int_0^1 \int_{\Qcuti} \shearmap \boldsymbol{U}(\xv, s)\cdot \overline{\shearmap \boldsymbol{V}(\xv, s)}\; d\xv ds = \int_{\mathbb{\Omega}_{\periodic, \cuti}} (\extper[2]{\boldsymbol{U}})\cdot (\overline{\extper[2]{\boldsymbol{V}}}).
  \]
  We then use the identity $(\extper[2]{\boldsymbol{U}})\cdot (\overline{\extper[2]{\boldsymbol{V}}}) = \extper[2]{(\boldsymbol{U}\cdot \overline{\boldsymbol{V}})}$, and we apply Lemma \ref{H2Dmod:lem:intg_Qtheta_Q_are_equal} to $\boldsymbol{U}\cdot \overline{\boldsymbol{V}}$, to deduce \eqref{H2Dmod:eq:oblique_plancherel} for $\boldsymbol{U}, \boldsymbol{V} \in [\mathscr{C}^\infty_{0, \periodic}(\overline{\mathbb{\Omega}_\periodic})]^d$. Finally, the density of $[\mathscr{C}^\infty_{0, \periodic}(\overline{\mathbb{\Omega}_\periodic})]^d$ in $[L^2 (\mathbb{\Omega}_\periodic)]^d$ leads to \eqref{H2Dmod:eq:oblique_plancherel}. Moreover, by choosing $\boldsymbol{U} = \boldsymbol{V}$ in \eqref{H2Dmod:eq:oblique_plancherel}, it follows that $\shearmap$ is bounded from $[L^2 (\mathbb{\Omega}_\periodic)]^d$ to $L^2(0, 1; [L^2(\Qcuti)]^d)$.

  \vspace{1\baselineskip} \noindent
  The bijectivity of $\shearmap$ results directly from the inverse of $\operatorname{T}_\cuti$, which leads to the expression \eqref{H2Dmod:eq:inverse_shear} of $\invshearmap$. The continuity of $\invshearmap$ is then implied by \eqref{H2Dmod:eq:oblique_plancherel}.
\end{dem}

\subsection{Proof of Proposition \ref{H2Dmod:prop:link_H1Cper_H1C_and_HdivCper_HdivC}: Action of the shear map \texorpdfstring{$\shearmap$}{Stheta} on differential operators}\label{H2Dmod:sec:proof_link_H1Cper_H1C_and_HdivCper_HdivC}
\begin{dem}[of Proposition \ref{H2Dmod:prop:link_H1Cper_H1C_and_HdivCper_HdivC}]
  We only prove \eqref{H2Dmod:eq:shear_commutes_with_grad} since the proof of \eqref{H2Dmod:eq:shear_commutes_with_dive} is very similar. The proof is merely a weak adaptation of the chain rule. Let $V \in \Hgradper{\mathbb{\Omega}_\periodic}$. To obtain the expression of $\transp{\cutmat} \bsnabla V$ in the sense of distributions, consider a test function $\boldsymbol{W} \in [\mathscr{C}^\infty_0(\mathbb{\Omega}_\periodic)]^2$ and define $\boldsymbol{W}_\cutmat := \shearmap \boldsymbol{W} \in [\mathscr{C}^\infty_0(\Qcuti \times (0, 1))]^2$. Then $\extper[2]{\boldsymbol{W}} \in [\mathscr{C}^\infty(\mathbb{\Omega}_\periodic)]^2$ since $\boldsymbol{W} \in [\mathscr{C}^\infty_0(\mathbb{\Omega}_\periodic)]^2$ and one easily computes using the chain rule that
  \begin{equation}\label{H2Dmod:eq:link_H1Cper_H1C_and_HdivCper_HdivC:0}
    \displaystyle
    \spforall (\xv, s) \in \Qcuti \times (0, 1), \quad \shearmap (\transp{\bsnabla} \cutmat \boldsymbol{W})(\xv, s) = \transp{\nabla}_\xv \boldsymbol{W}_\cutmat (\xv, s).
  \end{equation}
  Let also $V_\cutmat := \shearmap V$. Using derivation in the sense of distributions with $V \in L^2(\mathbb{\Omega}_\periodic)$ leads to:
  \begin{align}
    \big\langle \transp{\cutmat} \bsnabla V,\; \boldsymbol{W} \big\rangle_{[\mathscr{C}^\infty_0(\mathbb{\Omega}_\periodic)^2]', \mathscr{C}^\infty_0(\mathbb{\Omega}_\periodic)^2} &
    = - \int_{\mathbb{\Omega}_\periodic} V(\xva) \; \overline{\transp{\bsnabla} \cutmat \boldsymbol{W} (\xva)} \; d\xva \nonumber 
    \\
    &= - \frac{1}{\cuti_1} \int_0^1 \int_{\Qcuti} \shearmap V (\xv, s) \  \overline{\shearmap (\transp{\bsnabla} \cutmat \boldsymbol{W}) (\xv, s)}\; d\xv ds \quad \textnormal{from \eqref{H2Dmod:eq:oblique_plancherel}} \nonumber
    \\
    &= - \frac{1}{\cuti_1} \int_0^1 \int_{\Qcuti} V_\cutmat (\xv, s) \; \overline{\transp{\nabla}_\xv \boldsymbol{W}_\cutmat (\xv, s)}\; d\xv ds \quad \textnormal{from \eqref{H2Dmod:eq:link_H1Cper_H1C_and_HdivCper_HdivC:0}}. \label{H2Dmod:eq:link_H1Cper_H1C_and_HdivCper_HdivC:1_0}
  \end{align}
  But the definition of $\Hgradper{\mathbb{\Omega}_\periodic}$ implies that $V_\cutmat(\cdot, s) \in H^1(\Qcuti)$, $\aeforall s \in (0, 1)$. Thus in \eqref{H2Dmod:eq:link_H1Cper_H1C_and_HdivCper_HdivC:1_0}, we can apply the usual Green formula to $V_\cutmat(\cdot, s)$ and $\boldsymbol{W}_\cutmat (\cdot, s) \in [\mathscr{C}^\infty_0(\Qcuti)]^2$ (with no boundary term since $\boldsymbol{W}_\cutmat (\cdot, s)$ is compactly supported) to obtain
  \begin{align}
    \big\langle \transp{\cutmat} \bsnabla V,\; \boldsymbol{W} \big\rangle_{[\mathscr{C}^\infty_0(\mathbb{\Omega}_\periodic)^2]', \mathscr{C}^\infty_0(\mathbb{\Omega}_\periodic)^2} &=  \surligner{\frac{1}{\cuti_1}} \int_0^1 \int_{\Qcuti} \nabla_\xv V_\cutmat (\xv, s) \cdot \overline{\boldsymbol{W}_\cutmat (\xv, s)}\; d\xv ds \nonumber
    \\
    &= \int_{\mathbb{\Omega}_\periodic} \invshearmap \nabla_\xv V_\cutmat (\xva) \cdot \overline{\boldsymbol{W}(\xva)} \; d\xva \quad \textnormal{from \eqref{H2Dmod:eq:oblique_plancherel}} \nonumber
    \\
    &= \big\langle \invshearmap \nabla_\xv V_\cutmat,\; \boldsymbol{W} \big\rangle_{[\mathscr{C}^\infty_0(\mathbb{\Omega}_\periodic)^2]', \mathscr{C}^\infty_0(\mathbb{\Omega}_\periodic)^2}.\label{H2Dmod:eq:link_H1Cper_H1C_and_HdivCper_HdivC:1}
  \end{align}
  The identity \eqref{H2Dmod:eq:link_H1Cper_H1C_and_HdivCper_HdivC:1} being true for any $\boldsymbol{W} \in [\mathscr{C}^\infty_0(\mathbb{\Omega}_\periodic)]^2$ implies that $\transp{\cutmat} \bsnabla V = \invshearmap \nabla_\xv V_\cutmat$ in the sense of distributions. But thanks to the bijectivity of $\shearmap$ and to the fact that $\nabla_\xv V_\cutmat \in L^2(0, 1; [L^2(\Qcuti)]^2)$, it follows that $\invshearmap \nabla_\xv V_\cutmat \in [L^2(\mathbb{\Omega}_\periodic)]^2$. Therefore, $\transp{\cutmat} \bsnabla V$ also belongs to $[L^2(\mathbb{\Omega}_\periodic)]^2$, and one has
  \begin{equation}
    \displaystyle
    \aeforall \xva \in \mathbb{\Omega}_\periodic, \quad \transp{\cutmat} \bsnabla V (\xva) = \invshearmap \nabla_\xv V_\cutmat (\xva).
    \label{H2Dmod:eq:link_H1Cper_H1C_and_HdivCper_HdivC:2}
  \end{equation}
  Consequently, $V \in \Hgradper{\mathbb{\Omega}_\periodic}$ and \eqref{H2Dmod:eq:shear_commutes_with_grad} follows by applying the transform $\shearmap$ to \eqref{H2Dmod:eq:link_H1Cper_H1C_and_HdivCper_HdivC:2}.
\end{dem}

\subsection{Proof of Proposition \ref{H2Dmod:prop:Cinftyper_dense_H1Cper}: density of smooth functions in \texorpdfstring{$\Hgradper{\mathbb{\Omega}_\periodic}, \Hdiveper{\mathbb{\Omega}_\periodic}$}{H1theta, Hdivtheta}}\label{H2Dmod:sec:Cinftyper_dense_H1Cper}
\begin{dem}[of Proposition \ref{H2Dmod:prop:Cinftyper_dense_H1Cper}]
  We prove the first item, since the proof of the second one is very similar. For $V \in \Hgradper{\mathbb{\Omega}_\periodic}$, we construct a sequence of smooth functions that converges to $\shearmap V \in L^2(0, 1; H^1(\Qcuti))$. For almost any $\xv \in \Qcuti$, let $\zeta_n(\xv, \cdot)$ denote the $n$--th partial Fourier sum of $\shearmap V (\xv, \cdot)$:
  \begin{equation*}
    \displaystyle
    \aeforall (\xv, s) \in \Qcuti \times (0, 1), \quad \zeta_n (\xv, s) := \sum_{\ell = -n}^n \widehat{V}^{(\ell)}_\cutmat (\xv)\, \euler^{2\icplx \pi j s}, \quad \widehat{V}^{(\ell)}_\cutmat (\xv) := \int_0^1 \shearmap V (\xv, t)\, \euler^{- 2\icplx \pi \ell t}\, dt.
  \end{equation*}
  It can be checked easily that $\zeta(\xv, \cdot) \in \mathscr{C}^\infty([0, 1])$ is $1$--periodic, and that $\zeta(\cdot, s) \in H^1(\Qcuti)$ (as a \emph{finite} linear combination of $\widehat{V}^{(j)}_\cutmat \in H^1(\Qcuti)$). Moreover, using dominated convergence theorem, one shows that $\zeta_n \to \shearmap V$ in $L^2(0, 1; H^1(\Qcuti))$, $n \to +\infty$.

  \vspace{1\baselineskip} \noindent
  Now, consider a sequence of functions $\widehat{V}^{(\ell, m)}_\cutmat \in \mathscr{C}^\infty_0(\overline{\Qcuti})$ converging to $\widehat{V}^{(\ell)}_\cutmat$ in $H^1(\Qcuti)$, and define
  \begin{equation*}
    \displaystyle
    \spforall (\xv, s) \in \Qcuti \times (0, 1), \quad \zeta_{n, m} (\xv, s) := \sum_{\ell = -n}^n \widehat{V}^{(\ell, m)}_\cutmat (\xv)\, \euler^{2\icplx \pi \ell s}.
  \end{equation*}
  For any $(m, n) \in \N^2$, $\zeta_{n, m} \in \mathscr{C}^\infty(\overline{\Qcuti \times (0, 1)})$ is $1$--periodic with respect to $s$. Moreover, from the above, we deduce by the triangle inequality that $\zeta_{n, n} \to \shearmap V$ in $L^2(0, 1; H^1(\Qcuti))$, $n \to +\infty$. As a consequence, $\invshearmap \zeta_{n, n} \in \mathscr{C}^\infty_{0, \periodic}(\overline{\mathbb{\Omega}_\periodic}) \to V$ in $\Hgradper{\mathbb{\Omega}_\periodic}$, $n \to +\infty$. 
\end{dem}
  \section{Proof of Proposition \ref{H2Dmod:prop:TFB_properties}: properties of the Floquet-Bloch transform}\label{H2Dmod:sec:proof_FB_properties}
We prove the properties of the Floquet-Bloch transform stated in Proposition \ref{H2Dmod:prop:TFB_properties} by adapting the ideas of the isotropic case \cite{kuchment1993floquet}.

\begin{dem}
  \textbf{\textit{Item $(a)$}, Step $1$}-- We first show that $\fbtransform[]$ is bounded from $\Hgradper{\Omegaper}$ to  $L^2(-\pi, \pi; \Hgradperper{\Omegaperper})$, and that \eqref{H2Dmod:eq:fb_grad} holds. First, one computes directly that
  \begin{equation*}
    \spforall \zeta \in \mathscr{C}^\infty_{0, \periodic}(\Omegaper), \quad \aeforall \fbvar \in (-\pi, \pi), \quad \left\{
      \begin{array}{l}
        \fbtransform[] \zeta(\cdot, \fbvar) \in \mathscr{C}^\infty_{0, \perperiodic}(\overline{\Omegaperper}),
        \ret
        \fbtransform[] (\gradcut \zeta)\vts (\cdot, \fbvar) = \gradcutFB\, \fbtransform[] \zeta (\cdot, \fbvar).
      \end{array}
    \right.
  \end{equation*}
  This implies from the Parseval-like formula \eqref{H2Dmod:eq:partial_TFB_plancherel} (with $\boldsymbol{U} = \boldsymbol{V} = \zeta, \gradcut \zeta$) that $\fbtransform[]$ is continuous from $\mathscr{C}^\infty_{0, \periodic}(\Omegaper)$ to $L^2(-\pi, \pi; \mathscr{C}^\infty_{0, \perperiodic}(\overline{\Omegaperper}))$ (which are equipped respectively with the $H^1_\cutmat(\Omegaper)$--norm and the $L^2(-\pi, \pi; H^1_\cutmat(\Omegaperper))$--norm). By \emph{density} (Proposition \ref{H2Dmod:prop:Cinftyper_dense_H1Cper}), $\fbtransform[]$ extends to a continuous map from $\Hgradper{\Omegaper}$ to $L^2(-\pi, \pi; \Hgradperper{\Omegaperper})$, which satisfies \eqref{H2Dmod:eq:fb_grad}.

  \vspace{1\baselineskip} \noindent
  \textbf{\textit{Item $(a)$}, Step $2$}-- It remains to show that $\invfbtransform[]$ (see \eqref{H2Dmod:eq:partial_TFB_inversion}) is bounded from $L^2(-\pi, \pi; \Hgradperper{\Omegaperper})$ to $\Hgradper{\Omegaper}$. To do so, it is sufficient to prove that 
  \begin{equation}\label{H2Dmod:eq:fb_grad:dem_1}
    \displaystyle
    \spforall \widehat{V} \in L^2(-\pi, \pi; \Hgradperper{\Omegaperper}), \quad \gradcut \invfbtransform[] \widehat{V} = \invfbtransform[] [\gradcutFB\, \widehat{V}].
  \end{equation}
  For isotropic Sobolev spaces ($\cutmat = \mathbb{I}_3$), the classical way to prove \eqref{H2Dmod:eq:fb_grad:dem_1} is by using the jump rule, which involves traces on the faces $\zvi_1 \in \{0, 1\}$. But since we have not defined these traces in this paper, we shall resort to the Green's formula \eqref{H2Dmod:eq:Green_cylinder}, obtained by density. Given $\boldsymbol{W} \in \mathscr{C}^\infty_0(\Omegaper)$, using derivation in the sense of distributions with $\invfbtransform[] \widehat{V} \in L^2(\Omegaper)$ gives
  \begin{align}
    \big\langle \gradcut \invfbtransform[] \widehat{V}, \boldsymbol{W} \big\rangle_{[\mathscr{C}^\infty_0(\Omegaper)^2]', \mathscr{C}^\infty_0(\Omegaper)^2} &= \int_{\Omegaper} \invfbtransform[] \widehat{V}(\xva)\, \overline{\divecut \boldsymbol{W}(\xva)} \; d\xva  \nonumber
    \\
    &= \int_{-\pi}^\pi \int_{\Omegaperper} \widehat{V}(\xva, \fbvar)\, \overline{\fbtransform[] \divecut \boldsymbol{W}(\xva, \fbvar)} \; d\xva\, d\fbvar \quad \textnormal{from \eqref{H2Dmod:eq:partial_TFB_plancherel}} \nonumber
    \\
    &= \int_{-\pi}^\pi \int_{\Omegaperper} \widehat{V}(\xva, \fbvar)\, \overline{\divecutFB \fbtransform[] \boldsymbol{W}(\xva, \fbvar)} \; d\xva d\fbvar,\label{H2Dmod:eq:fb_grad:dem_2}
  \end{align}
  where the last equality is obtained from direct computations. Since $\widehat{V}(\cdot, \fbvar) \in \Hgradperper{\Omegaperper}$ and $\boldsymbol{W}$ is smooth, we can apply Green's formula \eqref{H2Dmod:eq:Green_cylinder} to \eqref{H2Dmod:eq:fb_grad:dem_2}, to deduce that
  \begin{align*}
    \big\langle \gradcut \invfbtransform[] \widehat{V}, \boldsymbol{W} \big\rangle_{[\mathscr{C}^\infty_0(\Omegaper)^2]', \mathscr{C}^\infty_0(\Omegaper)^2} &= \int_{-\pi}^\pi \int_{\Omegaperper} \gradcutFB\vts \widehat{V}(\xva, \fbvar)\, \overline{\fbtransform[] \boldsymbol{W}(\xva, \fbvar)} \; d\xva d\fbvar
    \\
    &= \int_{\Omegaper} \invfbtransform[] [\gradcutFB \widehat{V}](\xva)\, \overline{\boldsymbol{W}(\xva)}\; d\xva \quad \textnormal{from \eqref{H2Dmod:eq:partial_TFB_plancherel}}
    \\
    &= \big\langle \invfbtransform[] [\gradcutFB \widehat{V}], \boldsymbol{W} \big\rangle_{[\mathscr{C}^\infty_0(\Omegaper)^2]', \mathscr{C}^\infty_0(\Omegaper)^2}.
  \end{align*}
  The above implies \eqref{H2Dmod:eq:fb_grad:dem_1} in the sense of distributions. But since $\invfbtransform[] [\gradcutFB \widehat{V}] \in L^2(\Omegaper)$, we deduce that $\gradcut \invfbtransform[] \widehat{V} \in L^2(\Omegaper)$, so that \eqref{H2Dmod:eq:fb_grad:dem_1} holds.

  \vspace{1\baselineskip} \noindent
  \textbf{\textit{Item $(b)$}}-- The proof is very similar to the one of Item $(a)$, and therefore is omitted.

  \vspace{1\baselineskip} \noindent
  \textbf{\textit{Item $(c)$}}--  For the sake of clarity in this part, we highlight the distinction between the \emph{volume} Floquet-Bloch transform $\fbtransform[]_v$ for functions in $\Omegaper$ and the $\emph{surface}$ Floquet-Bloch transform $\fbtransform[]_s$ for functions on $\Sigmaper$. In addition, let $\gamma_{0, \periodic}$ (\emph{resp.} $\gamma_{0, \perperiodic}$) denote the trace operator on $\Sigmaper$ (\emph{resp.} $\Sigmaperper$). We begin by proving that the Floquet-Bloch transform commutes with the trace operator in the following sense:
  \begin{equation}\label{H2Dmod:eq:fbtransform_gamma0}
    \displaystyle
    \spforall V \in \Hgradper{\Omegaper}, \quad \fbtransform[]_s \, (\gamma_{0, \periodic} V) = \gamma_{0, \perperiodic}\, (\fbtransform[]_v V).
  \end{equation}
  Note that \eqref{H2Dmod:eq:fbtransform_gamma0} is straightforward for smooth functions $V \in \mathscr{C}^\infty_{0, \periodic}(\overline{\Omegaper})$. It extends to $V \in \Hgradper{\Omegaper}$ using $(i)$ the density of $\smash{\mathscr{C}^\infty_{0, \periodic}(\overline{\Omegaper})}$ in $\smash{\Hgradper{\Omegaper}}$ (Proposition \ref{H2Dmod:prop:Cinftyper_dense_H1Cper}), $(ii)$ the continuity of $\gamma_{0, \periodic}$ and $\gamma_{0, \perperiodic}$ (Proposition \ref{H2Dmod:prop:gamma0_continuous_surjective}), $(iii)$ the continuity of $\fbtransform[]_v$ from $\Hgradper{\Omegaper}$ to $L^2(-\pi, \pi; \Hgradperper{\Omegaperper})$ (proved in Item $(a)$), and finally $(iv)$ the continuity of $\fbtransform[]_s$ from $L^2(\Sigmaper)$ to $L^2(-\pi, \pi; L^2(\Sigmaperper))$.

  \vspace{1\baselineskip} \noindent
   Now we can show Item $(c)$, namely that the surface Floquet-Bloch transform $\fbtransform[]_s$ is an isomorphism from $\Honehalfper{\Sigmaper}$ to $L^2(-\pi, \pi; \Honehalfperper{\Sigmaperper})$. Since $\gamma_{0, \periodic}$ and $\gamma_{0, \perperiodic}$ are surjective, there exist two bounded operators $\mathcal{R}_\periodic \in \mathscr{L}(\smash{\Honehalfper{\Sigmaper}}, \Hgradper{\Omegaper})$ and $\mathcal{R}_\perperiodic \in \mathscr{L}(\smash{\Honehalfperper{\Sigmaperper}}, \Hgradperper{\Omegaperper})$ such that $\gamma_{0, \periodic}\, \mathcal{R}_\periodic = I$ and $\gamma_{0, \perperiodic}\, \mathcal{R}_\perperiodic = I$. Thus \eqref{H2Dmod:eq:fbtransform_gamma0} leads to
  \[\displaystyle
    \fbtransform[]_s = \gamma_{0, \perperiodic}\, \fbtransform[]_v\, \mathcal{R}_\periodic \quad \textnormal{and} \quad [\fbtransform[]_s]^{-1} = \gamma_{0, \periodic}\, [\fbtransform[]_v]^{-1}\, \mathcal{R}_\perperiodic.
  \]
  Using Item $(a)$ then allows to conclude.
\end{dem}

}{%
}

\newpage
\printbibliography

@book{kittel1996introduction,
  title     = {Introduction to solid state physics},
  author    = {Kittel, Charles and McEuen, Paul},
  volume    = {8},
  year      = {1996},
  publisher = {Wiley New York}
}

@article{joannopoulos2008molding,
  title   = {Molding the flow of light},
  author  = {Joannopoulos, John and Johnson, Steven Glenn and Winn, Joshua N. and Meade, Robert D.},
  journal = {Princeton Univ. Press, Princeton, NJ},
  year    = {2008}
}

@book{johnson2001photonic,
  title     = {Photonic crystals: the road from theory to practice},
  author    = {Johnson, Steven Glenn and Joannopoulos, John},
  year      = {2001},
  publisher = {Springer Science \& Business Media}
}

@incollection{kuchment2001mathematics,
  title     = {The mathematics of photonic crystals},
  author    = {Kuchment, Peter},
  booktitle = {Mathematical modeling in optical science},
  pages     = {207--272},
  year      = {2001},
  publisher = {SIAM},
  chapter = 7,
}

@book{sakoda2004optical,
  title     = {Optical properties of photonic crystals},
  author    = {Sakoda, Kazuaki},
  volume    = {80},
  year      = {2004},
  publisher = {Springer Science \& Business Media}
}

@article{joly2006exact,
  title   = {Exact boundary conditions for periodic waveguides containing a local perturbation},
  author  = {Joly, Patrick and Li, Jing-Rebecca and Fliss, Sonia},
  journal = {Commun. Comput. Phys},
  volume  = {1},
  number  = {6},
  pages   = {945--973},
  year    = {2006}
}

@phdthesis{fliss2009analyse,
  title       = {{Analyse math{\'e}matique et num{\'e}rique de probl{\`e}mes de propagation des ondes dans des milieux p{\'e}riodiques infinis localement perturb{\'e}s}},
  author      = {Fliss, Sonia},
  url         = {https://pastel.archives-ouvertes.fr/pastel-00005464},
  school      = {{Ecole Polytechnique X}},
  year        = {2009},
  month       = May,
  type        = {Theses},
  hal_id      = {pastel-00005464},
  hal_version = {v1}
}

@article{fliss2009exact,
  title     = {Exact boundary conditions for time-harmonic wave propagation in locally perturbed periodic media},
  author    = {Fliss, Sonia and Joly, Patrick},
  journal   = {Applied Numerical Mathematics},
  volume    = {59},
  number    = {9},
  pages     = {2155--2178},
  year      = {2009},
  publisher = {Elsevier}
}

@article{flisscassanbernier2010,
  title       = {{Computation of light refraction at the surface of a photonic crystal using DtN approach}},
  author      = {Fliss, Sonia and Cassan, Eric and Bernier, Damien},
  journal     = {{Journal of the Optical Society of America B}},
  publisher   = {{Optical Society of America}},
  volume      = {27},
  number      = {7},
  pages       = {1492-1503},
  year        = {2010},
  hal_id      = {hal-00873060},
  hal_version = {v1}
}

@article{fliss2012wave,
  title     = {Wave propagation in locally perturbed periodic media (case with absorption): Numerical aspects},
  author    = {Fliss, Sonia and Joly, Patrick},
  journal   = {Journal of Computational Physics},
  volume    = {231},
  number    = {4},
  pages     = {1244--1271},
  year      = {2012},
  publisher = {Elsevier}
}

@article{fliss2016solutions,
  title     = {Solutions of the time-harmonic wave equation in periodic waveguides: asymptotic behaviour and radiation condition},
  author    = {Fliss, Sonia and Joly, Patrick},
  journal   = {Archive for Rational Mechanics and Analysis},
  volume    = {219},
  pages     = {349--386},
  year      = {2016},
  publisher = {Springer}
}

@article{fliss2021dirichlet,
  title     = {A Dirichlet-to-Neumann approach to the mathematical and numerical analysis in waveguides with periodic outlets at infinity},
  author    = {Fliss, Sonia and Joly, Patrick and Lescarret, Vincent},
  journal   = {Pure and Applied Analysis},
  volume    = {3},
  number    = {3},
  pages     = {487--526},
  year      = {2021},
  publisher = {Mathematical Sciences Publishers}
}

@article{besse:hal-00698916,
  title       = {{Transparent boundary conditions for locally perturbed infinite hexagonal periodic media.}},
  author      = {Besse, Christophe and Coatleven, Julien and Fliss, Sonia and Lacroix-Violet, Ingrid and Ramdani, Karim},
  url         = {https://hal.archives-ouvertes.fr/hal-00698916},
  journal     = {{Communications in Mathematical Sciences}},
  publisher   = {{International Press}},
  volume      = {11},
  number      = {4},
  pages       = {907-938},
  year        = {2013},
  hal_id      = {hal-00698916},
  hal_version = {v1}
}

@article{coatleven2012helmholtz,
  title     = {Helmholtz equation in periodic media with a line defect},
  author    = {Coatl{\'e}ven, Julien},
  journal   = {Journal of Computational Physics},
  volume    = {231},
  number    = {4},
  pages     = {1675--1704},
  year      = {2012},
  publisher = {Elsevier}
}

@article{kuchment2004some,
  title     = {On some spectral problems of mathematical physics},
  author    = {Kuchment, Peter},
  journal   = {Partial differential equations and inverse problems},
  volume    = {362},
  pages     = {241--276},
  year      = {2004},
  publisher = {Amer. Math. Soc. Providence, RI}
}

@book{kuchment1993floquet,
  title     = {Floquet theory for partial differential equations},
  author    = {Kuchment, Peter A},
  volume    = {60},
  year      = {1993},
  publisher = {Springer Science \& Business Media}
}

@article{kozlov1979averaging,
  title     = {Averaging differential operators with almost periodic, rapidly oscillating coefficients},
  author    = {Kozlov, Sergei M},
  journal   = {Mathematics of the USSR-Sbornik},
  volume    = {35},
  number    = {4},
  pages     = {481},
  year      = {1979},
  publisher = {IOP Publishing}
}

@article{gerard2011homogenization,
  title   = {Homogenization in polygonal domains},
  author  = {G{\'e}rard-Varet, David and Masmoudi, Nader},
  journal = {Journal of the European Mathematical society},
  volume  = {13},
  number  = {5},
  pages   = {1477--1503},
  year    = {2011}
}

@article{gerard2012homogenization,
  title     = {Homogenization and boundary layers},
  author    = {{G{\'e}rard-Varet}, David and {Masmoudi}, Nader},
  journal   = {Acta mathematica},
  volume    = {209},
  number    = {1},
  pages     = {133--178},
  year      = {2012},
  publisher = {Springer}
}

@article{blanc2015local,
  title={Local profiles for elliptic problems at different scales: defects in, and interfaces between periodic structures},
  author={Blanc, Xavier and Le Bris, Claude and Lions, P-L},
  journal={Communications in Partial Differential Equations},
  volume={40},
  number={12},
  pages={2173--2236},
  year={2015},
  publisher={Taylor \& Francis}
}

@article{bouchitte2010homogenization,
  title     = {Homogenization of dielectric photonic quasi crystals},
  author    = {Bouchitt{\'e}, Guy and Guenneau, S{\'e}bastien and Zolla, Fr{\'e}d{\'e}ric},
  journal   = {Multiscale Modeling \& Simulation},
  volume    = {8},
  number    = {5},
  pages     = {1862--1881},
  year      = {2010},
  publisher = {SIAM}
}

@article{wellander2019homogenization,
  title   = {Homogenization of quasiperiodic structures and two-scale cut-and-projection convergence},
  author  = {Wellander, Niklas and Guenneau, S{\'e}bastien and Cherkaev, Elena},
  journal = {arXiv preprint arXiv:1911.03560},
  year    = {2019}
}

@article{combes1973asymptotic,
  title     = {Asymptotic behaviour of eigenfunctions for multiparticle Schr{\"o}dinger operators},
  author    = {Combes, Jean-Michel and Thomas, Lyn Carey},
  journal   = {Communications in Mathematical Physics},
  volume    = {34},
  number    = {4},
  pages     = {251--270},
  year      = {1973},
  publisher = {Springer}
}

@article{eidus1986limiting,
  title     = {The limiting absorption and amplitude principles for the diffraction problem with two unbounded media},
  author    = {Eidus, Daniel},
  journal   = {Communications in mathematical physics},
  volume    = {107},
  number    = {1},
  pages     = {29--38},
  year      = {1986},
  publisher = {Springer}
}

@book{weder1990spectral,
  title     = {Spectral and scattering theory for wave propagation in perturbed stratified media},
  author    = {Weder, Ricardo},
  volume    = {87},
  year      = {1990},
  publisher = {Springer Science \& Business Media}
}

@article{wilcox1966wave,
  title     = {Wave operators and asymptotic solutions of wave propagation problems of classical physics},
  author    = {Wilcox, Calvin H},
  journal   = {Archive for Rational Mechanics and Analysis},
  volume    = {22},
  number    = {1},
  pages     = {37--76},
  year      = {1966},
  publisher = {Springer}
}

@article{agmon1975spectral,
  title   = {Spectral properties of Schr{\"o}dinger operators and scattering theory},
  author  = {Agmon, Shmuel},
  journal = {Annali della Scuola Normale Superiore di Pisa-Classe di Scienze},
  volume  = {2},
  number  = {2},
  pages   = {151--218},
  year    = {1975}
}

@article{kirsch2018radiation,
  title     = {A radiation condition arising from the limiting absorption principle for a closed full-or half-waveguide problem},
  author    = {Kirsch, Andreas and Lechleiter, Armin},
  journal   = {Mathematical Methods in the Applied Sciences},
  volume    = {41},
  number    = {10},
  pages     = {3955--3975},
  year      = {2018},
  publisher = {Wiley Online Library}
}

@article{hoang2011limiting,
  title     = {The limiting absorption principle for a periodic semi-infinite waveguide},
  author    = {Hoang, Vu},
  journal   = {SIAM Journal on Applied Mathematics},
  volume    = {71},
  number    = {3},
  pages     = {791--810},
  year      = {2011},
  publisher = {SIAM}
}

@article{radosz2015new,
  title     = {New limiting absorption and limit amplitude principles for periodic operators},
  author    = {Radosz, Maria},
  journal   = {Zeitschrift f{\"u}r angewandte Mathematik und Physik},
  volume    = {66},
  number    = {2},
  pages     = {253--275},
  year      = {2015},
  publisher = {Springer}
}

@article{kirsch2018limiting,
  title     = {The limiting absorption principle and a radiation condition for the scattering by a periodic layer},
  author    = {Kirsch, Andreas and Lechleiter, Armin},
  journal   = {SIAM Journal on Mathematical Analysis},
  volume    = {50},
  number    = {3},
  pages     = {2536--2565},
  year      = {2018},
  publisher = {SIAM}
}

@article{kirsch2022scattering,
  title     = {A scattering problem for a local perturbation of an open periodic waveguide},
  author    = {Kirsch, Andreas},
  journal   = {Mathematical Methods in the Applied Sciences},
  year      = {2022},
  publisher = {Wiley Online Library}
}

@article{sommerfeld1912,
  author  = {Sommerfeld, A.},
  journal = {Jahresbericht der Deutschen Mathematiker-Vereinigung},
  pages   = {309-352},
  title   = {Die Greensche Funktion der Schwingungslgleichung.},
  url     = {http://eudml.org/doc/145344},
  volume  = {21},
  year    = {1912}
}

@article{murata2006asymptotics,
  title     = {Asymptotics of Green functions and the limiting absorption principle for elliptic operators with periodic coefficients},
  author    = {Murata, Minoru and Tsuchida, Tetsuo},
  journal   = {Journal of Mathematics of Kyoto University},
  volume    = {46},
  number    = {4},
  pages     = {713--754},
  year      = {2006},
  publisher = {Duke University Press}
}

@article{bonnet2002diffraction,
  title     = {Diffraction by an acoustic grating perturbed by a bounded obstacle},
  author    = {Bonnet-Bendhia, Anne-Sophie and Ramdani, Karim},
  journal   = {Advances in Computational Mathematics},
  volume    = {16},
  pages     = {113--138},
  year      = {2002},
  publisher = {Springer}
}

@article{mandel2019limiting,
  title     = {The limiting absorption principle for periodic differential operators and applications to nonlinear Helmholtz equations},
  author    = {Mandel, Rainer},
  journal   = {Communications in Mathematical Physics},
  volume    = {368},
  number    = {2},
  pages     = {799--842},
  year      = {2019},
  publisher = {Springer}
}

@book{hardy1979introduction,
  title     = {An introduction to the theory of numbers},
  author    = {Hardy, Godfrey Harold and Wright, Edward Maitland and others},
  year      = {1979},
  publisher = {Oxford university press}
}

@book{diestel_vector_1977,
  title     = {Vector {Measures}},
  isbn      = {978-0-8218-1515-1},
  language  = {en},
  publisher = {American Mathematical Soc.},
  author    = {Diestel, Joseph and Uhl, John Jerry},
  month     = jun,
  year      = {1977}
}

@article{trefethen2014exponentially,
  title={The exponentially convergent trapezoidal rule},
  author={Trefethen, Lloyd N and Weideman, JAC},
  journal={SIAM review},
  volume={56},
  number={3},
  pages={385--458},
  year={2014},
  publisher={SIAM}
}

@article{amenoagbadji2023wave,
  title   = {Wave propagation in one-dimensional quasiperiodic media},
  author  = {Amenoagbadji, Pierre and Fliss, Sonia and Joly, Patrick},
  journal = {Communications in Optimization Theory},
  doi     = {10.23952/cot.2023.17},
  year    = {2023},
}

\end{document}